\documentclass[12pt]{book}

\usepackage{yhmath}
\usepackage{tabularx}

\usepackage{amstext    }
\usepackage{amsthm    }
\usepackage{a4}
\usepackage[mathscr]{eucal} 
\usepackage{mathrsfs} 

\usepackage{amsmath}
\usepackage{amssymb}
\usepackage{amscd}

\newtheorem{theorem}{Theorem}[section]
\newtheorem{definition}[theorem]{Definition}

\newtheorem{proposition}[theorem]{Proposition}
\newtheorem{corollary}[theorem]{Corollary}
\newtheorem{lemma}[theorem]{Lemma}
\newtheorem{remark}[theorem]{Remark}

\newtheorem{example}[theorem]{Example}
\newtheorem{examples}[theorem]{Examples}
\newtheorem{exercise}{Exercise}[section]

\newcommand{\cali}[1]{\mathscr{#1}}

\newcommand{\GL}{{\rm GL}}

\newcommand{\PGL}{{\rm PGL}}
\newcommand{\Jac}{{\rm Jac}}
\newcommand{\area}{{\rm area}}
\newcommand{\volume}{{\rm volume}}
\newcommand{\vol}{{\rm volume}}

\newcommand{\lov}{{\rm lov}}
\newcommand{\sing}{{\rm sing}}
\newcommand{\reg}{{\rm reg}}
\newcommand{\trans}{{\ \!\!\sp t\ \!\!}}
\newcommand{\supp}{{\rm supp}}
\newcommand{\const}{{\rm const}}
\newcommand{\dist}{{\rm dist}}
\newcommand{\diam}{{\rm diam}}
\newcommand{\mult}{{\rm mult}}
\newcommand{\codim}{{\rm codim}}

\newcommand{\Grass}{{\rm Grass}}

\newcommand{\loc}{{loc}}
\newcommand{\ddc}{{dd^c}}
\newcommand{\dc}{{d^c}}
\newcommand{\dbar}{{\overline\partial}}
\newcommand{\ddbar}{{\partial\overline\partial}}
\newcommand{\PSH}{{\rm PSH}}
\newcommand{\DSH}{{\rm DSH}}
\newcommand{\FS}{{\rm FS}}

\newcommand{\ind}{{\bf 1}}
\newcommand{\id}{{\rm id}}

\newcommand{\capacity}{{\rm cap}}

\renewcommand{\L}{{\cal L}}

\newcommand{\Bt}{{\tt B}}
\newcommand{\Et}{{\tt E}}
\newcommand{\St}{{\tt S}}

\newcommand{\gr}{{\tt g}}

\newcommand{\Ac}{\cali{A}}
\newcommand{\Bc}{\cali{B}}
\newcommand{\Cc}{\cali{C}}
\newcommand{\Dc}{\cali{D}}
\newcommand{\Ec}{\cali{E}}
\newcommand{\Fc}{\cali{F}}

\newcommand{\Hc}{\cali{H}}
\newcommand{\Jc}{\cali{J}}
\newcommand{\Kc}{\cali{K}}

\newcommand{\Mc}{\cali{M}}

\newcommand{\Oc}{\cali{O}}
\newcommand{\Pc}{\cali{P}}
\newcommand{\Rc}{\cali{R}}
\newcommand{\Sc}{\cali{S}}
\newcommand{\Uc}{\cali{U}}
\newcommand{\Vc}{\cali{V}}

\newcommand{\C}{\mathbb{C}}

\newcommand{\Z}{\mathbb{Z}}
\newcommand{\R}{\mathbb{R}}

\newcommand{\T}{\mathbb{T}}
\renewcommand{\P}{\mathbb{P}}

\title{Dynamics in several complex variables: endomorphisms of
  projective spaces and polynomial-like mappings}
\author{Tien-Cuong Dinh and Nessim Sibony}
\date{}

\begin{document}

\maketitle

\newpage 
\thispagestyle{empty}
\ \newpage

\ 

\bigskip\bigskip\bigskip

{\small

\begin{center}
{\bf Abstract}
\end{center}

\medskip

The emphasis of this introductory course is on pluripotential methods in
complex dynamics in higher dimension. They are based on the
compactness properties of plurisubharmonic (p.s.h.) functions and on the 
theory of positive closed currents. Applications of these methods are not limited
to the dynamical systems that we consider
here. Nervertheless, we choose to show their effectiveness 
and to describe the theory for two large families of maps: the endomorphisms of
projective spaces and the polynomial-like mappings.

The first chapter deals with holomorphic endomorphisms of 
the projective space $\P^k$. We establish the first
properties and give several constructions for the Green currents $T^p$
and the
equilibrium measure $\mu=T^k$. The emphasis is on quantitative properties and
speed of convergence. We then treat equidistribution problems. We show
the existence of a proper algebraic set
$\Ec$, totally invariant, i.e. $f^{-1}(\Ec)=f(\Ec)=\Ec$, such that
when $a\not \in \Ec$, the probability measures, equidistributed on
the fibers $f^{-n}(a)$, converge towards the equilibrium measure $\mu$,
as $n$ goes to infinity. A similar result holds for the restriction of
$f$ to invariant subvarieties.
We survey the equidistribution problem when
points are replaced by varieties of arbitrary dimension, and discuss the 
equidistribution of
periodic points. We then
establish ergodic properties of $\mu$: 
K-mixing, exponential decay of
correlations for various classes of observables, central limit theorem
and large deviations theorem. We heavily use the compactness of the
space $\DSH(\P^k)$ of differences of quasi-p.s.h. functions. In
particular, we show that the measure $\mu$ is moderate, i.e. $\langle
\mu, e^{\alpha|\varphi|}\rangle \leq c$, on bounded sets of $\varphi$ in
$\DSH(\P^k)$, for suitable positive constants $\alpha,c$. 
Finally, we study the entropy, the Lyapounov exponents and the
dimension of $\mu$.

The second chapter develops the theory of polynomial-like maps,
i.e. proper holomorphic maps $f:U\rightarrow V$ where $U,V$ are open
subsets of $\C^k$ with $V$ convex and $U\Subset V$. We introduce the
dynamical degrees for such maps and construct the equilibrium measure
$\mu$ of maximal entropy. Then,  under a natural assumption on the dynamical
degrees, we prove equidistribution
properties of points and various statistical properties of the measure
$\mu$. The assumption is stable under small pertubations on the map. 
We also study the dimension of $\mu$, the Lyapounov exponents
and their variation.

Our aim is to get a self-contained text that requires only a minimal background. 
In order to help the reader, 
an appendix gives the basics on p.s.h. functions,
positive closed currents and super-potentials on projective spaces. Some exercises
are proposed and an extensive bibliography is given.
}

\bigskip



\newpage
\thispagestyle{empty}

\tableofcontents

\newpage
\thispagestyle{empty}

\chapter*{Introduction}
\addcontentsline{toc}{chapter}{Introduction}

These  notes are based on a series of lectures given by
the authors at IHP in 2003, Luminy in 2007, Cetraro in 2008 and
Bedlewo 2008. 
The purpose is to provide an introduction to some
developments in dynamics of several complex variables. 
We have chosen to treat here only two chapters of the theory: the dynamics
of endomorphisms of the projective space $\P^k$ and the dynamics of
polynomial-like mappings in higher dimension. Besides the basic notions and results, we
describe the recent developments and the new tools introduced in the
theory. These tools are useful in other settings.  
We tried to give a complete picture of the theory for the above families
of dynamical systems.
Meromorphic maps on compact K\"ahler manifolds, in particular
polynomial automorphisms of
$\C^k$, will be studied in a forthcoming survey.

Let us comment on how complex dynamics fits in the general theory of dynamical
systems. The abstract ergodic theory is well-developed with
remarkable achievements like the Oseledec-Pesin theory. It is however
difficult to show in concrete examples that an invariant measure 
satisfies exponential decay of correlations for smooth
observables or is hyperbolic, i.e. has only non-zero Lyapounov
exponents, see e.g.  Benedicks-Carleson \cite{BenedicksCarleson}, Viana
\cite{Viana}, L.S. Young
\cite{Young,Young1}.
One of our goals is to show that holomorphic dynamics in several
variables provides remarkable examples of non-uniformly hyperbolic
systems where the abstract theory can be applied. Powerful tools from 
the theory of several complex variables permit to avoid delicate
combinatorial estimates. 
Complex dynamics also require a development of new tools like the calculus on
currents and the introduction of new spaces of observables, which are
of independent interest.

Complex dynamics in dimension one, i.e. dynamics of
rational maps on $\P^1$,
is well-developed and has in some sense reached maturity. The
main tools there are Montel's theorem on normal families, the Riemann
measurable mapping theorem and the theory of quasi-conformal
maps, see e.g. Beardon, Carleson-Gamelin \cite{Beardon,
  CarlesonGamelin}. 
When dealing with maps in several variables such tools are not
available: the Kobayashi hyperbolicity of a manifold and the
possibility to apply normal family arguments, are more difficult to
check. Holomorphic maps in several variables are not conformal and
there is no Riemann measurable mapping theorem.

The theory in higher dimension is developed using mostly pluripotential theory, i.e. the
theory of plurisubharmonic (p.s.h. for short) functions and positive closed currents. The
Montel's compactness property is replaced by the compactness properties
of p.s.h. or quasi-p.s.h. functions. Another crucial tool is the use
of good estimates for the $\ddc$-equation in various settings. One
of the main ideas is: in order to study the statistical behavior of
orbits of a holomorphic map, we consider its action on some
appropriate functional spaces. We then decompose the action into the
``harmonic'' part and the ``non-harmonic'' one. This is done solving a
$\ddc$-equation with estimates. The non-harmonic part of the dynamical action may be
controled thanks to good estimates for the solutions of a
$\ddc$-equation. The harmonic part can be treated using either 
Harnack's inequality in the local setting or the linear action of maps
on cohomology groups in the case of dynamics on compact K{\"a}hler manifolds.
This approach has permitted to
give a satisfactory theory of the ergodic properties of holomorphic
and meromorphic dynamical systems: construction of the measure of
maximal entropy, decay of correlations, central limit theorem, large
deviations theorem, etc. with respect to that measure. 

In order to use the pluripotential methods, we are led to develop
the calculus on positive closed currents. 
Readers not familiar with these theories may
start with the appendix at the end of these notes where 
we have gathered some notions and 
results on currents and pluripotential theory. 
A large part in the appendix is classical but there are
also some recent results, mostly on 
new spaces of currents and on the notion of super-potential associated
to positive closed currents in higher bidegree. 
Since we only deal here with projective spaces
and open sets in $\C^k$, this is easier and the background is limited.

The main problem in the dynamical study of a map is to understand the behavior of
the orbits of points under the action of the map.
Simple examples show that in general there is a set (Julia set) where
the dynamics is unstable: the orbits may diverge exponentially. 
Moreover, the geometry of the Julia set is in general very
wild. In order to study complex dynamical systems, we follow the
classical concepts. We introduce and establish basic
properties of some invariants
associated to the system, like the topological entropy and the
dynamical degrees which are
the analogues of volume growth indicators in the real dynamical
setting. These invariants give a
rough classification of the system. The remarkable fact in
complex dynamics is that they can be computed or estimated in many
non-trivial situations.

A central question in dynamics is to construct
interesting invariant measures, in particular, measures with positive
entropy. Metric entropy is an indicator of the complexity of the
system with respect to an invariant measure. We focus our
study on the measure of maximal entropy. Its support is 
in some sense the most chaotic part of the
system. For the maps we consider here, measures of maximal entropy are
constructed using pluripotential methods. For endomorphisms in $\P^k$,
they can be obtained as self-intersections of some invariant positive
closed $(1,1)$-currents (Green currents). 
We give estimates on the
Hausdorff dimension and on Lyapounov exponents of these measures.
The results give the behavior on the most chaotic part. 
Lyapounov exponents are shown to be strictly positive. This means in some
sense that the system is expansive in all directions, despite of the
existence of a critical set.

Once, the measure of maximal entropy is constructed, we 
study its fine dynamical properties. Typical orbits can be observed using
test functions. Under the action of the map, each observable provides a
sequence of functions that can be seen as dependent random variables. 
The aim is to show that the dependence is weak and then 
to establish  stochastic properties
which are known for 
independent random variables in probability theory.
Mixing, decay of correlations, central limit theorem, large deviations
theorems, etc. are proved for the measure of maximal entropy. It is crucial here that 
the Green currents and the measures of maximal entropy are obtained using an
iterative process with estimates; we can then bound the speed of convergence.

Another problem, we consider in these notes, is the
equidistribution of periodic points or of preimages of points with
respect to the measure of maximal entropy. For endomorphisms of
$\P^k$, we also study the equidistribution of
varieties with respect to the Green currents. Results in this
direction give some informations on the rigidity of the system and
also some strong ergodic properties that the Green currents or the 
measure of maximal entropy satisfy. 
The results we obtain are in spirit similar to a second main theorem in
value distribution theory and 
should be useful in order to study the
arithmetic analogues. 
We give complete proofs for most results,
but we only survey the equidistribution of hypersurfaces and results
using super-potentials, in particular, the equidistribution of
subvarieties of higher codimension. We have given exercises,
basically in each section, some of them are not straightforward.

The text is organized as follows. In the first chapter, 
we study holomorphic endomorphisms of $\P^k$. We
introduce several methods in order to construct and to study the Green currents and the
Green measure, i.e. equilibrium measure or measure of maximal entropy. 
These methods were not originally introduced in this
setting but here they are simple and very effective. The reader will find a
description and the references of the earlier approach in the ten
years old survey by the second author \cite{Sibony}. 
The second chapter deals with a very large family of maps:
polynomial-like maps. In this case, $f:U\rightarrow V$ is proper and
defined on an open set $U$, strictly contained in a convex domain $V$
of $\C^k$. Holomorphic endomorphisms of $\P^k$ can be lifted to a
polynomial-like maps on some open set in $\C^{k+1}$. So, we can
consider polynomial-like maps as a semi-local version of the
endomorphisms studied in the first chapter. They can appear in the
study of meromorphic maps or in the dynamics of transcendental
maps. The reader will find in the end of these notes an appendix on
the theory of currents and an extensive bibliography.


\chapter{Endomorphisms of projective spaces} \label{chapter_endo}

In this chapter, we give the main results on the dynamics of holomorphic
maps on the projective space $\P^k$. Several results are recent and
some of them are new even in dimension 1. The reader will find here 
an introduction to methods that can be developed
in other situations, in particular, in the study of meromorphic maps
on arbitrary compact K\"ahler manifolds. The main references for this
chapter are \cite{BriendDuval1, BriendDuval2,
  DinhNguyenSibony3, DinhSibony9, DinhSibony10, FornaessSibony1, Sibony}.

\section{Basic properties and examples}

Let $f:\P^k\rightarrow\P^k$ be a holomorphic endomorphism.
Such a map is always induced by a  polynomial self-map $F=(F_0,\ldots,F_k)$
on $\C^{k+1}$ such that $F^{-1}(0)=\{0\}$ and 
the components $F_i$ are homogeneous polynomials of the same
degree $d\geq 1$. Given an endomorphism $f$, the associated map $F$ is unique up to a
multiplicative constant and is called {\it a lift} of $f$ to $\C^{k+1}$.
From now on, assume that $f$ is non-invertible, i.e. the {\it
  algebraic degree} $d$ is at least  2. Dynamics of an invertible
map is simple to study. If
$\pi:\C^{k+1}\setminus\{0\}\rightarrow\P^k$ is the natural projection,
we have $f\circ\pi=\pi\circ F$. Therefore, dynamics of holomorphic
maps on
$\P^k$ can be deduced from the polynomial case in $\C^{k+1}$. 
We will count preimages of points, periodic points, introduce Fatou and
Julia sets and give some examples. 

It is easy to construct examples of holomorphic maps in $\P^k$. The family of
homogeneous polynomial maps $F$ of a given degree $d$ is parametrized
by a complex vector space of dimension $N_{k,d}:=(k+1)(d+k)!/(d!k!)$. The maps
satisfying $F^{-1}(0)=\{0\}$ define a Zariski dense open
set. Therefore, the parameter space $\Hc_d(\P^k)$, of holomorphic
endomorphisms of algebraic degree $d$, is a Zariski dense open set in
$\P^{N_{k,d}-1}$, in particular, it is connected.

If $f:\C^k\rightarrow\C^k$ is a polynomial map, we can extend $f$ to
$\P^k$ but the extension is not always holomorphic. The extension is
holomorphic 
when the dominant homogeneous part $f^+$ of $f$, satisfies
$(f^+)^{-1}(0)=\{0\}$. Here, if $d$ is the maximal degree in the
polynomial expression of $f$,  then $f^+$ is composed by the
monomials of degree $d$ in the components of $f$. So, it is easy to
construct examples using products of one dimensional polynomials or
their pertubations.

A general meromorphic map $f:\P^k\rightarrow\P^k$ of algebraic
degree $d$ is given in homogeneous coordinates by
$$f[z_0:\cdots:z_k]=[F_0:\cdots:F_k],$$
where the components $F_i$ are homogeneous polynomials of degree $d$
without common factor, except
constants. The map $F:=(F_0,\ldots,F_k)$ on $\C^{k+1}$ is still called
{\it a lift} of $f$. In general, $f$ is not defined on the analytic set 
$I=\{[z]\in\P^k, F(z)=0\}$ which is of codimension $\geq 2$ since the
$F_i$'s have no common factor. This is the {\it indeterminacy set of
  $f$} which is empty when $f$ is holomorphic.

It is easy to check that if $f$ is in $\Hc_d(\P^k)$ and $g$ is in
$\Hc_{d'}(\P^k)$, the composition $f\circ g$  belongs to $\Hc_{dd'}(\P^k)$. This is in
general false for meromorphic maps: the
algebraic degree of the composition is not necessarily equal to the
product of the algebraic degrees. It is enough to consider the
meromorphic involution of algebraic degree $k$
$$f[z_0:\cdots:z_k]:=\Big[{1\over z_0}:\cdots:{1\over z_k}\Big]
=\Big[{z_0\ldots z_k\over z_0}:\cdots:{z_0\ldots z_k\over z_k}\Big].$$ 
The composition $f\circ f$ is the identity map.

We say that $f$ is {\it dominant} if $f(\P^k\setminus I)$ contains a
non-empty open set. The space of dominant meromorphic maps of
algebraic degree $d$, is denoted by $\Mc_d(\P^k)$. It is also a 
Zariski dense open set in $\P^{N_{k,d}-1}$. A result by 
Guelfand, Kapranov and Zelevinsky
shows that $\Mc_d(\P^k)\setminus \Hc_d(\P^k)$ is an irreducible
algebraic variety \cite{GuelfandKapranov}. We will be concerned in this chapter mostly with
holomorphic maps. We show that they are open
and their topological degree, i.e. the number of points in a generic
fiber, is equal to $d^k$. We recall here the classical B\'ezout's
theorem which is a central tool for the dynamics in $\P^k$. 

\begin{theorem}[B\'ezout] Let $P_1,\ldots,P_k$ be homogeneous
  polynomials in $\C^{k+1}$ of degrees $d_1,\ldots,d_k$
  respectively. Let $Z$ denote the set of common zeros of $P_i$, in
  $\P^k$, i.e. the set of points $[z]$ such that $P_i(z)=0$ for $1\leq
  i\leq k$. If $Z$ is discrete, then the number of points in $Z$,
  counted with multiplicity, is $d_1\ldots d_k$. 
\end{theorem}

The multiplicity of a point $a$ in $Z$ can be defined in several
ways. For instance, if $U$ is a small neighbourhood of $a$ and
if $P_i'$ are generic homogeneous polynomials of degrees $d_i$ close
enough to $P_i$, then the hypersurfaces $\{P_i'=0\}$ in $\P^k$
intersect transversally. The number of points of the intersection in
$U$ does not depend on the choice of $P_i'$ and is {\it the multiplicity} of
$a$ in $Z$.

\begin{proposition}
Let $f$ be an endomorphism of algebraic degree $d$ of
$\P^k$. Then for every $a$ in $\P^k$, the fiber $f^{-1}(a)$ contains
exactly $d^k$ points, counted with multiplicity. In particular, $f$
is open and defines a ramified covering of degree $d^k$. 
\end{proposition}
\proof
For the multiplicity of $f$ and the notion of ramified covering, we
refer to Appendix \ref{section_pk}. Let
$f=[F_0:\cdots :F_k]$ be an expression of $f$ in homogeneous
coordinates. Consider a point $a=[a_0:\cdots:a_k]$ in $\P^k$. Without
loss of generality, we can assume $a_0=1$, hence
$a=[1:a_1:\cdots:a_k]$. The points in $f^{-1}(a)$ are the common
zeros, in $\P^k$,
of the polynomials $F_i-a_iF_0$ for $i=1,\ldots,k$. 

We have to check that the common zero
set is discrete, then B\'ezout's theorem asserts that the cardinality
of this set is equal to the product of the degrees of $F_i-a_iF_0$,
i.e. to $d^k$. If the set were not discrete, then the common zero set of
$F_i-a_iF_0$ in $\C^{k+1}$ is analytic of dimension $\geq 2$. This
implies that the set of common zeros of the $F_i$'s, $0\leq i\leq k$, in $\C^{k+1}$ is of positive
dimension. This is impossible when $f$ is holomorphic. So, $f$ is a
ramified covering of degree $d^k$. In particular, it is open. 

Note that when $f$ is a map in $\Mc_d(\P^k)\setminus \Hc_d(\P^k)$ with
indeterminacy set $I$, we can prove that the generic fibers of 
$f:\P^k\setminus I\rightarrow\P^k$ contains at most $d^k-1$
points. Indeed, for every $a$, the hypersurfaces $\{F_i-a_iF_0=0\}$ in $\P^k$ contain $I$.
\endproof

Periodic points of order $n$, i.e. points which satisfy $f^n(z)=z$,
play an important role in dynamics. Here, $f^n:=f\circ \cdots\circ f$,
$n$ times, is {\it the iterate of order $n$} of $f$. Periodic points
of order $n$ of $f$ are fixed points of $f^n$ which is an endomorphism
of algebraic degree $d^n$. In the present case, their
counting is simple. We have the following result.

\begin{proposition} \label{prop_number_periodic_pk}
Let $f$ be an endomorphism of algebraic
  degree $d\geq 2$ in $\P^k$. Then the number of fixed points of $f$, counted
  with multiplicity, is equal to $(d^{k+1}-1)/(d-1)$. In particular, the number
  of periodic points of order $n$ of $f$ is $d^{kn}+o(d^{kn})$. 
\end{proposition}
\proof
There are several methods to count the periodic points. 
In $\P^{k+1}$, with homogeneous coordinates $[z:t]=[z_0:\cdots:z_k:t]$,
we consider the system of equations 
$F_i(z)-t^{d-1}z_i=0$. The set is
discrete since it is analytic and does not intersect the hyperplane
$\{t=0\}$.  
So, we can count the solutions of the above system using B\'ezout's theorem and we find $d^{k+1}$
points counting with multiplicity. Each point $[z:t]$ in this set, except
$[0:\cdots:0:1]$, corresponds to a fixed point $[z]$ of $f$. The
correspondence is $d-1$ to $1$. Indeed, if we multiply $t$ by a
$(d-1)$-th root of unity, we get the same fixed point. Hence, the
number of fixed points of $f$ counted
  with multiplicity is $(d^{k+1}-1)/(d-1)$.

The number of fixed points of $f$ is also the number of points in the
intersection of the graph of $f$ with the diagonal of
$\P^k\times\P^k$. So, we can count these points using the cohomology
classes associated to the above analytic sets, i.e. using the
Lefschetz fixed point formula, see \cite{GriffithsHarris}. We can also observe that
this number depends continuously on $f$. So, it is constant for $f$ in
$\Hc_d(\P^k)$ which is connected. We obtain the result by counting the fixed points of an
example, e.g. for $f[z]=[z_0^d:\cdots:z_k^d]$. 
\endproof

Note that the periodic points of period $n$ are isolated. If $p$ is
such a point, a theorem due to Shub-Sullivan
\cite[p.323]{KatokHasselblatt} implies that the multiplicity at $p$ of
the equation $f^{mn}(p)=p$ is bounded independently on $m$. The result
holds for $\Cc^1$ maps. We deduce from the above result that $f$ admits infinitely many 
distinct periodic points.

The set of fixed points of a meromorphic map could be
empty or infinite.  One checks
easily that the map $(z_1,z_2)\mapsto (z_1^2,z_2)$  in $\C^2$ admits $\{z_1=0\}$
as a curve of fixed points. 

\begin{example}\rm
Consider the following map:
$$f(z_1,z_2):=(z_1+1,P(z_1,z_2)),$$ 
where $P$ is a homogeneous polynomial of degree $d\geq 2$ such that
$P(0,1)=0$. It is clear that $f$ has no periodic point in $\C^2$.
The meromorphic extension of $f$ is given in homogeneous
coordinates $[z_0:z_1:z_2]$ by 
$$f[z]=[z_0^d:z_0^{d-1}z_1+z_0^d:P(z_1,z_2)].$$ 
Here, $\C^2$ is identified to the open set $\{z_0=1\}$ of $\P^2$. The
indeterminacy set $I$ of $f$ is defined by $z_0=P(z_1,z_2)=0$ and 
is contained in the line at infinity $L_\infty:=\{z_0=0\}$. We have
$f(L_\infty\setminus I)= [0:0:1]$ which is an indeterminacy point. So,
$f:\P^2\setminus I\rightarrow\P^2$ has no periodic point.
\end{example}

\begin{example}\rm
Consider the holomorphic map $f$ on $\P^2$ given by
$$f[z]:=[z_0^d+P(z_1,z_2),z_2^d+\lambda z_0^{d-1}z_1:z_1^d],$$
with $P$ homogeneous of degree $d\geq 2$. Let $p:=[1:0:0]$, then
$f^{-1}(p)=p$. Such a point is called {\it totally invariant}. In
general, $p$ is not necessarily an attractive point. Indeed, the
eigenvalues of the differential of $f$ at $p$ are $0$ and
$\lambda$. When $|\lambda|>1$, there is an expansive direction for $f$
in a neighbourhood of $p$. 
In dimension one, totally invariant points are always
attractive.   
\end{example}

For a holomorphic map $f$ on $\P^k$, a point $a$ in $\P^k$ is {\it critical} if $f$ is not
injective in a neighbourhood of $a$ or equivalently the multiplicity
of $f$ at $a$ in the fiber $f^{-1}(f(a))$ is strictly larger than 1, see
Theorem \ref{th_anal_set_local}.
We say
that $a$ is a critical point of {\it multiplicity $m$} if  the multiplicity
of $f$ at $a$ in the fiber $f^{-1}(f(a))$ is equal to $m+1$. 

\begin{proposition}
Let $f$ be a holomorphic endomorphism of algebraic degree $d\geq 2$ of
$\P^k$. Then, the critical set of $f$ is an algebraic  hypersurface of degree
$(k+1)(d-1)$ counted with multiplicity. 
\end{proposition}
\proof
If $F$ is a lift of $f$ to $\C^{k+1}$, the Jacobian $\Jac(F)$ is a
homogeneous polynomial of degree $(k+1)(d-1)$. The zero set of
$\Jac(F)$ in $\P^k$ is exactly the critical set of $f$. The result follows. 
\endproof

Let $\Cc$ denote the critical set of $f$. The orbit $\Cc,
f(\Cc)$, $f^2(\Cc),\ldots$ is either a hypersurface or a countable union of
hypersurfaces. We say that $f$ is {\it postcritically finite} if this
orbit is a hypersurface, i.e. has only finitely many irreducible components. 
Besides very simple examples, postcritically finite maps are 
difficult to construct, because the image of a variety
is in general a variety of larger degree.
We give few examples of postcritically finite maps, see
\cite{FornaessSibony4, FornaessSibony7}. 

\begin{examples} \rm
We can check that for $d\geq 2$ and $(1-2\lambda)^d=1$
$$f[z_0:\cdots:z_k]:=[z_0^d:\lambda(z_0-2z_1)^d:\cdots:\lambda(z_0-2z_k)^d]$$ 
is postcritically finite. For some parameters
$\alpha\in\C$ and $0\leq l\leq d$, the map 
$$f_\alpha[z]:=[z_0^d:z_1^d: z_2^d+\alpha z_1^{d-l}z_2^l]$$ 
is also postcritically finite. In particular, for
$f_0[z]=[z_0^d:z_1^d:z_2^d]$, the associated critical set
is equal to $\{z_0z_1z_2=0\}$ which is invariant under $f_0$. So, $f_0$ is postcritically
 finite.
\end{examples}

Arguing as above, using B\'ezout's theorem, we can prove that if $Y$ is
an analytic set of pure codimension $p$ in $\P^k$ then $f^{-1}(Y)$ is
an analytic set of pure codimension $p$. Its degree, counting with
multiplicity, is equal to $d^p\deg(Y)$. Recall that the degree
$\deg(Y)$ of $Y$ is the number of points in the intersection of $Y$
with a generic projective subspace of dimension $p$. We deduce 
that the pull-back operator $f^*$ on the Hodge cohomology group
$H^{p,p}(\P^k,\C)$ is simply a multiplication by $d^p$. Since $f$ is a
ramified covering of degree $d^k$, $f_*\circ f^*$ is the multiplication
by $d^k$. Therefore, the push-forward operator $f_*$ acting on
$H^{p,p}(\P^k,\C)$ is the multiplication by $d^{k-p}$. In particular,
the image $f(Y)$ of $Y$ by $f$ is an analytic set of pure codimension $p$ 
and of degree $d^{k-p}\deg(Y)$, counted with multiplicity.

We now introduce the Fatou and Julia sets associated to 
an endomorphism. The following definition is analogous to
the one variable case.

\begin{definition} \rm 
{\it The Fatou set} of
$f$ is the largest open set $\Fc$ in $\P^k$ where the sequence of iterates
$(f^n)_{n\geq 1}$ is locally equicontinuous. 
The complement $\Jc$ of $\Fc$ is called {\it the Julia set} of $f$.
\end{definition}

Fatou and Julia sets are totally invariant by $\Fc$, that is,
$f^{-1}(\Fc)=f(\Fc)=\Fc$ and the same property holds for $\Jc$. Julia and
Fatou sets associated to $f^n$ are also equal to $\Jc$ and $\Fc$. We see
here that the space $\P^k$ is divided into two parts: on $\Fc$ the dynamics
is stable and tame while the dynamics on $\Jc$ is a priori chaotic. If $x$ is a
point in $\Fc$ and $y$ is
close enough to $x$, the orbit of $y$ is close to the orbit of $x$ when
the time $n$ goes to infinity. On the Julia set, this property is not
true. Attractive fixed points and their basins are examples of
points in the Fatou set. Siegel domains, i.e. invariant domains on
which $f$ is conjugated to a rotation, are also in the Fatou set. Repelling periodic points are always in the
Julia set. Another important notion in dynamics is non-wandering set.

\begin{definition} \rm
A point $a$ in $\P^k$ is {\it non-wandering} with respect to $f$ if
for every neighbourhood $U$ of $a$, there is an $n\geq 1$ such that
$f^n(U)\cap U\not=\varnothing$. 
\end{definition}

The study of the Julia and Fatou sets is a fundamental problem in
dynamics. It is quite well-understood in the one variable case
where the Riemann measurable theorem is a basic tool. The help of
computers is also important there. In higher dimension, Riemann measurable
theorem is not valid and the use of computers is more delicate. 
The most important tool in higher dimension is 
pluripotential theory.

For instance, Fatou and Julia sets for a general map are far from being
understood.
Many fundamental questions are still open. We
do not know if wandering Fatou components exist in higher
dimension. In dimension one, a theorem due to Sullivan
\cite{Sullivan} says that such a domain does not exist. 
The classification of Fatou components is not known,
see \cite{FornaessSibony8} for a partial answer in dimension 2
and \cite{FornaessSibony4, Sibony,Ueda} for the case of postcritically finite maps. 
The reader will find
in the survey \cite{Sibony} some results on local dynamics
near a fixed point, related to the Fatou-Julia theory. We now give  few examples.

The following construction is due to Ueda \cite{Ueda}. It is useful in
order to obtain interesting examples, in particular, to show that some
properties hold for generic maps. The strategy is to check that the set
of maps satisfying these
properties is a Zariski open set in the space of parameters and then
to produce an example using Ueda's construction.

\begin{examples} \rm \label{example_ueda}
Let $h:\P^1\rightarrow\P^1$ be a rational map of degree $d\geq
2$. Consider the multi-projective space $\P^1\times\cdots\times\P^1$,
$k$ times. The permutations of coordinates define a finite group $\Gamma$
acting on this space and the quotient of  $\P^1\times\cdots\times\P^1$ by $\Gamma$ is equal to $\P^k$. Let
$\Pi:\P^1\times\cdots\times\P^1\rightarrow\P^k$ denote the canonical
projection. Let $\widetilde f$ be the endomorphism of
$\P^1\times\cdots\times\P^1$ defined by $\widetilde
f(z_1,\ldots,z_k):=(h(z_1),\ldots, h(z_k))$. 
If $\sigma$ is a permutation of coordinates $(z_1,\ldots,z_k)$, then
$\sigma\circ \widetilde f= \widetilde f\circ\sigma$.
It is not difficult to
deduce that there is an endomorphism $f$ on $\P^k$ of algebraic degree
$d$ semi-conjugated to $\widetilde f$, that is, $f\circ\Pi=\Pi\circ\widetilde
f$. One can deduce dynamical
properties of $f$ from properties of $h$. For example, if $h$ is
chaotic, i.e. has a dense orbit, then $f$ is also chaotic. The first chaotic maps on $\P^1$ were
constructed by Latt{\`e}s. Ueda's construction gives Latt{\`e}s maps in
higher dimension. A {\it Latt{\`e}s map} $f$ on $\P^k$ is a map
semi-conjugated to an affine map on a torus. More precisely, there is
an open holomorphic map $\Psi:\T\rightarrow\P^k$ from a $k$-dimensional
torus $\T$ onto
$\P^k$ and an affine map $A:\T\rightarrow \T$ such that
$f\circ\Psi=\Psi\circ A$. We refer to \cite{BertelootDupont, BertelootLoeb,
Dinh4, DinhSibony0, Milnor} for a discussion of
Latt{\`e}s maps.
\end{examples}

The following map is the simplest in our context. Its iterates
can be explicitely computed. The reader may use this map and its
pertubations as basic examples 
in order to get a picture on the objects we will introduce latter.

\begin{example} \rm \label{example_power_map}
Let $f:\C^k\rightarrow \C^k$ be the polynomial map defined by 
$$f(z_1,\ldots,z_k):=(z_1^d,\ldots,z_k^d),\quad d\geq 2.$$
We can extend $f$ holomorphically to $\P^k$. 
Let $[z_0:\cdots:z_k]$ denote the homogeneous coordinates on $\P^k$
such that $\C^k$ is identified to the chart $\{z_0\not=0\}$. Then, the
extension of $f$ to $\P^k$ is
$$f[z_0:\cdots:z_k]=[z_0^d:\cdots:z_k^d].$$
The Fatou set is the union of the basins of the $k+1$ attractive fixed points
$[0:\cdots:0:1:0:\cdots:0]$.
These components are defined by
$$\Fc_i:=\big\{z\in\P^k,\quad |z_j|<|z_i|\quad \mbox{for every } j\not=i\big\}.$$ 
The Julia set of $f$ is the union of the following sets $\Jc_{ij}$ with
$0\leq i<j\leq k$, where 
$$\Jc_{ij}:=\big\{z\in\P^k,\quad |z_i|=|z_j|\quad \mbox{and} \quad |z_l|\leq
|z_i| \quad \mbox{for every } l \big\}.$$ 
We have
$f^n(z)=(z_1^{d^n},\ldots, z_k^{d^n})$ for $n\geq 1$. 
\end{example}

\bigskip\bigskip

\begin{exercise}
Let $h:\P^1\rightarrow \P^1$ be a rational map. Discuss Fatou components for the
associated map $f$ defined in Example \ref{example_ueda}. Prove
in particular that there exist Fatou components which are bi-holomorphic to a disc
cross an annulus. Describe the set of non-wandering points of $f$.
\end{exercise}

\begin{exercise}
Let $a$ be a fixed point of $f$. Show that the eigenvalues of the
differential $Df$ of $f$ at $a$ do not depend on the local
coordinates. Assume that $a$ is in the Fatou set. Show that these eigenvalues
are of modulus $\leq 1$. If all the eigenvalues are of modulus $1$, show
that $Df(a)$ is diagonalizable.
\end{exercise}

\begin{exercise}
Let $f$ be a Latt{\`e}s map associated to an affine map $A$ as in
Example \ref{example_ueda}. Show that $f$ is postcritically
finite. Show that $d^{-1/2} DA$ is an unitary matrix where $DA$ is the
differential of $A$. Deduce that 
 the orbit of $a$ is
dense in $\P^k$ for almost every $a$ in $\P^k$. Show that the periodic points of $f$ are dense in $\P^k$.
\end{exercise}

\begin{exercise} 
Let $f:\P^k\rightarrow\P^k$ be a dominant meromorphic map. Let $I$ be
the indeterminacy set of $f$, defined as above. Show that $f$ cannot be
extended to a holomorphic map on any open set which intersects $I$.
\end{exercise}


\section{Green currents and Julia sets} \label{section_green_pk}

Let $f$ be an endomorphism of algebraic degree $d\geq 2$ as above.
In this paragraph, we give the first construction of canonical invariant
currents $T^p$ associated to $f$ (Green currents).
The construction is now classical and is used
in most of the references, see \cite{FornaessSibony1,FornaessSibony3,
  HubbardPapadopol, Sibony}.
We will show that the support of the Green $(1,1)$-current is exactly
the Julia set of $f$ \cite{FornaessSibony1}. In some examples, Green currents describe the
distribution of stable varieties but in general their geometric
structure is not yet well-understood. We will see later that  
$\mu:=T^k$ is the invariant measure of maximal entropy.

\begin{theorem} \label{th_green_11_pk}
Let $S$ be a positive closed $(1,1)$-current of mass $1$ on
$\P^k$. Assume that $S$ has bounded local potentials. Then
$d^{-n}(f^n)^*(S)$ converge weakly to a positive closed
$(1,1)$-current $T$ of mass $1$. This current has continuous local
potentials and does not depend on
$S$. Moreover, it is totally invariant: $f^*(T)=dT$ and
$f_*(T)=d^{k-1}T$. We also have for a smooth $(k-1,k-1)$-form $\Phi$
$$\big|\langle d^{-n} (f^n)^*(S)-T,\Phi\rangle\big| \leq cd^{-n}
\|\Phi\|_\DSH,$$
where $c>0$ is a constant independent of $\Phi$ and of $n$.  
\end{theorem}
\proof
We refer to Appendix for the basic properties of
quasi-p.s.h. functions, positive closed currents and DSH currents.
Since $S$ has mass 1, it is cohomologous to $\omega_\FS$. Therefore, we can write
$S=\omega_\FS+\ddc u$ where $u$ is a quasi-p.s.h. function. By
hypothesis, this function is bounded.
The current $d^{-1}f^*(\omega_\FS)$ is smooth and of mass 1 since
$f^*:H^{1,1}(\P^k,\C)\rightarrow H^{1,1}(\P^k,\C)$ is the multiplication
by $d$ and the mass of a positive closed current can be computed cohomologically. So, we
can also write $d^{-1}f^*(\omega_\FS)=\omega_\FS+\ddc v$ where $v$ is a
 quasi-p.s.h. function. Here, $v$ is smooth since $\omega_\FS$ and
 $f^*(\omega_\FS)$ are smooth.
We have
\begin{eqnarray*}
d^{-1}f^*(S) & = & d^{-1}f^*(\omega_\FS)+\ddc\big(d^{-1}u\circ
f\big) \\
& = & \omega_\FS+\ddc v+\ddc\big(d^{-1}u\circ f\big).
\end{eqnarray*}
By induction, we obtain
$$d^{-n}(f^n)^*(S)=\omega_\FS+\ddc \big(v+\cdots+d^{-n+1}v\circ f^{n-1}\big)+
 \ddc\big(d^{-n}u\circ f\big).$$
Observe that, since $v$ is smooth, the sequence of smooth functions
$v+\cdots+d^{-n+1}v\circ f^{n-1}$ converges uniformly to a
continuous function $\gr$. 
Since $u$ is bounded,
the functions  $d^{-n}u\circ f$ tend to $0$. It
follows that $d^{-n} (f^n)^*(S)$ converge weakly to a current $T$
which satisfies 
$$T=\omega_\FS+\ddc \gr.$$
Clearly, this current does not depend on $S$ since $\gr$ does not depend
on $S$. Moreover, the currents $d^{-n}(f^n)^*(S)$ are positive closed of mass
1. So, $T$ is also a positive closed current of mass 1. We deduce that $\gr$ is
quasi-p.s.h. since it is continuous and satisfies $\ddc \gr\geq -\omega_\FS$.

Applying the above computation to $T$ instead of $S$, we obtain that
$$d^{-1}f^*(T)=\omega_\FS+\ddc v+\ddc\big(d^{-1}\gr\circ
f\big)=\omega_\FS+\ddc \gr.$$
Hence, $f^*(T)=dT$. On smooth forms $f_*\circ f^*$ is equal to $d^k$
times the identity; this holds by
continuity for positive closed currents. Therefore,
$$f_*(T)=f_*(f^*(d^{-1}T))=d^{k-1}T.$$
It remains to prove the estimate in the theorem. 
Recall that we can write
$\ddc\Phi=R^+-R^-$ where $R^\pm$ are positive measures
such that $\|R^\pm\|\leq \|\Phi\|_\DSH$. 
We have
\begin{eqnarray*}
\big|\langle d^{-n} (f^n)^*(S)-T,\Phi\rangle\big| & = &
\big|\langle \ddc (v+\cdots+d^{-n+1} v\circ f^{n-1}+ d^{-n} u\circ
f^n -\gr),\Phi\rangle \big| \\
& = & \big|\langle v+\cdots+d^{-n+1} v\circ f^{n-1}+ d^{-n} u\circ
f^n -\gr,\ddc \Phi\rangle \big| \\
& = & \big|\big\langle d^{-n} u\circ
f^n - \sum_{i\geq n} d^{-i} v\circ f^i, R^+-R^-\big\rangle \big|.
\end{eqnarray*}
Since $u$ and $v$ are bounded, the mass estimate for $R^\pm$ implies 
that the last integral is $\lesssim
d^{-n}\|\Phi\|_\DSH$. The result follows.
\endproof

Theorem \ref{th_green_11_pk} gives a convergence result for $S$ quite
diffuse (with bounded potentials). It is
like the first main theorem in value distribution theory. The question
that we will address is the convergence for singular $S$,
e.g. hypersurfaces. 

\begin{definition}\rm
We call $T$ {\it the Green $(1,1)$-current} and $\gr$ {\it the Green
  function} of $f$. The power $T^p:=T\wedge\ldots\wedge T$,
$p$ factors, is 
{\it the Green  $(p,p)$-current} of $f$, and its support $\Jc_p$
is called {\it the Julia set of order $p$}.
\end{definition}

Note that the Green function
is defined up to an additive constant and since $T$ has a continuous quasi-potential, $T^p$ is
well-defined. Green currents are totally invariant: we have $f^*(T^p)=d^pT^p$ and
$f_*(T^p)=d^{k-p}T^p$. 
The Green
$(k,k)$-current $\mu:=T^k$ is also called {\it the Green measure, 
the  equilibrium measure or the measure of maximal entropy}. 
We will give in the next paragraphs results which justify the terminologies.
The iterates $f^n$, $n\geq 1$, have the same Green currents and Green
function. We have the following result. 

\begin{proposition} \label{prop_green_holder}
The local potentials of the Green current $T$ are $\gamma$-H{\"o}lder continuous
for every  $\gamma$ such that $0<\gamma < \min(1,\log
  d/\log d_\infty)$, where $d_\infty:=\lim \|Df^n\|^{1/n}_\infty$. In
  particular, the Hausdorff dimension of $T^p$ is strictly larger than
  $2(k-p)$ and $T^p$ has no mass on pluripolar sets and on proper analytic
  sets of $\P^k$. 
\end{proposition}

Since $Df^{n+m}(x)=Df^m(f^n(x))\circ Df^n(x)$, it is not difficult to check that the sequence
$\|Df^n\|^{1/n}_\infty$ is decreasing. So,
$d_\infty=\inf\|Df^n\|^{1/n}_\infty$.
The last assertion of the proposition is deduced from Corollary \ref{cor_cln} and
Proposition \ref{prop_hol_dim_wedge} in Appendix.
The first assertion is equivalent to the H{\"o}lder continuity of the Green
function $\gr$, it was obtained by Sibony \cite{Sibony_O} for
one variable  polynomials and by
Briend \cite{Briend} and  Kosek \cite{Kosek} in higher
dimension.

The following lemma, due to Dinh-Sibony
\cite{DinhSibony4, DinhSibony10}, implies the above proposition and can be applied in a
more general setting. Here, we apply it to $\Lambda:=f^m$ with $m$
large enough and to the
above smooth function $v$. We choose $\alpha:=1$, $A:=\|Df^m\|_\infty$
and $d$ is replaced by $d^m$.

\begin{lemma} \label{lemma_holder}
Let $K$ be a metric space with finite diameter and  $\Lambda:K\rightarrow K$ be a
Lipschitz map: $\|\Lambda(a)-\Lambda(b)\|\leq A\|a-b\|$ with $A>0$. Here, $\|a-b\|$ denotes the distance between
two points $a$, $b$ in $K$. 
Let $v$ be an $\alpha$-H{\"o}lder continuous function on $K$ with
$0<\alpha\leq 1$. Then, 
$\sum_{n\geq 0} d^{-n}v\circ\Lambda^n$ converges pointwise to a
function which is $\beta$-H{\"o}lder continuous on $K$ for every $\beta$
such that $0<\beta<\alpha$ and $\beta\leq \log d/\log A$. 
\end{lemma}
\proof
By hypothesis, there is a constant  $A'>0$ such that
 $|v(a)-v(b)|\leq A'\|a-b\|^\alpha$.
Define $A'':=\|v\|_\infty$. Since $K$ has finite diameter, $A''$ is
finite and we only have to consider the case where $\|a-b\|\ll 1$.
If $N$ is an integer, we have
\begin{eqnarray*}
\lefteqn{\big|\sum_{n\geq 0} d^{-n}v\circ\Lambda^n(a)-\sum_{n\geq 0}
d^{-n}v\circ\Lambda^n(b)\big|}\\
& \leq & \sum_{0\leq n\leq N} d^{-n}|v\circ\Lambda^n(a)-
v\circ\Lambda^n(b)| + \sum_{n>N} d^{-n}|v\circ\Lambda^n(a)-
v\circ\Lambda^n(b)|\\
& \leq &  A' \sum_{0\leq n \leq N} d^{-n}\|\Lambda^n(a)-\Lambda^n(b)\|^\alpha
+ 2A''\sum_{n>N} d^{-n} \\
&\lesssim & \|a-b\|^\alpha\sum_{0\leq n \leq N} d^{-n}A^{n\alpha} + d^{-N}.
\end{eqnarray*}
If $A^\alpha\leq d$, the last sum is of order at most equal to
$N\|a-b\|^\alpha+d^{-N}$. For a given $0<\beta<\alpha$, choose
$N\simeq -\beta\log\|a-b\|/\log d$. So, the last expression is
$\lesssim \|a-b\|^\beta$. In this case, the function is
$\beta$-H{\"o}lder continuous for every $0<\beta<\alpha$.
When $A^\alpha>d$, the sum is $\lesssim
d^{-N}A^{N\alpha}\|a-b\|^\alpha+d^{-N}$.
For $N\simeq -\log\|a-b\|/\log A$, the last expression is $\lesssim
\|a-b\|^\beta$ with $\beta:=\log d/\log A$. Therefore, the
function is $\beta$-H{\"o}lder continuous. 
\endproof

\begin{remark}\rm
Lemma  \ref{lemma_holder} still holds for $K$ with infinite diameter
if $v$ is H{\"o}lder continuous and 
bounded. We can also replace the distance on $K$ by any positive
symmetric function on $K\times K$ which vanishes on the diagonal.
Consider a family $(f_s)$ of endomorphisms of $\P^k$
depending holomorphically on $s$ in a space of parameters $\Sigma$. In the
above construction of the Green current, we can locally on $\Sigma$, choose $v_s(z)$ smooth
such that $dd^c_{s,z} v_s(z)\geq -\omega_\FS(z)$.
Lemma \ref{lemma_holder} implies that the Green function $\gr_s(z)$ of $f_s$ 
is locally H{\"o}lder continuous on $(s,z)$ in $\Sigma\times \P^k$.
Then, $\omega_\FS(z)+dd^c_{s,z} \gr_s(z)$ is a positive closed $(1,1)$-current on
$\Sigma\times \P^k$. 
Its slices by  $\{s\}\times\P^k$ are 
the Green currents $T_s$ of $f_s$.
\end{remark}

We want to use the properties of the Green currents in order to
establish some properties of the Fatou and Julia sets. We will show
that the Julia set coincides with the Julia set of order 1. 
We recall the notion of Kobayashi hyperbolicity on a complex
manifold $M$. Let $p$ be a point in $M$ and $\xi$ a tangent vector of
$M$ at $p$. Consider the holomorphic maps $\tau:\Delta\rightarrow M$ on the
unit disc $\Delta$ in $\C$ such that $\tau(0)=p$ and $D\tau(0)=c\xi$ where
$D\tau$ is the differential of $\tau$ and $c$ is a constant. {\it The
  Kobayashi-Royden pseudo-metric} is defined by
$$K_M(p,\xi):=\inf_\tau |c|^{-1}.$$
It measures the size of a disc that one can immerse in $M$. In
particular, if $M$ contains an image of $\C$ passing through $p$ in
the direction $\xi$, we have $K_M(p,\xi)=0$.  

Kobayashi-Royden pseudo-metric is contracting for holomorphic maps: if
$\Psi:N\rightarrow M$ is a holomorphic map between complex manifolds, 
we have 
$$K_M(\Psi(p),D\Psi(p)\cdot\xi)\leq K_N(p,\xi).$$
The Kobayashi-Royden pseudo-metric on $\Delta$ coincides
with the Poincar\'e metric.
A complex manifold $M$ is {\it Kobayashi hyperbolic} if $K_M$ is a
metric \cite{Kobayashi}. In which case, holomorphic self-maps of $M$, form a locally 
equicontinuous family of maps. We have the following result where
the norm of $\xi$ is with respect to a smooth metric on $X$.

\begin{proposition} \label{prop_mfd_fct_psh}
Let  $M$ be a relatively compact open set of a compact complex manifold
$X$. Assume that there is a bounded function $\rho$ on $M$ which
is strictly p.s.h., i.e. $\ddc\rho\geq\omega$ on $M$ for some positive
Hermitian form
$\omega$ on $X$.
Then $M$ is Kobayashi hyperbolic and hyperbolically embedded in
$X$. More precisely, there is a constant $\lambda>0$ such that  
$K_M(p,\xi)\geq \lambda \|\xi\|$ for every $p\in M$ and every tangent
vector $\xi$ of $M$ at $p$.
\end{proposition}
\proof
If not, we can find holomorphic discs $\tau_n:\Delta\rightarrow M$ such that 
$\|D\tau_n(0)\|\geq n$ for $n\geq 1$. So, this family is not equicontinuous. A
lemma due to Brody \cite{Kobayashi} says that, after reparametrization, there is a subsequence
converging to an image of $\C$ in $\overline M$. More precisely, up to
extracting a subsequence, there
are holomorphic maps $\Psi_n:\Delta_n\rightarrow\Delta$ on discs
$\Delta_n$ centered at 0, of radius $n$,
such that $\tau_n\circ\Psi_n$ converge locally
uniformly to a non-constant map $\tau_\infty:\C\rightarrow \overline
M$. Since $\rho$ is bounded, up to extracting a subsequence, 
the subharmonic functions $\rho_n:=\rho\circ \tau_n\circ\Psi_n$
converge in $L^1_\loc(\C)$ to some subharmonic function 
$\rho_\infty$. Since the function $\rho_\infty$ is bounded, it should be constant. 

For simplicity, we use here the metric on $X$ induced by $\omega$.
Let $L,K$ be arbitrary 
compact subsets of $\C$ such that $L\Subset K$. For $n$ large enough, the area of
$\tau_n(\Psi_n(L))$ counted with multiplicity, satisfies
$$\area(\tau_n(\Psi_n(L)))=\int_L (\tau_n\circ \Psi_n)^*(\omega)\leq
\int_L \ddc\rho_n.$$
We deduce that 
$$\area(\tau_\infty(L))=\lim_{n\rightarrow\infty} \area(\tau_n(\Psi_n(L)))
\leq \int_K \ddc\rho_\infty=0.$$
This is a contradiction.
\endproof

The following result was obtained by Forn\ae ss-Sibony in
\cite{FornaessSibony3, FornaessSibony2} and by Ueda for the assertion on the
Kobayashi hyperbolicity of the Fatou set \cite{Ueda2}.

\begin{theorem}
Let $f$ be an endomorphism of algebraic degree $d\geq 2$ of $\P^k$. Then,
the Julia set of order $1$ of $f$, i.e. the support $\Jc_1$ of the Green
$(1,1)$-current $T$, coincides with the Julia
set $\Jc$. The Fatou set $\Fc$ is Kobayashi hyperbolic and hyperbolically embedded in $\P^k$. 
Moreover, for $p\leq k/2$, the Julia set of order $p$ of
$f$ is connected.
\end{theorem}
\proof
The sequence $(f^n)$ is equicontinuous on the Fatou set $\Fc$ and $f^n$ are
holomorphic, hence the differential $Df^n$ are locally uniformly bounded on
$\Fc$. Therefore, $(f^n)^*(\omega_\FS)$ are locally uniformly bounded on
$\Fc$. We deduce that $d^{-n}(f^n)^*(\omega_\FS)$ converge to 0 on
$U$. Hence, $T$ is supported on the Julia set $\Jc$.

Let $\Fc'$ denote the complement of the support of $T$ in
$\P^k$. Observe that $\Fc'$ is invariant under $f^n$ and 
that $-\gr$ is a smooth function which is strictly p.s.h. on
$\Fc'$. 
Therefore, by Proposition \ref{prop_mfd_fct_psh}, $\Fc'$ is Kobayashi
hyperbolic and  hyperbolically embedded in $\P^k$. Therefore, the maps
$f^n$, which are self-maps of $\Fc'$, are
equicontinuous with respect to the Kobayashi-Royden metric. On the other
hand, the fact that $\Fc'$ is hyperbolically embedded implies that
the Kobayashi-Royden metric is bounded from below by a constant times
the Fubini-Study metric.
It follows that $(f^n)$ is locally
equicontinuous on $\Fc'$ with respect to the Fubini-Study metric. We
conclude that $\Fc'\subset \Fc$, hence $\Fc=\Fc'$ and $\Jc=\supp(T)=\Jc_1$.

In order to show that $\Jc_p$ are connected, it is enough to prove that if
$S$ is a positive closed current of bidegree $(p,p)$ with $p\leq k/2$
then the support of $S$ is connected. Assume that the support of $S$
is not connected, then we can write $S=S_1+S_2$ with $S_1$ and $S_2$
non-zero, positive closed
with disjoint supports. Using a convolution on the automorphism group
of $\P^k$, we can construct smooth positive closed
$(p,p)$-forms $S_1',S_2'$ with disjoint supports. So, we have $S_1'\wedge
S_2'=0$. This contradicts that the cup-product of the classes $[S_1']$ and
$[S_2']$ is non-zero in $H^{2p,2p}(\P^k,\R)\simeq\R$: we have
$[S_1']=\|S_1'\|[\omega_\FS^p]$, $[S_2']=\|S_2'\|[\omega_\FS^p]$ and 
$[S_1']\smile [S_2'] = \|S_1'\|\|S_2'\|[\omega_\FS^{2p}]$, a contradiction. 
Therefore, the support of
$S$ is connected.
\endproof

\begin{example}\rm
Let $f$ be a polynomial map of algebraic degree $d\geq 2$ on $\C^k$
which extends holomorphically to $\P^k$. If $B$ is a ball large enough
centered at 0, then $f^{-1}(B)\Subset B$.
Define $G_n:=d^{-n} \log^+\|f^n\|$, where $\log^+:=\max(\log,0)$. As
in Theorem \ref{th_green_11_pk}, we
can show that $G_n$ converge uniformly to a continuous p.s.h. function
$G$ such that $G\circ f=dG$. On $\C^k$, the Green current $T$ of $f$ is equal to $\ddc G$ and
$T^p=(\ddc G)^p$. The Green measure is equal to $(\ddc
G)^k$. If $\Kc$ denotes the set of points in $\C^k$ with bounded orbit,
then $\mu$ is supported on $\Kc$. Indeed, outside $\Kc$ we have $G=\lim
d^{-n} \log\|f^n\|$ and the convergence is locally uniform. It
follows that $(\ddc G)^k=\lim d^{-kn} (\ddc \log\|f^n\|)^k$ on
$\C^k\setminus \Kc$. One 
easily check that $(\ddc \log\|f^n\|)^k=0$ out of $f^{-n}(0)$. Therefore,
$(\ddc G)^k=0$ on $\C^k\setminus \Kc$. The set $\Kc$ is called {\it
  the filled Julia set}. We can show that $\Kc$ is the zero set of
$G$.
In particular, if $f(z)=(z_1^d,\ldots,z_k^d)$, then $G(z)=\sup_i
\log^+|z_i|$. One can check that the support of $T^p$ is foliated
(except for a set of zero measure with respect to the trace of $T^p$)
by stable manifolds of dimension $k-p$ and that $\mu=T^k$ is the
Lebesgue measure on the torus $\{|z_i|=1,\ i=1,\ldots,k\}$.  
\end{example}

\begin{example}\rm
We consider Example \ref{example_ueda}. Let $\nu$ be the Green measure of $h$ on
$\P^1$, i.e. $\nu=\lim d^{-n} (h^n)^*(\omega_\FS)$. Here, $\omega_\FS$
denotes also the Fubini-Study form on $\P^1$. Let $\pi_i$ denote the projections of
$\P^1\times\cdots\times\P^1$ on the factors. Then, the Green current
of $f$ is equal to
$$T={1\over k!}\pi_*\big(\pi_1^*(\nu)+\cdots+\pi_k^*(\nu)\big),$$ 
as can be easily checked.
\end{example}

\begin{example} \rm
The following family of maps on $\P^2$ was studied in \cite{FornaessSibony6}:
$$f[z_0:z_1:z_2]:=[z_0^d+\lambda z_0z_1^{d-1}, \nu
(z_1-2z_2)^d+cz_0^d:z_1^d+cz_0^d].$$
For appropriate choices of the parameters $c$ and $\lambda$, one can
show that $\supp(T)$ and $\supp(\mu)$ coincide and have non-empty
interior. Moreover, $f$ has an attracting fixed point, so the Fatou
set is not empty. The situation is then quite different from the one
variable case, where either the Julia set is equal to $\P^1$ or it has
empty interior. Observe that the restriction of $f$ to the projective
line $\{z_0=0\}$, for appropriate $\nu$, is chaotic, i.e. has dense
orbits. One shows that $\{z_0=0\}$ is in the support of $\mu$ and that
$\P^2\setminus \cup_{i=0}^m f^{-i}\{z_0=0\}$ is Kobayashi
hyperbolic. Hence using the total invariance of $\supp(\mu)$, we get
that the complement of $\supp(\mu)$ is in the Fatou set. It is
possible to choose the parameters so that $\P^2\setminus\supp(\mu)$
contains an attractive fixed point. Several other examples are
discussed in \cite{FornaessSibony6}. 
\end{example}

We now give a characterization of the Julia sets in term of volume
growth. There is an interesting gap in the possible volume growth.

\begin{proposition}
Let $f$ be a holomorphic endomorphism of algebraic degree $d\geq 2$ of
$\P^k$. Let $T$ be its Green $(1,1)$-current. Then the following
properties are equivalent:
\begin{enumerate}
\item $x$ is a point in the Julia set of order $p$, i.e. $x\in
      \Jc_p:=\supp(T^p)$;
\item For every neighbourhood $U$ of $x$, we have 
$$\liminf_{n\rightarrow\infty} d^{-pn}\int_U (f^n)^*(\omega_\FS^p)\wedge
  \omega_\FS^{k-p} \not =0;$$
\item For every neighbourhood $U$ of $x$, we have 
$$\limsup_{n\rightarrow\infty} d^{-(p-1)n}\int_U (f^n)^*(\omega_\FS^p)\wedge
  \omega_\FS^{k-p}  =+\infty.$$
\end{enumerate}
\end{proposition}
\proof
We have seen in Theorem \ref{th_green_11_pk} that $d^{-n} (f^n)^*(\omega_\FS)$ converges
to $T$ when $n$ goes to infinity. Moreover, $d^{-n} (f^n)^*(\omega_\FS)$
admits a quasi-potential which converges uniformly to a quasi-potential of
$T$. It follows that $\lim d^{-pn}(f^n)^*(\omega_\FS^p)=T^p$. We deduce
that Properties 1) and 2) are equivalent. Since 2) implies 3), it
remains to show that 3) implies 1). For this purpose, it is enough to
show that  for any 
open set $V$ with $\overline V\cap \Jc_p=\varnothing$, 
$$\int_V(f^n)^*(\omega_\FS^p)\wedge
  \omega_\FS^{k-p}= O(d^{(p-1)n}).$$
This is a consequence of a more general inequality in Theorem \ref{th_entropy} below. 
We give here a direct proof. 

Since $(\omega_\FS+\ddc\gr)^p=0$ on
$\P^k\setminus\Jc_p$, we can write there $\omega_\FS^p=\ddc \gr\wedge (S^+-S^-)$
where $S^\pm$ are positive closed $(p-1,p-1)$-currents on $\P^k$. Let $\chi$ be
a cut-off function with compact support in $\P^k\setminus\Jc_p$ and equal to 1 on
$V$. The above integral is bounded by
$$\int_{\P^k} \chi (f^n)^*\big(\ddc \gr\wedge (S^+-S^-)\big)
\wedge \omega_\FS^{k-p}=\int_{\P^k} \ddc\chi \wedge (\gr\circ f^n)(f^n)^*(S^+-S^-)\wedge \omega_\FS^{k-p}.$$
Since $\gr$ is bounded, the last integral is bounded by a constant
times $\|(f^n)^*(S^+)\|+\|(f^n)^*(S^-)\|$. We conclude using the
identity $\|(f^n)^*(S^\pm)\|=d^{(p-1)n}\|S^\pm\|$.
\endproof

The previous proposition suggests a notion of {\it local dynamical
degree}. Define
$$\delta_p(x,r):=\limsup_{n\rightarrow\infty}
\Big(\int_{B(x,r)}(f^n)^*(\omega_\FS^p)\wedge \omega_\FS^{k-p}\Big)^{1/n}$$
and
$$\delta_p(x):=\inf_{r>0} \delta_p(x,r)=\lim_{r\rightarrow 0}
\delta_p(x,r).$$
It follows from the above proposition that $\delta_p(x)=d^p$ for
$x\in\Jc_p$ and $\delta_p(x)=0$ for $x\not\in\Jc_p$.
This notion can be extended to general
meromorphic maps or correspondences and the sub-level sets
$\{\delta_p(x)\geq c\}$ can be seen
as a kind of Julia sets.

\bigskip\bigskip

\begin{exercise}
Let $f$ be an endomorphism of algebraic degree $d\geq 2$ of $\P^k$. Suppose a
subsequence $(f^{n_i})$ is equicontinuous on an open set $U$. Show
that $U$ is contained in the Fatou set.
\end{exercise}

\begin{exercise}
Let $f$ and $g$ be two commuting holomorphic endomorphisms of $\P^k$, 
i.e. $f\circ g=g\circ f$. Show that $f$ and $g$
have the same Green currents. Deduce that they have the same Julia and
Fatou sets. 
\end{exercise}

\begin{exercise} \label{exo_tp_laminar}
Determine the Green $(1,1)$-current and the Green measure 
for the map $f$ in Example \ref{example_power_map}. Study the
lamination on $\supp(T^p)\setminus\supp(T^{p+1})$. Express the current
$T^p$ on that set as an integral on appropriate manifolds. 
\end{exercise}

\begin{exercise}
Let $f$ be an endomorphism of algebraic degree $d\geq 2$ of
$\P^k$ and $T$ its Green $(1,1)$-current. 
Consider the family of maps $\tau:\Delta\rightarrow\P^k$ such that
$\tau^*(T)=0$. The last equation means that if $u$ is a local
potential of $T$, i.e. $\ddc u=T$ on some open set, then $\ddc
u\circ\tau=0$ on its domain of definition. Show that the sequence 
$(f^n_{|\tau(\Delta)})_{n\geq 1}$ is equicontinuous. Prove that
there is a constant $c>0$ such that $\|D\tau(0)\|\leq c$ for every
$\tau$ as above (this property holds for any positive closed
$(1,1)$-current $T$ with continuous potentials). Find the
corresponding discs for $f$ as in Exercise \ref{exo_tp_laminar}.
\end{exercise} 

\begin{exercise}
Let $f$ be an endomorphism of algebraic degree $d\geq 2$ of $\P^k$.
Let $X$ be an analytic set of pure dimension $p$ in an open set
$U\subset\P^k$. Show that for every compact $K\subset U$
$$\limsup_{n\rightarrow\infty} {1\over n} \log\vol(f^n(X\cap K))\leq
p\log d.$$
Hint. For an appropriate cut-off function $\chi$, estimate $\int_X\chi
(f^n)^*(\omega_\FS^p)$. 
\end{exercise}


\section{Other constructions of the Green currents}

In this paragraph, we give other methods, introduced and developped
by the authors, in order to construct the
Green currents and Green measures for meromorphic maps. We obtain estimates on the
Perron-Frobenius operator and on the thickness of the Green measure,
that will be applied 
in the stochastic study of the dynamical system. A key point
here is the use of d.s.h. functions as observables.

We  first present 
a recent direct construction of Green $(p,p)$-currents
using super-potentials\footnote{These super-potentials correspond to
  super-potentials of mean 0 in \cite{DinhSibony10}.}. 
Super-potentials are a tool in order to compute with positive closed
$(p,p)$-currents. They play the same role as potentials for bidegree
$(1,1)$ currents.
In dynamics, they are used in particular in the equidistribution
problem for algebraic sets of arbitrary dimension and allow to get
some estimates on the speed of convergence. 

\begin{theorem} \label{th_green_pp}
Let $S$ be a positive closed $(p,p)$-current of mass $1$ on
$\P^k$. Assume that the super-potential of $S$ is bounded. Then
$d^{-pn} (f^n)^*(S)$ converge to the Green $(p,p)$-current $T^p$ of
$f$. Moreover, $T^p$ has a H{\"o}lder continuous super-potential.
\end{theorem}
\noindent
{\bf Sketch of proof.}
We refer to Appendix \ref{section_positive} and \ref{section_current_pk} for an introduction to super-potentials 
and to the action of holomorphic maps on positive closed currents. Recall
that $f^*$ and $f_*$ act on $H^{p,p}(\P^k,\C)$ as the multiplications
by $d^p$ and $d^{k-p}$ respectively. So, if $S$ is a positive closed
$(p,p)$-current of mass 1, then $\|f^*(S)\|=d^p$ and
$\|f_*(S)\|=d^{k-p}$ since the mass can be computed cohomologically. 
Let $\Lambda$ denote the operator $d^{-p+1}f_*$ acting on
$\Cc_{k-p+1}(\P^k)$, the
convex set of positive closed currents of bidegree $(k-p+1,k-p+1)$ and of mass 1. 
It is continuous and it takes values also in $\Cc_{k-p+1}(\P^k)$. 
Let $\Vc$, $\Uc$, $\Uc_n$ denote the super-potentials of $d^{-p}f^*(\omega_\FS^p)$, 
$S$ and $d^{-pn}(f^n)^*(S)$ respectively. Consider a quasi-potential
$U$ of mean 0 of $S$ which is a DSH current satisfying $\ddc
U=S-\omega_\FS^p$. 
The following computations are valid for $S$ smooth and can be
extended to all currents $S$ using a suitable regularization
procedure.

By Theorem \ref{th_pullback_local} in the Appendix, the current $d^{-p}f^*(U)$ is DSH and satisfies 
$$\ddc \big(d^{-p}f^*(U) \big)=d^{-p} f^*(S) - d^{-p}f^*(\omega_\FS^p).$$
If $V$ is a smooth quasi-potential of mean 0 of
$d^{-p}f^*(\omega_\FS^p)$, i.e. a smooth real $(p-1,p-1)$-form such
that
$$\ddc V= d^{-p} f^*(\omega_\FS^p)-\omega_\FS^p\quad \mbox{and}\quad
\langle \omega_\FS^{k-p+1},V\rangle=0,$$
then $V+d^{-p}f^*(U)$ is a quasi-potential
of $d^{-p}f^*(S)$.  Let $m$ be the real number such that 
$V+d^{-p}f^*(U)+m\omega_\FS^{p-1}$ is
a quasi-potential of mean 0 of $d^{-1}f^*(S)$. 
We have
\begin{eqnarray*}
\Uc_1(R) & = & \langle V+d^{-p}f^*(U)+m\omega_\FS^{p-1},R\rangle \\
& = & \langle V,R\rangle+d^{-1} \langle U,\Lambda(R)\rangle  +m\\
& = & \Vc(R)+d^{-1} \Uc(\Lambda(R)) +m.
\end{eqnarray*}
By induction, we obtain
\begin{eqnarray*}
\Uc_n(R) & = &
\Vc(R)+d^{-1}\Vc(\Lambda(R))+\cdots+d^{-n+1}\Vc(\Lambda^{n-1}(R))\\
& & +d^{-n} \Uc(\Lambda^n(R))+m_n,
\end{eqnarray*}
where $m_n$ is a constant depending on $n$ and on $S$.

Since $d^{-p} f^*(\omega_\FS^p)$ is smooth, $\Vc$ is a H{\"o}lder continuous
function. It is not difficult to show that $\Lambda$ is Lipschitz with respect to the distance
$\dist_\alpha$ on $\Cc_{k-p+1}(\P^k)$.
Therefore, by Lemma \ref{lemma_holder}, the sum
$$\Vc(R)+d^{-1}\Vc(\Lambda(R))+\cdots+d^{-n+1}\Vc(\Lambda^{n-1}(R))$$
converges uniformly to a H{\"o}lder continuous function $\Vc_\infty$
which does not depend on $S$. 
Recall that super-potentials vanish at $\omega_\FS^{k-p+1}$, in
particular, $\Uc_n(\omega_\FS^{k-p+1})=0$. Since
$\Uc$ is bounded,
the above expression of $\Uc_n(R)$  for $R=\omega_\FS^{k-p+1}$ implies that $m_n$ converge to
$m_\infty:=-\Vc_\infty(\omega_\FS^{k-p+1})$ which is independent of
$S$. So, $\Uc_n$ converge uniformly to $\Vc_\infty+m_\infty$. We
deduce that $d^{-pn}(f^n)^*(S)$ converge to a current $T_p$ which does
not depend on $S$. Moreover, the super-potential of $T_p$ is 
the H{\"o}lder continuous function $\Vc_\infty+m_\infty$.

We deduce from the above discussion that
$d^{-pn}(f^n)^*(\omega_\FS^p)$ converge in the Hartogs' sense to
$T_p$, see Appendix \ref{section_current_pk}. Theorem \ref{th_wedge_product_pp} implies that
$T_{p+q}=T_p\wedge T_q$ if $p+q\leq k$. Therefore,
if $T$ is the Green $(1,1)$-current, $T_p$ is equal to $T^p$ the Green
$(p,p)$-current of $f$. \hfill $\square$

\begin{remark}\rm
Let $S_n$ be positive closed $(p,p)$-currents of mass 1 on
$\P^k$. Assume that their super-potentials $\Uc_{S_n}$ satisfy
$\|\Uc_{S_n}\|_\infty =o(d^n)$. Then $d^{-pn} (f^n)^*(S_n)$ converge
to $T^p$. If $(f_s)$ is a family of maps depending holomorphically on
$s$ in a space of parameters $\Sigma$, then the Green super-functions are also
locally H{\"o}lder continuous with respect to $s$ and define a
positive closed $(p,p)$-current on $\Sigma\times \P^k$. Its slice by
$\{s\}\times\P^k$ is the Green $(p,p)$-current of $f_s$.   
\end{remark}

We now introduce two other constructions of the Green measure.
The main point is the use of appropriate spaces of
test functions adapted to complex analysis. Their norms take into account
the complex structure of $\P^k$. The reason to
introduce these spaces  is that they are invariant under the
push-forward by a holomorphic map. This is not the case for spaces of
smooth forms because of the critical set.
Moreover, we will see that there is a spectral gap for the action of
endomorphisms of $\P^k$ which is a useful property in the stochastic
study of the dynamical system.
The first method, called {\it the $\ddc$-method}, was introduced in
\cite{DinhSibony1} and developped in \cite{DinhSibony6}.
It can be extended to Green currents of any
bidegree. We show a convergence result for PB measures $\nu$ towards
the Green measure. PB measures are diffuse in some sense; we will
study equidistribution of Dirac masses in the next paragraph.

Recall that $f$ is an endomorphism of $\P^k$ of algebraic degree
$d\geq 2$. Define the {\it Perron-Frobenius operator}
$\Lambda$ on test functions $\varphi$ by
$\Lambda(\varphi):=d^{-k}f_*(\varphi)$. More precisely, we have
$$\Lambda(\varphi)(z):=d^{-k}\sum_{w\in f^{-1}(z)} \varphi(w),$$
where the points in $f^{-1}(z)$ are counted with multiplicity.
The following proposition is crucial. 

\begin{proposition} \label{prop_Perron_dsh}
The operator $\Lambda:\DSH(\P^k)\rightarrow\DSH(\P^k)$ is
well-defined, bounded and continuous with respect to the weak topology
on $\DSH(\P^k)$. 
The operator
$\widetilde\Lambda:\DSH(\P^k)\rightarrow\DSH(\P^k)$ defined by
$$\widetilde\Lambda(\varphi):=\Lambda(\varphi)-\langle\omega_\FS^k,\Lambda(\varphi)\rangle$$
is contracting and satisfies the estimate
$$\|\widetilde\Lambda(\varphi)\|_\DSH\leq d^{-1}\|\varphi\|_\DSH.$$
\end{proposition}
\proof
We prove the first assertion. Let $\varphi$ be a quasi-p.s.h. function
such that $\ddc\varphi\geq -\omega_\FS$. We show that
$\Lambda(\varphi)$ is d.s.h. Since $\varphi$ is strongly upper semi-continuous,
$\Lambda(\varphi)$ is strongly upper semi-continuous, see Appendix
\ref{section_positive}. 
If $\ddc\varphi=S-\omega_\FS$ with
$S$ positive closed, we have $\ddc \Lambda(\varphi)=d^{-k}
f_*(S)-d^{-k}f_*(\omega_\FS)$. Therefore, if $u$ is a quasi-potential
of $d^{-k}f_*(\omega_\FS)$, then $u+\Lambda(\varphi)$ is strongly
semi-continuous and is a quasi-potential of $d^{-k} f_*(S)$. So, this
function is 
quasi-p.s.h. We deduce that $\Lambda(\varphi)$ is d.s.h., and hence $\Lambda:\DSH(\P^k)\rightarrow\DSH(\P^k)$ is
well-defined.

Observe that $\Lambda:L^1(\P^k)\rightarrow L^1(\P^k)$ is
continuous. Indeed, if $\varphi$ is in $L^1(\P^k)$, we have 
$$\|\Lambda(\varphi)\|_{L^1} = \langle \omega^k_\FS, d^{-k}|
f_*(\varphi)|\rangle \leq  \langle \omega^k_\FS, d^{-k}
f_*(|\varphi|)\rangle = d^{-k}\langle f^*(\omega_\FS^k),
|\varphi|\rangle\lesssim \|\varphi\|_{L^1}.$$
Therefore,  $\Lambda:\DSH(\P^k)\rightarrow\DSH(\P^k)$ is continuous
with respect to the weak topology. This and the estimates below imply
that $\Lambda$ is a bounded operator.

We prove now the last estimate in the proposition. Write $\ddc\varphi=S^+-S^-$ with $S^\pm$
positive closed. We have
$$\ddc\widetilde\Lambda(\varphi)=
\ddc\Lambda(\varphi)=d^{-k}f_*(S^+-S^-)=
d^{-k}f_*(S^+)-d^{-k}f_*(S^-).$$
Since $\|f_*(S^\pm)\|= d^{k-1}\|S^\pm\|$ and 
$\langle\omega_\FS^k,\widetilde\Lambda(\varphi)\rangle=0$, we obtain that 
$\|\widetilde\Lambda(\varphi)\|_\DSH\leq d^{-1}\|S^\pm\|$.
The result follows.
\endproof

Recall that if $\nu$ is a positive measure on $\P^k$, the pull-back $f^*(\nu)$ is defined by the formula $\langle
f^*(\nu),\varphi\rangle =\langle \nu, f_*(\varphi)\rangle$ for
$\varphi$ continuous on $\P^k$.  
Observe that since $f$ is a ramified
covering, $f_*(\varphi)$ is continuous when $\varphi$ is
continuous, see Exercise \ref{Exo_push_c0} in Appendix. So, the above definition makes sense. For $\varphi=1$, we
obtain that $\|f^*(\nu)\|=d^k\|\nu\|$, since the covering is of degree
$d^k$. If $\nu$ is the Dirac mass at a point $a$, $f^*(\nu)$ is the sum of
Dirac masses at the points in $f^{-1}(a)$.

Recall that a measure $\nu$ is PB if quasi-p.s.h. are
$\nu$-integrable and $\nu$ is PC if it is PB and acts
continuously on $\DSH(\P^k)$ with respect to the weak topology on this
space, see Appendix \ref{section_current_pk}.
We deduce from Proposition \ref{prop_Perron_dsh} 
the following result where 
the norm $\|\cdot\|_\mu$ on $\DSH(\P^k)$ is defined by 
$$\|\varphi\|_\mu:=|\langle \mu,\varphi\rangle|+\inf\|S^\pm\|,$$
with $S^\pm$ positive closed such
that $\ddc\varphi=S^+-S^-$. We will see that $\mu$ is PB, hence this norm is equivalent to
$\|\cdot\|_\DSH$, see Proposition \ref{prop_dsh_norm_eq}.

\begin{theorem} \label{th_green_measure_dsh}
Let $f$ be as above. If $\nu$ is a PB probability
  measure, then $d^{-kn} (f^n)^*(\nu)$ converge to a PC probability 
measure $\mu$ which is independent of $\nu$ and totally invariant
under $f$. Moreover, if $\varphi$
is a d.s.h. function and $c_\varphi:=\langle \mu,\varphi\rangle$, then
$$\|\Lambda^n(\varphi)-c_\varphi\|_\mu\leq d^{-n} \|\varphi\|_\mu
\quad \mbox{and}\quad \|\Lambda^n(\varphi)-c_\varphi\|_\DSH\leq
Ad^{-n} \|\varphi\|_\DSH,$$
where $A>0$ is a constant independent of $\varphi$ and $n$. In
particular, there is a constant $c>0$ depending on $\nu$ such that
$$|\langle d^{-kn}(f^n)^*(\nu)-\mu,\varphi\rangle|\leq cd^{-n} \|\varphi\|_\DSH.$$
\end{theorem}
\proof
Since $\nu$ is PB, d.s.h. functions are 
$\nu$ integrable. It follows that
there is a constant $c>0$ such that $|\langle\nu,\varphi\rangle|\leq
c\|\varphi\|_\DSH$. Otherwise, there are d.s.h. functions $\varphi_n$
with $\|\varphi_n\|_\DSH\leq 1$ and $\langle \nu,\varphi_n\rangle \geq
2^n$, hence the d.s.h. function $\sum 2^{-n}\varphi_n$
is not $\nu$-integrable.

It follows from
Proposition \ref{prop_Perron_dsh} that $f^*(\nu)$ is PB. So, 
$d^{-kn} (f^n)^*(\nu)$ is PB for every $n$. Define for $\varphi$ d.s.h.,
$$c_0:=\langle \omega_\FS^k,\varphi\rangle \quad \mbox{and} \quad
\varphi_0:=\varphi-c_0$$ 
and inductively
$$c_{n+1}:=\langle \omega_\FS^k,\Lambda(\varphi_n)\rangle
\quad \mbox{and} \quad
\varphi_{n+1}:=\Lambda(\varphi_n)-c_{n+1}=\widetilde\Lambda(\varphi_n).$$
A straighforward computation gives
$$\Lambda^n(\varphi)=c_0+\cdots+c_n+\varphi_n.$$
Therefore,
$$\langle d^{-kn} (f^n)^*(\nu),\varphi\rangle= \langle
\nu,\Lambda^n(\varphi)\rangle= c_0+\cdots+c_n
+\langle \nu,\varphi_n\rangle.$$
Proposition \ref{prop_Perron_dsh} applied inductively on $n$ implies that
$\|\varphi_n\|_\DSH\leq d^{-n}\|\varphi\|_\DSH$. Since $\Lambda$ is bounded,
it follows that $|c_n|\leq Ad^{-n}\|\varphi\|_\DSH$, where $A>0$ is a constant.
The property that $\nu$ is PB and the above estimate on $\varphi_n$ imply that
$\langle\nu,\varphi_n\rangle$ converge to 0.

We deduce that $\langle  d^{-kn} (f^n)^*(\nu),\varphi\rangle$ converge
to $c_\varphi:=\sum_{n\geq 0} c_n$  and $|c_\varphi|\lesssim
\|\varphi\|_\DSH$. Therefore, $d^{-kn} (f^n)^*(\nu)$ converge
to a PB measure $\mu$ defined by $\langle
\mu,\varphi\rangle:=c_\varphi$.
The constant $c_\varphi$ does not depend on $\nu$, hence the measure $\mu$ is
independent of $\nu$. The above convergence implies that $\mu$ is
totally invariant, i.e. $f^*(\mu)=d^k\mu$.
Finally, since $c_n$ depends continuously on the d.s.h. function
$\varphi$, the constant $c_\varphi$, which is defined by a normally convergent
series, 
depends also continuously on $\varphi$. It
follows that $\mu$ is PC.

We prove now the estimates in the theorem. The total invariance of
$\mu$ implies that $\langle
\mu,\Lambda^n(\varphi)\rangle = \langle \mu,\varphi\rangle=c_\varphi$.
If $\ddc\varphi=S^+-S^-$ with $S^\pm$ positive closed, we have
$\ddc\Lambda^n(\varphi)=d^{-kn}(f^n)_*(S^+)- d^{-kn}(f^n)_*(S^-)$, hence 
$$\|\Lambda^n(\varphi)-c_\varphi\|_\mu\leq d^{-kn}\|(f^n)_*(S^\pm)\|=d^{-n}
\|S^\pm\|.$$
It follows that 
$$\|\Lambda^n(\varphi)-c_\varphi\|_\mu\leq d^{-n}\|\varphi\|_\mu.$$
For the second estimate,
 we have
$$\|\Lambda^n(\varphi)-c_\varphi\|_\DSH=\|\varphi_n\|_\DSH+\sum_{i\geq
  n} c_i.$$
The last sum is clearly bounded by a constant times
$d^{-n}\|\varphi\|_\DSH$. This together with the inequality
$\|\varphi_n\|_\DSH\lesssim d^{-n}\|\varphi\|_\DSH$ implies  
$\|\Lambda^n(\varphi)-c_\varphi\|_\DSH\lesssim
d^{-n}\|\varphi\|_\DSH$. We can also use that $\|\ \|_\mu$ and $\|\
\|_\DSH$ are equivalent.

The last inequality in
the theorem is then deduced from the identity
$$\langle d^{-kn}(f^n)^*(\nu)-\mu,\varphi\rangle = \langle
\nu,\Lambda^n(\varphi)-c_\varphi\rangle$$
and the fact that $\nu$ is PB. 
\endproof

\begin{remark} \rm
In the present case, the $\ddc$-method is quite simple. The function
$\varphi_n$ is the normalized solution of the equation $\ddc\psi=\ddc
\Lambda^n(\varphi)$. It satisfies automatically good estimates. The
other solutions differ from $\varphi_n$ by constants. 
We will see that for polynomial-like maps, the solutions differ by
pluriharmonic functions and the estimates are less straightforward.
In the construction of Green $(p,p)$-currents with $p>1$, $\varphi$ is
replaced by a test form of bidegree $(k-p,k-p)$ and $\varphi_n$ is a solution of
an appropriate $\ddc$-equation. The constants $c_n$ will be
replaced by $\ddc$-closed currents with a control of their cohomology
class.  
\end{remark}

The second construction of the Green measure follows the same
lines but we use the complex Sobolev space $W^*(\P^k)$ instead of
$\DSH(\P^k)$. We obtain that the Green measure $\mu$ is WPB, see
Appendix \ref{section_current_pk} for the terminology.
We only mention here the result which
replaces Proposition \ref{prop_Perron_dsh}.

\begin{proposition}
The operator $\Lambda:W^*(\P^k)\rightarrow W^*(\P^k)$ is
well-defined, bounded and continuous with respect to the weak topology
on $W^*(\P^k)$. 
The operator
$\widetilde\Lambda:W^*(\P^k)\rightarrow W^*(\P^k)$ defined by
$$\widetilde\Lambda(\varphi):=\Lambda(\varphi)-\langle\omega_\FS^k,\Lambda(\varphi)\rangle$$
is contracting and satisfies the estimate
$$\|\widetilde\Lambda(\varphi)\|_{W^*}\leq d^{-1/2}\|\varphi\|_{W^*}.$$
\end{proposition}
\noindent
{\bf Sketch of proof.}
As in Proposition \ref{prop_Perron_dsh}, since $\varphi$ is in $L^1(\P^k)$, 
$\Lambda(\varphi)$ is also in $L^1(\P^k)$ and the main point here is to estimate
$\partial \varphi$. Let $S$ be a positive closed $(1,1)$-current on $\P^k$
such that $\sqrt{-1}\partial \varphi\wedge \dbar \varphi\leq S$. We show
that $\sqrt{-1}\partial f_*(\varphi)\wedge \dbar f_*(\varphi)\leq d^k f_*(S)$,
in particular,  the Poincar\'e differential $d\Lambda(\varphi)$ of
$\Lambda(\varphi)$ is in $L^2(\P^k)$.

If $a$ is not a critical value of $f$ and $U$ a small neighbourhood of
$a$, then $f^{-1}(U)$ is the union of $d^k$ open sets $U_i$ which are
sent bi-holomorphically on $U$. Let $g_i:U\rightarrow U_i$ be the inverse
branches of $f$. On $U$, we obtain using Schwarz's inequality that
\begin{eqnarray*}
\sqrt{-1}\partial f_*(\varphi)\wedge \dbar f_*(\varphi) & = &  
\sqrt{-1}\Big(\sum \partial g_i^*(\varphi)\Big)\wedge \Big(\sum 
\dbar g_i^*(\varphi)\Big) \\
& \leq & d^k \sum \sqrt{-1}  \partial g_i^*(\varphi)\wedge \dbar
g_i^*(\varphi) \\
& = & d^k f_*\big(\sqrt{-1}\partial\varphi\wedge \dbar \varphi\big).
\end{eqnarray*}
Therefore, we have $\sqrt{-1}\partial f_*(\varphi)\wedge \dbar
f_*(\varphi)\leq d^k f_*(S)$ out of the critical values of $f$ which
is a manifold of real codimension 2. 

Recall that $f_*(\varphi)$ is in
$L^1(\P^k)$. It is a classical result in Sobolev spaces theory
 that an $L^1$ function whose gradient out of a
 submanifold of codimension 2 is in $L^2$, is in fact in the Sobolev space
$W^{1,2}(\P^k)$. We deduce that the inequality $\sqrt{-1}\partial f_*(\varphi)\wedge \dbar
f_*(\varphi)\leq d^k f_*(S)$ holds on $\P^k$, because the left hand side
term is an $L^1$ form and has no mass on critical values of $f$.
Finally, we have 
$$\sqrt{-1}\partial\Lambda(\varphi)\wedge
\dbar\Lambda(\varphi)\leq d^{-k}f_*(S).$$ 
This, together with the
identity $\|f_*(S)\|=d^{k-1}\|S\|$, implies that 
$\|\widetilde\Lambda(\varphi)\|_{W^*}\leq d^{-1/2}\|S\|$. The
proposition follows.
\hfill $\square$

\medskip

In the rest of this paragraph, we show that the Green measure $\mu$ is
moderate, see Appendix \ref{section_current_pk}. Recall that
a positive measure $\nu$ on $\P^k$ is {\it moderate} if there are constants
$\alpha>0$ and $c>0$ such that 
$$\|e^{-\alpha\varphi}\|_{L^1(\nu)}\leq c$$
for  $\varphi$ quasi-p.s.h. such that $\ddc\varphi\geq -\omega_\FS$
and $\langle\omega_\FS^k,\varphi\rangle=0$. Moderate measures are
PB and by linearity, if $\nu$ is
moderate and $\Dc$ is a bounded set of d.s.h. functions then
there are constants
$\alpha>0$ and $c>0$ such that 
$$\|e^{\alpha|\varphi|}\|_{L^1(\nu)}\leq c \quad \mbox{for}\quad \varphi\in\Dc.$$
Moderate measures were introduced in \cite{DinhSibony1}. The
fundamental estimate in Theorem \ref{th_hormander} in Appendix implies that smooth measures are
moderate. So, when we use test d.s.h. functions, several estimates for
the Lebesgue measure can be extended to moderate measures. For
example, we will see
that quasi-p.s.h. functions have comparable repartition functions with
respect to the Lebesgue measure $\omega_\FS^k$ and to the equilibrium
measure $\mu$.

It is shown in
\cite{DinhNguyenSibony3} that measures which are wedge-products of
positive closed $(1,1)$-currents with H\"older continuous potentials,
are moderate. In particular, the Green measure $\mu$ is moderate. We
will give here another proof of this result using the following
criterion. Since $\DSH(\P^k)$ is a subspace of $L^1(\P^k)$, the
$L^1$-norm induces a distance on $\DSH(\P^k)$ that we denote by
$\dist_{L^1}$. 

\begin{proposition}
Let $\nu$ be a PB positive measure on $\P^k$. Assume that $\nu$
restricted to any bounded subset of $\DSH(\P^k)$ is H\"older
continuous\footnote{This property is close to the property that $\nu$
  has a H\"older continuous super-potential.} with respect to $\dist_{L^1}$. Then $\nu$ is moderate.
\end{proposition}
\proof
Let $\varphi$ be a quasi-p.s.h. function such that $\ddc\varphi\geq
-\omega_\FS$ and $\langle \omega_\FS^k,\varphi\rangle=0$. 
We want to
prove that $\langle \nu, e^{-\alpha\varphi}\rangle \leq c$ for some
positive constants $\alpha,c$. For this purpose, we only have to show
that $\nu\{\varphi\leq -M\}\lesssim e^{-\alpha M}$ for some
constant $\alpha>0$ and for $M\geq 1$. Define for $M>0$,
$\varphi_M:=\max(\varphi,-M)$. These functions $\varphi_M$ belong to a
compact family $\Dc$ of d.s.h. functions. Observe that
$\varphi_{M-1}-\varphi_M$ is positive with support in $\{\varphi\leq
-M+1\}$. It is bounded by 1 and equal to 1 on
$\{\varphi<-M\}$. Therefore,  the H\"older continuity
of $\nu$ on $\Dc$ implies that there is a constant $\lambda>0$ such that
\begin{eqnarray*}
\nu\{\varphi<-M\} & \leq &  \langle \nu,\varphi_{M-1}-\varphi_M\rangle
=\nu(\varphi_{M-1})-\nu(\varphi_M)\\
& \lesssim  & \dist_{L^1}(\varphi_{M-1},\varphi_M)^\lambda \leq \vol\{\varphi<-M+1\}^\lambda.
\end{eqnarray*}
Since the Lebesgue measure is moderate, the last expression is $\lesssim e^{-\alpha
  M}$ for some positive constant $\alpha$. The proposition follows.
\endproof

We have the following result obtained in \cite{DinhNguyenSibony3}. It
will be used to establish several stochastic properties of d.s.h. observables for
the equilibrium measure.

\begin{theorem} \label{th_end_moderate}
Let $f$ be an endomorphism of algebraic degree $d\geq 2$ of
$\P^k$. Then, the Green measure $\mu$ of $f$ is H\"older continuous on
bounded subsets of $\DSH(\P^k)$. In particular, it is moderate.
\end{theorem}
\proof
Let $\Dc$ be a bounded set of d.s.h. functions. We have to show that
$\mu$ is H\"older continuous on $\Dc$ with respect to
$\dist_{L^1}$. By linearity, since $\mu$ is PC, it is enough to
consider the case where 
$\Dc$ is the set of d.s.h. functions $\varphi$ such that
$\langle\mu,\varphi\rangle\geq 0$ and $\|\varphi\|_\mu\leq 1$. Let
$\widetilde\Dc$ denote the set of d.s.h. functions
$\varphi-\Lambda(\varphi)$ with $\varphi\in\Dc$. By Proposition
\ref{prop_Perron_dsh}, $\widetilde\Dc$ is a bounded family of d.s.h. functions. We claim that
$\widetilde\Dc$ is invariant under
$\widetilde\Lambda:=d\Lambda$. Observe that if
$\varphi$ is in $\Dc$, then
$\widetilde\varphi:=\varphi-\langle\mu,\varphi\rangle$
is also in $\Dc$. Since $\langle\mu,\varphi\rangle =
\langle\mu,\Lambda(\varphi)\rangle$, we have
$$\widetilde\Lambda(\varphi-\Lambda(\varphi))=\widetilde\Lambda(\widetilde
\varphi-\Lambda(\widetilde \varphi))=
\widetilde\Lambda(\widetilde \varphi)
-\Lambda(\widetilde\Lambda(\widetilde \varphi)).$$
By Theorem \ref{th_green_measure_dsh}, $\widetilde\Lambda(\widetilde
\varphi)$ belongs to $\Dc$. This proves the claim. So, the crucial
point is that $\Lambda$ is contracting on an appropriate hyperplane.

For $\varphi,\psi$ in $L^1(\P^k)$ we have
$$\|\widetilde\Lambda (\varphi)-\widetilde\Lambda (\psi)\|_{L^1} \leq 
\int \widetilde \Lambda (|\varphi-\psi|)\omega_\FS^k = d^{1-k}\int
|\varphi-\psi|f^*(\omega_\FS^k) \lesssim \|\varphi-\psi\|_{L^1}.$$
So, the map $\widetilde\Lambda$ is Lipschitz with respect to
$\dist_{L^1}$. In particular, the map $\varphi\mapsto
\varphi-\Lambda(\varphi)$ is Lipschitz with respect to this distance.
Now, we have  for
$\varphi\in\Dc$
\begin{eqnarray*}
\langle\mu,\varphi\rangle & = & \lim_{n\rightarrow\infty} \langle
d^{-kn} (f^n)^*(\omega_\FS^k),\varphi\rangle =
\lim_{n\rightarrow\infty}
\langle\omega_\FS^k,\Lambda^n(\varphi)\rangle\\
& = &  
\langle\omega_\FS^k,\varphi\rangle -\sum_{n\geq 0} \langle
\omega_\FS^k,\Lambda^n(\varphi)-\Lambda^{n+1}(\varphi)\rangle \\
& = &  
\langle\omega_\FS^k,\varphi\rangle -\sum_{n\geq 0} d^{-n} \langle
\omega_\FS^k,\widetilde\Lambda^n(\varphi-\Lambda(\varphi))\rangle.
\end{eqnarray*}
By Lemma \ref{lemma_holder}, the last series defines a function on
$\widetilde\Dc$ which is H\"older continuous with respect
to $\dist_{L^1}$. Therefore, $\varphi\mapsto \langle \mu,\varphi\rangle$ is
H\"older continuous on $\Dc$.
\endproof

\begin{remark}\rm
Let $f_s$ be a family of endomorphisms of
algebraic degree $d\geq 2$, depending holomorphically on a parameter
$s\in\Sigma$. Let $\mu_s$ denote its equilibrium measure. 
We get that $(s,\varphi)\mapsto \mu_s(\varphi)$ is H\"older
continuous on bounded subsets of $\Sigma\times\DSH(\P^k)$. 
\end{remark}

The following results are  useful in the stochastic study of the
dynamical system.

\begin{corollary} \label{cor_exp_dsh}
Let $f$, $\mu$ and $\Lambda$ be as above. There are constants
$c>0$ and $\alpha>0$ such that if $\psi$ is d.s.h. with $\|\psi\|_\DSH\leq 1$, then 
$$\big\langle \mu, e^{\alpha d^n |\Lambda^n(\psi)-\langle\mu,\psi\rangle|}\big\rangle \leq c \quad
\mbox{and}\quad \|\Lambda^n(\psi)-\langle\mu,\psi\rangle\|_{L^q(\mu)}\leq cq d^{-n}$$
for every $n\geq 0$ and every  $1\leq q <+\infty$. 
\end{corollary}
\proof
By Theorem \ref{th_green_measure_dsh}, $d^n
(\Lambda^n(\psi)-\langle\mu,\psi\rangle)$ 
belong to a compact family in $\DSH(\P^k)$.
The first inequality in the corollary follows from Theorem \ref{th_end_moderate}.
For the second one, we can assume that $q$ is integer and we easily
deduce the result from the first inequality and the inequalities
 $x^q\leq q!e^x\leq q^qe^x$ for $x\geq 0$.
\endproof

\begin{corollary} \label{cor_exp_holder}
Let $0<\nu\leq 2$ be a constant. There are constants $c>0$ and
$\alpha>0$ such that if $\psi$ is a $\nu$-H{\"o}lder continuous function
with $\|\psi\|_{\Cc^\nu}\leq 1$, then 
$$\big\langle \mu, e^{\alpha d^{n\nu/2} |\Lambda^n(\psi)-\langle \mu,\psi\rangle|}\big\rangle \leq c
\quad \mbox{and}\quad 
\|\Lambda^n(\psi)-\langle\mu,\psi\rangle\|_{L^q(\mu)}\leq cq^{\nu/2} d^{-n\nu/2}$$
for every $n\geq 0$ and every  $1\leq q <+\infty$.
\end{corollary}
\proof
By Corollary \ref{cor_exp_dsh}, since $\|\cdot\|_\DSH\lesssim
\|\cdot\|_{\Cc^2}$, we have 
$$\|\Lambda^n(\psi)-\langle\mu,\psi\rangle\|_{L^q(\mu)}\leq cq d^{-n}\|\psi\|_{\Cc^2},$$
with $c>0$ independent of $q$ and of $\psi$. On the other hand, by
definition of $\Lambda$, we have 
$$\|\Lambda^n(\psi)-\langle\mu,\psi\rangle\|_{L^q(\mu)}\leq
\|\Lambda^n(\psi)-\langle \mu,\psi\rangle\|_{L^\infty(\mu)}\leq 2\|\psi\|_{\Cc^0}.$$
The theory of interpolation between the Banach spaces $\Cc^0$ and
$\Cc^2$ \cite{Triebel}, applied to the linear operator
$\psi\mapsto\Lambda^n(\psi)-\langle\mu,\psi\rangle$, implies that
$$\|\Lambda^n(\psi)-\langle\mu,\psi\rangle\|_{L^q(\mu)}\leq A_\nu 2^{1-\nu/2}[cq
d^{-n}]^{\nu/2}\|\psi\|_{\Cc^\nu},$$
for some constant $A_\nu>0$ depending only on $\nu$ and on $\P^k$. This
gives the second inequality in the corollary.

Recall that if $L$ is a
linear continuous functional on the space $\Cc^0$ of continuous
functions, then we have for every $0< \nu < 2$
$$\|L\|_{\Cc^\nu}\leq A_\nu
\|L\|_{\Cc^0}^{1-\nu/2}\|L\|_{\Cc^2}^{\nu/2}$$
for some constant $A_\nu>0$ independent of $L$ (in our case, the
functional is with values in $L^q(\mu)$).   

For the first inequality, we have for a fixed
constant $\alpha>0$ small enough, 
$$\big\langle\mu, e^{\alpha d^{n\nu/2} |\Lambda^n(\psi)-\langle\mu,\psi\rangle|}\big\rangle  = 
\sum_{q\geq 0} {1\over q!} \big\langle \mu,
|\alpha d^{n\nu/2}(\Lambda^n(\psi)-\langle\mu,\psi\rangle)|^q\big\rangle 
\leq \sum_{q\geq 0} {1\over q!} \alpha^q c^qq^q.$$
By Stirling's formula, the last sum converges. The result follows.
\endproof

\bigskip\bigskip

\begin{exercise} \label{exercise_equi_measure_pb}
Let $\varphi$ be a smooth function and $\varphi_n$ as in Theorem
\ref{th_green_measure_dsh}. Show that we can write $\varphi_n=\varphi_n^+-\varphi_n^-$
with $\varphi_n^\pm$ quasi-p.s.h. such that
$\|\varphi_n^\pm\|_\DSH\lesssim d^{-n}$ and $\ddc\varphi_n^\pm\gtrsim
-d^{-n}\omega_\FS$. Prove that $\varphi_n$ converge pointwise to $0$
out of a
pluripolar set. Deduce that if $\nu$ is a probability measure with no
mass on pluripolar sets, then $d^{-kn}(f^n)^*(\nu)$ converge to $\mu$.
\end{exercise}

\begin{exercise}
Let $\DSH_0(\P^k)$ be the space of d.s.h. functions $\varphi$ such that $\langle
\mu,\varphi\rangle=0$. Show that $\DSH_0(\P^k)$ is a closed subspace
of $\DSH(\P^k)$, invariant under $\Lambda$, and that the spectral
radius of $\Lambda$ on this space is equal to $1/d$. Note that $1$ is an eigenvalue of $\Lambda$
on $\DSH(\P^k)$, so, $\Lambda$ has a spectral gap on $\DSH(\P^k)$. Prove a similar
result for $W^*(\P^k)$. 
\end{exercise}


\section{Equidistribution of points}

In this paragraph, we show that the preimages of a generic point by
$f^n$ are equidistributed with respect to the Green measure $\mu$ when
$n$ goes to infinity. 
The proof splits in two parts. First, we prove that there is a maximal
proper algebraic 
set $\Ec$ which is totally invariant, then we show that for $a\not\in\Ec$, the
preimages of $a$ are equidistributed. We will also prove that the
convex set of probability measures $\nu$, which are totally invariant,
i.e. $f^*(\nu)=d^k\nu$, is finite dimensional. The equidistribution
for $a$ out of an algebraic set is reminiscent of the main questions in
value distribution theory (we will see in the next paragraph that
using super-potentials we can get
an estimate on the speed of convergence towards $\mu$, at least for
generic maps). 
Finally, we prove a theorem due to
Briend-Duval on the equidistribution of the repelling periodic points.
The following result was obtained by the authors in \cite{DinhSibony1},
see also \cite{DinhSibony9} and \cite{Brolin, FreireLopesMane,Lyubich}
for the case of dimension 1.

\begin{theorem} \label{th_equi_k}
Let $f$ be an endomorphism of $\P^k$ of algebraic
  degree $d\geq 2$ and $\mu$ its Green measure. Then there is a proper algebraic set
  $\Ec$ of $\P^k$, possibly empty, such that $d^{-kn} (f^n)^*(\delta_a)$ converge to
  $\mu$ if and only if $a\not\in\Ec$. Here, $\delta_a$ denotes the
  Dirac mass at $a$. Moreover, $\Ec$ is totally
  invariant: $f^{-1}(\Ec)=f(\Ec)=\Ec$ and is maximal in the sense that
  if $E$ is a proper algebraic set of $\P^k$ such that
  $f^{-n}(E)\subset E$ for some $n\geq 1$, then $E$ is contained in
  $\Ec$. 
\end{theorem}

Briend-Duval proved in \cite{BriendDuval2} the above convergence for $a$ outside
the orbit of the critical set. 
They announced the property for $a$ out of an algebraic set, but
there is a problem with
the counting of multiplicity in their lemma in \cite[p.149]{BriendDuval2}.

We also have the following earlier
result due to Forn\ae ss-Sibony \cite{FornaessSibony1}.

\begin{proposition} \label{prop_except_FS}
There is a pluripolar set $\Ec'$ such that if $a$
  is out of $\Ec'$, then $d^{-kn} (f^n)^*(\delta_a)$ converge to
  $\mu$.
\end{proposition}
\noindent
{\bf Sketch of proof.} We use here a version of the
above $\ddc$-method which is given in \cite{DinhSibony6} in a more general
setting. Let $\varphi$ be a smooth function and $\varphi_n$ as in
Theorem \ref{th_green_measure_dsh}. Then, the functions $\varphi_n$
are continuous. The estimates on $\varphi_n$
imply that 
the series $\sum \varphi_n$ converges in $\DSH(\P^k)$, hence converges
pointwise out of a pluripolar set. Therefore,
$\varphi_n(a)$ converge to $0$ for $a$ out of some
pluripolar set $E_\varphi$, see Exercise \ref{exercise_equi_measure_pb}. 
If $c_n:=\langle \omega_\FS^k,\Lambda(\varphi_n)\rangle$, we have as
in Theorem \ref{th_green_measure_dsh}
$$\langle d^{-kn}(f^n)^*(\delta_a),\varphi\rangle =
c_0+\cdots+c_n+\langle \delta_a,\varphi_n\rangle = c_0+\cdots+c_n+\varphi_n(a).$$
Therefore, $\langle d^{-kn}(f^n)^*(\delta_a),\varphi\rangle$ converge
to $c_\varphi=\langle\mu,\varphi\rangle$, for $a$ out of $E_\varphi$. 

Now, consider $\varphi$ in a countable family $\Fc$ which
is dense in the space of smooth functions. If $a$ is not in the union
$\Ec'$ of the pluripolar sets $E_\varphi$, the above convergence 
of $\langle d^{-kn}(f^n)^*(\delta_a),\varphi\rangle$ together with 
the density of $\Fc$ implies that $d^{-kn}(f^n)^*(\delta_a)$ converge
to $\mu$. Finally, $\Ec'$ is pluripolar since it is a countable union
of such sets.
\hfill $\square$ 

\medskip

For the rest of the proof, we follow a geometric method introduced by
Lyubich \cite{Lyubich} in dimension one  and developped in
higher dimension by Briend-Duval and Dinh-Sibony. We first prove the
existence of the exceptional set and give
several characterizations in the following general situation.
Let $X$ be an analytic set of pure dimension $p$ in $\P^k$
invariant under $f$, i.e. $f(X)=X$. Let
$g:X\rightarrow X$ denote the restriction of $f$ to $X$. 
The following result can be deduced from Section 3.4
in \cite{DinhSibony1}, see also \cite{Dinh5,DinhSibony9}.

\begin{theorem} \label{th_exceptional}
There is a proper analytic set $\Ec_X$ of $X$, possibly empty,
totally invariant under $g$, which is maximal in the following
sense. If $E$ is an analytic set of $X$, of dimension $<\dim X$, such
that $g^{-s}(E)\subset E$ for some $s\geq 1$, then $E\subset \Ec_X$.
Moreover, there are at most finitely many  analytic sets in
$X$ which are totally invariant under $g$.
\end{theorem}

Since $g$ permutes the irreducible components of $X$, we can find an
integer $m\geq 1$ such that $g^m$ fixes the components of $X$. 

\begin{lemma} \label{lemma_topological}
The topological degree of  $g^m$ is equal to 
$d^{mp}$. More precisely, there is a hypersurface $Y$ of $X$ containing
$\sing(X)\cup g^m(\sing(X))$ such that for $x\in X$ out of $Y$, the
fiber $g^{-m}(x)$ has exactly $d^{mp}$ points.
\end{lemma}
\proof
Since $g^m$ fixes the components of $X$, 
we can assume that $X$ is irreducible. It follows
that $g^m$ defines a covering over some Zariski dense open set of $X$.  We want to
prove that $\delta$, the degree of this covering, is equal to $d^{mp}$.
Consider the positive measure
$(f^m)^*(\omega_\FS^p)\wedge[X]$. Since $(f^m)^*(\omega_\FS^p)$ is cohomologous to
$d^{mp}\omega_\FS^p$, this measure is of mass 
$d^{mp}\deg(X)$. Observe that  $(f^m)_*$
preserves the mass of positive measures and that we have 
$(f^m)_*[X]=\delta[X]$. Hence,
\begin{eqnarray*}
d^{mp}\deg(X) & = & \|(f^m)^*(\omega_\FS^p)\wedge [X]\| =
\|(f^m)_*((f^m)^*(\omega_\FS^p)\wedge [X])\| \\
& = &
\|\omega_\FS^p\wedge (f^m)_*[X]\|  
=  \delta\|\omega_\FS^p\wedge [X]\|=\delta\deg(X).
\end{eqnarray*}
It follows that $\delta=d^{mp}$.  So, we can take for $Y$, a hypersurface
which contains the ramification values of $f^m$ and the set $\sing(X)\cup g^m(\sing(X))$.
\endproof

Let $Y$ be as above. Observe that if $g^m(x)\not\in Y$ then $g^m$ defines a bi-holomorphic
map between a neighbourhood of $x$ and a neighbourhood of $g^m(x)$ in
$X$. 
Let $[Y]$ denote the $(k-p+1,k-p+1)$-current of integration
on $Y$ in $\P^k$. Since $(f^{mn})_*[Y]$ is a positive closed
$(k-p+1,k-p+1)$-current of mass $d^{mn(p-1)}\deg(Y)$, we can define
the following ramification current
$$R=\sum_{n\geq 0} R_n:=\sum_{n\geq 0} d^{-mnp} (f^{mn})_*[Y].$$  
Let $\nu(R,x)$ denote the Lelong number of $R$ at $x$.
By Theorem \ref{th_siu}, for $c>0$,
$E_c:=\{\nu(R,x)\geq c\}$ 
is an analytic set of dimension $\leq p-1$ contained in $X$. Observe
that $E_1$  contains $Y$. We will see that $R$ measures the obstruction for
constructing good backwards orbits.

For any point $x\in X$ let $\lambda'_n(x)$ denote the number of  distinct orbits
$$x_{-n},x_{-n+1},\ldots,x_{-1}, x_0$$
such that $g^m(x_{-i-1})=x_{-i}$, $x_0=x$ and $x_{-i}\in X\setminus Y$ for
$0\leq i\leq n-1$. These are the ``good" orbits. Define $\lambda_n:=d^{-mpn}\lambda'_n$.
The function $\lambda_n$ is lower semi-continuous with respect to
the Zariski topology on $X$. Moreover, by Lemma
\ref{lemma_topological}, 
we have $0\leq \lambda_n\leq 1$  and $\lambda_n=1$ 
out of the analytic set $\cup_{i=0}^{n-1}g^{mi}(Y)$. The sequence 
$(\lambda_n)$ decreases to a function $\lambda$, which represents the
asymptotic proportion of backwards orbits in $X\setminus Y$.

\begin{lemma} \label{lemma_exceptional}
There is a constant $\gamma>0$ such that
  $\lambda\geq\gamma$ on   $X\setminus E_1$.
\end{lemma}
\proof
We deduce from Theorem \ref{th_siu},  the existence of a constant
$0<\gamma<1$ satisfying $\{\nu(R,x)> 1-\gamma\}=E_1$. 
Indeed, the sequence of analytic sets $\{\nu(R,x)\geq 1-1/i\}$ is
decreasing, hence stationary.
Consider a
point $x\in X\setminus E_1$. We have $x\not\in Y$ and if
$\nu_n:=\nu(R_n,x)$, then $\sum
\nu_n\leq 1-\gamma$. 
Since $E_1$ contains $Y$, $\nu_0=0$ and $F_1:=g^{-m}(x)$
contains exactly $d^{mp}$ points. The definition of $\nu_1$, which is
``the multiplicity'' of $d^{-mp} (f^m)_*[Y]$ at $x$, implies
that 
$g^{-m}(x)$ contains at most $\nu_1d^{mp}$ points in $Y$. 
Then 
$$\#g^{-m}(F_1\setminus Y)=d^{mp}\#(F_1\setminus Y) \geq
(1-\nu_1)d^{2mp}.$$
Define $F_2:=g^{-m}(F_1\setminus Y)$. The definition of $\nu_2$
implies that $F_2$ contains at most 
$\nu_2 d^{2mp}$ points in $Y$. Hence, $F_3:=g^{-m}(F_2\setminus Y)$
contains at least $(1-\nu_1-\nu_2)d^{3mp}$ points. In the same way, we
define $F_4$, $\ldots$, $F_n$ with $\#F_n\geq (1-\sum
\nu_i)d^{mpn}$. Hence, for every $n$ we get the following estimate: 
$$\lambda_n(x)\geq
d^{-mpn}\#F_n\geq 1-\sum\nu_i\geq \gamma.$$ 
This proves the lemma.
\endproof

\medskip

\noindent
{\bf End of the proof of Theorem \ref{th_exceptional}.}
Let $\Ec_X^n$ denote the set of  $x\in X$ such that
$g^{-ml}(x)\subset E_1$ for $0\leq l\leq n$ and define
$\Ec_X:=\cap_{n\geq 0 }\Ec_X^n$.
Then, $(\Ec_X^n)$ is a decreasing sequence of analytic subsets of $E_1$.
It should be stationary. So, there is $n_0\geq 0$ such that
$\Ec^n_X=\Ec_X$ for $n\geq n_0$.

By definition, $\Ec_X$ is the set of $x\in X$ such that
$g^{-mn}(x)\subset E_1$ for every $n\geq 0$. Hence,
$g^{-m}(\Ec_X)\subset \Ec_X$. It follows that the sequence of analytic
sets $g^{-mn}(\Ec_X)$ is decreasing and there is $n\geq 0$ such that
$g^{-m(n+1)}(\Ec_X)=g^{-mn} (\Ec_X)$. Since $g^{mn}$  is surjective,
we deduce that $g^{-m}(\Ec_X)=\Ec_X$ and hence $\Ec_X=g^m(\Ec_X)$.

Assume as in the theorem that $E$ is analytic with $g^{-s}(E)\subset
E$. Define $E':=g^{-s+1}(E)\cup\ldots\cup E$. We have $g^{-1}(E')\subset
E'$ which implies $g^{-n-1}(E') \subset g^{-n}(E')$ for every $n\geq
0$. Hence, $g^{-n-1}(E')=g^{-n}(E')$ for $n$ large enough. This and
the surjectivity of $g$ imply that
$g^{-1}(E')=g(E')=E'$. 
By Lemma \ref{lemma_topological}, the topological degree of
$(g^{m'})_{|E'}$ is at most $d^{m'(p-1)}$ for some $m'\geq 1$.
This, the identity $g^{-1}(E')=g(E')=E'$ together with Lemma \ref{lemma_exceptional}
imply that $E'\subset E_1$. Hence,
$E'\subset \Ec_X$ and $E\subset \Ec_X$. 

Define $\Ec_X':=g^{-m+1}(\Ec_X)\cup\ldots\cup \Ec_X$. We have
$g^{-1}(\Ec_X')=g(\Ec_X')=\Ec_X'$. Applying the previous assertion to
$E:=\Ec'_X$ yields $\Ec'_X\subset \Ec_X$. Therefore,
$\Ec_X'=\Ec_X$ and $g^{-1}(\Ec_X)=g(\Ec_X)=\Ec_X$. So, $\Ec_X$ is the
maximal proper analytic set in $X$ which is totally invariant under $g$.

We prove now that there are only finitely many totally invariant
algebraic sets.
We only have to consider totally invariant sets $E$ of pure dimension
$q$. The proof is by induction on the dimension $p$ of $X$. The case
$p=0$ is clear. Assume
that the assertion is true for $X$ of dimension $\leq p-1$ and
consider the case of dimension $p$. If $q=p$ then $E$ is a union of
components of $X$. There are only a finite number of such analytic
sets. If $q< p$, we have seen that $E$ is contained in
$\Ec_X$. Applying the induction hypothesis to the restriction of
$f$ to $\Ec_X$ gives the result.
\hfill $\square$

\medskip

We now give another characterization of $\Ec_X$.
Observe that if $X$ is not locally irreducible at a point $x$ then
$g^{-m}(x)$ may contain more than $d^{mp}$ points. 
Let $\pi:\widetilde X\rightarrow X$ be the normalization of $X$, see
Appendix \ref{section_pk}. By Theorem \ref{th_normalization} applied to $g\circ\pi$, $g$ can be lifted to a
map $\widetilde g:\widetilde X\rightarrow \widetilde X$ such that
$g\circ \pi=\pi\circ\widetilde g$. Since $g$ is finite, $\widetilde
g$ is also finite. We deduce that $\widetilde g^m$ defines ramified
coverings of degree $d^{mp}$ on each component of $\widetilde X$. In
particular, any fiber of $\widetilde g^m$ contains at most $d^{mp}$
points. Observe that if $g^{-1}(E)\subset E$ then
$\widetilde g^{-1}(\widetilde E)\subset \widetilde E$ where
$\widetilde E:=\pi^{-1}(E)$. Theorem
\ref{th_exceptional} can be extended to $\widetilde g$.
For simplicity, we consider the case where $X$ is itself a normal analytic
space. If $X$ is not normal, one should work with its normalization.

Let $Z$ be a hypersurface of $X$ containing $E_1$.
Let $N_n(a)$ denote the number of orbits of $g^m$
$$a_{-n},\ldots,a_{-1},a_0$$
with $g^m(a_{-i-1})=a_{-i}$ and $a_0=a$ such that $a_{-i}\in Z$
for every $i$. Here, the orbits are counted with multiplicity. So,
$N_n(a)$ is the number of
negative orbits of order $n$ of $a$ which stay in $Z$. Observe
that the sequence of functions $\tau_n:=d^{-pmn} N_n$ decreases to some
function $\tau$. Since $\tau_n$ are upper semi-continuous with respect
to the Zariski topology and $0\leq \tau_n\leq 1$ (we use here the
assumption that $X$
is normal), the function $\tau$
satisfies the same properties. 
Note that $\tau(a)$ is the probability that an infinite negative
orbit of $a$ stays in $Z$.

\begin{proposition} \label{prop_tau_except}
Assume that $X$ is normal.
Then, $\tau$ is the characteristic function
  of $\Ec_X$, that is, $\tau=1$ on $\Ec_X$ and $\tau=0$ on $X\setminus
  \Ec_X$.
\end{proposition}
\proof
Since $\Ec_X\subset Z$ and $\Ec_X$ is totally invariant by
$g$, we have $\Ec_X\subset \{\tau=1\}$. Let $\theta\geq 0$ denote the
maximal value of $\tau$ on $X\setminus \Ec_X$. This value exists since
$\tau$ is upper semi-continuous with respect to the Zariski topology
(indeed, it is enough to consider the algebraic subset
$\{\tau\geq\theta'\}$ of $X$ which decreases when $\theta'$
increases; the family is stationary).
We have to check that $\theta=0$. Assume in order to obtain a
contradiction that $\theta>0$. Since $\tau\leq 1$, we always have
$\theta\leq 1$. 
Consider the non-empty analytic set
$E:=\tau^{-1}(\theta)\setminus\Ec_X$ in $Z\setminus \Ec_X$. Let $a$ be a point in $E$. Since
$\Ec_X$ is totally invariant, we have
$g^{-m}(a)\cap\Ec_X=\varnothing$. 
Hence, $\tau(b)\leq\theta$ for every $b\in g^{-m}(a)$.
We deduce
from the definition of $\tau$ and $\theta$ that
$$\theta=\tau(a)\leq d^{-pm}\sum_{b\in g^{-m}(a)} \tau(b)\leq
\theta.$$
It follows that $g^{-m}(a)\subset E$. Therefore, the analytic subset $\overline E$ 
of $Z$ satisfies $g^{-m}(\overline E)\subset\overline
E$. This contradicts the maximality of $\Ec_X$. 
\endproof

We continue the proof of Theorem \ref{th_equi_k}. We will use the
above results for $X=\P^k$, $Y$ the set of critical values of $f$.
Let $R$ be the ramification current defined as above by
$$R=\sum_{n\geq 0} R_n:=\sum_{n\geq 0} d^{-kn} (f^n)_*[Y].$$  
The following proposition was obtained in \cite{DinhSibony1},
a weaker version was independently obtained by Briend-Duval 
\cite{BriendDuval3}. Here, an inverse branch on $B$ for $f^n$ is a
bi-holomorphic map $g_i:B\rightarrow U_i$ such that 
$g_i\circ f^n$ is identity on $U_i$.

\begin{proposition} \label{prop_inverse_ball}
Let $\nu$ be a strictly positive constant. Let $a$ be a point in $\P^k$ such that
the Lelong number $\nu(R,a)$ of $R$ at
$a$ is strictly smaller than $\nu$. Then, there is a ball $B$ centered at $a$
such that $f^n$ admits at least $(1-\sqrt{\nu})d^{kn}$ inverse
branches $g_i:B\rightarrow U_i$ where $U_i$ are open sets in $\P^k$
of diameter $\leq d^{-n/2}$. 
In particular, if $\mu'$ is a limit value
of the measures $d^{-kn}(f^n)^*(\delta_a)$ then $\|\mu'-\mu\|\leq 2\sqrt{\nu(R,a)}$. 
\end{proposition}

Given a local coordinate system at $a$, let $\Fc$ denote the family of
complex lines passing through $a$. 
For such a line $\Delta$ denote by $\Delta_r$ the disc of center $a$
and of radius $r$.
The family $\Fc$ is parametrized by
$\P^{k-1}$ where the probability measure (the volume form) 
associated to the Fubini-Study metric is denoted by $\L$. Let $B_r$
denote the ball of center $a$ and of radius $r$.

\begin{lemma} \label{lemma_slice_lelong}
Let $S$ be a positive closed $(1,1)$-current on a neighbourhood of
$a$. Then for any $\delta>0$ there is an $r>0$ and a family
$\Fc'\subset\Fc$, such that $\L(\Fc')\geq 1-\delta$ and for every
$\Delta$ in $\Fc'$, the measure $S\wedge [\Delta_r]$ is well-defined
and of mass $\leq \nu(S,a)+\delta$, where $\nu(S,a)$ is the Lelong
number of $S$ at $a$. 
\end{lemma}
\proof
Let $\pi:\widehat \P^k\rightarrow \P^k$ be the blow-up of $\P^k$ at $a$ and $E$
the exceptional hypersurface. Then, we can write
$\pi^*(S)=\nu(S,a)[E]+S'$ with $S'$ a current having no mass on $E$,
see Exercise \ref{exercise_lelong_slice}. It is clear that for almost
every $\Delta_r$, the restriction of the potentials of $S$ to
$\Delta_r$ is not identically $-\infty$, so, the measure $S\wedge [\Delta_r]$ is well-defined. Let
$\widehat \Delta_r$ denote the strict transform of $\Delta_r$ by $\pi$,
i.e. the closure of $\pi^{-1}(\Delta_r\setminus\{a\})$. Then, the
$\widehat\Delta_r$ define a smooth holomorphic fibration over $E$.
The measure $S\wedge [\Delta_r]$ is equal to the push-forward of
$\pi^*(S)\wedge [\widehat\Delta_r]$ by $\pi$. 
Observe that $\pi^*(S)\wedge [\widehat\Delta_r]$ is equal to $S'\wedge
[\widehat\Delta_r]$ plus $\nu(S,a)$ times the Dirac mass at
$\widehat\Delta_r\cap E$. Therefore, we only have to consider the
$\Delta_r$ such that $S'\wedge [\widehat\Delta_r]$ are of mass $\leq\delta$.

Since $S'$ have no mass on $E$, its mass on $\pi^{-1}(B_r)$ tends
to 0 when $r$ tends to 0.
It follows from Fubini's theorem that when $r$ is small enough the mass of the slices
$S'\wedge [\widehat\Delta_r]$ is $\leq\delta$ except for a small
family of $\Delta$. This proves the lemma. 
\endproof

\begin{lemma} \label{lemma_slice_mass}
Let $U$ be a neighbourhood of $\overline B_r$.
Let $S$ be
a positive closed $(1,1)$-current on $U$. Then, for every $\delta>0$,
there is a family $\Fc'\subset \Fc$ with $\L(\Fc')>1-\delta$, such that for $\Delta$ in $\Fc'$,
the measure $S\wedge [\Delta_r]$ is well-defined and of mass $\leq
A\|S\|$, where 
$A>0$ is a constant depending on $\delta$ but independent of $S$.
\end{lemma}
\proof
We can assume that $\|S\|=1$. 
Let $\pi$ be as in Lemma \ref{lemma_slice_lelong}.
Then, by continuity of $\pi^*$, the mass of $\pi^*(S)$ on $\pi^{-1}(B_r)$
is bounded by a constant. It is enough to apply 
Fubini's theorem in order to estimate the mass of $\pi^*(S)\wedge
[\widehat\Delta_r]$. 
\endproof

Recall the following theorem due to Sibony-Wong \cite{SibonyWong}.

\begin{theorem} \label{th_sibony_wong}
Let $m>0$ be a positive constant.
Let $\Fc'\subset \Fc$ be such that $\L(\Fc')\geq m$ and
let $\Sigma$ denote the intersection of the family $\Fc'$ with $B_r$. Then any 
holomorphic function $h$ on a neighbourhood of $\Sigma$ can be
extended to a holomorphic function on $B_{\lambda r}$ where $\lambda>0$
is a constant depending on $m$ but independent of $\Fc'$ and $r$. Moreover, we have
$$\sup_{B_{\lambda r}}|h|\leq \sup_\Sigma |h|.$$ 
\end{theorem}

We will use the following version of a lemma due to Briend-Duval
\cite{BriendDuval2}. Their proof uses the theory of moduli of annuli.

\begin{lemma} \label{lemma_bd}
Let $g:\Delta_r\rightarrow\P^k$ be a holomorphic map from a disc of
center $0$ and of radius $r$ in $\C$. Assume that 
$\area(g(\Delta_r))$ counted with multiplicity, is smaller than $1/2$. Then for any
$\epsilon>0$ there is a constant $\lambda>0$ independent of $g,r$ such that
the diameter of $g(\Delta_{\lambda r})$ is smaller than
$\epsilon\sqrt{\area(g(\Delta_r))}$. 
\end{lemma}
\proof
Observe that the lemma is an easy consequence of the Cauchy formula if
$g$ has values in a compact set of $\C^k\subset\P^k$. In order to
reduce the problem to this case, it is enough to prove that given an
$\epsilon_0>0$, there is a
constant $\lambda_0>0$ such that $\diam(g(\Delta_{\lambda_0r}))\leq
\epsilon_0$. For $\epsilon_0$ small enough, we can apply the
above case to $g$ restricted to $\Delta_{\lambda_0r}$.

By hypothesis, the graphs $\Gamma_g$ of $g$ in $\Delta_r\times\P^k$ have bounded
area. So, according to Bishop's theorem \cite{Bishop}, these
graphs form a relatively compact family of analytic sets, that is, the limits of
these graphs in the Hausdorff sense, are analytic sets. Since
$\area(g(\Delta_r))$ is bounded by $1/2$, the limits have no
compact components. So, they are also graphs and the family of
the maps $g$ is compact. We deduce that  $\diam(g(\Delta_{\lambda_0r}))\leq
\epsilon_0$ for $\lambda_0$ small enough.
\endproof

\noindent
{\bf Sketch of the proof of Proposition \ref{prop_inverse_ball}.}
The last assertion in the proposition is deduced from the first one
and Proposition \ref{prop_except_FS} applied to a generic point in $B$. 
We obtain that $\|\mu'-\mu\|\leq 2\sqrt{\nu}$ for every $\nu$ strictly
larger than $\nu(R,a)$ which implies the result.

For the first assertion, the idea is to construct inverse branches for many discs passing
through $a$ and then to apply Theorem \ref{th_sibony_wong} in order to
construct inverse branches on balls. 
We can assume that $\nu$ is smaller than 1.
Choose constants $\delta>0$, $\epsilon>0$ 
small enough and then a constant $\kappa>0$ large enough; all 
independent of $n$. 
Fix now the integer $n$. 
Recall that $\|(f^n)_*(\omega_\FS)\|=d^{(k-1)n}$. 
By Lemmas \ref{lemma_slice_lelong} and \ref{lemma_slice_mass}, there is a family
$\Fc'\subset \Fc$ and a constant $r>0$ such that $\L(\Fc')>1-\delta$
and for any $\Delta$ in
$\Fc'$, the mass of $R\wedge [\Delta_{\kappa^2 r}]$ is smaller than
$\nu$ and the mass of $(f^n)_*(\omega_\FS)\wedge
[\Delta_{\kappa r}]$ is smaller than $Ad^{(k-1)n}$ with $A>0$. 

\medskip
\noindent
{\bf Claim.} For each $\Delta$ in $\Fc'$, $f^n$ admits at least 
$(1-2\nu)d^{kn}$ inverse branches $g_i:\Delta_{\kappa^2 r}\rightarrow
V_i$ with $\area(V_i)\leq A\nu^{-1}d^{-n}$. The inverse branches $g_i$
can be extended to a neighbourhood of $\Delta_{\kappa^2 r}$.
\medskip

Assuming the claim, we complete the proof of the proposition.
Let $a_1,\ldots,a_l$ be the points in $f^{-n}(a)$, with $l\leq
d^{kn}$, and $\Fc_s'\subset\Fc'$ the family of $\Delta$'s such that one
of the previous inverse branches $g_i:\Delta_{\kappa^2 r}\rightarrow
V_i$ passes through $a_s$, that is, $V_i$ contains $a_s$. The above
claim implies that $\sum \L(\Fc'_s)\geq
(1-\delta)(1-2\nu)d^{kn}$. There are at most $d^{kn}$ terms in this
sum. We only consider the family $\Sc$ of the indices $s$ such that $\L(\Fc'_s)\geq
1-3\sqrt{\nu}$. Since
$\L(\Fc'_s)\leq 1$ for every $s$, we have
$$\#\Sc + (d^{kn}-\#\Sc)(1-3\sqrt{\nu})\geq \sum \L(\Fc'_s)\geq
(1-\delta)(1-2\nu)d^{kn}.$$
Therefore, since $\delta$ is small, we have $\#\Sc\geq
(1-\sqrt{\nu})d^{kn}$. 
For any index $s\in\Sc$ and for $\Delta$ in
$\Fc'_s$, by Lemma \ref{lemma_bd}, the corresponding inverse branch on
$\Delta_{\kappa r}$, which passes through $a_s$, has diameter $\leq \epsilon d^{-n/2}$.
By Theorem \ref{th_sibony_wong}, $f^n$ admits an inverse branch
defined on the ball $B_r$ and 
passing through $a_s$, with diameter $\leq d^{-n/2}$. This implies the
result.

\medskip
\noindent
{\bf Proof of the claim.} Let $\nu_l$ denote the mass of $R_l\wedge
[\Delta_{\kappa^2r}]$. Then, $\sum \nu_l$ is the mass of $R\wedge
[\Delta_{\kappa^2r}]$. Recall that this mass is smaller than $\nu$. 
By definition,
$\nu_ld^{kl}$ is the number of points in $f^l(Y)\cap
\Delta_{\kappa^2r}$, counted with multiplicity. We only have to
consider the case $\nu<1$. So, we have $\nu_0=0$ and $\Delta_{\kappa^2r}$ does not
intersect $Y$, the critical values of $f$. It follows that
$\Delta_{\kappa^2r}$ admits $d^k$ inverse branches for $f$. By
definition of $\nu_1$, there are at most $\nu_1d^k$ such inverse
branches which intersect $Y$, i.e. the images intersect $Y$. 
So, $(1-\nu_1)d^k$ of them do not meet $Y$ and the image of such a
branch admits $d^k$ inverse branches for $f$. We conclude that
$\Delta_{\kappa^2r}$ admits at least $(1-\nu_1)d^{2k}$ inverse
branches for $f^2$. By induction, we construct for $f^n$ at least
$(1-\nu_1-\cdots-\nu_{n-1})d^{kn}$ inverse branches on
$\Delta_{\kappa^2r}$.

Now, observe that the mass of $(f^n)_*(\omega_\FS)\wedge
[\Delta_{\kappa r}]$ is exactly the area of $f^{-n}(\Delta_{\kappa
  r})$. We know that it is smaller than
$Ad^{(k-1)n}$. It is not difficult to see that $\Delta_{\kappa^2r}$
has at most $\nu d^{kn}$ inverse branches with area $\geq
A\nu^{-1}d^{-n}$.
This completes the proof.
\hfill $\square$

\medskip
\noindent
{\bf End of the proof of Theorem \ref{th_equi_k}.}
Let $a$ be a point out of the exceptional set $\Ec$ defined in
Theorem \ref{th_exceptional} for $X=\P^k$. Fix $\epsilon>0$ and a constant $\alpha>0$ small
enough. If $\mu'$ is a limit value of $d^{-kn}(f^n)^*(\delta_a)$, it is
enough to show that $\|\mu'-\mu\|\leq
2\alpha+2\epsilon$. 
Consider $Z:=\{\nu(R,z)>\epsilon\}$ and $\tau$ as in 
Proposition 
\ref{prop_tau_except} for $X=\P^k$.  We have
$\tau(a)=0$. So, for $r$ large enough we have $\tau_r(a)\leq\alpha$. 
Consider all the negative orbits $\Oc_j$ of order $r_j\leq r$ 
$$\Oc_j=\{a^{(j)}_{-r_j},\ldots,a^{(j)}_{-1},a^{(j)}_{0}\}$$
with  $f(a^{(j)}_{-i-1})=a^{(j)}_{-i}$ and $a^{(j)}_{0}=a$ such that 
$a^{(j)}_{-r_j}\not\in Z$ and $a^{(j)}_{-i}\in Z$ for
$i\not=r_j$. Each orbit is repeated according to its multiplicity. 
Let $S_r$ denote the family of points $b\in
f^{-r}(a)$ such that $f^{i}(b)\in Z$ for $0\leq i\leq r$.
Then $f^{-r}(a)\setminus S_r$ consists of the preimages of the points
$a^{(j)}_{-r_j}$. So, by definition of $\tau_r$, we have
$$d^{-kr}\#S_r=\tau_r(a)\leq \alpha$$
and
$$d^{-kr}\sum_j d^{k(r-r_j)}=d^{-kr}\#(f^{-r}(a)\setminus S_r)=1-\tau_r(a)\geq 1-\alpha.$$
We have for $n\geq r$
\begin{eqnarray*}
d^{-kn}(f^{n})^*(\delta_a) & = & d^{-kn}\sum_{b\in S_r} 
(f^{(n-r)})^*(\delta_b) + d^{-kn}\sum_j
(f^{(n-r_j)})^*(\delta_{a^{(j)}_{-r_j}}).
\end{eqnarray*}
Since $d^{-kn}(f^{n})^*$ preserves the mass of any measure, the first term
in the last sum is of mass $d^{-kr}\#S_r=\tau_r(a)\leq \alpha$ and the
second term is of mass $\geq 1-\alpha$. We apply Proposition \ref{prop_inverse_ball} to
the Dirac masses at $a^{(j)}_{-r_j}$. We deduce that if $\mu'$ is a limit value of
$d^{-kn}(f^{n})^*(\delta_a)$ then
$$\| \mu'-\mu\|\leq 2\alpha +(1-\alpha)2\epsilon\leq 2\alpha+2\epsilon.$$
This completes the proof of the theorem.
\hfill $\square$

\medskip
We have the following more general result. When $X$ is not normal, one
has an analogous result for the lift of $g$ to the normalization of $X$.

\begin{theorem} \label{th_equi_point_X}
Let $X$ be an irreducible analytic set of dimension $p$,
invariant under $f$. Let $g$ denote the restriction of $f$ to $X$ and
$\Ec_X$ the exceptional set of $g$. Assume that $X$ is a normal
analytic space.
Then
$d^{-pn}(g^{n})^*(\delta_a)$ converge to $\mu_X:=(\deg X)^{-1}
T^p\wedge [X]$ if and only if
$a$ is out of $\Ec_X$. 
Moreover, the convex set of probability measures on $X$ which are totally invariant under
$g$, is of finite dimension.
\end{theorem}
\proof
The proof of the first assertion follows the same lines as in Theorem
\ref{th_equi_k}. We use the fact that $g$ is the restriction of a
holomorphic map in $\P^k$ in order to define the ramification current $R$.
The assumption that $X$ is normal allows to define
$d^{-pn}(g^{n})^*(\delta_a)$. 
We prove the second  assertion.
Observe that an analytic set, totally
invariant by $g^n$, is not necessarily totally invariant by $g$, but
it is a union of components of such sets, see Theorem \ref{th_exceptional}. 
Therefore, we can replace $g$ by an iterate $g^n$ in
order to assume that $g$ fixes all the components of all 
the totally invariant analytic sets. Let $\mu'$ be an extremal element
in the convex set of totally invariant probability measures and $X'$ the smallest
totally invariant analytic set such that $\mu'(X')=1$. The first
assertion applied to $X'$ implies that $\mu'=\mu_{X'}$. Hence, the set
of such $\mu'$ is finite. We use a normalization of $X'$ if necessary.
\endproof

The following result due to Briend-Duval \cite{BriendDuval1}, shows that
repelling periodic points are equidistributed on the support of the
Green measure.

\begin{theorem} \label{th_periodic_pk}
Let $P_n$ denote the set of repelling periodic points of period
$n$ on the support of $\mu$. Then the sequence of measures
$$\mu_n:=d^{-kn}\sum_{a\in P_n}\delta_a$$
converges to $\mu$.
\end{theorem}
\proof
By Proposition \ref{prop_number_periodic_pk},
the number of periodic points of period $n$ of $f$, counted with
multiplicity, is equal to $(d^n-1)^{-1}(d^{k(n+1)}-1)$. Therefore, any
limit value $\mu'$ of $\mu_n$ is of mass $\leq 1$. Fix a small constant
$\epsilon>0$. It is enough to
check that for $\mu$-almost every point $a\in \P^k$, there is a ball
$B$ centered at $a$, arbitrarily small, such that $\#P_n\cap B\geq
(1-\epsilon)d^{kn}\mu(B)$ for large $n$. We will use in particular a trick due to
X. Buff, which simplifies the original proof.

Since $\mu$ is PC, it has no mass on analytic sets. So, it has no
mass on the orbit $\Oc_Y$ of $Y$, the set of critical
values of $f$. 
Fix a point $a$ on the support of
$\mu$ and out of $\Oc_Y$. We have $\nu(R,a)=0$. 
By Proposition \ref{prop_inverse_ball}, 
there is a ball $B$ of center $a$, with sufficiently small radius, which
admits $(1-\epsilon^2)d^{kn}$ inverse branches of
diameter $\leq d^{-n/2}$ for $f^n$ when $n$ is large enough. 
Choose a finite family of such balls $B_i$ of center
$b_i$ with $1\leq i\leq m$
such that $\mu(B_1\cup\ldots\cup B_m)> 1-\epsilon^2\mu(B)$ and each
$B_i$ admits $(1-\epsilon^2\mu(B))d^{kn}$ inverse branches of
diameter $\leq d^{-n/2}$ for $f^n$ when $n$ is large enough. Choose
balls $B_i'\Subset B_i$ such that $\mu(B'_1\cup\ldots\cup B'_m)>
1-\epsilon^2\mu(B)$. 

Fix a constant $N$ large enough. Since $d^{-kn} (f^n)^*(\delta_a)$ converge to
$\mu$, there are at least $(1-2\epsilon^2)d^{kN}$ inverse branches for
$f^N$ whose image intersects $\cup B_j'$ and
then with image contained in one of the $B_j$. In the same way, we show
that for $n$ large enough,
each $B_j$ admits $(1-2\epsilon^2)\mu(B)d^{k(n-N)}$ inverse branches 
for $f^{n-N}$ with images in $B$. Therefore, $B$ admits at least 
$(1-2\epsilon^2)^2\mu(B)d^{kn}$ inverse branches $g_i:B\rightarrow
U_i$ for $f^n$ with image $U_i\Subset B$. Observe that every 
holomorphic map $g:B\rightarrow U\Subset B$ contracts the Kobayashi
metric and then admits an attractive fixed point $z$. Moreover, $g^l$
converges uniformly to $z$ and $\cap_l g^l(\overline B)=\{z\}$. Therefore, each
$g_i$ admits a fixed attractive point $a_i$. This point is fixed
and repelling for $f^n$. They are different since the $U_i$ are disjoint. 
Finally, since $\mu$ is totally invariant, its support is also totally
invariant under $f$. Hence, $a_i$, which is equal to $\cap_l
g_i^l(\supp(\mu)\cap \overline B)$, is necessarily in
$\supp(\mu)$. 
We deduce that 
$$\#P_n\cap B\geq (1-2\epsilon^2)^2\mu(B)d^{kn}
\geq (1-\epsilon)d^{kn}\mu(B).$$ 
This completes the proof.
\endproof

Note that in the
previous theorem, one can replace $P_n$ by the set of all periodic
points counting with multiplicity or not.

\bigskip\bigskip

\begin{exercise}
Let $f$ be an endomorphism of algebraic degree $d\geq 2$ of
$\P^k$. Let $K$ be a compact set such that $f^{-1}(K)\subset K$. Show that
either $K$ contains $\Jc_k$, the Julia set of order $k$, or $K$ is contained in an analytic set. 
Let $U$ be an open set which intersects the Julia set $\Jc_k$. 
Show that $\cup_{n\geq 0} f^n(U)$ is a Zariski dense open set of
$\P^k$. 
Prove that
$a\not\in\Ec$ if and only if $\cup f^{-n}(a)$ is Zariski dense.
\end{exercise}

\begin{exercise} Assume that $p$ is a repelling fixed point in
  $\Jc_k$. 
If $g$ is another endomorphism close enough to $f$
  in $\Hc_d(\P^k)$ such that $g(p)=p$, show that $p$ belongs also to
  the Julia set of order $k$ of $g$. Hint: use that $g\mapsto \mu_g$
  is continuous.
\end{exercise}

\begin{exercise} 
Using Example \ref{example_ueda}, construct a map $f$ in $\Hc_d(\P^k)$, $d\geq
2$, such that for $n$ large enough, every fiber of $f^n$ contains more
than $d^{(k-1/2)n}$ points. 
Deduce that there is Zariski dense open set in
  $\Hc_d(\P^k)$ such that if $f$ is in that Zariski open set, its exceptional
  set is empty.
\end{exercise}

\begin{exercise}
Let $\epsilon$ be a fixed constant such that $0<\epsilon<1$.
Let $P_n'$ the set of repelling periodic points $a$ of prime period
$n$ on the support of $\mu$ such that all the eigenvalues of $Df^n$ at
$a$ are of modulus $\geq (d-\epsilon)^{n/2}$. Show that $d^{-kn}\sum_{a\in P'_n}\delta_a$
converges to $\mu$.
\end{exercise}

\begin{exercise} Let $g$ be as in Theorem \ref{th_equi_point_X}. Show that repelling
  periodic points on $\supp(\mu_X)$ are equidistributed with respect
  to $\mu_X$. In particular, they are Zariski dense.
\end{exercise}


\section{Equidistribution of varieties}

In this paragraph, we consider the inverse images by $f^n$ of
varieties in $\P^k$. The geometrical method in the last
paragraph is quite difficult to apply here. Indeed, the inverse image of
a generic variety of codimension $p<k$ is irreducible of degree
$O(d^{pn})$. The pluripotential method that we
introduce here is probably the right method for the equidistribution
problem. Moreover, it should give  some precise estimates on the convergence, see Remark \ref{rk_unif_cv_perron}.

The following result, due to the authors, gives a satisfactory
solution in the case of hypersurfaces. It was proved for Zariski generic maps
by Forn\ae ss-Sibony in \cite{FornaessSibony3, Sibony} and for maps in dimension 2 by
Favre-Jonsson in \cite{FavreJonsson}. More precise results
are given in \cite{DinhSibony9} and in \cite{FavreJonsson} when
$k=2$. The proof requires
some self-intersection estimates for currents, due to Demailly-M\'eo.

\begin{theorem} \label{th_equi_hyp}
Let $f$ be an endomorphism of algebraic degree $d\geq 2$ of $\P^k$. Let $\Ec_m$ denote the union of 
the totally invariant proper analytic sets in $\P^k$ which are minimal, i.e. do not
contain smaller ones. Let $S$ be a positive closed $(1,1)$-current of
mass $1$ on $\P^k$ whose local potentials are not identically $-\infty$
on any component of $\Ec_m$. Then, $d^{-n}(f^n)^*(S)$ converge weakly
to the Green $(1,1)$-current $T$ of $f$. 
\end{theorem}

The following corollary gives a solution to the equidistribution
problem for hypersurfaces: the exceptional hypersurfaces belong to a proper
analytic set in the parameter space of hypersurfaces of a given degree. 

\begin{corollary}
Let $f$, $T$ and $\Ec_m$ be as above. If $H$ is a hypersurface of degree
$s$ in $\P^k$, which does not contain any component of $\Ec_m$, then
$s^{-1}d^{-n} (f^n)^*[H]$ converge to $T$ in the sense of currents.
\end{corollary}

Note that $(f^n)^*[H]$ is the current of integration on $f^{-n}(H)$
where the components of $f^{-n}(H)$ are counted with multiplicity.

\medskip
\noindent
{\bf Sketch of the proof of Theorem \ref{th_equi_hyp}.}
We can write $S=T+\ddc u$ where $u$ is a p.s.h. function modulo $T$,
that is, the difference of quasi-potentials of $S$ and of
$T$. Subtracting from $u$ a constant allows to assume that $\langle
\mu,u\rangle=0$. We call $u$ {\it the dynamical quasi-potential} of $S$.
Since $T$ has continuous quasi-potentials, $u$ satisfies analogous
properties that quasi-p.s.h. functions do. We are mostly concerned
with the singularities of $u$.
 
The total invariance of $T$ and $\mu$ implies that the dynamical
quasi-potential of $d^{-n}(f^n)^*(S)$ is equal to $u_n:=d^{-n}u\circ
f^n$. We have to show that this sequence of functions converges to 0
in $L^1(\P^k)$. Since $u$ is bounded from above, we have $\limsup
u_n\leq 0$. Assume that $u_n$ do not converge to 0. By Hartogs' lemma,
see Proposition \ref{prop_hartogs_else}, there is a ball $B$ and a
constant $\lambda>0$ such that $u_n\leq -\lambda$ on $B$ for infinitely
many indices $n$. It follows
that $u\leq -\lambda d^n$ on $f^n(B)$ for such an index $n$. On the other hand, the
exponential estimate in Theorem \ref{th_hormander} implies that
$\|e^{\alpha|u|}\|_{L^1}\leq A$ for some positive constants $\alpha$
and $A$ independent of $u$. If the multiplicity of $f$ at
every point is $\leq d-1$, then a version of Lojasiewicz's theorem
implies that $f^n(B)$ contains a ball of radius $\simeq e^{-c(d-1)^n}$,
$c>0$. Therefore, we have
$$e^{-2kc(d-1)^n}e^{\lambda d^n\alpha}\lesssim \int_{f^n(B)}
e^{\lambda d^n\alpha}\omega_\FS^k\leq \int_{\P^k} e^{\alpha|u|}\omega_\FS^k.$$
This contradicts the above exponential estimate.

In general, by Lemma \ref{lemma_push_ball} below, $f^n(B)$ contains always a ball of
radius $\simeq e^{-cd^n}$.
So, a slightly stronger version of 
the above exponential estimate will be enough to get a
contradiction.
We may improve this exponential estimate: if the Lelong numbers of $S$ are small, we can
increase the constant $\alpha$ and get a contradiction; if the Lelong numbers of $S_n$ are
small, we replace $S$ by $S_n$. 

The assumption $u_n<-\lambda d^n$ on $f^n(B)$ allows to show that all the limit currents of the sequence 
$d^{-n}(f^n)^*(S)$ have Lelong numbers larger than some constant
$\nu>0$. If $S'$ is such a current, there are other
currents $S'_n$ such that $S'=d^{-n}(f^n)^*(S_n')$. Indeed, if $S'$ is
the limit of $d^{-n_i}(f^{n_i})^*(S)$ one can take $S_n'$ a limit
value of $d^{-n_i+n}(f^{n_i-n})^*(S)$.

Let $a$ be a point such that $\nu(S_n',a)\geq \nu$.
The assumption on the potentials of $S$ allows to prove by induction on the dimension of
the totally invariant analytic sets that  $u_n$ converge to 0 on the maximal
totally invariant set $\Ec$. So, $a$ is out of
$\Ec$.
Lemma \ref{lemma_exceptional} allows to construct many distinct
points in
$f^{-n}(a)$. The identity $S'=d^{-n}(f^n)^*(S_n')$ implies 
an estimate from below of the Lelong numbers of $S'$ on $f^{-n}(a)$. 
This holds for every $n$.
Finally, this permits to construct analytic sets of large degrees 
on which we have estimates on the Lelong numbers of $S'$. Therefore, $S'$ has
a too large self-intersection. This 
contradicts an inequality of Demailly-M{\'e}o \cite{Demailly1, Meo3} and completes the
proof. Note that the proof of Demailly-M\'eo inequality uses H\"ormander's
$L^2$ estimates for the $\dbar$-equation. 
\hfill $\square$ 

\medskip

The following lemma is proved in \cite{DinhSibony9}. It also holds for
meromorphic maps. Some earlier versions were given in
\cite{FornaessSibony3} and in terms of Lebesgue measure in \cite{FavreJonsson,Guedj2}.

\begin{lemma} \label{lemma_push_ball}
There is a constant $c>0$ such that if $B$ is a ball of radius $r$,
$0<r<1$, in $\P^k$, then $f^n(B)$ contains a ball $B_n$ of radius
$\exp(-cr^{-2k}d^n)$ for any $n\geq 0$. 
\end{lemma} 

The ball $B_n$ is centered at $f^n(a_n)$ for some point $a_n\in B$
which is not necessarily the center of $B$. The key point in the proof
of the lemma is to find a point $a_n$ with an estimate from below on the
Jacobian of $f^n$ at $a_n$. If $u$ is a quasi-potential of the current
of integration on the critical set, the logarithm  of this Jacobian is
essentially the 
value of $u+u\circ f+\cdots + u\circ f^{n-1}$ at $a_n$. So, in order
to prove the existence of a point $a_n$ with a good estimate, it
is enough to bound the $L^1$ norm of the last function. One easily
obtains the result using the operator $f^*:\DSH(\P^k)\rightarrow
\DSH(\P^k)$ and its iterates, as it is done for $f_*$.

\begin{remark} \rm
Let $\Cc$ denote the convex compact set of totally invariant
positive closed $(1,1)$-currents of mass 1 on $\P^k$, i.e. currents $S$ such that $f^*(S)=dS$. Define an operator $\vee$ on $\Cc$.
If $S_1$, $S_2$ are elements of $\Cc$ and $u_1,u_2$ their dynamical
quasi-potentials, then 
$u_i\leq 0$. Since  $\langle \mu,u_i\rangle=0$ and $u_i$ are upper
semi-continuous, we deduce that $u_i=0$ on $\supp(\mu)$. Define
$S_1\vee S_2:=T+\ddc \max(u_1,u_2)$. It is easy to check that $S_1\vee
S_2$ is an
element of $\Cc$. An element $S$ is said to be {\it minimal} if
$S=S_1\vee S_2$ implies $S_1=S_2=S$. It is clear that $T$ is not minimal
if $\Cc$ contains other currents.  In fact, for $S$ in $\Cc$, we have
$T\vee S=T$. A current of integration on a totally
invariant hypersurface is a minimal element. It is likely that 
$\Cc$ is generated by a finite number of currents, the
operation $\vee$, convex hulls and limits.
\end{remark}

\begin{example} \label{example_Green_cv} \rm
If $f$ is the map given in Example \ref{example_power_map}, the exceptional set
$\Ec_m$ is the union of the $k+1$ attractive fixed points 
$$[0:\cdots:0:1:0:\cdots:0].$$
The convergence of $s^{-1}d^{-n}(f^n)^*[H]$ towards $T$ holds for hypersurfaces $H$ of
degree $s$ which do not contain
these points. If $\pi:\C^{k+1}\setminus\{0\}\rightarrow\P^k$ is the
canonical projection,
the Green $(1,1)$-current $T$ of $f$ is
given by $\pi^*(T)=\ddc (\max_i\log|z_i|)$, or equivalently
$T=\omega_\FS+\ddc v$ where 
$$v[z_0:\cdots:z_k]:=\max_{0\leq i\leq k}\log |z_i|
-\frac{1}{2}\log(|z_0|^2+\cdots+|z_k|^2).$$
The currents $T_i$ of integration on $\{z_i=0\}$ belong to $\Cc$ and 
 $T_i=T+\ddc u_i$ with $u_i:=\log|z_i|-\max_j\log|z_j|$. These
 currents are minimal. If
$\alpha_0$, $\ldots$, $\alpha_k$ are positive real numbers such that 
$\alpha:=1-\sum\alpha_i$ is positive, then $S:=\alpha T+\sum \alpha_i
T_i$ is an element of $\Cc$. We have $S=T+\ddc u$ with $u:=\sum
\alpha_i u_i$. The current $S$ is minimal if and only if 
$\alpha=0$. One can obtain other elements of $\Cc$ using the
operator $\vee$. We show that $\Cc$ is infinite dimensional. 
Define for $A:=(\alpha_0,\ldots,\alpha_k)$ with $0\leq \alpha_i\leq 1$ and $\sum
\alpha_i=1$ the p.s.h. function $v_A$ by 
$$v_A:=\sum \alpha_i\log|z_i|.$$
If $\Ac$ is a family of such $(k+1)$-tuples $A$, define
$$v_\Ac:=\sup_{A\in\Ac} v_A.$$
Then, we can define a 
positive closed $(1,1)$-currents $S_\Ac$ on $\P^k$ by
$\pi^*(S_\Ac)=\ddc v_\Ac$. It is clear that $S_\Ac$ belongs to $\Cc$
and hence $\Cc$ is of infinite dimension. 

\end{example}

The equidistribution problem in higher codimension is much more delicate
and is still open for general maps. 
We first recall the following lemma.

\begin{lemma} \label{lemma_f_generic}
For every $\delta>1$, there
is a Zariski dense open set  $\Hc_d^*(\P^k)$ in
  $\Hc_d(\P^k)$ and a constant $A>0$ such that 
  for $f$  in $\Hc_d^*(\P^k)$,
  the maximal multiplicity $\delta_n$ of $f^n$ at a point in $\P^k$ is at most
  equal to $A\delta^n$. In particular, the exceptional set of such a
  map $f$ is empty when $\delta<d$.
\end{lemma}
\proof
Let $X$ be a component of a totally invariant analytic set $E$ of pure dimension
$p\leq k-1$. Then, $f$ permutes the components of $E$. We deduce that $X$ is totally
invariant under $f^n$ for some $n\geq 1$. Lemma
\ref{lemma_topological} implies that the maximal multiplicity of $f^n$ at a
point in $X$ is at least equal to $d^{(k-p)n}$. Therefore, the second
assertion in the lemma is a consequence of the first one.

Fix an $N$ large enough such that $\delta^N>2^kk!$. 
Let  $\Hc_d^*(\P^k)$ be the set of $f$ such that $\delta_N\leq 2^kk!$.
This set is  Zariski open  in $\Hc_d(\P^k)$. 
Since the sequence $(\delta_n)$ is sub-multiplicative, i.e.
$\delta_{n+m}\leq\delta_n\delta_m$ for $n,m\geq 0$, 
if $f$ is in  $\Hc_d^*(\P^k)$, we have $\delta_N<\delta^N$, hence $\delta_n\leq A\delta^n$ for
$A$ large enough and for all $n$. It remains to show that 
   $\Hc_d^*(\P^k)$ is not empty.
Choose a
rational map
$h:\P^1\rightarrow \P^1$ of degree $d$ whose critical points are
simple and have disjoint infinite orbits.  
Observe that the multiplicity of $h^N$ at every point is at most equal to 2.
We construct the map $f$ using the method described in Example \ref{example_ueda}.
We have $f^N\circ \Pi=\Pi\circ
\widehat f^N$. 
Consider a point $x$ in $\P^k$ and a
point $\widehat x$ in $\Pi^{-1}(x)$. 
The multiplicity of $\widehat f^N$ at $\widehat x$ is at most equal to
$2^k$.
It follows that the multiplicity of $f^N$ at $x$ is at most equal to
$2^kk!$ since  $\Pi$ has degree $k!$. Therefore, $f$ satisfies the
desired inequality.
\endproof

We have the following result due to the authors \cite{DinhSibony10}.

\begin{theorem} \label{th_endo_generic}
There is a Zariski dense open set $\Hc_d^*(\P^k)$ in
  $\Hc_d(\P^k)$ such that if $f$ is in
  $\Hc_d^*(\P^k)$, then $d^{-pn}(f^n)^*(S)\rightarrow T^p$ uniformly
  on positive closed $(p,p)$-currents $S$ of mass $1$ on $\P^k$. 
In particular, the Green $(p,p)$-current $T^p$ is the unique positive
closed $(p,p)$-current of mass $1$ which is totally invariant. If $V$
is an analytic set of pure codimension $p$ and of degree $s$ in
$\P^k$, then $s^{-1}d^{-pn} (f^n)^*[V]$ converge to $T^p$ in the sense
of currents. 
\end{theorem}

\noindent
{\bf Sketch of proof.} The proof uses the super-potentials of
currents. In order to simplify the notation, introduce the dynamical
super-potential $\Vc$ of $S$. Define $\Vc:=\Uc_S-\Uc_{T^p}+c$ where $\Uc_S,\Uc_{T^p}$ are
super-potentials of $S,T^p$ and the constant $c$ is chosen so that
$\Vc(T^{k-p+1})=0$. Using a computation as in Theorem \ref{th_green_pp}, we
obtain that the dynamical super-potential of $d^{-pn} (f^n)^*(S)$ is
equal to $d^{-n}\Vc\circ\Lambda^n$ where
$\Lambda:\Cc_{k-p+1}(\P^k)\rightarrow \Cc_{k-p+1}(\P^k)$ is the
operator $d^{-p+1}f_*$. Observe that the dynamical super-potential of
$T^p$ is identically 0. In order to prove the convergence $d^{-pn}
(f^n)^*(S)\rightarrow T^p$, we only have to check that
$d^{-n}\Vc(\Lambda^n(R))\rightarrow 0$ for $R$ smooth in
$\Cc_{k-p+1}(\P^k)$. Since $T^p$ has a continuous super-potential,
$\Vc$ is bounded from above. Therefore, $\limsup
d^{-n}\Vc(\Lambda^n(R))\leq 0$.

Recall that $\Uc_S(R)=\Uc_R(S)$. So, in order to prove that $\liminf
d^{-n}\Vc(\Lambda^n(R))\geq 0$, it is enough to estimate 
$\inf_S\Uc_S(\Lambda^n(R))$, or equivalently, to estimate the capacity
of $\Lambda^n(R)$ from below. 
Assume in order to explain the idea that the support of $R$ is contained in a compact set
$K$ such that $f(K)\subset K$ and $K$ does not intersect the critical
set of $f$ (this is possible when $p=1$). We easily obtain that $\|\Lambda^n(R)\|_{\infty}\lesssim
A^n$ for some constant $A>0$. The estimate in Theorem \ref{th_hormander_sp} implies
the result. In the general case, if  $\Hc_d^*(\P^k)$ is chosen as in
Lemma \ref{lemma_f_generic} for $\delta$ small enough and if $f$ is in
$\Hc_d^*(\P^k)$, we can prove the estimate
$\capacity(\Lambda^n(R))\lesssim d'^n$ for any fixed constant $d'$ such
that $1<d'<d$. This
implies the desired convergence of super-potentials. The main
technical difficulty is that when $R$ hits the critical set, then
$\Lambda(R)$ is not bounded. The estimates requires a smoothing and
precise evaluation of the error. 
\hfill $\square$

\begin{remark} \label{rk_unif_cv_perron}
\rm
The above estimate on $\capacity(\Lambda^n(R))$ can be seen as a version of
Lojasiewicz's inequality for currents. It is quite delicate to obtain.
We also have an explicit estimate on the speed of convergence. Indeed,
we have for an appropriate $d'<d$: 
$$\dist_2\big(d^{-pn}(f^n)^*(S),T^p\big):=\sup_{\|\Phi\|_{\Cc^2}\leq 1}
|\langle d^{-pn}(f^n)^*(S)-T^p, \Phi\rangle|\lesssim d'^n d^{-n}.$$
The theory of interpolation between Banach spaces \cite{Triebel} implies
a similar estimate for $\Phi$ H\"older continuous.
\end{remark}

\bigskip\bigskip

\begin{exercise}
If $f:=[z_0^d:\cdots:z_k^d]$, show that $\{z_1^p=z_2^q\}$, for
arbitrary $p,q$, is invariant
under $f$. Show that a curve invariant under an endomorphism is an
image of $\P^1$ or a torus, possibly singular.
\end{exercise}

\begin{exercise}
Let $f$ be as in Example \ref{example_Green_cv}. Let $S$ be a $(p,p)$-current with positive
Lelong number at $[1:0:\cdots:0]$. Show that  any limit of $d^{-pn}
(f^n)^*(S)$ has a strictly positive Lelong number at $[1:0:\cdots:0]$ and deduce that
$d^{-pn} (f^n)^*(S)$ 
do not converge to $T^p$.
\end{exercise}

\begin{exercise} Let $f$ be as in Theorem \ref{th_endo_generic} for
  $p=k$ and $\Lambda$ the associated Perron-Frobenius operator. 
If $\varphi$ is a $\Cc^2$ function on $\P^k$, show that
$$\|\Lambda^n(\varphi)-\langle\mu,\varphi\rangle\|_\infty\leq
cd'^nd^{-n}$$ 
for some constant $c>0$. Deduce that $\Lambda^n(\varphi)$ converge
uniformly to $\langle\mu,\varphi\rangle$. Give an estimate of
$\|\Lambda^n(\varphi)-\langle\mu,\varphi\rangle\|_\infty$ 
for $\varphi$ H\"older continuous.
\end{exercise}

\begin{exercise} 
Let $f$ be an endomorphism of algebraic degree $d\geq 2$. Assume that $V$
is a  totally invariant hypersurface,
i.e.
$f^{-1}(V)=V$. Let $V_i$  denote the irreducible components of $V$ and
$h_i$ minimal homogeneous polynomials such that 
$V_i=\{h_i=0\}$. Define $h=\prod h_i$.
Show that $h\circ f=ch^d$ where $c$ is a constant.
If $F$ is a lift of $f$ to $\C^{k+1}$, prove that $\Jac(F)$ contains $(\prod
      h_i)^{d-1}$ as a factor.
Show that $V$ is contained in the critical set of $f$ and deduce\footnote{It is known that in dimension $k=2$, $V$
        is a union of at most 3 lines, \cite{CerveauLinsNeto,FornaessSibony7,SSU}.} that
      $\deg V\leq k+1$.  
Assume now that $V$ is reducible. Find a totally invariant positive closed
$(1,1)$-current of mass $1$ which
      is not the Green current nor
the current associated to an analytic set.
\end{exercise}

\begin{exercise}
Let $u$ be a p.s.h. function in $\C^k$, such that for
$\lambda\in\C^*$, $u(\lambda z)=\log|\lambda|+u(z)$. If $\{u<0\}$ is
bounded in $\C^k$, show that
$\ddc u^+$ is a positive closed current on $\P^k$ which is extremal in the cone of positive closed
$(1,1)$-currents\footnote{Unpublished result by
  Berndtsson-Sibony.}. Deduce that the the Green $(1,1)$-current of a
polynomial map of $\C^k$ which extends holomorphically to $\P^k$, is extremal. 
\end{exercise} 

\begin{exercise} Let $v$ be a subharmonic function on $\C$. Suppose
  $v(e^{i\theta}z)=v(z)$ for every $z\in\C$ and for every $\theta\in\R$ such that
  $e^{i\theta d^n}=1$ for some integer $n$. Prove that
  $v(z)=v(|z|)$ for $z\in\C$. Hint: use the
  Laplacian of $v$. Let $f$ be as in Example \ref{example_Green_cv}, $R$ a current in
  $\C^{k+1}$, and $v$ a p.s.h. function on $\C^{k+1}$ such that $R=\ddc v$
  and $v(F(z))=dv(z)$, where $F(z):=(z_0^d,\ldots,z_k^d)$ is a lift of
  $f$ to $\C^{k+1}$. Show that $v$ is invariant under the action of
  the unit torus $\T^{k+1}$ in $\C^{k+1}$. Determine such functions
  $v$. Recall that $\T$ is the unit circle in $\C$ and $\T^{k+1}$ acts
  on $\C^{k+1}$ by multiplication.
\end{exercise}

\begin{exercise} 
Define the Desboves map $f_0$ in $\Mc_4(\P^2)$ as  
$$f_0[z_0:z_1:z_2]:=[z_0(z_1^3-z_2^3):z_1(z_2^3-z_0^3):z_2(z_0^3-z_1^3)].$$
Prove that $f_0$ has $12$ indeterminacy points. If $\sigma$ is a
      permutation of coordinates, compare $f_0\circ \sigma$ and
      $\sigma\circ f_0$.
Define  
$$\Phi_\lambda(z_0,z_1,z_2):=z_0^3+z_1^3+z_2^3-3\lambda z_0z_1z_2,\quad \lambda\in\C$$
and 
$$L[z_0:z_1:z_2]:=[az_0:bz_1:cz_2], \quad a,b,c\in\C.$$
Show that for Zariski  generic $L$, 
$f_L:=f_0+\Phi_\lambda L$ is in 
$\Hc_4(\P^2)$. Show that on the curve $\{\Phi_\lambda=0\}$ in $\P^2$,
$f_L$ coincides with $f_0$, and that  $f_0$ maps the cubic
$\{\Phi_\lambda=0\}$ onto itself.\footnote{This example was considered in
  \cite{BonifantDabijaMilnor}. It gives maps in  $\Hc_4(\P^2)$ which preserves a cubic.
The cubic is singular if $\lambda=1$, non singular if $\lambda\not=1$. 
In higher dimension, Beauville proved that a smooth hypersurface of
$\P^k$, $k\geq 3$, of degree $>1$ does not have an endomorphism with
$d_t>1$, unless the degree is 2, $k=3$ and the hypersurface is
isomorphic to $\P^1\times\P^1$ \cite{Beauville}.} 
\end{exercise}


\section{Stochastic properties of the Green measure}

In this paragraph, we are concerned with the stochastic properties of 
the equilibrium measure $\mu$
associated to an endomorphism $f$. 
If $\varphi$ is an observable, $(\varphi\circ f^n)_{n\geq 0}$ can
be seen as a sequence of dependent random variables. Since the measure
is invariant, these variables are identically distributed,
i.e. the Borel sets $\{\varphi\circ f^n<t\}$ have the same $\mu$
measure for any fixed constant $t$. The idea is to show that the
dependence is weak and then to extend classical results in probability
theory to our setting. One of the key point is
the spectral study of the Perron-Frobenius operator $\Lambda:=d^{-k}f_*$. It
allows to prove the exponential decay of correlations for d.s.h. and
H\"older continuous observables, the
central limit theorem, the large deviation theorem, etc. An important
point is to use the space of d.s.h. functions as a space of
observables. For the
reader's convenience, we  recall few general facts
from ergodic theory and probability theory. We refer to
\cite{KatokHasselblatt, Walters} 
for the general theory.

Consider a dynamical system associated to a map $g:X\rightarrow X$
which is measurable with respect to a $\sigma$-algebra $\Fc$ on $X$. 
The direct image of a probability measure $\nu$ by $g$ is the probability
measure $g_*(\nu)$ defined by
$$g_*(\nu)(A):=\nu(g^{-1}(A))$$
for every measurable set $A$. Equivalently, for any positive 
measurable function $\varphi$, we have
$$\langle g_*(\nu),\varphi\rangle:=\langle \nu, \varphi\circ g\rangle.$$
The measure $\nu$ is {\it invariant} if $g_*(\nu)=\nu$. 
When $X$ is a compact metric space and $g$ is continuous, the set $\Mc(g)$ of invariant probability
measures is convex, compact
and non-empty: for any sequence of probability measures $\nu_N$, the cluster points of 
$${1\over N}\sum_{j=0}^{N-1}(g^n)_*(\nu_N)$$
are invariant probability measures.

A measurable set $A$ is {\it totally invariant} if $\nu(A\setminus
g^{-1}(A))=\nu(g^{-1}(A)\setminus A)=0$. An invariant probability
measure $\nu$ is {\it ergodic} if any totally invariant set is of zero
or full $\nu$-measure. It is easy to show that $\nu$ is ergodic if and
only if when $\varphi\circ g=\varphi$, for $\varphi\in L^1(\nu)$, then
$\varphi$ is constant. Here, we can replace $L^1(\nu)$ by $L^p(\nu)$
with $1\leq p\leq+\infty$. 
The ergodicity of $\nu$ is also equivalent to the fact that it is
extremal in $\Mc(g)$.
We recall Birkhoff's ergodic theorem,
which is the analogue of the law of large
numbers for independent random variables \cite{Walters}.

\begin{theorem}[Birkhoff] Let $g:X\rightarrow X$ be a measurable map as above.
Assume that $\nu$ is an ergodic invariant
  probability measure. Let $\varphi$ be a function in $L^1(\nu)$. Then
$${1\over N} \sum_{n=0}^{N-1}\varphi(g^n(x))\rightarrow \langle\nu,\varphi\rangle$$
almost everywhere with respect to $\nu$.
\end{theorem}

When $X$ is a compact metric space, we can apply Birkhoff's theorem to continuous
functions $\varphi$ and deduce that 
for $\nu$ almost every $x$
$${1\over N} \sum_{n=0}^{N-1} \delta_{g^n(x)}\rightarrow \nu,$$
where $\delta_x$ denotes the Dirac mass at $x$. 
The sum 
$$\St_N(\varphi):=\sum_{n=0}^{N-1}\varphi\circ g^n$$ 
is called {\it Birkhoff's sum}. 
So, Birkhoff's theorem describes the behavior of ${1\over N} \St_N(\varphi)$
for an observable $\varphi$.  We will be concerned with
the precise behavior of $\St_N(\varphi)$ for various classes of
functions $\varphi$. 

A notion stronger than ergodicity is the notion of mixing. An invariant
probability measure $\nu$ is {\it mixing} if for every measurable sets $A,B$
$$\lim_{n\rightarrow\infty} \nu(g^{-n}(A)\cap B)=\nu(A)\nu(B).$$
Clearly, mixing implies ergodicity.
It is not difficult to see that $\nu$ is mixing if and only if for any
test functions $\varphi,\psi$ in $L^\infty(\nu)$ or in $L^2(\nu)$, we
have
$$\lim_{n\rightarrow\infty} \langle\nu,(\varphi\circ g^n)\psi\rangle
=\langle\nu,\varphi\rangle\langle\nu,\psi\rangle.$$
The quantity
$$I_n(\varphi,\psi):=|\langle\nu,(\varphi\circ g^n)\psi\rangle
-\langle\nu,\varphi\rangle\langle\nu,\psi\rangle|$$
is called {\it the correlation at time $n$} of $\varphi$ and
$\psi$. So, mixing is equivalent to the convergence of
$I_n(\varphi,\psi)$ to 0. We say that $\nu$ is {\it K-mixing} if for every $\psi\in L^2(\nu)$
$$\sup_{\|\varphi\|_{L^2(\nu)}\leq 1} I_n(\varphi,\psi)\rightarrow
0.$$
Note that K-mixing is equivalent to the
fact that the $\sigma$-algebra $\Fc_\infty:=\cap g^{-n}(\Fc)$ contains
only sets of zero and full measure. 
This is the strongest form of mixing for observables in $L^2(\nu)$.
It is however of interest to get a quantitative information on the
mixing speed for more regular
observables like smooth or H\"older continuous functions.

Consider now an endomorphism $f$ of algebraic degree $d\geq 2$ of
$\P^k$ as above and its equilibrium measure $\mu$. We know that $\mu$
is totally invariant: $f^*(\mu)=d^k\mu$. If $\varphi$ is a continuous
function, then 
$$\langle \mu,\varphi\circ f\rangle =\langle
d^{-k}f^*(\mu),\varphi\circ f\rangle = \langle \mu, d^{-k}
f_*(\varphi\circ f)\rangle =\langle\mu,\varphi\rangle.$$
We have used the obvious fact that
$f_*(\varphi\circ f)=d^k\varphi$. 
So, $\mu$ is invariant. We have the following proposition.

\begin{proposition}
The Perron-Frobenius
operator $\Lambda:=d^{-k}f_*$ has a continuous extension of norm $1$ to
$L^2(\mu)$. Moreover, the adjoint of $\Lambda$ satisfies 
$\trans\Lambda(\varphi)=\varphi\circ f$ and
$\Lambda\circ \trans\Lambda=\id$. Let $L_0^2(\mu)$ denote the
hyperplane of $L^2(\mu)$ defined by $\langle
\mu,\varphi\rangle=0$. Then, the spectral radius of $\Lambda$
on $L_0^2(\mu)$ is also
equal to $1$.
\end{proposition}
\proof
Schwarz's inequality implies that 
$$|f_*(\varphi)|^2\leq d^k f_*(|\varphi|^2).$$
Using the total invariance of $\mu$, we get
$$\langle  \mu, |\Lambda(\varphi)|^2\rangle \leq \langle
\mu,\Lambda(|\varphi|^2)\rangle = \langle \mu, |\varphi|^2\rangle.$$
Therefore, $\Lambda$ has a continuous extension to $L^2(\mu)$, with
norm $\leq 1$. Since $\Lambda(1)=1$, the norm of this operator is
equal to 1. The properties on the adjoint of $\Lambda$ are easily deduced from the total
invariance of $\mu$.

Let $\varphi$ be a function in $L_0^2(\mu)$ of norm 1. Then,
$\varphi\circ f^n$ is also in $L^2(\mu)$ and of norm 1. Moreover,
$\Lambda^n(\varphi\circ f^n)$, which is equal to $\varphi$, is of norm
1. So,  the spectral radius of $\Lambda$
on $L_0^2(\mu)$ is also
equal to $1$.
\endproof

Mixing for the measure $\mu$ was proved in \cite{FornaessSibony1}. We
give in this paragraph two proofs of K-mixing. The first one is from \cite{DinhSibony1}
and does not use that $\mu$ is moderate.

\begin{theorem} \label{th_K-mixing_end}
Let $f$ be an endomorphism of algebraic degree $d\geq 2$ of
$\P^k$. Then, its Green measure $\mu$ is K-mixing.
\end{theorem}
\proof
Let $c_\psi:=\langle \mu,\psi\rangle$. Since $\mu$ is totally
invariant, the correlations between two observables $\varphi$ and
$\psi$ satisfy
\begin{eqnarray*}
I_n(\varphi,\psi) & = & |\langle
\mu,(\varphi\circ f^n)\psi\rangle
-\langle\mu,\varphi\rangle\langle\mu,\psi\rangle|\\
& = &  |\langle
\mu,\varphi\Lambda^n(\psi)\rangle
-c_\psi\langle\mu,\varphi\rangle|\\
& = & |\langle
\mu,\varphi(\Lambda^n(\psi)-c_\psi)\rangle|.
\end{eqnarray*}
Hence,
$$\sup_{\|\varphi\|_{L^2(\mu)}\leq 1}I_n(\varphi,\psi)\leq
  \|\Lambda^n(\psi)-c_\psi\|_{L^2(\mu)}.$$
Since $\|\Lambda\|_{L^2(\mu)}\leq 1$, it is enough to show that
$\|\Lambda^n(\psi)-c_\psi\|_{L^2(\mu)}\rightarrow 0$ for a dense
family of functions $\psi\in L^2(\mu)$. 
So, we can assume that $\psi$ is a d.s.h. function such that $|\psi|\leq 1$. 
We have $|c_\psi|\leq 1$ and 
$\|\Lambda^n(\psi)-c_\psi\|_{L^\infty(\mu)} \leq 2$. Since $\mu$ is
PB, we deduce from Theorem \ref{th_green_measure_dsh} and
Cauchy-Schwarz's inequality that
$$\|\Lambda^n(\psi)-c_\psi\|_{L^2(\mu)}\lesssim
\|\Lambda^n(\psi)-c_\psi\|_{L^1(\mu)}^{1/2}
\lesssim \|\Lambda^n(\psi)-c_\psi\|_\DSH^{1/2}\lesssim d^{-n/2}.$$
This completes the proof. For the last argument, we can also use continuous test functions
$\psi$ and apply Proposition \ref{prop_except_FS}. We then obtain that
$\Lambda^n(\psi)-c_\psi$ converges to 0 pointwise out of a pluripolar set. 
Lebesgue's convergence theorem and the fact that $\mu$ has no mass on
pluripolar sets imply the result. 
\endproof

In what follows, we show that
the equilibrium measure $\mu$ satisfies remarkable stochastic properties
which are quite hard to obtain in the setting of real dynamical
systems. We will see the effectiveness of the pluripotential methods
which replace the delicate estimates, used in some real dynamical systems.
The following result was recently obtained by Nguyen and the
authors \cite{DinhNguyenSibony3}. It shows that the equilibrium measure is exponentially mixing
and generalizes earlier results of \cite{DinhSibony1, DinhSibony6,FornaessSibony1}.
Note that d.s.h. observables may be everywhere discontinuous.

\begin{theorem} \label{th_mixing}
Let $f$ be a holomorphic endomorphism of algebraic degree $d\geq 2$ on
$\P^k$. Let $\mu$ be the Green measure of $f$. Then for every $1< p\leq +\infty$
there is a constant $c>0$ such that
$$|\langle \mu, (\varphi\circ f^n)\psi\rangle
-\langle\mu,\varphi\rangle \langle\mu,\psi\rangle|\leq c
d^{-n}\|\varphi\|_{L^p(\mu)}\|\psi\|_\DSH$$
for $n\geq 0$, $\varphi$ in $L^p(\mu)$ and $\psi$
d.s.h. Moreover, for $0\leq \nu\leq 2$ there is a constant $c>0$
such that 
$$|\langle \mu, (\varphi\circ f^n)\psi\rangle
-\langle\mu,\varphi\rangle \langle\mu,\psi\rangle|\leq c
d^{-n\nu/2}\|\varphi\|_{L^p(\mu)}\|\psi\|_{\Cc^\nu}$$
for $n\geq 0$, $\varphi$ in $L^p(\mu)$ and $\psi$
of class $\Cc^\nu$.
\end{theorem}
\proof
We prove the first assertion. 
Observe that the correlations
$$I_n(\varphi,\psi):=|\langle \mu, (\varphi\circ f^n)\psi\rangle
-\langle\mu,\varphi\rangle \langle\mu,\psi\rangle|$$
vanish if $\psi$ is constant. Therefore, we can assume that $\langle
\mu,\psi\rangle=0$. In which case, we have 
$$I_n(\varphi,\psi)=|\langle \mu,\varphi\Lambda^n(\psi)\rangle|,$$
where  $\Lambda$ denotes the Perron-Frobenius operator associated to
$f$.

We can also assume that $\|\psi\|_\DSH\leq
1$. Corollary \ref{cor_exp_dsh} implies that for $1\leq
q<\infty$,  
$$\|\Lambda^n(\psi)\|_{L^q(\mu)}\leq cq d^{-n}$$ 
where $c>0$ is  a constant independent of $n,q$ and $\psi$.    
Now, if $q$ is chosen so that $p^{-1}+q^{-1}=1$, we obtain using 
H{\"o}lder's inequality that
$$I_n(\varphi,\psi)\leq \|\varphi\|_{L^p(\mu)}
\|\Lambda^n(\psi)\|_{L^q(\mu)} \leq cq \|\varphi\|_{L^p(\mu)}
d^{-n}.$$
This completes the proof of the first assertion. The second assertion
is proved in the same way using Corollary \ref{cor_exp_holder}. 
\endproof

Observe that the above estimates imply that for $\psi$ smooth
$$\lim_{n\rightarrow \infty}\sup_{\|\varphi\|_{L^2(\mu)}\leq 1}
  I_n(\varphi,\psi) = 0.$$
Since smooth functions are dense in $L^2(\mu)$, the convergence
holds for every $\psi$ in $L^2(\mu)$ and gives another proof of the K-mixing.
The following result \cite{DinhNguyenSibony3} gives estimates for the
exponential mixing of any order. It can be extended to 
H{\"o}lder continuous observables using the second assertion in
Theorem \ref{th_mixing}.

\begin{theorem} \label{th_mixing_order}
Let $f$, $d$, $\mu$ be as in Theorem \ref{th_mixing} and $r\geq 1$
an integer. Then there is a constant $c>0$ such that
$$\Big|\langle \mu, \psi_0(\psi_1\circ f^{n_1})\ldots (\psi_r\circ f^{n_r})\rangle
-\prod_{i=0}^r\langle\mu,\psi_i\rangle \Big|\leq c
d^{-n}\prod_{i=0}^r\|\psi_i\|_\DSH$$
for $0=n_0\leq n_1\leq \cdots\leq n_r$, $n:=\min_{0\leq i< r} (n_{i+1}-n_i)$ and $\psi_i$ 
d.s.h. 
\end{theorem}
\proof
The proof is by induction on $r$. The case $r=1$ is a consequence of Theorem
\ref{th_mixing}. Suppose the result is true for $r-1$, we have to
check it for $r$. Without loss of generality, assume that
$\|\psi_i\|_\DSH\leq 1$. This implies that $m:=\langle \mu,\psi_0\rangle$ is
bounded. The invariance of $\mu$ and the induction hypothesis imply that
\begin{eqnarray*}
\lefteqn{\Big|\langle \mu, m(\psi_1\circ f^{n_1})\ldots (\psi_r\circ f^{n_r})\rangle
-\prod_{i=0}^r\langle\mu,\psi_i\rangle \Big|}\\
&=&  \Big|\langle \mu, m\psi_1(\psi_2\circ f^{n_2-n_1})\ldots (\psi_r\circ f^{n_r-n_1})\rangle
-m\prod_{i=1}^r\langle\mu,\psi_i\rangle \Big|\leq  c d^{-n}
\end{eqnarray*}
for some constant $c>0$. In order to get the desired estimate, it is
enough to show that
$$\Big|\langle \mu, (\psi_0-m)(\psi_1\circ f^{n_1})\ldots 
(\psi_r\circ f^{n_r})\rangle\Big|\leq c d^{-n}.$$
Observe that the operator $(f^n)^*$ acts on $L^p(\mu)$ for $p\geq 1$
and its norm is bounded by 1. Using the invariance of $\mu$ and 
H{\"o}lder's inequality, we get for $p:=r+1$
\begin{eqnarray*}
\lefteqn{\Big|\langle \mu, (\psi_0-m)(\psi_1\circ f^{n_1})\ldots 
(\psi_r\circ f^{n_r})\rangle\Big| }\\
&=& \Big|\langle \mu, \Lambda^{n_1}(\psi_0-m)\psi_1\ldots 
(\psi_r\circ f^{n_r-n_1})\rangle\Big| \\
&\leq& \|\Lambda^{n_1}(\psi_0-m)\|_{L^p(\mu)}\|\psi_1\|_{L^p(\mu)}\ldots
\|\psi_r\circ f^{n_r-n_1}\|_{L^p(\mu)}\\
& \leq & cd^{-n_1}\|\psi_1\|_{L^p(\mu)}\ldots\|\psi_r\|_{L^p(\mu)},
\end{eqnarray*}
for some constant $c>0$.
Since $\|\psi_i\|_{L^p(\mu)}\lesssim \|\psi_i\|_\DSH$, the previous
estimates imply the result. Note that as in Theorem \ref{th_mixing},
it is enough to assume that $\psi_i$ is d.s.h. for $i\leq r-1$ and
$\psi_r$ is in $L^p(\mu)$ for some $p>1$. 
\endproof

The mixing of $\mu$ implies that for any measurable observable
$\varphi$,  
the times series $\varphi\circ f^n$, behaves like independent
random variables with the same distribution. For example, the dependence of $\varphi\circ f^n$ and
$\varphi$ is weak when $n$ is large: if
$a,b$ are real numbers, then the measure of $\{\varphi\circ
f^n\leq a\ \mbox{and}\ \varphi\leq b\}$ is almost equal to 
$\mu\{\varphi\circ f^n\leq a\}\mu\{\varphi\leq b\}$. Indeed, it is equal to
$$\big\langle \mu, (\ind_{]-\infty,a]}\circ \varphi\circ f^n)(\ind_{]-\infty,b]}\circ\varphi)\big\rangle,$$
and when $n$ is large, mixing implies that the last integral is approximatively equal to 
$$\langle \mu, \ind_{]-\infty,a]}\circ \varphi\rangle
\langle \mu, \ind_{]-\infty,b]}\circ\varphi\rangle=\mu\{\varphi\leq
a\}\mu\{\varphi\leq b\}=\mu\{\varphi\circ f^n\leq a\}\mu\{\varphi\leq
b\}.$$

The estimates on the decay of correlations obtained in the above
results, give at
which speed the observables become ``almost independent''. 
We are going to
show that under weak assumptions on the regularity of observables
$\varphi$, the times series $\varphi\circ f^n$, satisfies the Central
Limit Theorem (CLT for short). We recall the classical CLT for
independent random variables. In what follows, $\Et(\cdot)$ denotes
expectation, i.e. the mean, of a
random variable.

\begin{theorem}
Let $(X,\Fc,\nu)$ be a probability space. Let $Z_1,Z_2,\ldots$
be independent  identically distributed  
(i.i.d. for short) random variables with values in
$\R$, and of mean zero, i.e. $\Et(Z_n)=0$. Assume also that 
$0<\Et(Z_n^2)<\infty$. Then for any
open interval $I\subset \R$ and for $\sigma:=\Et(Z_n^2)^{1/2}$, we have
$$\lim_{N\rightarrow\infty} \nu\Big\{{1\over\sqrt{N}} \sum_{n=0}^{N-1}
Z_n\in I\Big\}={1\over \sqrt{2\pi}\sigma} \int _I e^{-{t^2\over
    2\sigma^2}}dt.$$
\end{theorem}

The important hypothesis here is that the variables 
have the same distribution (this means that for every interval
$I\subset \R$, $\langle\nu, \ind_I\circ Z_n\rangle$ is independent of $n$, where $\ind_I$ is the
characteristic function of $I$) and that they are independent. The result can be phrased as
follows. If we define the random variables $\widehat Z_N$ by 
$$\widehat Z_N:={1\over\sqrt{N}} \sum_{n=0}^{N-1} Z_n,$$
then the sequence of probability measures $(\widehat Z_N)_*(\nu)$
on $\R$ converges to the probability measure of density ${1\over \sqrt{2\pi}\sigma}e^{-{t^2\over
    2\sigma^2}}$. This is also called {\it the convergence in law}. 

We want to replace the random variables $Z_n$ by the functions
$\varphi\circ f^n$ on 
the probability space $(\P^k,\Bc,\mu)$ where $\Bc$ is the Borel $\sigma$-algebra.
The fact that $\mu$ is
invariant means exactly that $\varphi\circ f^n$ are identically
distributed.
We state first a central limit theorem due to Gordin
\cite{Gordin}, see also \cite{Viana}. 
For simplicity, we consider a measurable map $g:(X,\Fc)\rightarrow (X,\Fc)$ 
as above. 
Define $\Fc_n:=g^{-n}(\Fc)$, $n\geq 0$, the $\sigma$-algebra of
$g^{-n}(A)$, with $A\in\Fc$. This sequence is non-increasing.
Denote by $\Et(\varphi|\Fc')$
the conditional expectation of $\varphi$ with respect to a
$\sigma$-algebra $\Fc'\subset\Fc$. We say that $\varphi$ is {\it a
  coboundary} if  $\varphi=\psi\circ g-\psi$ for some function $\psi\in L^2(\nu)$.

\begin{theorem}[Gordin] \label{th_gordin} 
Let $\nu$ be an ergodic invariant probability
measure on $X$. Let $\varphi$ be an observable in $L^1(\nu)$ such that
$\langle \nu,\varphi\rangle=0$. Suppose
$$\sum_{n\geq 0} \|\Et(\varphi|\Fc_n)\|_{L^2(\nu)}^2<\infty.$$ 
Then $\langle \nu,\varphi^2\rangle+2\sum_{n\geq 1}
\langle \nu,\varphi(\varphi\circ g^n)\rangle$ is a finite positive
number which vanishes if and
only if $\varphi$ is a coboundary. Moreover, if
$$\sigma:=\Big[\langle \nu,\varphi^2\rangle+2\sum_{n\geq 1}
\langle \nu,\varphi(\varphi\circ g^n)\rangle\Big]^{1/2}$$
is strictly positive, then $\varphi$ satisfies the central limit
theorem 
with variance $\sigma$: for any interval $I\subset \R$
$$\lim_{N\rightarrow\infty} \nu\Big\{{1\over\sqrt{N}} \sum_{n=0}^{N-1}
\varphi\circ g^n\in I\Big\}={1\over \sqrt{2\pi}\sigma} \int _I e^{-{t^2\over
    2\sigma^2}}dt.$$
\end{theorem}

It is not difficult to see that a function $u$ is
$\Fc_n$-measurable if and only if  $u=v\circ g^n$ with $v$
$\Fc$-measurable. Let $L^2(\nu,\Fc_n)$ denote the space of
$\Fc_n$-measurable functions which are in $L^2(\nu)$. Then,
$\Et(\varphi|\Fc_n)$ is the orthogonal projection
of $\varphi\in L^2(\nu)$ into $L^2(\nu,\Fc_n)$. 

A straighforward
computation using the invariance of $\nu$ gives that the variance
$\sigma$ in the above theorem is equal to 
$$\sigma=\lim_{n\rightarrow\infty}
n^{-1/2}\|\varphi+\cdots+\varphi\circ g^{n-1}\|_{L^2(\nu)}.$$
When $\varphi$ is orthogonal to all $\varphi\circ g^n$, we find that
$\sigma=\|\varphi\|_{L^2(\mu)}$. So, Gordin's theorem assumes a weak
dependence and concludes that the observables satisfy the central
limit theorem.

Consider now the dynamical system associated to an endomorphism $f$ of
$\P^k$ as above. 
Let $\Bc$ denote the Borel $\sigma$-algebra on $\P^k$ and define
$\Bc_n:=f^{-n}(\Bc)$. 
Since the measure $\mu$ satisfies $f^*(\mu)=d^k\mu$,
the norms
$\|\Et(\cdot|\Bc_n)\|_{L^2(\mu)}$ can be expressed in terms of the
operator $\Lambda$. We have the following lemma.

\begin{lemma} \label{lemma_pk_experance}
Let $\varphi$ be an observable in $L^2(\mu)$. Then
$$\Et(\varphi|\Fc_n)=\Lambda^n(\varphi)\circ f^n \quad\mbox{and}\quad 
\|\Et(\varphi|\Bc_n)\|_{L^p(\mu)}=\|\Lambda^n(\varphi)\|_{L^p(\mu)},$$
for $1\leq p\leq 2$.
\end{lemma}
\proof
We have 
\begin{eqnarray*}
\big\langle \mu,\varphi(\psi\circ f^n)\big\rangle  & = &  \big\langle d^{-kn}(f^n)^*(\mu),
\varphi(\psi\circ f^n)\big\rangle 
 =  \big\langle\mu,
d^{-kn} (f^n)_*[\varphi(\psi\circ f^n)]\big\rangle \\
& = & 
\big\langle\mu,\Lambda^n(\varphi)\psi\big\rangle 
 =  \big\langle\mu, [\Lambda^n(\varphi)\circ f^n][\psi\circ f^n]\big\rangle.
\end{eqnarray*}
This proves the first assertion.
The invariance of $\mu$ implies that 
$\|\psi\circ f^n\|_{L^p(\mu)}=\|\psi\|_{L^p(\mu)}$.
Therefore, the second assertion is a consequence of the first one.
\endproof

Gordin's Theorem \ref{th_gordin}, Corollaries
\ref{cor_exp_dsh} and \ref{cor_exp_holder}, applied to $q=2$, give the following result. 

\begin{corollary} \label{cor_clt_end}
Let $f$ be an endomorphism of algebraic degree $d\geq 2$ of $\P^k$ and
$\mu$ its equilibrium measure. Let $\varphi$ be
a d.s.h. function or a H\"older continuous function on
$\P^k$, such that $\langle\mu, \varphi\rangle=0$. Assume that $\varphi$ is not a
coboundary. Then $\varphi$ satisfies the central limit theorem with
the variance $\sigma>0$ given by 
$$\sigma^2:=
\langle\mu,\varphi^2\rangle +2 \sum_{n\geq 1} \langle \mu,\varphi(\varphi\circ f^n)\rangle.$$
\end{corollary}

We give an interesting decomposition of the space $L_0^2(\mu)$ which
shows that $\Lambda$, acts like a ``generalized shift''. Recall that
$L_0^2(\mu)$ is the space of functions $\psi\in L^2(\mu)$ such that
$\langle\mu,\psi\rangle=0$. Corollary \ref{cor_clt_end} can also be deduced from the
following result.

\begin{proposition}
Let $f$ be an endomorphism of algebraic degree $d\geq 2$ of $\P^k$ and $\mu$ 
the corresponding equilibrium measure. Define 
$$V_n:=\{\psi\in L_0^2(\mu),\ \Lambda^n(\psi)=0\}.$$ 
Then, we have
$V_{n+1}=V_n\oplus V_1\circ f^n$ as an orthogonal sum and
$L_0^2(\mu)=\oplus_{n=0}^\infty V_1\circ f^n$ as a Hilbert sum. Let
$\psi=\sum \psi_n\circ f^n$, with $\psi_n\in V_1$, be a function
in $L_0^2(\mu)$. Then, $\psi$
satisfies the Gordin's condition, see Theorem \ref{th_gordin}, if and
only if the sum $\sum_{n\geq 1}
n\|\psi_n\|_{L^2(\mu)}^2$ is finite. Moreover, if $\psi$ is d.s.h. (resp. of class $\Cc^\nu$,
$0<\nu\leq 2$) with $\langle\mu,\psi\rangle=0$, then
$\|\psi_n\|_{L^2(\mu)}\lesssim d^{-n}$
(resp. $\|\psi_n\|_{L^2(\mu)}\lesssim d^{-n\nu/2}$). 
\end{proposition}
\proof
It is easy to check that $V_n^\perp=\{\theta\circ f^n,\ \theta\in
L^2(\mu)\}$. Let $W_{n+1}$ denote the orthogonal complement of $V_n$
in $V_{n+1}$. Suppose $\theta\circ f^n$ is in $W_{n+1}$. Then,
$\Lambda(\theta)=0$. This gives the first decomposition in the
proposition.

For the second decomposition, observe that $\oplus_{n=0}^\infty
V_1\circ f^n$ is a direct orthogonal sum. We only have to show that $\cup V_n$ is
dense in $L_0^2(\mu)$. Let $\theta$ be an element in $\cap
V_n^\perp$. We have to show that $\theta=0$. For every $n$, 
$\theta=\theta_n\circ f^n$ for appropriate
$\theta_n$. Hence, $\theta$ is measurable with respect to the
$\sigma$-algebra $\Bc_\infty:=\cap_{n\geq 0}\Bc_n$.
We show that $\Bc_\infty$ is the
trivial algebra. Let $A$ be an element of $\Bc_\infty$. Define
$A_n=f^n{(A)}$. Since $A$ is in $\Bc_\infty$,  $\ind_A=\ind_{A_n}\circ f^n$ and
  $\Lambda^n(\ind_A)=\ind_{A_n}$. K-mixing implies that
  $\Lambda^n(\ind_A)$ converges in $L^2(\mu)$ to a constant, see
  Theorem \ref{th_K-mixing_end}. So,
  $\ind_{A_n}$ converges to a constant which is necessarily 0 or
  1. We deduce that $\mu(A_n)$ converges  to 0 or 1.
On the other hand, we have
$$\mu(A_n)= \langle \mu,\ind_{A_n}\rangle = \langle
\mu,\ind_{A_n}\circ f^n\rangle = \langle \mu,\ind_A\rangle=\mu(A).$$
Therefore, $A$ is of measure 0 or 1. This implies the decomposition of
$L_0^2(\mu)$. 

Suppose now that $\psi:=\sum \psi_n\circ f^n$ with
$\Lambda(\psi_n)=0$ is an element of $L^2_0(\mu)$. We have
$\Et(\psi|\Bc_n)=\sum_{i\geq n} \psi_i\circ f^i$. So, 
$$\sum_{n\geq 0} \|\Et(\psi|\Bc_n)\|_{L^2(\mu)}^2=\sum_{n\geq 0} (n+1)
\|\psi_n\|_{L^2(\mu)}^2,$$
and $\psi$ satisfies Gordin's condition if and only if the last sum is
finite.

Let $\psi$ be a d.s.h. function with $\langle
\mu,\psi\rangle=0$. 
It follows from
Theorem \ref{th_mixing} that 
$$\sup_{\|\varphi\|_{L^2(\mu)}\leq 1} |\langle \mu,(\varphi\circ f^n)
\psi\rangle|\lesssim d^{-n}.$$
Choose $\varphi=\psi_n/\|\psi_n\|_{L^2(\mu)}$. 
The above estimate implies that
$$\|\psi_n\|_{L^2(\mu)}= |\langle \mu,(\varphi\circ f^n)
\psi\rangle|\lesssim d^{-n}.$$ 
The case of $\Cc^\nu$ observables is proved in the same way. 

Observe that if $(\psi_n\circ f^n)_{n\geq 0}$ is the sequence of
projections of $\psi$ on the factors of the direct sum
$\oplus_{n=0}^\infty V_1\circ f^n$, then the coordinates of
$\Lambda(\psi)$ are $(\psi_n\circ f^{n-1})_{n\geq 1}$.
\endproof

We continue the study with other types of convergence.
Let us recall the almost sure version of the central limit theorem in
probability theory. Let             
$Z_n$ be random variables, identically distributed in
$L^2(X,\Fc,\nu)$, such that $\Et(Z_n)=0$ and $\Et(Z_n^2)=\sigma^2$,
$\sigma>0$. We say that {\it the almost
sure central limit theorem} holds if at $\nu$-almost every point in $X$, the sequence of measures 
$${1\over \log N}\sum_{n=1}^N {1\over n}
\delta_{n^{-1/2}\sum_{i=0}^{n-1}Z_i}$$
converges in law to the normal distribution of mean 0 and variance
$\sigma$. In particular, $\nu$-almost surely
$${1\over \log N}\sum_{n=1}^N {1\over n}
\ind_{\big\{n^{-1/2}\sum_{i=0}^{n-1}Z_i\leq t_0\big\}}\rightarrow {1\over
  \sqrt{2\pi}\sigma}\int_{-\infty}^{t_0} e^{-{t^2\over 2\sigma^2}} dt,$$
for any $t_0\in\R$. 
In the central limit theorem, we only get the 
$\nu$-measure of the set 
$\{N^{-1/2}\sum_{n=0}^{N-1}Z_n<t_0\}$ when $N$ goes to infinity. Here,
we get an information at $\nu$-almost every point for the
logarithmic averages.

The almost sure central limit theorem can be deduced from the
so-called almost sure invariance principle (ASIP for short). In the
case of i.i.d. random variables as above, this principle 
compares the variables  $\widehat Z_N$ with Brownian motions and 
gives some information about the fluctuations of $\widehat Z_N$ around
0. 

\begin{theorem}
Let $(X,\Fc,\nu)$ be a probability space. Let $(Z_n)$ be a sequence of
i.i.d. random variables with mean $0$ and variance $\sigma>0$. Assume
that there is an $\alpha>0$ such that $Z_n$ is in $L^{2+\alpha}(\nu)$. Then,
there is another probability space $(X',\Fc',\nu')$ with a sequence of
random variables 
$\St_N'$ on $X'$ which has the same joint distribution as 
$\St_N:=\sum_{n=0}^{N-1} Z_n$,
and a Brownian motion $\Bt$ of variance $\sigma$ on $X'$ such that 
$$|\St_N'-\Bt(N)|\leq c N^{1/2-\delta},$$
for some positive constants $c,\delta$. It follows that
$$|N^{-1/2}\St_N'-\Bt(1)|\leq c N^{-\delta}.$$
\end{theorem}

For weakly dependent
variables, this type of result is a consequence of a 
theorem due to Philipp-Stout 
\cite{PhilippStout}. 
It gives conditions which imply that the ASIP holds.  
Lacey-Philipp proved in  \cite{LaceyPhilipp} 
that the ASIP implies the almost sure
central limit theorem.  For holomorphic endomorphisms of $\P^k$, we
have the following result due to Dupont which holds in particular for
H\"older continuous observables \cite{Dupont2}.

\begin{theorem} \label{th_asip_pk}
Let $f$ be an endomorphism of algebraic degree $d\geq 2$ as above and
$\mu$ its equilibrium measure. 
Let $\varphi$  be an observable with values in $\R\cup\{-\infty\}$ such that
$e^\varphi$ is H\"older continuous, $H:=\{\varphi=-\infty\}$ is an
analytic set and $|\varphi|\lesssim |\log\dist(\cdot,H)|^\rho$ near $H$ for some
$\rho>0$.
If $\langle\mu,\varphi\rangle =0$ and $\varphi$ is not a coboundary, then 
the almost sure invariance principle holds for
$\varphi$. In particular, the almost sure central
limit theorem holds for such observables.
\end{theorem}

The ASIP in the above setting says that 
there is a probability space $(X',\Fc',\nu')$ with a sequence of
random variables 
$\St_N'$ on $X$ which has the same joint distribution as 
$\St_N:=\sum_{n=0}^{N-1} \varphi\circ f^n$,
and a Brownian motion $\Bt$ of variance $\sigma$ on $X'$ such that 
$$|\St_N'-\Bt(N)|\leq c N^{1/2-\delta},$$
for some positive constants $c,\delta$.

The ASIP implies other stochastic results, see \cite{PhilippStout}, in
particular, {\it the law of the iterated logarithm}. With our
notations, it implies  that for $\varphi$ as above 
$$\limsup_{N\rightarrow\infty} {\St_N(\varphi)\over \sigma
  \sqrt{N\log\log (N\sigma^2)}}=1 \quad \mu\mbox{-a.e.}$$
Dupont's approach is based on the Philipp-Stout's result applied to a
Bernoulli system and a quantitative Bernoulli property of the
equilibrium measure of $f$, i.e. a construction of coding tree. We
refer to Dupont and Przytycki-Urbanski-Zdunik
\cite{PrzytyckiUrbanskiZdunik} for these results.
Note that Bernoulli property of this measure was proved by
Briend in \cite{Briend1}. 
It says that outside sets of measure zero, the
system is conjugated to a shift. The dimension one case is due to Heicklen-Hoffman \cite{HeicklenHoffman}.

The last stochastic property we consider here is the large deviations
theorem. As above, we first recall the classical result in probability
theory.

\begin{theorem}
Let $Z_1,Z_2,\ldots$
be independent random variables on $(X,\Fc,\nu)$, identically distributed with values in
$\R$, and of mean zero, i.e. $\Et(Z_1)=0$. Assume also that for
$t\in\R$, $\exp(tZ_n)$ is integrable. Then, the limit
$$I(\epsilon):=-\lim_{N\rightarrow\infty}
\log\nu\Big\{\Big|{Z_1+\cdots+Z_N\over N}\Big|>\epsilon\Big\}.$$
exists and $I(\epsilon)>0$ for $\epsilon>0$.
\end{theorem}

The theorem estimates the size of the set where the average is away
from zero, the expected value. We have 
$$\nu\Big\{\Big|{Z_1+\cdots+Z_N\over N}\Big|>\epsilon\Big\}\sim
e^{-NI(\epsilon)}.$$ 
Our goal is to give an analogue for the
equilibrium measure of endomorphisms of $\P^k$. 
We first prove an abstract result corresponding to the above Gordin's
result for the central limit theorem.

Consider a dynamical system $g:(X,\Fc,\nu)\rightarrow (X,\Fc,\nu)$ as
above where $\nu$ is an invariant probability measure. 
So, $g^*$ defines a linear operator of norm 1 from $L^2(\nu)$
into itself.   
We say that $g$ has {\it bounded Jacobian} if
there is a constant $\kappa>0$ such that $\nu(g(A))\leq \kappa \nu(A)$
for every $A\in\Fc$. 
The following result was obtained in \cite{DinhNguyenSibony3}.

\begin{theorem} \label{th_deviation_real}
Let $g:(X,\Fc,\nu)\rightarrow (X,\Fc,\nu)$ be a map with bounded Jacobian
as above. Define $\Fc_n:=g^{-n}(\Fc)$.
Let $\psi$ be a bounded real-valued measurable function. Assume there
are constants $\delta>1$ and $c>0$ such that 
$$\big\langle \nu,
e^{\delta^n|\Et(\psi|\Fc_n)-\langle \nu,\psi\rangle|}\big\rangle \leq c\quad \mbox{for every}\quad
n\geq 0.$$ 
Then $\psi$ satisfies a weak large deviations theorem. More precisely, for every
$\epsilon>0$, there exists a constant $h_\epsilon>0$ such that 
\begin{equation*}
\nu\Big\{ x\in X:\   \big\vert\frac{1}{N}\sum_{n=0}^{N-1}\psi\circ
g^n(x)- \langle \nu,\psi\rangle \big\vert >\epsilon \Big\}\leq
e^{-N (\log N)^{-2}h_\epsilon}
\end{equation*}
for all $N$ large enough\footnote{In the LDT 
for independent random variables, there is no factor $(\log
  N)^{-2}$ in the estimate.}.
\end{theorem} 

We first prove some preliminary lemmas. The following one is a version
of the classical Bennett's inequality see \cite[Lemma 2.4.1]{DemboZeitouni}.

\begin{lemma}\label{lemma_Bennett}
Let $\psi$ be an
observable such that $\|\psi\|_{L^\infty(\nu)}\leq b$ for some
constant $b\geq 0$, and $\Et(\psi)=0$.  
Then  
$$\Et(e^{\lambda\psi}) \leq {e^{-\lambda b}+ e^{\lambda b} \over 2}$$ 
for every $\lambda\geq 0$. 
\end{lemma}
\proof
We can assume $\lambda=1$. Consider first the case where there is a
measurable set $A$ such that $\nu(A)=1/2$. 
Let $\psi_0$ be the function which is equal to
$-b$ on $A$ and to 
$b$ on $X\setminus A$. We have
$\psi_0^2= b^2\geq \psi^2$. Since $\nu(A)=1/2$, we have $\Et(\psi_0)=0$. 
Let  $g(t)=a_0t^2+a_1t+a_2$,  be the unique quadratic function such that
$h(t):=g(t)-e^t$ satisfies $h(b)=0$ and $h(-b)=h'(-b)=0$.
We have $g(\psi_0)=e^{\psi_0}$.

Since $h''(t)=2a_0-e^t$ admits at most one zero, $h'$ admits at most two
zeros. The fact that $h(-b)=h(b)=0$ implies that $h'$ vanishes in
$]-b,b[$. Hence $h'$ admits exactly one zero at $-b$ and another one in
$]-b,b[$. We deduce that $h''$ admits a zero. This implies that
$a_0>0$. Moreover, $h$ vanishes only at $-b$, $b$ and
$h'(b)\not=0$. It follows that $h(t)\geq 0$ on $[-b,b]$ because
$h$ is negative near $+\infty$.
Thus, $e^t\leq g(t)$ on $[-b,b]$ and then $e^\psi\leq  g(\psi)$.

Since  $a_0>0$, if an observable $\phi$ satisfies $\Et(\phi)=0$, 
then $\Et(g(\phi))$ is an increasing function of $\Et(\phi^2)$. 
Now, using the properties of $\psi$ and $\psi_0$, 
we obtain 
$$\Et(e^{\psi})  \leq    \Et(g(\psi))\leq  
\Et(g(\psi_0))  =  \Et(e^{\psi_0})
 =  { e^{-b}+e^b\over 2}.$$
This completes the proof under the assumption that $\nu(A)=1/2$ for
some measurable set $A$.

The general case is deduced from the previous particular case. Indeed,
it is enough to apply the first case to the disjoint union of $(X,\Fc,\nu)$ with a
copy $(X',\Fc',\nu')$ of this space, i.e. to the space $(X\cup
X',\Fc\cup\Fc',{\nu\over 2}+{\nu'\over 2})$, and to the function equal to
$\psi$ on $X$ and on $X'$.
\endproof

\begin{lemma} \label{lemma_Bennett_bis}
Let $\psi$ be an observable such that
$\|\psi\|_{L^\infty(\nu)}\leq b$ for some constant $b\geq 0$, and $\Et(\psi|\Fc_1)=0$.  
Then  
$$\Et(e^{\lambda\psi}|\Fc_1) \leq { e^{-\lambda b}+e^{\lambda b} \over 2}$$ 
for every $\lambda\geq 0$. 
\end{lemma}
\proof
We consider the desintegration of $\nu$ with respect to $g$. For
$\nu$-almost every $x\in X$, there is a
positive measure $\nu_x$ on $g^{-1}(x)$ such that if $\varphi$ is a
function in $L^1(\nu)$ then
$$\langle \nu,\varphi\rangle = \int_X \langle \nu_x,\varphi\rangle
d\nu(x).$$
Since $\nu$ is $g$-invariant, we have
$$\langle \nu,\varphi\rangle = \langle \nu,\varphi\circ g\rangle =\int_X
\langle \nu_x,\varphi\circ g\rangle d\nu(x)=\int_X \|\nu_x\|\varphi(x) d\nu(x).$$
Therefore,  $\nu_x$ is a probability measure for $\nu$-almost every $x$.
Using also the invariance of $\nu$, we obtain for $\varphi$ and $\phi$ in $L^2(\nu)$ that
\begin{eqnarray*}
\langle \nu ,\varphi(\phi\circ g)\rangle & = &  \int_X \langle
\nu_x,\varphi(\phi\circ g)\rangle
d\nu(x)= \int_X \langle
\nu_x,\varphi\rangle \phi(x) d\nu(x)\\
& = & \int_X \langle
\nu_{g(x)},\varphi\rangle \phi(g(x)) d\nu(x).
\end{eqnarray*}
We deduce that
$$\Et(\varphi|\Fc_1)(x)=\langle \nu_{g(x)},\varphi\rangle.$$
So, the hypothesis in the lemma is  that $\langle \nu_x,\psi\rangle=0$ for $\nu$-almost
every $x$. It suffices to check that 
$$\langle \nu_x,e^{\lambda\psi}\rangle \leq { e^{-\lambda
    b}+e^{\lambda b} \over 2}.$$
But this is a consequence of Lemma \ref{lemma_Bennett} applied to
$\nu_x$ instead of $\nu$. 
\endproof

We continue the proof of Theorem \ref{th_deviation_real}. 
Without loss of generality we can assume  that  $\langle
\nu,\psi\rangle=0$ and $|\psi|\leq 1$.
The  general idea is to  write $\psi=\psi'+(\psi''-\psi''\circ g)$
for  functions $\psi'$ and $\psi''$ in  $L^2(\nu)$  such that 
\begin{equation*}
\Et(\psi'\circ g^n|\Fc_{n+1})=0,\qquad  n\geq 0.
\end{equation*}
In the  language of probability theory, these identities
mean that  $(\psi'\circ g^n)_{n\geq 0}$ is  a {\it reversed martingale
  difference} as in Gordin's approach, see also \cite{Viana}.  
The strategy  is  to  prove the weak LDT  for  $\psi'$ and
for the coboundary $\psi''-\psi''\circ g$. Theorem \ref{th_deviation_real} is
then a consequence of Lemmas \ref{lemma_ldt} and \ref{lemma_ldt_bis} below.  

Let $\Lambda_g$ denote the adjoint of the operator $\varphi\mapsto
\varphi\circ g$ on $L^2(\nu)$. These operators are of norm 1. The
computation in Lemma \ref{lemma_Bennett_bis} shows that
$\Et(\varphi|\Fc_1)=\Lambda_g(\varphi)\circ g$. 
We obtain in the same way that $\Et(\varphi|\Fc_n)=\Lambda_g^n(\varphi)\circ g^n$. 
Define
\begin{equation*}
\psi'':=-\sum_{n=1}^{\infty} \Lambda_g^n(\psi),\qquad  \psi':=\psi-(\psi''-\psi''\circ g).
\end{equation*}  
Using the hypotheses in Theorem \ref{th_deviation_real},
we  see that  $\psi'$ and $\psi''$ are in $L^2(\nu)$ with norms
bounded by some constant. However, we loose the uniform boundedness:
these functions are not necessarily in $L^\infty(\nu)$.

\begin{lemma} \label{lemma_martingal}
We have $\Lambda_g^n(\psi')=0$ for $n\geq 1$ and $\Et(\psi'\circ g^n|\Fc_m)=0$ for $m>n\geq 0$.
\end{lemma}
\proof
Clearly 
$\Lambda_g(\psi''\circ g)=\psi''$. 
We deduce from the definition of $\psi''$ that
$$\Lambda_g(\psi')=\Lambda_g(\psi) -\Lambda_g (\psi'')+\Lambda_g(\psi''\circ g)
=\Lambda_g(\psi)-\Lambda_g(\psi'')+\psi''=0.$$
Hence, $\Lambda_g^n(\psi')=0$ for $n\geq 1$. 
For every function $\phi$ in $L^2(\nu)$, since $\nu$ is invariant, we
have for $m>n$
$$\langle \nu, (\psi'\circ g^n)(\phi\circ g^m) \rangle  = 
\langle \nu, \psi'(\phi\circ g^{m-n})\rangle = \langle \nu, \Lambda_g^{m-n}
(\psi')\phi\rangle =0.$$
It follows that $\Et(\psi'\circ g^n|\Fc_m)=0$.
\endproof

\begin{lemma} \label{lemma_level_estimate}
There are constants $\delta_0>1$ and $c>0$ such that
$$\nu\{|\psi'|>b\}\leq ce^{-\delta_0^b} \quad \mbox{and}\quad 
\nu\{|\psi''|>b\}\leq ce^{-\delta_0^b}$$ 
for any $b\geq 0$. In particular,
 $t\psi'$ and
  $t\psi''$ are $\nu$-integrable for every $t\geq 0$. 
\end{lemma}
\proof
Since $\psi':=\psi-(\psi''-\psi''\circ g)$ and $\psi$ is bounded, it
is enough to prove the estimate on 
$\psi''$. Indeed, the invariance of $\nu$ implies that $\psi''\circ g$
satisfies a similar inequality.

Fix a positive constant $\delta_1$ such that $1<\delta_1^2<\delta$, where $\delta$ is the constant in
Theorem \ref{th_deviation_real}. Define $\varphi:=\sum_{n\geq
  1}\delta_1^{2n}|\Lambda_g^n(\psi)|$. We first show that there is a
constant $\alpha>0$ such that $\nu\{\varphi\geq b\} \lesssim e^{-\alpha b}$
for every $b\geq 0$.
Recall that 
$\Et(\psi|\Fc_n)=\Lambda_g^n(\psi)\circ g^n$.
Using  
the hypothesis of Theorem \ref{th_deviation_real}, the inequality 
$\sum {1\over 2n^2}\leq 1$ and the invariance of $\nu$, we obtain for $b\geq 0$
\begin{eqnarray*}
\nu\{\varphi\geq b\}  & \leq & \sum_{n\geq 1} \nu\Big\{|\Lambda_g^n(\psi)|\geq
{\delta_1^{-2n}b\over 2n^2}\Big\}
 \leq   \sum_{n\geq 1} \nu\Big\{|\Et(\psi|\Fc_n)|\geq {\delta_1^{-2n}b\over 2n^2}\Big\}\\
& = & \sum_{n\geq 1} \nu\Big\{\delta^n|\Et(\psi|\Fc_n)|\geq
{\delta^n\delta_1^{-2n}b\over 2n^2}\Big\}\lesssim \sum_{n\geq 1}
\exp\Big({-\delta^n\delta_1^{-2n}b\over 2n^2}\Big). 
\end{eqnarray*}
It follows that $\nu\{\varphi\geq b\} \lesssim e^{-\alpha b}$ for some
constant $\alpha>0$.  

We prove now the estimate $\nu\{|\psi''|>b\}\leq ce^{-\delta_0^b}$. 
It is enough to consider the case where $b=2l$ for some positive
integer $l$. Recall that for simplicity we assumed 
$|\psi|\leq 1$. It follows that $|\Et(\psi|\Fc_n)|\leq 1$ and hence
$|\Lambda_g^n(\psi)|\leq 1$. 
We have
$$|\psi''| \leq  \sum_{n\geq 1} |\Lambda_g^n(\psi)| \leq  
\delta_1^{-2l} \sum_{n\geq 1} \delta_1^{2n}| \Lambda_g^n(\psi)|
+ \sum_{1\leq n\leq l} |\Lambda_g^n(\psi)|\\
\leq \delta_1^{-2l} \varphi+ l.$$
Consequently, 
$$\nu\big\{|\psi''|>2l\big\}\leq \nu\big\{\varphi> \delta_1^{2l} \big\}
\lesssim  e^{-\alpha \delta_1^{2l}}.$$
It is enough to choose $\delta_0<\delta_1$ and $c$ large enough.
\endproof

\begin{lemma} \label{lemma_ldt}
The coboundary $\psi''-\psi''\circ g$ satisfies the LDT.
\end{lemma}
\proof
Given a function $\phi\in L^1(\mu)$, recall that Birkhoff's sum 
$\St_N(\phi)$ is defined by  
$$ \St_0(\phi):=0\qquad \mbox{and} \qquad \St_N(\phi):= \sum\limits_{n=0}^{N-1}
\phi\circ g^n \quad \mbox{for } N\geq 1.$$
Observe that $\St_N(\psi''-\psi'' \circ g) =\psi''-\psi''\circ g^N$. 
Consequently, for  a given $\epsilon >0$, 
using the invariance of $\nu$, we have 
\begin{eqnarray*}
\nu\big\{ | \St_N(\psi''-\psi''\circ g)|  >N\epsilon \big\}& \leq&
\nu \Big\{|  \psi''\circ g^N|> {N\epsilon \over 2}\Big\}+ \mu
\Big\{|\psi''|>{N\epsilon\over 2}\Big\}\\
&=& 2\nu \Big\{|\psi''|>\frac{N\epsilon}{2}\Big\}.
\end{eqnarray*}
Lemma \ref{lemma_level_estimate} implies that the last expression is
smaller than $e^{-Nh_\epsilon}$ for some $h_\epsilon>0$ and for $N$
large enough. This completes the proof.
\endproof

It remains to show that  $\psi'$ satisfies the weak LDT. We use the
following lemma.

\begin{lemma} \label{lemma_bounded_part}
For every $b\geq 1$, there are Borel sets $W_N$ such that
$\nu(W_N)\leq c Ne^{-\delta_0^b}$ and 
$$\int_{X\setminus W_N} e^{\lambda \St_N(\psi')}d\nu\leq
2\Big [{e^{-\lambda b}+e^{\lambda b}\over 2}\Big]^N,$$
where $c>0$ is a constant independent of $b$.
\end{lemma}
\proof
For $N=1$, define $W:=\{|\psi'|>b\}$, $W':=g(W)$ and
$W_1:=g^{-1}(W')$.
Recall that the Jacobian of $\nu$ is of bounded by some constant
$\kappa$. This 
and Lemma \ref{lemma_level_estimate} imply that
$$\nu(W_1)=\nu(W')=\nu(g(W))\leq \kappa \nu(W)\leq ce^{-\delta_0^b}$$
for some constant $c>0$.
We also have
$$\int_{X\setminus W_1} e^{\lambda \St_1(\psi')}d\nu
=\int_{X\setminus W_1} e^{\lambda \psi'}d\nu\leq e^{\lambda b} \leq
2 \Big[{e^{-\lambda b}+e^{\lambda b}\over 2}\Big].$$
So, the lemma holds for $N=1$. 

Suppose the lemma for $N\geq 1$, we prove  it for $N+1$.
Define  
$$W_{N+1}:=g^{-1}(W_N) \cup W_1=g^{-1}(W_N\cup W').$$ 
We have
$$\nu(W_{N+1}) \leq  \nu(g^{-1}(W_N))+ \nu(W_1) = \nu(W_N)+
\nu(W_1)\leq c(N+1)e^{-\delta_0^b}.$$ 
We will apply Lemma \ref{lemma_Bennett_bis} 
to the function $\psi^*$ such that 
$\psi^*=\psi'$ on $X\setminus W_1$ and $\psi^*=0$ on $W_1$. 
By Lemma \ref{lemma_martingal}, we have
$\Et(\psi^*|\Fc_1)=0$ since $W_1$ is an element of $\Fc_1$. 
The choice of $W_1$ gives that $|\psi^*|\leq b$. 
By Lemma \ref{lemma_Bennett_bis}, we have 
$$\Et(e^{\lambda\psi^*}|\Fc_1)\leq  {e^{-\lambda b}+e^{\lambda
    b}\over 2} \quad \mbox{on}\quad X\quad \mbox{for}\quad \lambda\geq
0.$$ 
It follows that 
$$\Et(e^{\lambda\psi'}|\Fc_1)\leq  {e^{-\lambda b}+e^{\lambda
    b}\over 2} \quad \mbox{on}\quad X\setminus W_1\quad \mbox{for}\quad \lambda\geq
0.$$ 
Now, using the fact that $W_{N+1}$ and 
$e^{\lambda  \St_{N}(\psi'\circ g)}$ are $\Fc_1$-measurable, we can write
\begin{eqnarray*}
\int_{X\setminus W_{N+1}} e^{\lambda  \St_{N+1}(\psi')}d\nu
&=&\int_{X\setminus W_{N+1}} 
e^{\lambda\psi'} e^{\lambda  \St_N(\psi'\circ g)}d\nu\\
&=&\int_{X\setminus W_{N+1}} 
\Et( e^{\lambda\psi'}|\Fc_1) e^{\lambda  \St_N(\psi'\circ g)}d\nu.
\end{eqnarray*}
Since $W_{N+1}=g^{-1}(W_N)\cup W_1$, the last integral is bounded by
\begin{eqnarray*}
\lefteqn{\sup_{X\setminus W_1} \Et( e^{\lambda\psi'}|\Fc_1)\int_{X\setminus g^{-1}(W_N)}  
e^{\lambda \St_N(\psi'\circ g)}d\nu}\\ 
&\leq &  \Big[{e^{-\lambda b}+e^{\lambda b}\over 2}\Big] 
\int_{X\setminus W_N} 
e^{\lambda  \St_N(\psi')}d\nu \\
&\leq &  2\Big[{e^{-\lambda b}+e^{\lambda b}\over 2}\Big]^{N+1},
\end{eqnarray*}
where the last inequality  follows 
from the induction hypothesis. So, the lemma holds for $N+1$.
\endproof

The following lemma, together with Lemma \ref{lemma_ldt}, implies Theorem \ref{th_deviation}.

\begin{lemma} \label{lemma_ldt_bis}
The function $\psi'$ satisfies the weak LDT.
\end{lemma}
\proof
Fix an $\epsilon>0$. By Lemma \ref{lemma_bounded_part},
we  have, for  every $\lambda\geq 0$
\begin{eqnarray*}
\nu\big\{\St_N(\psi')\geq  N\epsilon\big\}&\leq&  
\nu (W_N)+e^{-\lambda N\epsilon}
\int_{X\setminus W_N} e^{\lambda  \St_N(\psi')}d\nu\\
&\leq &  cNe^{-\delta_0^b}+2e^{-\lambda N\epsilon} \Big 
[{e^{-\lambda b}
+e^{\lambda b}\over 2}\Big]^N.
\end{eqnarray*}
Let $b:=\log N (\log \delta_0)^{-1}$. We have 
$$cNe^{-\delta_0^b} = cN e^{-N} \leq e^{-N/2}$$
for $N$ large. We also have
$${e^{-\lambda b}+e^{\lambda b}\over 2} =   \sum_{n\geq 0}
{\lambda^{2n}b^{2n}\over (2n)!}\leq e^{\lambda^2b^2}.$$
Therefore,  if $\lambda:=u\epsilon b^{-2}$
with a fixed $u>0$ small enough
$$2e^{-\lambda N\epsilon} \Big [{e^{-\lambda b}+e^{\lambda
    b}\over 2}\Big]^N \leq 2e^{-\epsilon^2 b^{-2}(1-u)Nu}=
2e^{-2N(\log N)^{-2}h_\epsilon}$$
for some constant $h_\epsilon>0$. We deduce from the previous
estimates that
$$\nu\big\{\St_N(\psi')\geq  N\epsilon\big\}\leq e^{-N(\log
  N)^{-2}h_\epsilon}$$
for $N$ large. A similar estimate holds for $-\psi'$. So, $\psi'$
satisfies the weak LDT. 
\endproof

We deduce from Theorem \ref{th_deviation_real}, Corollaries
\ref{cor_exp_dsh} and \ref{cor_exp_holder} the following result
\cite{DinhNguyenSibony3}.

\begin{theorem} \label{th_deviation}
Let $f$ be a holomorphic endomorphism of $\P^k$ of algebraic degree
$d\geq 2$. Then the equilibrium measure $\mu$ of $f$ satisfies the
weak large deviations theorem for bounded
d.s.h. observables and also for H\"older continuous observables. More
precisely, if a function $\psi$ is bounded d.s.h. or H{\"o}lder continuous, then for every
$\epsilon>0$ there is a constant $h_\epsilon>0$ such that 
\begin{equation*}
\mu\Big\{ z\in \P^k:\   \big\vert\frac{1}{N}\sum_{n=0}^{N-1}\psi\circ
f^n(z)- \langle \mu,\psi\rangle \big\vert >\epsilon \Big\}\leq
e^{-N (\log N)^{-2}h_\epsilon}
\end{equation*}
for all $N$ large enough. 
\end{theorem}

The exponential estimate on $\Lambda^n(\psi)$ is crucial in the
proofs of the previous results. It is nearly an
estimate in sup-norm. Note that if $\|\Lambda^n(\psi)\|_{L^\infty(\mu)}$
converge exponentially fast to 0 then $\psi$ satisfies the LDT. This is
the case for H{\"o}lder continuous observables in dimension 1, following a result by Drasin-Okuyama
\cite{DrasinOkuyama}, and when $f$ is a  generic map 
in higher dimension, see Remark \ref{rk_unif_cv_perron}.
The LDT was recently obtained in dimension 1 by Xia-Fu in
\cite{XiaFu} for Lipschitz observables.

\bigskip\bigskip

\begin{exercise}
Show that if a Borel set $A$ satisfies $\mu(A)>0$, then
$\mu(f^n(A))$ converges to $1$. 
\end{exercise}

\begin{exercise}
Show that $\sup_\psi I_n(\varphi,\psi)$ with $\psi$ smooth
$\|\psi\|_\infty\leq 1$ is equal to $\|\varphi\|_{L^1(\mu)}$. Deduce
that there is no decay of correlations which is uniform on
$\|\psi\|_\infty$.
\end{exercise}

\begin{exercise} Let $V_1:=\{\psi\in L^2_0(\mu),\
  \Lambda(\psi)=0\}$. Show that $V_1$ is infinite dimensional and that
  bounded functions in $V_1$ are dense in $V_1$ with respect to the $L^2(\mu)$-topology.
\end{exercise}

\begin{exercise}
Let $\varphi$ be a d.s.h. function as in Corollary \ref{cor_clt_end}. 
Show that 
$$\|\varphi+\cdots+\varphi\circ
f^{n-1}\|^2_{L^2(\mu)}-n\sigma^2+\gamma= O(d^{-n}),$$
where $\gamma:=2\sum_{n\geq 1} n\langle \mu,\varphi(\varphi\circ
f^n)\rangle$ is a finite constant. Prove an analogous property for
$\varphi$ H\"older continuous. 
\end{exercise}


\section{Entropy, hyperbolicity and dimension}

There are various ways to describe the complexity of a dynamical
system. A basic measurement is the entropy which is closely related to
the volume growth of the images of subvarieties. We will compute the topological entropy and the metric
entropy of holomorphic endomorphisms of 
$\P^k$. We will also estimate the Lyapounov exponents with
respect to the measure of maximal entropy and the Hausdorff dimension
of this measure. 

We recall few notions.
Let $(X,\dist)$ be a compact metric space where $\dist$ is a distance
on $X$. Let $g:X\rightarrow X$ be a continuous map. We introduce {\it the
Bowen metric} associated to $g$. For a positive integer $n$, define
the distance $\dist_n$  on $X$ by
$$\dist_n(x,y):=\sup_{0\leq i\leq n-1}\dist(g^i(x),g^i(y)).$$
We have $\dist_n(x,y)>\epsilon$ if the orbits
$x,g(x),g^2(x),\ldots$ of $x$ and $y,g(y),g^2(y),\ldots$ of
$y$ are
distant by more than $\epsilon$ at a time $i$ less than $n$. In which case, we
say that $x,y$ are {\it $(n,\epsilon)$-separated}. 

The topological entropy
measures the rate of growth in function of time $n$, of the number of
orbits that can be distinguished at $\epsilon$-resolution. In other
words, it measures the divergence of the
orbits. More precisely, for $K\subset X$, not necessarily invariant,
let $N(K,n,\epsilon)$ denote the maximal number of points in $K$ which
are pairwise $(n,\epsilon)$-separated. This number increases as
$\epsilon$ decreases.
{\it The topological entropy} of $g$ on $K$ is 
$$h_t(g,K):=\sup_{\epsilon>0} \limsup_{n\rightarrow\infty} {1\over n}
\log N(K,n,\epsilon).$$
{\it The topological entropy} of $g$ is the entropy on $X$ and is
denoted by $h_t(g)$.
The reader can check that if $g$ is an isometry, then $h_t(g)=0$. 
In complex dynamics, we often have that for $\epsilon$ small enough, ${1\over n}
\log N(X,n,\epsilon)$ converge to $h_t(g)$.

Let $f$ be an endomorphism of algebraic degree $d\geq 2$ of $\P^k$as above. 
As we have seen, the iterate $f^n$ of $f$ has algebraic degree $d^n$.
If $Z$ is an algebraic set in $\P^k$ of codimension
$p$ then the degree of $f^{-n}(Z)$, counted with multiplicity, is
equal to $d^p\deg(Z)$ and the degree of $f^n(Z)$, counting with
multiplicity, is equal to $d^{k-p}\deg(Z)$.
This is a consequence of 
B{\'e}zout's theorem. Recall that the degree of an algebraic set of
codimension $p$ in $\P^k$
is the number of points in the intersection with a generic projective
subspace of dimension $p$. 

The pull-back 
by $f$ induces a linear map
$f^*:H^{p,p}(\P^k,\C)\rightarrow H^{p,p}(\P^k,\C)$ which is just the
multiplication by $d^p$.
The constant $d^p$ is {\it the dynamical degree of order $p$} of
$f$. 
Dynamical degrees were considered by Gromov in \cite{Gromov}
where he introduced a method to bound the topological entropy from
above. We will see that they measure the volume growth of the graphs.  
The degree of maximal order $d^k$ is also called {\it the topological
  degree}. It is equal to the number of points in a fiber counting
with multiplicity.
The push-forward
by $f^n$ induces a linear map
$f_*:H^{p,p}(\P^k,\C)\rightarrow H^{p,p}(\P^k,\C)$ which is the
multiplication by $d^{k-p}$. These operations act continuously on
positive closed currents and hence, the actions are compatible with
cohomology, see Appendix \ref{section_pk}.

We have the following result due to Gromov \cite{Gromov} for the upper
bound and to
Misiurewicz-Przytycky \cite{MisiurewiczPrzytycky}
for the lower bound of the entropy.

\begin{theorem} \label{th_gromov_Pk}
Let $f$ be a holomorphic endomorphism of algebraic
  degree $d$ on $\P^k$. Then the topological entropy $h_t(f)$ of $f$ is equal
  to $k\log d$, i.e. to the logarithm of the maximal dynamical
  degree. 
\end{theorem}

The inequality $h_t(f)\geq k\log d$ is a consequence of the following result
which is valid for arbitrary $\Cc^1$ maps \cite{MisiurewiczPrzytycky}.

\begin{theorem}[Misuriewicz-Przytycki] \label{th_MisiurewiczPrzytycky}
Let $X$ be a compact smooth orientable
manifold and $g:X\rightarrow X$ a $\Cc^1$ map. Then 
$$h_t(g)\geq \log |\deg (g)|.$$
\end{theorem}

Recall that {\it the degree} of $g$ is defined as follows. Let $\Omega$ be a
continuous form of maximal degree on $X$ such that
$\int_X\Omega\not=0$. Then
$$\deg(g):={\int_X g^*(\Omega)\over \int_X \Omega}\cdot$$
The number is independent of the choice of $\Omega$. When $X$ is a
complex manifold, it is necessarily orientable and $\deg(g)$ is just
the generic number of preimages of a point, i.e. the
topological degree of $g$. In our case, the topological degree of $f$
is equal to $d^k$. So, $h_t(f)\geq k\log d$. 

Instead of using Misuriewicz-Przytycki theorem, it is also possible to
apply the following important result due to Yomdin \cite{Yomdin}.

\begin{theorem}[Yomdin] \label{th_yomdin}
Let $X$ be a compact smooth manifold and
  $g:X\rightarrow X$ a smooth map. Let $Y$ be a manifold in $X$ smooth
  up to the boundary, then  
$$\limsup_{n\rightarrow\infty}{1\over n} \log \vol(g^n(Y))\leq h_t(g),$$
where the volume of $g^n(Y)$ is counted with multiplicity.
\end{theorem}

In our situation, when $Y=\P^k$, we have $\vol(g^n(Y))\simeq
d^{kn}$. Therefore, $h_t(f)\geq k\log d$. We can also deduce 
this inequality from  Theorem \ref{th_unique_entropy}
and the variational principle below. 

\medskip
\noindent
{\bf End of the proof of Theorem \ref{th_gromov_Pk}.}
It remains to prove that $h_t(f)\leq k\log
d$. Let $\Gamma_n$ denote the graph of $(f,f^2,\ldots, f^{n-1})$ in
$(\P^k)^n$, i.e. the set of points 
$$(z,f(z), f^2(z),\ldots,f^{n-1}(z))$$
 with $z$ in $\P^k$. 
This is a manifold of dimension $k$.
Let $\Pi_i$, $i=0,\ldots,n-1$, denote
the projections from $(\P^k)^n$ onto the factors $\P^k$. 
We use on $(\P^k)^n$ the metric and the distance associated
to the K{\"a}hler form
$\omega_n:=\sum\Pi_i^*(\omega_\FS)$
induced by the Fubini-Study metrics $\omega_\FS$ on the factors
$\P^k$, see Appendix \ref{section_pk}.
The following indicator $\lov(f)$ was introduced by Gromov, it
measures the growth rate of the volume of $\Gamma_n$,
$$\lov(f):=\lim_{n\rightarrow\infty} {1\over n}\log\vol(\Gamma_n).$$
The rest of the proof splits into two parts. We first show that 
the previous limit exists and is equal to $k\log d$ and then 
we prove the inequality $h_t(f)\leq \lov(f)$.

Using that $\Pi_0:\Gamma_n\rightarrow\P^k$ is
a bi-holomorphic map and that $f^i=\Pi_i\circ(\Pi_{0|\Gamma_n})^{-1}$,
we obtain 
\begin{eqnarray*}
k!\vol(\Gamma_n) &=& \int_{\Gamma_n}\omega_n^k
= \sum_{0\leq i_s\leq n-1}\int_{\Gamma_n}
\Pi_{i_1}^*(\omega_\FS)\wedge\ldots \wedge \Pi_{i_k}^*(\omega_\FS)\\
&=&  \sum_{0\leq i_s\leq n-1}\int_{\P^k}
(f^{i_1})^*(\omega_\FS)\wedge\ldots \wedge (f^{i_k})^*(\omega_\FS).
\end{eqnarray*}
The last sum contains $n^k$ integrals that we can compute
cohomologically. The above discussion on the action of $f^n$ on
cohomology implies that the last integral is equal to
$d^{i_1+\cdots+i_k}\leq d^{kn}$. 
So, the sum is bounded from above by $n^kd^{kn}$. When
$i_1=\cdots=i_k=n-1$, we see that $k!\vol(\Gamma_n)\geq
d^{(n-1)k}$. Therefore, the limit in the definition of $\lov(f)$
exists and is equal to $k\log d$. 

For the second step, we need the following classical estimate due to
Lelong \cite{Lelong}, see also
Appendix \ref{section_positive}.

\begin{lemma}[Lelong]
Let $A$ be an analytic set of pure dimension $k$ in a ball $B_r$ of
radius $r$ in $\C^N$. Assume that $A$ contains the center of
$B_r$. Then the $2k$-dimensional volume of $A$ is at least equal to
the volume of a ball of radius $r$ in $\C^k$. In particular, we have 
$$\vol(A)\geq c_k r^{2k},$$ 
where $c_k>0$ is a constant independent of $N$ and of $r$.
\end{lemma}

We prove now the inequality $h_t(f)\leq\lov(f)$. 
Consider an $(n,\epsilon)$-separated set $\Fc$ in $\P^k$. 
For each
point $a\in\Fc$, let $a^{(n)}$ denote the corresponding point 
$(a,f(a),\ldots, f^{n-1}(a))$ in $\Gamma_n$ and $B_{a,n}$ the ball of
center $a^{(n)}$ and of radius $\epsilon/2$ in $(\P^k)^n$. Since $\Fc$
is $(n,\epsilon)$-separated, these balls are disjoint. On the other
hand, by Lelong's inequality, 
$\vol(\Gamma_n\cap B_{a,n})\geq c_k'\epsilon^{2k}$, $c_k'>0$. 
Note that Lelong's inequality is stated in the Euclidean
metric. We can apply it using a fixed atlas of $\P^k$ and the
corresponding product atlas of $(\P^k)^n$, the distortion is bounded. 
So, $\#\Fc\leq c_k'^{-1}\epsilon^{-2k} \vol(\Gamma_n)$ and hence, 
$${1\over n} \log\#\Fc\leq {1\over n}\log(\vol(\Gamma_n))+O\Big({1\over
  n}\Big).$$
It follows that  $h_t(f)\leq \lov(f)=k\log d$.
\hfill $\square$

\medskip

We study the entropy of $f$ on some subsets of $\P^k$. 
The following result is due to de Th{\'e}lin and Dinh \cite{deThelin3,
  Dinh3}.

\begin{theorem} \label{th_entropy}
Let $f$ be a holomorphic endomorphism of $\P^k$ of algebraic degree
$d\geq 2$ and $\Jc_p$ its Julia set of order $p$, $1\leq p\leq k$. If $K$ is a subset of
$\P^k$ such that $\overline K\cap \Jc_p=\varnothing$, then 
$h_t(f,K)\leq (p-1)\log d$.
\end{theorem}
\proof
The proof is based on Gromov's idea as in 
Theorem \ref{th_gromov_Pk} and on the
speed of convergence towards the Green current.
Recall that $\Jc_p$ is the support of the Green $(p,p)$-current $T^p$ of $f$.
Fix an open neighbourhood $W$ of $\overline K$ such that  
$W\Subset \P^k\setminus \supp(T^p)$. 
Using the notations in Theorem \ref{th_gromov_Pk},
we only have to prove that
$$\lov(f,W):=\limsup_{n\rightarrow\infty}{1\over n} \log
\vol(\Pi_0^{-1}(W)\cap \Gamma_n)\leq (p-1)\log d.$$

It is enough to
show that $\volume(\Pi_0^{-1}(W)\cap \Gamma_n)\lesssim
n^kd^{(p-1)n}$. As in Theorem \ref{th_gromov_Pk}, it is sufficient to
check  that for $0\leq n_i\leq n$
$$\int_{W}(f^{n_1})^*(\omega_\FS)\wedge \ldots \wedge
(f^{n_k})^*(\omega_\FS) \lesssim d^{(p-1)n}.$$
To this end, we prove by induction on $(r,s)$, $0\leq r\leq p$ and $0\leq s\leq 
k-p+r$, that 
$$\|T^{p-r}\wedge  (f^{n_1})^*(\omega_\FS)\wedge\ldots\wedge
(f^{n_s})^*(\omega_\FS)\|_{W_{r,s}}\leq c_{r,s}d^{n(r-1)},$$
where $W_{r,s}$ is a neighbourhood of $\overline W$
and $c_{r,s}\geq 0$ is a constant independent of $n$ and of $n_i$. We obtain the
result by taking $r=p$ and $s=k$.

It is clear that the previous inequality holds when $r=0$
and also when $s=0$. In both cases, we can take 
$W_{r,s}=\P^k\setminus\supp(T^p)$ and $c_{r,s}=1$. 
Assume the inequality for $(r-1,s-1)$ and $(r,s-1)$.
Let $W_{r,s}$ be a neighbourhood of $\overline W$ strictly contained
in $W_{r-1,s-1}$ and $W_{r,s-1}$. Let $\chi\geq 0$ be a smooth cut-off
function with support in $W_{r-1,s-1}\cap W_{r,s-1}$ which is equal to 1
on $W_{r,s}$. We only have to prove that
$$\int T^{p-r}\wedge (f^{n_1})^*(\omega_\FS)\wedge\ldots\wedge
(f^{n_s})^*(\omega_\FS) \wedge \chi\omega_\FS^{k-p+r-s}\leq c_{r,s}d^{n(r-1)}.$$

If $\gr$ is the Green function
of $f$, we have
$$(f^{n_1})^*(\omega_\FS)=d^{n_1}T-\ddc (\gr\circ f^{n_1}).$$
The above integral is equal to the sum of the following integrals 
$$ d^{n_1}\int T^{p-r+1}\wedge (f^{n_2})^*(\omega_\FS)\wedge\ldots\wedge
(f^{n_s})^*(\omega_\FS) \wedge \chi\omega_\FS^{k-p+r-s}$$
and 
$$ -\int T^{p-r}\wedge \ddc(\gr\circ f^{n_1})\wedge (f^{n_2})^*(\omega_\FS)\wedge\ldots\wedge
(f^{n_s})^*(\omega_\FS) \wedge \chi\omega_\FS^{k-p+r-s}.$$
Using the case of $(r-1,s-1)$ we can bound the first integral by $c
d^{n(r-1)}$. Stokes' theorem implies that the second integral is equal to 
$$ -\int T^{p-r}\wedge (f^{n_2})^*(\omega_\FS)\wedge\ldots\wedge
(f^{n_s})^*(\omega_\FS) \wedge (\gr\circ f^{n_1})\ddc\chi\wedge
\omega_\FS^{k-p+r-s}$$
which is bounded by
$$ \|\gr\|_\infty\|\chi\|_{\Cc^2}\|T^{p-r}\wedge (f^{n_2})^*(\omega_\FS)\wedge\ldots\wedge
(f^{n_s})^*(\omega_\FS)\|_{W_{r,s-1}}$$
since $\chi$ has support in $W_{r,s-1}$.
We obtain the result using the $(r,s-1)$ case.
\endproof

The above result suggests a local indicator of volume growth. Define
for $a\in\P^k$
$$\lov(f,a):=\inf_{r>0}\limsup_{n\rightarrow\infty} {1\over
  n} \log\vol(\Pi_0^{-1}(B_r)\cap\Gamma_n),$$
where $B_r$ is the ball of center $a$ and of radius $r$. We can show
that if $a\in\Jc_p\setminus\Jc_{p+1}$ and if $\overline B_r$ does not
intersect $\Jc_{p+1}$, the above limsup is in fact a
limit and is equal to $p\log d$. One can also
consider the graph of $f^n$ instead of $\Gamma_n$. The notion can be
extended to meromorphic maps and its sub-level sets are analogues of
Julia sets.

\medskip

We discuss now the metric entropy, i.e. the entropy of an invariant
measure, a notion due to Kolgomorov-Sinai.
Let $g:X\rightarrow X$ be map on a space $X$ which is measurable with
respect to a $\sigma$-algebra $\Fc$. Let $\nu$ be an invariant probability measure for
$g$. Let $\xi=\{\xi_1,\ldots,\xi_m\}$ be a measurable partition of $X$. The entropy of $\nu$
with respect to $\xi$ is a measurement of the information we gain when
we know that a point $x$ belongs to a member of the partition
generated by $g^{-i}(\xi)$ with $0\leq i\leq n-1$.

The information we
gain when we know that a point $x$ belongs to  $\xi_i$ is
a positive function
$I(x)$ which depends only on $\nu(\xi_i)$,
i.e. $I(x)=\varphi(\nu(\xi_i))$. The information given by independent events should be
additive. In other words, we have 
$$\varphi(\nu(\xi_i)\nu(\xi_j))=\varphi(\nu(\xi_i))+\varphi(\nu(\xi_j))$$
for $i\not=j$. Hence, $\varphi(t)=-c\log t$
with $c>0$. With the normalization $c=1$, the information function for the
partition $\xi$ is defined by
$$I_\xi(x):=\sum -\log\nu(\xi_i) \ind_{\xi_i}(x).$$
The entropy of $\xi$ is the average of $I_\xi$:
$$H(\xi):=\int I_\xi(x)d\nu(x)=-\sum \nu(\xi_i)\log
\nu(\xi_i).$$
It is useful to observe that the function $t\mapsto -t\log t$ is
concave on $]0,1]$ and has the maximal value $e^{-1}$ at $e^{-1}$. 

Consider now the information obtained if we measure the position of
the orbit
$x,g(x),\ldots, g^{n-1}(x)$ relatively to $\xi$.
By definition, this is the measure of the entropy of the partition generated by
$\xi, g^{-1}(\xi),\ldots, g^{-n+1}(\xi)$, which we denote by
$\bigvee_{i=0}^{n-1}g^{-i}(\xi)$. The elements of this partition are
$\xi_{i_1}\cap g^{-1}(\xi_{i_2})\cap\ldots\cap g^{-n+1}(\xi_{i_{n-1}})$.
It can be shown \cite{Walters} that 
$$h_\nu(g,\xi):=\lim_{n\rightarrow\infty} {1\over n} H\big(\bigvee_{i=0}^{n-1}
g^{-i}(\xi)\big)$$
exists.
{\it The entropy of the measure} $\nu$ is defined as
$$h_\nu(g):=\sup_\xi h_\nu(g,\xi).$$

Two measurable dynamical systems $g$ on $(X,\Fc,\nu)$ and $g'$ on
$(X',\Fc',\nu')$ are said to be {\it measurably conjugate} if there is a
measurable invertible map $\pi:X\rightarrow X'$ such that $\pi\circ
g=g'\circ \pi$ and $\pi_*(\nu)=\nu'$. In that case, we have
$h_\nu(g)=h_{\nu'}(g')$. So, entropy is a conjugacy invariant. 
Note also that $h_\nu(g^n)=n h_\nu(g)$ and if $g$ is invertible, 
$h_\nu(g^n)=|n|h_\nu(g)$ for $n\in\Z$. Moreover, if $g$ is a
continuous map of a compact metric space, then  $\nu\mapsto
h_\nu(g)$ is affine function on the convex set of
$g$-invariant probability measures
\cite[p.164]{KatokHasselblatt}. 

We say that a measurable partition $\xi$  is {\it a
  generator} if up to sets of measure zero, $\Fc$ is the
smallest $\sigma$-algebra containing $\xi$ which is invariant under $g^{-1}$.
A finite partition $\xi$
is called {\it a strong generator} for a measure preserving dynamical
system $(X,\Fc,\nu,g)$ as above, if $\bigvee_{n=0}^\infty g^{-n}(\xi)=\Fc$ up to sets
of zero $\nu$-measure. The
following result of Kolmogorov-Sinai is useful in order to compute the
entropy \cite{Walters}. 

\begin{theorem}[Kolmogorov-Sinai] Let $\xi$ be a strong generator for
the dynamical system  $(X,\Fc,\nu,g)$ as above. Then 
$$h_\nu(g)=h_\nu(g,\xi).$$
\end{theorem}

We recall another useful theorem due to Brin-Katok \cite{BrinKatok}
which is valid for continuous maps $g:X\rightarrow X$ on a compact
metric space. 
Let $B_n^g(x,\delta)$ denote the ball of center $x$
and of radius $\delta$ with respect to the Bowen distance $\dist_n$.
We call $B^g_n(x,\delta)$ {\it the Bowen $(n,\delta)$-ball}. Define
local entropies of an invariant probability measure $\nu$ by
$$h_\nu(g,x):=\sup_{\delta>0} \limsup_{n\rightarrow\infty} -{1\over n}
\log \nu(B_n^g(x,\delta))$$
and
$$h_\nu^-(g,x):=\sup_{\delta>0} \liminf_{n\rightarrow\infty} -{1\over n}
\log \nu(B_n^g(x,\delta)).$$

\begin{theorem}[Brin-Katok] \label{th_brin_katok}
Let $g:X\rightarrow X$ be a continuous map
  on a compact metric space. Let $\nu$ be an invariant probability
  measure of finite entropy. Then, $h_\nu(g,x)=h^-_\nu(g,x)$ and $h_\nu(g,g(x))=
  h_\nu(g,x)$ for $\nu$-almost every $x$. Moreover,
$\langle\nu, h_\nu(g,\cdot)\rangle$ is equal to the entropy $h_\nu(g)$
of $\nu$. In particular, if $\nu$ is ergodic, we have
$h_\nu(g,x)=h_\nu(g)$ $\nu$-almost everywhere.
\end{theorem}

One can roughly say that $\nu(B_n^g(x,\delta))$ goes to zero at the
exponential rate $e^{-h_\nu(g)}$ for $\delta$ small. We can deduce
from the above theorem that if $Y\subset X$ is a Borel set with $\nu(Y)>0$, then
$h_t(g,Y)\geq h_\nu(g)$. The comparison with the topological entropy
is given by the variational principle \cite{KatokHasselblatt, Walters}.

\begin{theorem}[variational principle] Let $g:X\rightarrow X$ be a
  continuous map on a compact metric space. Then 
$$\sup h_\nu(g)=h_t(g),$$
where the supremum is taken over the invariant probability measures
$\nu$.
\end{theorem}

Newhouse proved in \cite{Newhouse} that if $g$ is a smooth map on a
smooth compact manifold, there is always a measure $\nu$ of maximal entropy,
i.e. $h_\nu(g)=h_t(g)$. One of the natural question in dynamics is to
find the measures which maximize entropy. Their supports are in some
sense the most chaotic parts of the system.
The notion of Jacobian of a
measure is useful in order to estimate the metric entropy.

Let $g:X\rightarrow X$ be a measurable map as above which preserves a probability
measure  
$\nu$. 
Assume there is a countable partition $(\xi_i)$ of $X$,  such that the map $g$
is injective on each $\xi_i$. {\it The Jacobian} $J_\nu(g)$ of $g$ with respect to $\nu$ is defined as
the Radon-Nikodym derivative of $g^*(\nu)$ with respect to $\nu$ on each
$\xi_i$. Observe that $g^*(\nu)$ is well-defined on $\xi_i$ since $g$
restricted to $\xi_i$ is injective.
We have the following theorem due to Parry \cite{Parry}. 

\begin{theorem}[Parry] \label{th_parry}
Let $g$, $\nu$ be  as above and $J_\nu(g)$ the Jacobian of $g$ with
respect to $\nu$. Then
$$h_\nu(g)\geq \int\log J_\nu(g) d\nu.$$
\end{theorem}

We now discuss the metric entropy of holomorphic maps on $\P^k$.
The
following result is a consequence of the variational principle and
Theorems \ref{th_gromov_Pk} and \ref{th_entropy}.

\begin{corollary} \label{cor_entropy_measure}
Let $f$ be an endomorphism of algebraic degree
  $d\geq 2$ of $\P^k$. Let  $\nu$ be 
an invariant probability measure. Then $h_\nu(f)\leq
  k\log d$. If the support of $\nu$ does not intersect the Julia set
  $\Jc_p$ of order $p$,
then $h_\nu(f)\leq (p-1)\log d$.
\end{corollary}

In the following result, the value of the metric entropy was
obtained in \cite{BriendDuval2,Sibony} 
and the uniqueness was obtained by Briend-Duval
in \cite{BriendDuval2}. The case of dimension 1 is due to
Freire-Lop\`es-Ma\~n\'e \cite{FreireLopesMane} and Lyubich
\cite{Lyubich}. 

\begin{theorem} \label{th_unique_entropy}
Let $f$ be an endomorphism of algebraic degree $d\geq
2$ of $\P^k$. Then the equilibrium measure $\mu$ of $f$ is the unique invariant
measure of maximal entropy $k\log d$.
\end{theorem}
\proof
We have seen in Corollary \ref{cor_entropy_measure} that $h_\mu(f)\leq k\log d$.
Moreover, $\mu$ has no mass on analytic sets, in particular on
the critical set of $f$. Therefore, if $f$ is injective on a Borel set
$K$, then  
$f_*(\ind_K)=\ind_{f(K)}$ and the total invariance of $\mu$
implies that
$\mu(f(K))=d^k\mu(K)$. So, $\mu$ is a measure of constant Jacobian
$d^k$. It follows from Theorem \ref{th_parry} that its entropy
is at least equal to $k\log d$. So, $h_\mu(f)=k\log d$.  

Assume now that there is another invariant probability measure $\nu$
of entropy $k\log d$. We are looking for a contradiction. Since
entropy is an affine function on $\nu$, we can
assume that $\nu$ is ergodic. This measure has no mass on proper
analytic sets of $\P^k$ since otherwise its entropy is at most equal
to $(k-1)\log d$, see Exercise \ref{exo_gromov} below. By Theorem \ref{th_equi_k},
$\nu$ is not totally invariant, so it is not of constant
Jacobian. Since $\mu$ has no mass on critical values of $f$, there is a simply connected open set $U$, not
necessarily connected, such that $f^{-1}(U)$ is a union
$U_1\cup\ldots\cup U_{d^k}$ of disjoint open sets such that $f:U_i\rightarrow
U$ is bi-holomorphic. One can choose $U$ and $U_i$ such that the $U_i$
do not have the same $\nu$-measure, otherwise $\mu=\nu$. So, we can assume that
$\nu(U_1)>d^{-k}$. This is possible since two ergodic measures are
multually singular.
Here, it is necessarily to chose $U$ so that
$\mu(\P^k\setminus U)$ is small.

Choose an open set $W\Subset U_1$ such that $\nu(W)>\sigma$ for
some constant $\sigma>d^{-k}$. Let
$m$ be a fixed integer and let $Y$ be the set of points $x$ such that for
every $n\geq m$, there
are at least $n\sigma$ points $f^i(x)$ with $0\leq i\leq n-1$ which belong
to $W$. If $m$ is large enough, Birkhoff's theorem implies that $Y$ has
positive $\nu$-measure. 
By Brin-Katok's theorem \ref{th_brin_katok}, we have
$h_t(f,Y)\geq h_\nu(f)= k\log d$. 

Consider the open sets $\Uc_\alpha:=U_{\alpha_0}\times\cdots \times U_{\alpha_{n-1}}$ in
$(\P^k)^n$ such that there are at least $n\sigma$ indices $\alpha_i$
equal to 1. 
A straighforward computation shows that the number of such open
sets is $\leq d^{k\rho n}$ for some constant $\rho<1$. 
Let $\Vc_n$ denote the union of these $\Uc_\alpha$. Using the same arguments as in 
Theorem \ref{th_gromov_Pk}, we get that
$$k\log d\leq h_t(f,Y)\leq \lim_{n\rightarrow\infty} {1\over n}
\log\vol(\Gamma_n\cap \Vc_n)$$
and
$$k!\vol(\Gamma_n\cap \Vc_n)=\sum_{0\leq i_s\leq n-1}\sum_\alpha 
\int_{\Gamma_n\cap \Uc_\alpha} \Pi_{i_1}^*(\omega_\FS)\wedge\ldots\wedge \Pi_{i_k}^*(\omega_\FS).$$ 

Fix a constant $\lambda$ such that $\rho<\lambda<1$. 
Let $I$ denote the set of multi-indices $i=(i_1,\ldots,i_k)$ in
$\{0,\ldots,n-1\}^k$ such that $i_s\geq n\lambda$ for every $s$. We
distinguish two cases where $i\not \in I$ or $i\in I$. 
In the first case, we have 
\begin{eqnarray*}
\sum_\alpha  \int_{\Gamma_n\cap \Uc_\alpha}
\Pi_{i_1}^*(\omega_\FS)\wedge\ldots\wedge \Pi_{i_k}^*(\omega_\FS)
& \leq &    \int_{\Gamma_n}
\Pi_{i_1}^*(\omega_\FS)\wedge\ldots\wedge \Pi_{i_k}^*(\omega_\FS)\\
& = & \int_{\P^k} (f^{i_1})^*(\omega_\FS)\wedge\ldots\wedge
(f^{i_k})^*(\omega_\FS)\\
& = & d^{i_1+\cdots+i_k}\leq d^{(k-1+\lambda)n}, 
\end{eqnarray*}
since $i_1+\cdots+i_k\leq (k-1+\lambda)n$. 

Consider the second case with multi-indices $i\in I$. 
Let $q$ denote the integer part of $\lambda n$ and
$W_\alpha$ the projection of $\Gamma_n\cap\Uc_\alpha$ on $\P^k$ by
$\Pi_0$. Observe that the choice of the open sets $U_i$ implies that
$f^q$ is injective on $W_\alpha$. Therefore,
\begin{eqnarray*}
\lefteqn{\sum_\alpha\int_{\Gamma_n\cap \Uc_\alpha}
  \Pi_{i_1}^*(\omega_\FS)\wedge\ldots\wedge
  \Pi_{i_k}^*(\omega_\FS)} \\
& \leq & \sum_\alpha \int_{W_\alpha} (f^{i_1})^*(\omega_\FS)\wedge\ldots\wedge
(f^{i_k})^*(\omega_\FS)\\
& = &  \sum_\alpha\int_{W_\alpha}
(f^q)^*\big[(f^{i_1-q})^*(\omega_\FS)\wedge\ldots\wedge
(f^{i_k-q})^*(\omega_\FS)\big]\\
& \leq &  \sum_\alpha\int_{\P^k}
(f^{i_1-q})^*(\omega_\FS)\wedge\ldots\wedge
(f^{i_k-q})^*(\omega_\FS).
\end{eqnarray*}  
Recall that the number of open sets $\Uc_\alpha$ is bounded by
$d^{k\rho n}$. So, the last sum is bounded by
$$d^{k\rho  n}d^{(i_1-q)+\cdots+(i_k-q)} \leq  d^{k\rho n}d^{k(n-q)}
 \lesssim  d^{k(1+\rho-\lambda)n}.$$

Finally, since the number of multi-indices $i$ is less than $n^k$, 
we deduce from the above estimates that
$$k!\vol(\Gamma_n\cap \Vc_n)\lesssim n^k
d^{(k-1+\lambda)n} +  n^k d^{k(1+\rho-\lambda)n}.$$
This contradicts the above bound from below of $\vol(\Gamma_n\cap \Vc_n)$.
\endproof

The remaining part of this paragraph deals with Lyapounov exponents
associated to the measure $\mu$ and their relations with the Hausdorff
dimension of $\mu$. Results in this direction give some 
information about the rough geometrical behaviour of the dynamical system on the
most chaotic locus. An abstract theory 
was developed by
Oseledec and Pesin, see e.g. \cite{KatokHasselblatt}. However, it is often difficult to show that a
given dynamical system has non-vanishing Lyapounov exponents. In
complex dynamics as we will see, the use of holomorphicity makes
the goal reachable. We first introduce few notions. 

Let $A$ be a linear endomorphism of $\R^k$. We can write 
$\R^k$ as the direct sum $\oplus E_i$ of invariant subspaces on which
all the complex eigenvalues of $A$ have the same modulus. This
decomposition of $\R^k$ describes clearly the
geometrical behaviour of the
dynamical system associated to $A$.
An important part in the 
dynamical study with respect to an invariant measure is to 
describe  geometrical aspects
following the directional dilation or contraction indicators.

Consider a smooth dynamical system $g:X\rightarrow X$ and an invariant
ergodic 
probability measure $\nu$. The map $g$ induces a linear map from the
tangent space at $x$ to the tangent space at $g(x)$. This linear map
is given by a square matrix when we fix local coordinates near $x$ and
$g(x)$. It is convenient in this setting to use local coordinates
depending smoothly on the point $x$. Then, we obtain a smooth function
on $X$ with values in $\GL(\R,k)$ where $k$ denotes the real dimension of $X$.
We will study the sequence of such functions associated to the
sequence of iterates $(g^n)$ of $g$.

Consider a more abstract setting. Let $g:X\rightarrow X$ be a
measurable map and $\nu$ an invariant probability measure. Let
$A:X\rightarrow \GL(\R,k)$ be a measurable function. Define for $n\geq
0$
$$A_n(x):=A(g^{n-1}(x))\ldots A(x).$$
These functions satisfy the identity
$$A_{n+m}(x)=A_n(g^m(x))A_m(x)$$
for $n,m\geq 0$. We say that the sequence $(A_n)$ is {\it the multiplicative cocycle} over $X$
{\it generated by $A$}. 

The following Oseledec's multiplicative ergodic theorem is related to
the Kingman's sub-multiplicative ergodic theorem
\cite{KatokHasselblatt, Walters}. It can be seen as a generalization
of the above property of a single square matrix $A$.

\begin{theorem}[Oseledec] \label{th_oseledec}
Let $g:X\rightarrow X$, $\nu$ and the
  cocycle $(A_n)$ be as above. Assume that $\nu$ is ergodic and that
  $\log^+\|A^{\pm 1}(x)\|$ are in
  $L^1(\nu)$. Then there is an integer
  $m$, real numbers $\chi_1<\cdots<\chi_m$, and for $\nu$-almost every $x$, a unique
  decomposition of $\R^k$ into a direct sum of linear subspaces
$$\R^k=\bigoplus_{i=1}^m E_i(x)$$ 
such that
\begin{enumerate}
\item The dimension of $E_i(x)$ does not depend on $x$.
\item The decomposition is invariant, that is, $A(x)$ sends $E_i(x)$
      to $E_i(g(x))$.
\item We have locally uniformly on vectors $v$ in $E_i(x)\setminus\{0\}$
$$\lim_{n\rightarrow\infty} {1\over n} \log \|A_n(x)\cdot v\|=\chi_i.$$
\item For $S\subset \{1,\ldots,m\}$, define $E_S(x):=\oplus_{i\in S}
      E_i(x)$. If $S,S'$ are disjoint, then the angle between $E_S(x)$
      and $E_{S'}(x)$ is a tempered function, that is,
$$\lim_{n\rightarrow\infty} {1\over n} \log\sin \big|\angle
\big(E_S(g^n(x)), E_{S'}(g^n(x))\big)\big|=0.$$ 
\end{enumerate}
\end{theorem}

The result is still valid for non-ergodic systems but the constants
$m$ and
$\chi_i$ should be replaced by invariant functions. If $g$ is
invertible, the previous decomposition is the same for $g^{-1}$ where
the exponents $\chi_i$ are replaced by $-\chi_i$. The result is also
valid in the complex setting where we replace $\R$ by $\C$ and
$\GL(\R,k)$ by $\GL(\C,k)$. In this case, the subspaces $E_i(x)$ are
complex.

We now come back to a smooth dynamical system $g:X\rightarrow X$ on a
compact manifold. We assume that the Jacobian $J(g)$ of $g$ associated
to a smooth volume form satisfies $\langle\nu,\log J(g)\rangle >-\infty$.
Under this hypothesis, we can
apply Oseledec's theorem to the cocycle induced by
$g$ on the tangent bundle of $X$; this allows to decompose, $\nu$-almost
everywhere, the tangent
bundle into invariant sub-bundles. The corresponding constants $\chi_i$
are called {\it Lyapounov exponents} of $g$ with respect to $\nu$. The
dimension of $E_i$ is {\it the multiplicity} of $\chi_i$. 
These notions do not depend on the choice of local coordinates on $X$.
The Lyapounov exponents of $g^n$ are equal to $n\chi_i$. 
We say that the measure $\nu$ is {\it hyperbolic} if no Lyapounov exponent is zero.
It is not difficult to deduce from the Oseledec's theorem that the sum
of Lyapounov exponents of $\nu$ is equal to $\langle\nu,\log
J(g)\rangle$. The reader will find in \cite{KatokHasselblatt} a theorem due
to Pesin, called {\it the $\epsilon$-reduction theorem}, which generalizes
Theorem \ref{th_oseledec}. It gives some coordinate changes on $\R^k$
which allow to write $A(x)$ in the form of a diagonal block
matrix with explicit estimates on the distortion.

The following result due to Briend-Duval \cite{BriendDuval1}, shows
that endomorphisms in $\P^k$ are expansive with
respect to the equilibrium measures. We give here a new proof using
Proposition \ref{prop_inverse_ball}. Note that there are $k$
Lyapounov exponents counted with multiplicity. If we consider these endomorphism as 
real maps, we should count twice the Lyapounov exponents.

\begin{theorem} \label{th_briend_duval}
Let $f$ be a holomorphic endomorphism of algebraic degree $d\geq
2$ of $\P^k$. Then the equilibrium measure $\mu$ of $f$ is
hyperbolic. More precisely, its Lyapounov exponents 
are larger or equal to ${1\over 2}\log d$.
\end{theorem}
\proof
Since the measure $\mu$ is PB, quasi-p.s.h. functions are
$\mu$-integrable. It is not difficult to check that if $J(f)$ is the Jacobian of
$f$ with respect to the Fubini-Study metric, then $\log J(f)$ is a
quasi-p.s.h. function. Therefore, we can apply Oseledec's theorem
\ref{th_oseledec}. We deduce from this result that the smallest Lyapounov
exponent of $\mu$ is equal to  
$$\chi:=\lim_{n\rightarrow \infty} -{1\over n} \log \|Df^n(x)^{-1}\|$$
for $\mu$-almost every $x$. By Proposition \ref{prop_inverse_ball}, there is a ball $B$
of positive $\mu$ measure which admits at least ${1\over 2}d^{kn}$ inverse
branches $g_i:B\rightarrow U_i$ for $f^n$ with $U_i$ of diameter $\leq
d^{-n/2}$. If we slightly reduce the ball $B$, we can assume that $\|Dg_i\|\leq
Ad^{-n/2}$ for some constant $A>0$. This is a simple consequence of
Cauchy's formula.
It follows that $\|(Df^n)^{-1}\|\leq
Ad^{-n/2}$ on $U_i$. The union $V_n$ of the $U_i$ is  
of measure at least equal to ${1\over 2}\mu(B)$. Therefore, by Fatou's lemma, 
$${1\over 2}\mu(B)\leq \limsup_{n\rightarrow \infty}
\langle\mu,\ind_{V_n}\rangle 
\leq \langle\mu,\limsup \ind_{V_n}\rangle = \langle
\mu, \ind_{\limsup V_n}\rangle.$$
Hence,
there is a set $K:=\limsup V_n$ of
positive measure such that if $x$ is in $K$, we have $\|Df^n(x)^{-1}\|\leq
Ad^{-n/2}$ for infinitely many of $n$. The result follows.
\endproof

Note that in a recent work \cite{deThelin2}, de Th\'elin proved that
for any invariant measure $\nu$ of entropy strictly larger than
$(k-1)\log d$, the Lyapounov exponents are strictly positive with some
explicit estimates from below.

\medskip

{\it The Hausdorff dimension} $\dim_H(\nu)$ of a probability measure
$\nu$ on $\P^k$ is the infimum  of the
numbers $\alpha\geq 0$ such that there is a Borel set $K$ of Hausdorff
dimension $\alpha$ of full measure, i.e. $\nu(K)=1$. Hausdorff
dimension says how
the measure fills out its support. The following result was obtained by
Binder-DeMarco \cite{BinderDeMarco} and Dinh-Dupont \cite{DinhDupont}.
The fact that $\mu$ has positive dimension has been proved
in \cite{Sibony}; indeed, a lower bound is given in terms of the
H\"older continuity exponent of the Green function $\gr$. 

\begin{theorem} \label{th_dinhdupont_pk}
Let $f$ be an endomorphism of algebraic degree $d\geq 2$ of $\P^k$ and
$\mu$ its equilibrium measure. Let $\chi_1,\ldots,\chi_k$ denote
the Lyapounov exponents of $\mu$ ordered by
$\chi_1\geq\cdots\geq\chi_k$ and $\Sigma$ their sum. Then
$${k\log d\over \chi_1}\leq \dim_H(\mu)\leq 2k-{2\Sigma-k\log d\over\chi_1}\cdot$$
\end{theorem}

The proof is quite technical. It is based on a delicate study
of the inverse branches of balls along a generic negative orbit. We
will not give the proof here. A better estimate in dimension 2,
probably the sharp one, was recently 
obtained by Dupont. Indeed, Binder-DeMarco conjecture
that the
Hausdorff dimension of $\mu$ satisfies 
$$\dim_H(\mu)={\log d\over\chi_1}+\cdots +{\log d\over\chi_k}\cdot$$
Dupont gives in \cite{Dupont3} results in this direction.

\bigskip\bigskip

\begin{exercise} \label{exo_gromov}
Let $X$ be an analytic subvariety of pure dimension $p$ in $\P^k$. Let
$f$ be an endomorphism of algebraic degree $d\geq 2$ of $\P^k$. Show
that $h_t(f,X)\leq p\log d$. 
\end{exercise}

\begin{exercise}
Let $f:X\rightarrow X$ be a smooth map and $K$ an invariant compact
subset of $X$. Assume that $K$ is hyperbolic, i.e. there is a
continuously varying decomposition $TX_{|K}=E\oplus F$ of the tangent bundle of $X$
restricted to $K$, into the sum of two invariant vector bundles such that
$\|Df\|<1$ on $E$ and $\|(Df)^{-1}\|<1$ on $F$ for some smooth metric
near $K$. Show that $f$ admits a
hyperbolic ergodic invariant measure supported on $K$.  
\end{exercise}

\begin{exercise}
Let $f:X\rightarrow X$ be a holomorphic map on a compact complex
manifold and let $\nu$ be an ergodic invariant measure. Show that 
 in Theorem \ref{th_oseledec} applied to the action of $f$ on the complex
tangent bundle, the spaces
$E_i(x)$ are complex.
\end{exercise}

\begin{exercise}
Let $\alpha>0$ be a constant. Show that there is an endomorphism $f$
of $\P^k$ such that the Hausdorff dimension of the equilibrium measure
of $f$ is smaller than $\alpha$. Show that there is an endomorphism
$f$ such that its Green function $\gr$ is not
$\alpha$-H{\"o}lder continuous.
\end{exercise}

\bigskip\bigskip

{\small

\addcontentsline{toc}{section}{Notes}
\noindent
{\bf Notes.} We do not give here results on local dynamics near a
fixed point. If this point is non-critical attractive or repelling, a
theorem of Poincar\'e says that the map is locally
conjugated to a polynomial map \cite{Sternberg}. Maps which are
tangent to the identity or are semi-attractive at a fixed point, were
studied by Abate and Hakim \cite{Abate,Hakim1,Hakim2, Hakim3}. Dynamics near a
super-attractive fixed point in dimension $k=2$ 
was studied by Favre-Jonsson using a
theory of valuations in \cite{FavreJonsson1}. 

The study of the dynamical system outside the support of the
equilibrium measure is not yet developped. Some results on attracting
sets, attracting currents, etc. were obtained by de Th\'elin, Dinh, Forn\ae ss, Jonsson,
Sibony, Weickert \cite{deThelin1, deThelin2, Dinh3,FornaessSibony6, FornaessWeickert,
  JonssonWeickert}, see also Mihailescu and Urba\'nski
\cite{Mihailescu,MihailescuUrbanski}.

In dimension 1, Fatou and Julia considered their theory as an
investigation to solve some functional equations. In particular, they found
all the commuting pairs of polynomials \cite{Fatou, Julia}, see also Ritt
\cite{Ritt} and Eremenko \cite{Eremenko} for the case of rational
maps. Commuting endomorphisms of $\P^k$
were studied by the authors in \cite{DinhSibony0}. A large family of
solutions are Latt\`es maps. We refer to Berteloot, Dinh, Dupont,
Loeb, Molino \cite{BertelootDupont, BertelootDupontMolino,
  BertelootLoeb,Dinh4, Dupont1} for a study of this class of maps, see also
Milnor \cite{Milnor} for the case of dimension 1.

We do not consider here bifurcation problems for families of maps and
refer to Bassanelli-Berteloot \cite{BassanelliBerteloot} and Pham \cite{Pham} for this subject. Some
results will be presented in the next chapter.

In \cite{Zhang}, Zhang considers some links between
complex dynamics and arithmetic questions. He is interested in polarized 
holomorphic maps on K\"ahler varieties, i.e. maps which multiply a
K\"ahler class by an integer. 
If the K\"ahler class is integral, the variety can be embedded into a
projective space
$\P^k$ and the maps extend to endomorphisms of $\P^k$. So, several
results stated above can be directely applied to that situation. In general,
most of the results for endomorphisms in $\P^k$ can be easily extended to
general polarized maps. In the unpublished preprint \cite{DinhSibony14}, the authors
considered the situation of smooth compact K\"ahler
manifolds. We recall here the main result. 

Let $(X,\omega)$ be an arbitrary compact K\"ahler manifold of dimension $k$. Let
$f$ be a holomorphic endomorphism of $X$. 
We assume that $f$ is open.
The spectral radius of $f^*$
acting on $H^{p,p}(X,\C)$ is called {\it the dynamical degree of order
  $p$} of $f$. It can be computed by the formula
$$d_p:=\lim_{n\rightarrow\infty}
\Big(\int_X(f^n)^*(\omega^p)\wedge\omega^{k-p}\Big)^{1/n}.$$
The last degree $d_k$ is the topological degree of $f$, i.e. equal to
the number of points in a generic fiber of $f$. We also denote it by
$d_t$. 

Assume that $d_t>d_p$ for $1\leq p\leq k-1$. Then, there is a maximal proper analytic
subset $\Ec$ of $X$ which is totally invariant by $f$,
i.e. $f^{-1}(\Ec)=f(\Ec)=\Ec$. If $\delta_a$ is a Dirac mass at 
$a\not\in \Ec$, then $d_t^{-n}(f^n)^*(\delta_a)$ converge to a
probability measure $\mu$, which does not depend on $a$. This is
{\it the equilibrium measure} of $f$. It satisfies $f^*(\mu)=d_t\mu$
and $f_*(\mu)=\mu$. If $J$ is the Jacobian of $f$ with respect to
$\omega^k$ then $\langle\mu,\log J\rangle \geq \log d_t$. The measure
$\mu$ is K-mixing and hyperbolic with Lyapounov exponents larger or
equal to ${1\over 2}\log(d_t/d_{k-1})$. Moreover, there are sets
$\Pc_n$ of repelling periodic points of order $n$, on $\supp(\mu)$ such that the
probability measures equidistributed on $\Pc_n$ converge to $\mu$, as
$n$ goes to infinity. If the periodic points of period $n$ are
isolated for every $n$, an estimate on the norm of $(f^n)^*$ on
$H^{p,q}(X,\C)$ obtained in \cite{Dinh2}, implies that the number of
these periodic points is $d_t^n+o(d_t^n)$. Therefore, periodic points
are equidistributed with respect to $\mu$. 
We can prove without difficulty that $\mu$ is the unique invariant measure of maximal entropy $\log
d_t$ and is moderate. Then, we can extend
the stochastic properties obtained for $\P^k$ to this more general setting.

When $f$ is polarized by the cohomology class $[\omega]$ of a K\"ahler
form $\omega$, there
is a constant $\lambda\geq 1$ such that $f^*[\omega]=\lambda[\omega]$.
It is not difficult to check that $d_p=\lambda^p$. The above results
can be applied for such a map when $\lambda>1$. In which case, periodic
points of a given period are isolated.
}


\chapter{Polynomial-like maps in higher dimension}

In this chapter we consider a large 
family of holomorphic maps in a semi-local setting: the polynomial-like
maps. They can appear as a
basic block in the study of some meromorphic maps on compact
manifolds. The main reference for this chapter is our article
\cite{DinhSibony1} where the $\ddc$-method in dynamics 
was introduced. Endomorphisms of $\P^k$ can be considered as a special
case of polynomial-like maps. However, in general, there is no Green
$(1,1)$-current for such maps. The notion of dynamical degrees for
polynomial-like maps replaces the algebraic degree.
Under natural assumptions on dynamical degrees, we prove that the
measure of maximal entropy is non-uniformly hyperbolic and we study
its sharp ergodic properties.

\section{Examples, degrees and entropy}

Let $V$ be a convex open set in $\C^k$ and $U\Subset V$ an open
subset. A proper holomorphic map $f:U\rightarrow V$ is called {\it
a  polynomial-like map}. Recall that a map $f:U\rightarrow V$ is proper
if $f^{-1}(K)\Subset U$ for every compact
subset $K$ of $V$.  
The map $f$ sends the boundary of $U$ to the boundary of $V$; more
precisely, the points near $\partial U$ are sent to points near
$\partial V$. So, polynomial-like maps are somehow expansive
in all directions, but the expansion is in the geometrical sense. In
general, they may have a non-empty critical set. 
A polynomial-like mapping
$f:U\rightarrow V$ defines a ramified covering over $V$. The degree
$d_t$ of this covering is also called {\it the topological degree}. It is
equal to the number of points in a generic fiber, or in any fiber if
we count the multiplicity. 

Polynomial-like maps are characterized by the property that their graph $\Gamma$
in $U\times V$ is in fact a submanifold of $V\times V$, that is,
$\Gamma$ is closed in $V\times V$. So, any small perturbation of $f$ is
polynomial-like of the same topological degree $d_t$, provided that we
reduce slightly the open set $V$. We will construct large families
of polynomial-like maps. In dimension one, it was proved by
Douady-Hubbard \cite{DouadyHubbard} that
such a map is conjugated to a polynomial via a H{\"o}lder continuous homeomorphism.
Many dynamical properties can be deduced from the corresponding
properties of polynomials.
In higher dimension, the analogous statement is not valid. 
Some new dynamical phenomena appear for polynomial-like mappings, that
do not exist for polynomial maps.
We use here an approach completely different from the
one dimensional case, where the basic tool is the Riemann measurable
mapping theorem.

Let $f:\C^k\rightarrow\C^k$ be a holomorphic map such that the
hyperplane at infinity is attracting in the sense that $\|f(z)\|\geq A\|z\|$
for some constant $A>1$ and for $\|z\|$ large enough. If $V$ is a
large ball centered at 0, then  $U:=f^{-1}(V)$ is strictly
contained in $V$. Therefore, $f:U\rightarrow V$ is a polynomial-like
map. Small transcendental perturbations of $f$, as we mentioned above,
give a large family of polynomial-like maps. Observe also that
the dynamical study of holomorphic endomorphisms on $\P^k$ can be
reduced to polynomial-like maps by lifting to a large ball in
$\C^{k+1}$. We give now other explicit examples.

\begin{example} \rm \label{example_poly_polylike}
Let $f=(f_1,\ldots,f_k)$ be a polynomial map in $\C^k$, with $\deg
f_i=d_i\geq 2$. Using a conjugacy with a permutation of coordinates, we
can assume that $d_1\geq \cdots\geq d_k$. Let $f_i^+$ denote the
homogeneous polynomial of highest degree in $f_i$. If
$\{f_1^+=\cdots=f_k^+=0\}$ is reduced to $\{0\}$, then $f$ is
polynomial-like in any large ball of center 0. Indeed, define $d:=d_1\ldots
d_k$ and $\pi(z_1,\ldots,z_k):=(z_1^{d/d_1},
\ldots,z_k^{d/d_k})$. Then, $\pi\circ f$ is a polynomial map of
algebraic degree $d$ which extends holomorphically at infinity to an
endomorphism of $\P^k$. Therefore, $\|\pi\circ f(z)\|\gtrsim \|z\|^d$
for $\|z\|$ large enough. The estimate $\|f(z)\|\gtrsim \|z\|^{d_k}$
near infinity
follows easily. If we consider the extension of $f$ to $\P^k$, we
obtain in general a meromorphic map which is not holomorphic. Small
pertubations $f_\epsilon$ of $f$ may have indeterminacy points in
$\C^k$ and a priori, indeterminacy points of the sequence
$(f_\epsilon^n)_{n\geq 1}$ may be dense in $\P^k$.
\end{example}

\begin{examples} \rm
The map $(z_1,z_2)\mapsto (z_1^2+az_2,z_1)$, $a\not=0$, is not
polynomial-like. It is invertible and the point $[0:0:1]$ at infinity,
in homogeneous coordinates $[z_0:z_1:z_2]$, is an
attractive fixed point for $f^{-1}$. Hence, the set $\Kc$ of points
$z\in\C^2$ with bounded orbit, clusters at $[0:0:1]$. 

The map
$f(z_1,z_2):=(z_2^d,2z_1)$, $d\geq 2$, is polynomial-like in any large
ball of center $0$. Considered as a map on $\P^2$, it is only meromorphic
with an indeterminacy point $[0:1:0]$. On a fixed large ball of
center 0, the perturbed maps $f_\epsilon:=(z_2^d+\epsilon
e^{z_1},2z_1+\epsilon e^{z_2})$ are polynomial-like, in an appropriate
open set $U$.
\end{examples}

Consider a general polynomial-like map $f:U\rightarrow V$ of
topological degree $d_t\geq 2$. 
We introduce several growth indicators of the action of $f$ on forms
or currents.
Define $f^n:=f\circ\cdots\circ f$, $n$ times, {\it the iterate of order $n$} of $f$.
This map is only defined on $U_{-n}:=f^{-n}(V)$. The sequence $(U_{-n})$
is decreasing: we have $U_{-n-1}=f^{-1}(U_{-n})\Subset U_{-n}$. Their intersection
$\Kc:=\cap_{n\geq 0} U_{-n}$ is a non-empty compact set that we call
{\it the filled Julia set} of $f$. The filled Julia set is totally
invariant: we have $f^{-1}(\Kc)=\Kc$ which implies that $f(\Kc)=\Kc$. 
Only for $x$ in $\Kc$, the infinite orbit
$x,f(x),f^2(x),\ldots$ is well-defined. The preimages $f^{-n}(x)$ by $f^n$ are
defined for every $n\geq 0$ and every $x$ in $V$. 

Let $\omega:=\ddc\|z\|^2$ denote the standard K{\"a}hler form on $\C^k$. 
Recall that the mass of a positive $(p,p)$-current $S$ on a Borel set
$K$ is  given by $\|S\|_K:=\int_K S\wedge\omega^{k-p}$. 
Define the {\it dynamical degree of order $p$} of $f$, for $0\leq
p\leq k$, by
$$d_p(f):=\limsup_{n\rightarrow\infty}\|(f^n)_*(\omega^{k-p})\|_W^{1/n}=\limsup_{n\rightarrow\infty}
\|(f^n)^*(\omega^p)\|_{f^{-n}(W)}^{1/n},$$
where $W\Subset V$ is a neighbourhood of $\Kc$. For simplicity, when there is no
confusion, this degree is also denoted by $d_p$. We have the following lemma.

\begin{lemma} \label{lemma_degree_well_def}
The degrees $d_p$ do not depend on the choice of $W$. Moreover, we
have $d_0\leq 1$, $d_k=d_t$ and the dynamical degree of order $p$ of
$f^m$ is equal to $d_p^m$. 
\end{lemma}
\proof
Let $W'\subset W$ be another neighbourhood of $\Kc$. For the first
assertion, we only have to show that 
$$\limsup_{n\rightarrow\infty}
\|(f^n)^*(\omega^p)\|_{f^{-n}(W)}^{1/n} \leq \limsup_{n\rightarrow\infty}
\|(f^n)^*(\omega^p)\|_{f^{-n}(W')}^{1/n}.$$
By definition of $\Kc$, there is an integer $N$ such that
$f^{-N}(V)\Subset W'$. Since $(f^N)^*(\omega^p)$ is smooth on
$f^{-N}(V)$, we can find a constant $A>0$ such that
$(f^N)^*(\omega^p)\leq A\omega^p$ on $f^{-N}(W)$. We have
\begin{eqnarray*}
\limsup_{n\rightarrow\infty}
\|(f^n)^*(\omega^p)\|_{f^{-n}(W)}^{1/n} & = &  \limsup_{n\rightarrow\infty}
\|(f^{n-N})^*\big((f^N)^*(\omega^p)\big)\|_{f^{-n+N}(f^{-N}(W))}^{1/n}\\
&\leq &  \limsup_{n\rightarrow\infty}
\|A(f^{n-N})^*(\omega^p)\|_{f^{-n+N}(W')}^{1/n}\\
& = & \limsup_{n\rightarrow\infty}
\|(f^n)^*(\omega^p)\|_{f^{-n}(W')}^{1/n}.
\end{eqnarray*}
This proves the first assertion. 

It is clear from the definition that $d_0\leq 1$. Since $f$ has topological
degree $d_t$ the pull-back of a positive measure multiplies the
mass by $d_t$. Therefore, $d_k=d_t$. For the last assertion of the
lemma, we only have to check that 
$$\limsup_{n\rightarrow\infty}
\|(f^n)^*(\omega^p)\|_{f^{-n}(W)}^{1/n}\leq \limsup_{s\rightarrow\infty}
\|(f^{ms})^*(\omega^p)\|_{f^{-ms}(W)}^{1/ms}.$$
To this end, we proceed as above. 
Write $n=ms+r$ with $0\leq r\leq
m-1$. We obtain the result using that $(f^r)^*(\omega^p)\leq A\omega^p$ on a
fixed neighbourhood of $\Kc$ for $0\leq r\leq m-1$.
\endproof

The main result of this paragraph is the following formula for the entropy.

\begin{theorem} \label{th_entropy_top_pol_like}
Let $f:U\rightarrow V$ be a polynomial-like map of topological degree
$d_t\geq 2$. Let $\Kc$ be the filled Julia set of $f$. Then, the topological
entropy of $f$ on $\Kc$ is equal to $h_t(f,\Kc)=\log d_t$. Moreover,
all the dynamical degrees $d_p$ of $f$ are smaller or equal to
$d_t$. 
\end{theorem}

We need the following lemma where we use standard metrics on Euclidean
spaces.

\begin{lemma} \label{lemma_volum_cn}
Let $V$ be an open set of $\C^k$, $U$ a relatively compact subset
of $V$ and $L$ a compact subset of $\C$.   Let $\pi$ denote the canonical projection
from $\C^m\times V$ onto $V$.
Suppose $\Gamma$ is an analytic subset of pure dimension $k$ of $\C^m\times V$
contained in $L^m\times V$. 
Assume also that $\pi:\Gamma\rightarrow V$
defines a ramified covering  of degree 
$d_\Gamma$.
Then, there exist constants $c>0$, $s>0$, independent of $\Gamma$ and  $m$, such that
$$\vol(\Gamma\cap \C^m\times U)\leq cm^sd_\Gamma.$$ 
\end{lemma}
\proof
Since the problem is local on $V$, we can assume that $V$ is the unit
ball of $\C^k$ and $U$ is the closed ball of
center 0 and of radius $1/2$. We can also assume that $L$ is the
closed unit disc in $\C$. 
Denote by $x=(x_1,\ldots , x_m)$ and $y=(y_1,\ldots,y_k)$ the
coordinates on $\C^m$ and on $\C^k$. Let $\epsilon$ be a 
$k\times m$ matrix whose entries have modulus bounded by 
$1/8mk$. Define $\pi_\epsilon(x,y):=y+\epsilon x$,
$\Gamma_\epsilon:=\Gamma\cap \{\|\pi_\epsilon\|< 3/4\}$ and
$\Gamma^*:=\Gamma\cap (L^m\times U)$.

We first show that $\Gamma^*\subset\Gamma_\epsilon$. Consider a point 
$(x,y)\in \Gamma^*$. We have $|x_i|< 1$ and $\|y\|\leq 1/2$. Hence,
$$\|\pi_{\epsilon}(x,y)\|\leq \|y\| +\|\epsilon x\|< 3/4.$$
This implies that $(x,y)\in\Gamma_\epsilon$.

Now, we prove that for every 
$a\in\C^k$ with
$\|a\|<3/4$, we have $\#\pi_\epsilon^{-1}(a)\cap \Gamma =d_\Gamma$,
where we count the multiplicities of points.
To this end, we show that $\#\pi_{t\epsilon}^{-1}(a)\cap \Gamma$
does not depend on $t\in [0,1]$. So, it is sufficient to 
check that the union of the sets $\pi_{t\epsilon}^{-1}(a)\cap \Gamma$
is contained in the compact subset
$\Gamma\cap\{\|\pi\|\leq 7/8\}$ of 
$\Gamma$. Let $(x,y)\in\Gamma$ and
$t\in[0,1]$ such that $\pi_{t\epsilon}(x,y)=a$. We have
$$3/4> \|a\| =\|\pi_{t\epsilon}(x,y)\| \geq\|y\| -
t\|\epsilon x\|.$$
It follows that $\|y\|<7/8$ and hence $(x,y)\in 
\Gamma\cap\{\|\pi\|\leq 7/8\}$.

Let $B$ denote the ball of center 0 and of radius 
$3/4$ in $\C^k$. We have for some constant
$c'>0$
$$\int_{\Gamma^*} \pi_\epsilon^*(\omega^k) \leq
\int_{\Gamma_\epsilon} \pi_\epsilon^*(\omega^k) = d_\Gamma \int_B
\omega^k =c'd_\Gamma.$$
Let $\Theta:=\ddc\|x\|^2+\ddc\|y\|^2$ be the standard K{\"a}hler $(1,1)$-form on $\C^m\times
\C^k$. We have
$$\vol(\Gamma\cap\C^m\times U)=\int_{\Gamma\cap\C^m\times U} \Theta^k.$$
It suffices to
bound $\Theta$ by a linear combination of
$2m+1$ forms of type $\pi_\epsilon^*(\omega)$ with coefficients of
order $\simeq m^2$ and then to use the previous estimates. Recall that $\omega=\ddc\|y\|^2$. So, we
only have to bound
$\sqrt{-1}d x_i \wedge d \overline x_i$ by a combination of 
$(1,1)$-forms of type 
$\pi_\epsilon^*(\omega)$. Consider
$\delta:=1/8mk$ and $\pi_\epsilon(x,y):=(y_1+\delta x_i,
y_2,\ldots,y_k)$. We have
\begin{eqnarray*}
\sqrt{-1}d x_i \wedge d \overline x_i & = & \frac{4\sqrt{-1}}{3\delta^2} 
\Big[3d y_1\wedge d \overline y_1 + d (y_1+\delta x_i) \wedge
d (\overline y_1+\delta \overline x_i) \\
& & - d (2y_1+\delta x_i/2) \wedge
d (2\overline y_1+\delta \overline x_i/2)\Big] \\
& \leq &  \frac{4\sqrt{-1}}{3\delta^2} 
\Big[3dy_1\wedge d \overline y_1 + d (y_1+\delta x_i) \wedge
d (\overline y_1+\delta \overline x_i)\Big] 
\end{eqnarray*}
The last form can be bounded by a combination of $\pi_0^*(\omega)$ and
$\pi_\epsilon^*(\omega)$.
This completes the proof.
\endproof

\noindent
{\bf Proof of Theorem \ref{th_entropy_top_pol_like}.} 
We prove that $h_t(f,\Kc)\leq \log d_t$. We will prove in Paragraphs
\ref{section_poly-like-construction} and \ref{section_poly-like-entropy} 
that $f$ admits a totally invariant measure of maximal entropy
$\log d_t$ with support in the boundary of $\Kc$. The variational
principle then implies that $h_t(f,\Kc)=h_t(f,\partial\Kc)=\log d_t$.
We can also conclude using Misiurewicz-Przytycky's theorem 
\ref{th_MisiurewiczPrzytycky} or 
Yomdin's theorem \ref{th_yomdin} which can be extended to this case. 

Let $\Gamma_n$ denote the graph of
$(f,\ldots,f^{n-1})$ in $V^n\subset (\C^k)^{n-1}\times V$. 
Let $\pi:(\C^k)^{n-1}\times V\rightarrow V$ be the canonical
projection. Since $f:U\rightarrow V$ is polynomial-like, it is easy to
see that $\Gamma_n\subset U^{n-1}\times V$ and that
$\pi:\Gamma_n\rightarrow V$ defines a ramified covering of degree $d_t^n$.
As in Theorem
\ref{th_gromov_Pk}, we have
$$h_t(f,\Kc)\leq \lov(f):=\limsup_{n\rightarrow \infty} {1\over n}
\log\vol(\Gamma_n\cap \pi^{-1}(U)).$$
But, it follows from Lemma
\ref{lemma_volum_cn} that 
$$\vol(\Gamma_n\cap \pi^{-1}(U))\leq c (kn)^s d_t^n.$$
Hence, $\lov(f)\leq \log d_t$. This implies the
inequality $h_t(f,\Kc)\leq \log d_t$. Note that the limit in the definition of $\lov(f)$
exists and we have $\lov(f)=\log d_t$. Indeed, since $\Gamma_n$ is a
covering of degree $d_t^n$ over $V$, we always have 
$$\vol(\Gamma_n\cap
\pi^{-1}(U))\geq d_t^n\vol(U).$$ 

We show that $d_p\leq d_t$. 
Let $\Pi_i$, $0\leq i\leq n-1$, denote the projection of $V^n$ onto its
factors. We have
$$\vol(\Gamma_n\cap\pi^{-1}(U))=\sum_{0\leq i_s\leq n-1}
\int_{\Gamma_n\cap\pi^{-1}(U)} \Pi_{i_1}^*(\omega)\wedge \ldots\wedge
\Pi_{i_k}^*(\omega).$$
The last sum contains the term
$$\int_{\Gamma_n\cap\pi^{-1}(U)} \Pi_0^*(\omega^{k-p})\wedge
\Pi_{n-1}^*(\omega^p)
=\int_{f^{-n+1}(U)} \omega^{k-p}\wedge (f^{n-1})^*(\omega^p).$$
We deduce from the estimate on
$\vol(\Gamma_n\cap\pi^{-1}(U))$ and from the definition of $d_p$ that $d_p\leq \lov(f)=d_t$. 
\hfill $\square$

\medskip

We introduce now others useful dynamical degrees. We call {\it dynamical
  $\ast$-degree of order $p$} of $f$ the following limit
$$d_p^*:=\limsup_{n\rightarrow\infty}  \sup_S\|(f^n)_*(S)\|_W^{1/n},$$
where $W\Subset V$ is a neighbourhood of $\Kc$ and the supremum is
taken over positive closed $(k-p,k-p)$-current of mass $\leq 1$ on a fixed
neighbourhood  $W'\Subset V$ of $\Kc$. Clearly, $d_p^*\geq d_p$,
since we can take $S=c\omega^k$ with $c>0$ small enough. 

\begin{lemma} \label{lemma_deg_star_def}
The above definition does not depend on the choice of
  $W$, $W'$. Moreover, we have $d_0^*=1$, $d_k^*=d_t$ and the dynamical
  $\ast$-degree of order $p$ of $f^n$ is equal to ${d_p^*}^n$.
\end{lemma}
\proof
If $N$ is an integer large enough, the operator $(f^N)_*$ sends
continuously positive closed currents on $W'$ to the ones on
$V$. Therefore, the independence of the definition on $W'$ is clear.
If $S$
is a probability measure on $\Kc$, then $(f^n)_*$ is also a
probability measure on $\Kc$. Therefore, $d_0^*=1$.
Observe that
$$\|(f^n)_*(S)\|_W^{1/n}=\Big[\int_{f^{-n}(W)} S\wedge
(f^n)^*(\omega^p)\Big]^{1/n}.$$
So, for the other properties, it is enough to follow the arguments given 
in Lemma \ref{lemma_degree_well_def}.
\endproof

Many results below are proved  under the hypothesis $d^*_{k-1} <
d_t$. The following proposition shows that this condition is stable under
small perturbations on $f$. This gives large families of maps
satisfying the hypothesis. Indeed, the condition is satisfied for
polynomial maps in $\C^k$ which extend at infinity as endomorphisms
of $\P^k$. For such maps, if $d$ is the algebraic degree, 
one can check that $d_p^*\leq d^p$. 

\begin{proposition} \label{prop_polylike_variation}
Let $f:U\rightarrow V$ be a polynomial-like map of
  topological degree $d_t$. Let $V'$ be a convex open set such that
  $U\Subset V'\Subset V$. If $g:U\rightarrow\C^k$ is a holomorphic
  map, close enough to $f$ and $U':=g^{-1}(V')$, then $g:U'\rightarrow
  V'$ is a polynomial-like map of topological degree $d_t$. If
  moreover, $f$ satisfies the condition $d^*_p<d_t$ for some $1\leq
  p\leq k-1$, then
  $g$ satisfies the same property.
\end{proposition}
\proof
The first assertion is clear, using the characterization of
polynomial-maps by their graphs. We prove the second one. Fix a
constant $\delta$ with $d_p^*<\delta<d_t$ and an open set $W$ such
that $U\Subset W\Subset V$.
Fix an integer $N$ large enough such that $\|(f^N)_*(S)\|_W\leq
\delta^N$ for any positive closed $(k-p,k-p)$-current $S$ of mass 1 on $U$. If $g$ is close enough
to $f$, we have $g^{-N}(U)\Subset f^{-N}(W)$ and 
$$\|(g^N)^*(\omega^p)-(f^N)^*(\omega^p)\|_{L^\infty(g^{-N}(U))}\leq
\epsilon$$ 
with $\epsilon>0$ a small constant. We have
\begin{eqnarray*}
\|(g^N)_*(S)\|_U & = & \int_{g^{-N}(U)} S\wedge (g^N)^*(\omega^p)\\
& \leq &  \int_{f^{-N}(W)} S\wedge (f^N)^*(\omega^p) +  \int_{g^{-N}(U)}
S\wedge \big[(g^N)^*(\omega^p)-(f^N)^*(\omega^p)\big]\\
& \leq & \|(f^N)_*(S)\|_W+\epsilon\leq \delta^N+\epsilon<d_t^N.
\end{eqnarray*}
Therefore, the dynamical $\ast$-degree $d_p^*(g^N)$ of $g^N$ is
strictly smaller than $d_t^N$. Lemma \ref{lemma_deg_star_def} implies
that $d_p^*(g)<d_t$. 
\endproof

\begin{remark}\rm
The proof gives that $g\mapsto d_p^*(g)$ is upper semi-continuous on
$g$.
\end{remark}

Consider a simple example.
Let $f:\C^2\rightarrow \C^2$ be the polynomial map
$f(z_1,z_2)=(2z_1,z_2^2)$. The restriction of $f$ to
$V:=\{|z_1|<2,|z_2|<2\}$ is polynomial-like and using the current $S=[z_1=0]$, 
it is not difficult to check that
$d_1=d_1^*=d_t=2$. The example shows that in general one may have $d_{k-1}^*=d_t$.

\bigskip\bigskip

\begin{exercise}
Let $f:\C^2\rightarrow\C^2$ be the polynomial map defined by
$f(z_1,z_2):=(3z_2,z_1^2+z_2)$. Show that the hyperplane at infinity
is attracting. Compute the topological degree of $f$. Compute the
topological degree of the map in Example \ref{example_poly_polylike}.
\end{exercise}

\begin{exercise}
Let $f$ be a polynomial map on $\C^k$ of algebraic degree $d\geq
2$, which extends to a holomorphic endomorphism of $\P^k$. Let $V$ be
a ball large enough centered at $0$ and $U:=f^{-1}(V)$. Prove that
 the polynomial-like map $f:U\rightarrow V$ satisfies
 $d_p^*= d^p$ and $d_t=d^k$. Hint: use the Green function and Green
 currents. 
\end{exercise}


\section{Construction of the Green measure} \label{section_poly-like-construction}

In this paragraph, we introduce the first version of the
$\ddc$-method. It allows to 
construct for a polynomial-like map $f$ a canonical measure which
is totally invariant. 
As we have seen  in the case of endomorphisms of
$\P^k$, the method gives good estimates and allows to obtain precise
stochastic properties. Here, we will see that it 
applies under a very weak hypothesis. The construction of the measure
does not require any hypothesis on the dynamical degrees and give useful convergence results.

Consider a polynomial-like map $f:U\rightarrow V$ of topological
degree $d_t>1$ as above. Define the Perron-Frobenius operator
$\Lambda$ acting on test functions $\varphi$ by
$$\Lambda(\varphi)(z):=d_t^{-1}f_*(\varphi)(z):=d_t^{-1}\sum_{w\in
  f^{-1}(z)} \varphi(w),$$
where the points in $f^{-1}(z)$ are counted with multiplicity.
Since $f$ is a ramified covering, $\Lambda(\varphi)$ is continuous
when $\varphi$ is continuous. If $\nu$ is a probability measure on $V$,
define the measure $f^*(\nu)$ by
$$\langle f^*(\nu),\varphi\rangle := \langle \nu,
f_*(\varphi)\rangle.$$
This is a positive measure of mass $d_t$ supported on
$f^{-1}(\supp(\nu))$. Observe that the operator $\nu\mapsto d_t^{-1}
f^*(\nu)$ is continuous on positive measures, see Exercise \ref{Exo_push_c0}.

\begin{theorem} \label{th_polylike_measure}
Let $f:U\rightarrow V$ be a polynomial-like map as above. Let
$\nu$ be a probability measure supported on $V$ which is defined by an
$L^1$ form. Then $d_t^{-n}(f^n)^*(\nu)$ converge to a probability
measure $\mu$ which does not depend on $\nu$. The measure $\mu$ is
supported on the boundary of the filled Julia set $\Kc$ and is totally
invariant: $d_t^{-1}f^*(\mu)=f_*(\mu)=\mu$. Moreover, if  $\Lambda$ is 
the Perron-Frobenius operator associated to $f$ and $\varphi$ is a
p.s.h. function on a neighbourhood of $\Kc$, then $\Lambda^n(\varphi)$
converge to $\langle \mu,\varphi\rangle$. 
\end{theorem}

Note that in general $\langle\mu,\varphi\rangle$ may be $-\infty$. 
If $\langle \mu,\varphi\rangle
=-\infty$, the above convergence means that $\Lambda^n(\varphi)$ tend locally uniformly to
$-\infty$; otherwise, the convergence is in $L^p_\loc$ for $1\leq
p<+\infty$, see Appendix \ref{section_positive}.
The
above result still holds for measures $\nu$ which have no mass on
pluripolar sets. The
proof in that case is more delicate. We have the following lemma.

\begin{lemma} \label{lemma_cv_perron_polylike}
If $\varphi$ is p.s.h. on a neighbourhood of $\Kc$, then
  $\Lambda^n(\varphi)$ converge to a constant $c_\varphi$ in $\R\cup\{-\infty\}$.
\end{lemma}
\proof
Observe that $\Lambda^n(\varphi)$ is defined on $V$ for $n$ large
enough. It is not difficult to check that these functions are
p.s.h. Indeed, when $\varphi$ is a continuous p.s.h. function,
$\Lambda^n(\varphi)$ is a continuous function, see Exercise \ref{Exo_push_c0}, and
$\ddc\Lambda^n(\varphi)=d_t^{-n} (f^n)_*(\ddc\varphi)\geq 0$. So,
$\Lambda^n(\varphi)$ is p.s.h. The general case is obtained using an
approximation of $\varphi$ by a decreasing sequence of smooth p.s.h. functions.

Consider $\psi$ the upper semi-continuous regularization of
$\limsup \Lambda^n(\varphi)$. We deduce from Proposition
\ref{prop_hartogs_else} that $\psi$ is a p.s.h. function. We first prove
that $\psi$ is constant.
Assume not. By maximum principle, there is a
constant $\delta$ such that 
$\sup_{\overline U} \psi<\delta<\sup_V\psi$. By Hartogs' lemma \ref{prop_hartogs_else}, for $n$ large enough, we
have $\Lambda^n(\varphi)<\delta$ on $U$. Since the fibers of $f$ are
contained in $U$, we deduce from the definition of
$\Lambda$ that
$$\sup_V\Lambda^{n+1}(\varphi)=\sup_V\Lambda(\Lambda^n(\varphi))\leq \sup_U\Lambda^n(\varphi)<\delta.$$
This implies that $\psi\leq \delta$ which contradicts the choice of $\delta$. So $\psi$ is constant.

Denote by $c_\varphi$ this constant.
If $c_\varphi=-\infty$, it is clear that $\Lambda^n(\varphi)$ converge
to $-\infty$ uniformly on compact sets.
Assume that $c_\varphi$ is finite and $\Lambda^{n_i}(\varphi)$ does not converge to $c_\varphi$
for some sequence $(n_i)$. By Hartogs' lemma, we have
$\Lambda^{n_i}(\varphi)\leq c_\varphi-\epsilon$ for some constant
$\epsilon>0$ and for $i$ large enough. We deduce as above that
$\Lambda^n(\varphi)\leq c_\varphi-\epsilon$ for $n\geq n_i$. This
contradicts the definition of $c_\varphi$. 
\endproof

\noindent
{\bf Proof of Theorem \ref{th_polylike_measure}.} We can replace $\nu$
by $d_t^{-1}f^*(\nu)$ in order to assume that $\nu$ is supported on $U$.
The measure $\nu$ can be written as a finite or countable sum of bounded
positive forms, we can assume that $\nu$ is a bounded form.

Consider a smooth p.s.h. function $\varphi$ on a neighbourhood of $\Kc$. It is clear that
$\Lambda^n(\varphi)$ are uniformly bounded for $n$ large enough. Therefore,
the constant $c_\varphi$ is finite. We deduce from Lemma \ref{lemma_cv_perron_polylike} that
$\Lambda^n(\varphi)$ converge in $L^1_\loc(V)$ to $c_\varphi$. It
follows that 
$$\langle d_t^{-n}(f^n)^*(\nu),\varphi\rangle = \langle
\nu,\Lambda^n(\varphi)\rangle \rightarrow c_\varphi.$$

Let $\phi$ be a general smooth function on $V$. We can always write
$\phi$ as a difference of p.s.h. functions on $U$. Therefore,
$\langle d_t^{-n}(f^n)^*(\nu),\phi\rangle$ converge. It follows that
the sequence of probability measures $d_t^{-n}(f^n)^*(\nu)$
converges to some probability measure $\mu$. Since $c_\varphi$ does
not depend on $\nu$, the measure $\mu$ does not depend on
$\nu$. Consider a measure $\nu$ supported on $U\setminus\Kc$. So, 
the limit $\mu$ of $d_t^{-n}(f^n)^*(\nu)$ is supported on
$\partial\Kc$. The total invariance is a direct consequence of the
above convergence.

For the rest of the theorem, assume that $\varphi$ is a general
p.s.h. function on a neighbourhood of $\Kc$. Since  $\limsup\Lambda^n(\varphi)\leq c_\varphi$, 
Fatou's lemma implies that
$$\langle \mu,\varphi\rangle = \langle
d_t^{-n}(f^n)^*(\mu),\varphi\rangle
= \langle\mu,\Lambda^n(\varphi)\rangle \leq
\langle\mu,\limsup_{n\rightarrow\infty}\Lambda^n(\varphi)\rangle= c_\varphi.$$
On the other hand, for $\nu$ smooth on $U$, we have since $\varphi$ is
upper semi-continuous
$$c_\varphi= \lim_{n\rightarrow\infty} \langle \nu,\Lambda^n(\varphi)\rangle = \lim_{n\rightarrow\infty} \langle
d_t^{-n}(f^n)^*(\nu),\varphi\rangle \leq \langle \lim_{n\rightarrow\infty}
d_t^{-n}(f^n)^*(\nu),\varphi\rangle =\langle\mu,\varphi\rangle.$$
Therefore, $c_\varphi=\langle \mu,\varphi\rangle$. Hence,
$\Lambda^n(\varphi)$ converge to $\langle\mu,\varphi\rangle$ for an
arbitrary p.s.h. function $\varphi$.
\hfill $\square$

\medskip

The measure $\mu$ is called {\it the equilibrium measure} of $f$.
We deduce from the above arguments the following result.

\begin{proposition} \label{prop_test_ph_polylike}
Let $\nu$ be a totally invariant probability measure. 
Then $\nu$
is supported on $\Kc$. Moreover, $\langle \nu,\varphi\rangle\leq
\langle\mu,\varphi\rangle$ for every function $\varphi$ which is
p.s.h. in a neighbourhood of $\Kc$ and $\langle \nu,\varphi\rangle =
\langle\mu,\varphi\rangle$ if $\varphi$ is pluriharmonic in a
neighbourhood of $\Kc$. 
\end{proposition} 
\proof
Since $\nu=d_t^{-n}(f^n)^*(\nu)$, it is supported on $f^{-n}(V)$ for
every $n\geq 0$. So, $\nu$ is supported on $\Kc$. 
We know that
$\limsup\Lambda^n(\varphi)\leq c_\varphi$, then 
Fatou's lemma implies that
$$\langle \nu,\varphi\rangle =\lim_{n\rightarrow\infty} \langle
d_t^{-n}(f^n)^*(\nu),\varphi\rangle
=\lim_{n\rightarrow\infty} \langle\nu,\Lambda^n(\varphi)\rangle \leq c_\varphi.$$
When $\varphi$ is pluriharmonic, the inequality holds for
$-\varphi$; we then deduce that $\langle\nu,\varphi\rangle \geq
c_\varphi$. The proposition follows.
\endproof

\begin{corollary}
Let $X_1,X_2$ be two analytic subsets of $V$ such that
$f^{-1}(X_1)\subset X_1$ and $f^{-1}(X_2)\subset X_2$. Then $X_1\cap
X_2\not=\varnothing$. In particular, $f$ admits at most one 
point $a$ such that $f^{-1}(a)=\{a\}$. 
\end{corollary}
\proof
Observe that there are totally invariant probability measures $\nu_1,\nu_2$
supported on $X_1,X_2$. Indeed, if $\nu$ is a probability measure
supported on $X_i$, then any limit value of 
$${1\over N}\sum_{n=0}^{N-1}d_t^{-n} (f^n)^*(\nu)$$
is supported on $X_i$ and is totally invariant. We are using the
continuity of the operator $\nu\mapsto d_t^{-1}f^*(\nu)$, for the weak
topology on measures. 

On the other hand, if $X_1$ and $X_2$ are disjoint, we can find a
holomorphic function $h$ on $U$ such that $h=c_1$ on $X_1$ and $h=c_2$ on
$X_2$, where $c_1,c_2$ are distinct constants. We consider the function defined on $X_1\cup X_2$ as claimed and
extend it as a holomorphic function in $U$. This is possible since $V$
is convex \cite{Hormander2}. 
Adding to $h$ a constant
allows to assume that $h$ does not vanish on $\Kc$. Therefore,
$\varphi:=\log|h|$ is pluriharmonic on a neighbourhood of $\Kc$. We have 
$\langle\nu_1,\varphi\rangle\not =\langle\nu_2,\varphi\rangle$.
This contradicts Proposition \ref{prop_test_ph_polylike}.
\endproof

When the test function is pluriharmonic, we have the following
exponential convergence.

\begin{proposition} \label{prop_ph_test_polylike}
Let $W$ be a neighbourhood of $\Kc$ and $\Fc$ a bounded family of
pluriharmonic functions on $W$. There are constants $N\geq 0$,
$c>0$ and $0<\lambda<1$ such
that if $\varphi$ is a function in $\Fc$, then
$$|\langle \Lambda^n(\varphi)-\langle\mu,\varphi\rangle|\leq
c\lambda^n \quad \mbox{on}\quad V$$
for $n\geq N$.
\end{proposition}
\proof
Observe that if $N$ is large enough, the functions $\Lambda^N(\varphi)$ are
pluriharmonic and they absolute values are bounded on
$V$ by the same constant. We can replace $\varphi$ by
$\Lambda^N(\varphi)$ in order to assume
that $W=V$, $N=0$ and that $|\varphi|$ is bounded by some constant $A$. Then,
$|\Lambda^n(\varphi)|\leq A$ for every $n$. Subtracting from $\varphi$
the constant $\langle \mu,\varphi\rangle$ allows to assume that
$\langle\mu,\varphi\rangle=0$. 

Let $\Fc_\alpha$ denote the family of pluriharmonic functions
$\varphi$ on $V$ such that $|\varphi|\leq \alpha$ and $\langle
\mu,\varphi\rangle=0$. It is enough to show that $\Lambda$ sends
$\Fc_\alpha$ into $\Fc_{\lambda\alpha}$ for some constant
$0<\lambda<1$. We can assume $\alpha=1$. 
Since $\mu$ is totally invariant, $\Lambda$ preserves the subspace
$\{\varphi,\ \langle
\mu,\varphi\rangle=0\}$. The family $\Fc_1$ is compact and does not
contain the function identically equal to 1. By maximum
principle applied to $\pm\varphi$, there is a constant $0<\lambda<1$
such that $\sup_U|\varphi|\leq\lambda$ for $\varphi$ in $\Fc_1$. 
We deduce that $\sup_V|\Lambda(\varphi)|\leq\lambda$. The result follows.
\endproof

The following result shows that the equilibrium measure $\mu$
satisfies the Oseledec's theorem hypothesis.
It can be extended to a class of orientation
preserving smooth maps on
Riemannian manifolds \cite{DinhSibony1}.

\begin{theorem} \label{th_log_jac_polylike}
Let $f:U\rightarrow V$ be a polynomial-like map as above. Let $\mu$ be
the equilibrium measure and $J$ 
the Jacobian of $f$ with respect to the standard volume form on
$\C^k$. Then 
$$\langle \mu,\log J\rangle \geq \log d_t.$$
In particular, $\mu$ has no mass on the critical set of $f$.
\end{theorem}
\proof
Let $\nu$ be the restriction of the Lebesgue measure to $U$ multiplied
by a constant so that $\|\nu\|=1$.
Define 
$$\nu_n:=d_t^{-n}(f^n)^*(\nu)\quad \mbox{and}\quad \mu_N:={1\over
  N}\sum_{n=1}^N\nu_n.$$
By Theorem \ref{th_polylike_measure}, $\mu_N$ converge to $\mu$. 
Choose a constant $M>0$ such that $J\leq M$ on $U$. 
For any constant $m>0$, define
$$g_m(x):=\min\Big(\log\frac{M}{J(x)}, m+\log M\Big)=
\min\Big(\log\frac{M}{J(x)}, m'\Big)$$
with $m':=m+\log M$. 
This is a family of continuous functions which are positive, bounded
on $U$ and which converge to  $\log M/J$  when $m$ goes to
infinity. Define 
$$s_N(x):=\frac{1}{N}\sum_{q=0}^{N-1} g_m(f^q(x)).$$ 
Using the definition of $f^*$ on measures, we obtain
\begin{eqnarray*}
\langle \nu_N, s_N \rangle & = & \frac{1}{N}\sum_{q=0}^{N-1} d_t^{-N}
\big\langle (f^N)^*(\nu),g_m\circ f^q \big\rangle \\
& = & \frac{1}{N}\sum_{q=0}^{N-1} d_t^{-N+q} \big \langle (f^{N-q})^*(\nu), 
g_m \big\rangle\\
& = & \frac{1}{N}\sum_{q=0}^{N-1} \langle \nu_{N-q},g_m\rangle
=\langle \mu_N, g_m\rangle.
\end{eqnarray*}
In order to bound $\langle \mu, \log J \rangle$ from below,
we will bound $\langle \mu_N, g_m\rangle$ from above.

For $\alpha>0$, let $U^\alpha_N$ denote the set of points $x\in U$ such that 
$s_N(x)>\alpha$. Since $s_N(x)\leq m'$, we have
\begin{eqnarray*}
\langle \mu_N, g_m \rangle =\langle \nu_N, s_N \rangle & \leq &
m'\nu_N(U^\alpha_N)+\alpha(1-\nu_N(U^\alpha_N))\\
& = & \alpha+(m'-\alpha)\nu_N(U^\alpha_N).
\end{eqnarray*}
If $\nu_N(U^\alpha_N)$ converge to 0 when $N\rightarrow\infty$, then
$$\langle\mu,g_m\rangle=\lim_{N\rightarrow\infty}\langle \mu_N, g_m\rangle
\leq \alpha \quad \mbox{ and hence }\quad \langle \mu, \log M/J \rangle\leq \alpha.$$
We determine a value of $\alpha$ such that
$\nu_N(U^\alpha_N)$ tend to 0.

By definition
of  $\nu_N$, we have
$$\nu_N(U^\alpha_N) = \int_{U^\alpha_N} d_t^{-N} (f^N)^*(\nu)  =  \int_{U^\alpha_N} d_t^{-N}
\Big(\prod_{q=0}^{N-1}J\circ f^q\Big) d\nu.$$
Define for a given $\delta>0$ and for any integer $j$,
$$W_j:=\big\{\exp(-j\delta)<J\leq
\exp(-(j-1)\delta)\big\}$$ and
$$\tau_j(x):=\frac{1}{N} \# \big\{q, \quad 
f^q(x)\in W_j \mbox{ and } 0\leq q\leq N-1\big\}.$$
We have $\sum\tau_j=1$ and
\begin{eqnarray*}
\nu_N(U^\alpha_N) & \leq & 
\int_{U^\alpha_N}\left[\frac{1}{d_t}\exp\left(\sum -(j-1)\delta
\tau_j\right)\right]^N d\nu.
\end{eqnarray*}
Using the inequality $g_m\leq\log M/J$, we have on $U^\alpha_N$ 
$$\alpha < s_N< \sum \tau_j (\log M + j\delta) =\sum
j\delta\tau_j +\log M.$$
Therefore,
$$ -\sum(j-1)\delta\tau_j <-\alpha +(\log M +\delta).$$
We deduce from the above estimate on $\nu_N(U^\alpha_N)$ that
$$\nu_N(U^\alpha_N)\leq \int_{U^\alpha_N}\left[\frac{\exp(-\alpha) M
\exp(\delta)}{d_t}\right]^N d\nu.$$
So, for every $\alpha>\log(M/d_t)+\delta$, we have
$\nu_N(U^\alpha_N)\rightarrow 0$. 

Choosing $\delta$ 
arbitrarily small, we deduce from the above discussion that
$$\lim_{N\rightarrow\infty}
\langle \mu_N,g_m\rangle \leq \log (M/d_t).$$ 
Since $g_m$ is continuous and $\mu_N$ converge to $\mu$, we have
$\langle\mu, g_m \rangle \leq \log (M/d_t)$. Letting $m$ go to infinity
gives 
$\langle \mu, \log J\rangle\geq \log d_t$.
\endproof

\bigskip\bigskip

\begin{exercise} Let $\varphi$ be a strictly p.s.h. function on a
  neighbourhood of $\Kc$, i.e. a p.s.h. function satisfying
  $\ddc\varphi\geq c\ddc\|z\|^2$, with $c>0$, in a neighbourhood of
  $\Kc$. Let $\nu$ be a probability measure such that
  $\langle d_t^{-n}(f^n)^*(\nu),\varphi\rangle$ converge to $\langle
  \mu,\varphi\rangle$ and that $\langle \mu,\varphi\rangle$ is finite.
Show that $d_t^{-n}(f^n)^*(\nu)$ converge to $\mu$. 
\end{exercise}

\begin{exercise} \label{exo_d_k-1}
Using the test function $\varphi=\|z\|^2$, show that
 $$\int_{f^{-n}(U)} (f^n)^*(\omega^{k-1})\wedge\omega = o(d_t^n),$$
when $n$ goes to infinity.
\end{exercise}

\begin{exercise}
Let $Y$ denote the set of critical values of $f$. Show that the volume
of $f^n(Y)$ in $U$ satisfies $\vol(f^n(Y)\cap
U)=o(d_t^n)$ when $n$ goes to infinity.
\end{exercise}

\begin{exercise}
Let $f$ be a polynomial endomorphism of $\C^2$ of algebraic degree
$d\geq 2$. Assume that $f$ extends at infinity to an endomorphism of
$\P^2$. Show that $f$ admits at most three totally invariant 
points\footnote{This result was proved in \cite{DinhSibony1}. Amerik and Campana proved in
\cite{AmerikCampana} for a general endomorphism of $\P^2$ that the
number of totally invariant points is at most equal to 9. The sharp
bound (probably 3) is unknown.}. 
\end{exercise}


\section{Equidistribution problems}

In this paragraph, we consider polynomial-like maps $f$ satisfying the
hypothesis that the dynamical degree $d_{k-1}^*$ is strictly smaller than
$d_t$. We say that $f$ has a {\it large topological
  degree}. We have seen that this property is stable under small pertubations
of $f$. 
Let $Y$ denote the hypersurface of critical values of $f$. As in the
case of endomorphisms of $\P^k$, define the
ramification current $R$ by
$$R:=\sum_{n\geq 0} d_t^{-n} (f^n)_*[Y].$$
The following result is a version of Proposition \ref{prop_inverse_ball}.

\begin{proposition} \label{prop_inverse_ball_polylike}
Let $f:U\rightarrow V$ be a polynomial-like map as
  above with large topological degree. 
Let $\nu$ be a strictly positive constant and let $a$ be a point in $V$ such that
the Lelong number $\nu(R,a)$ is strictly smaller than $\nu$. Let
$\delta$ be a constant such that $d_{k-1}<\delta<d_t$. 
Then, there is a ball $B$ centered at $a$
such that $f^n$ admits at least $(1-\sqrt{\nu})d_t^n$ inverse
branches $g_i:B\rightarrow W_i$ where $W_i$ are open sets in $V$
of diameter $\leq \delta^{n/2}d_t^{-n/2}$. 
In particular, if $\mu'$ is a limit value
of the measures $d_t^{-n}(f^n)^*(\delta_a)$ then $\|\mu'-\mu\|\leq 2\sqrt{\nu(R,a)}$. 
\end{proposition}
\proof
Since $d^*_{k-1}<d_t$, the current $R$ is well-defined and has locally
finite mass. If $\omega$ is the standard K{\"a}hler form on $\C^k$, we
also have $\|(f^n)_*(\omega)\|_{V'}\lesssim \delta^n$ for every open set
$V'\Subset V$. So, for the first part of the proposition, it is enough
to follow the arguments in Proposition \ref{prop_inverse_ball}. The
proof there is written in such way that the
estimates are local and can be extended without difficulty to the
present situation. In particular, we did not use B\'ezout's theorem.

For the second assertion, we do not have yet the analogue of Proposition
\ref{prop_except_FS}, but it is enough to compare $d_t^{-n}(f^n)^*(\delta_a)$ with
the pull-backs of a smooth measure supported on $B$ and to apply
Theorem \ref{th_polylike_measure}. 
\endproof

We deduce the following result as in the case of endomorphisms of
$\P^k$.

\begin{theorem} \label{th_periodic_polylike}
Let $f:U\rightarrow V$ be a polynomial-like map as
  above with large topological degree.
Let $P_n$ denote the set of repelling periodic points of period
$n$ on the support of $\mu$. Then the sequence of measures
$$\mu_n:=d_t^{-n}\sum_{a\in P_n}\delta_a$$
converges to $\mu$.
\end{theorem}
\proof
It is enough to repeat the proof of Theorem \ref{th_periodic_pk} and to
show that $f^n$ admits exactly $d_t^n$ fixed points counted
with multiplicity. We can assume $n=1$. Let $\Gamma$ denote the graph
of $f$ in $U\times V\subset V\times V$. The number of fixed points of $f$ is the number
of points in the intersection of $\Gamma$ with the diagonal $\Delta$
of $\C^k\times \C^k$. Observe that this intersection is contained in the
compact set $\Kc\times\Kc$.  
For simplicity, assume that $V$ contains the point
0 in $\C^k$. Let $(z,w)$ denote the standard coordinates on
$\C^k\times\C^k$ where the diagonal is given by the equation
$z=w$. Consider the deformations $\Delta_t:=\{w=tz\}$ with $0\leq t\leq 1$ of
$\Delta$. Since $V$ is convex, it is not difficult to see that the intersection of $\Gamma$
with this family stays in a compact subset of $V\times V$. Therefore,
the number of points in $\Delta_t\cap\Gamma$, counted with multiplicity, does not depend on
$t$. For $t=0$, this is just the number of points in the fiber
$f^{-1}(0)$. The result follows.
\endproof

The equidistribution of negative orbits of points is more delicate 
than in the case of endomorphisms of $\P^k$. It turns out that the
exceptional set $\Ec$ does not satisfy in general $f^{-1}(\Ec)
=\Ec\cap U$.
We have the following
result.

\begin{theorem} \label{th_equi_preimage_poly}
Let $f:U\rightarrow V$ be a polynomial-like map as
  above with large topological degree. Then there is a proper analytic
  subset $\Ec$ of $V$, possibly empty, such that $d_t^{-n}(f^n)^*(\delta_a)$ converge
  to the equilibrium measure $\mu$ if and only if $a$ does not belong
  to the orbit of $\Ec$. Moreover, $\Ec$ satisfies
  $f^{-1}(\Ec)\subset\Ec\subset f(\Ec)$ and is maximal in the sense
  that if $E$ is a proper analytic subset of $V$ contained in the
  orbit of critical values such that
  $f^{-n}(E)\subset E$ for some $n\geq 1$ then $E\subset \Ec$.
\end{theorem}

The proof follows the main lines of the case of endomorphisms of
$\P^k$ using the following proposition applied to $Z$ the set of
critical values of $f$. The set $\Ec$ will be defined as $\Ec:=\Ec_Z$. 
Observe that unlike in the case of endomorphisms of $\P^k$, we need to
assume that $E$ is in the orbit of the critical values.

Let $Z$ be an arbitrary
analytic subset of $V$ not necessarily of pure dimension. 
Let $N_n(a)$ denote the number of orbits 
$$a_{-n},\ldots,a_{-1},a_0$$
with $f(a_{-i-1})=a_{-i}$ and $a_0=a$ such that $a_{-i}\in Z$
for every $i$. Here, the orbits are counted with multiplicity,
i.e. we count the multiplicity of $f^n$ at $a_{-n}$. So,
$N_n(a)$ is the number of
negative orbits of order $n$ of $a$ which stay in $Z$. Observe
that the sequence of functions $\tau_n:=d_t^{-n} N_n$ decreases to some
function $\tau$. Since $\tau_n$ are upper semi-continuous with respect
to the Zariski topology and $0\leq \tau_n\leq 1$, the function $\tau$
satisfies the same properties. 
Note that $\tau(a)$ is the proportion of infinite negative
orbits of $a$ staying in $Z$. Define $\Ec_Z:=\{\tau=1\}$. 
The Zariski upper semi-continuity of $\tau$ implies that $\Ec_Z$ is
analytic. 
It is clear
that $f^{-1}(\Ec_Z)\subset \Ec_Z$ which implies that $\Ec_Z\subset
f(\Ec_Z)$. 

\begin{proposition}
If a point $a\in V$ does not belong to the orbit $\cup_{n\geq 0} f^n(\Ec_Z)$ of $\Ec_Z$, then $\tau(a)=0$.
\end{proposition}
\proof
Assume there is $\theta_0>0$ such that $\{\tau\geq \theta_0\}$ is not
contained in the orbit of $\Ec_Z$.  We claim that there is a maximal value $\theta_0$ satisfying
the above property. Indeed, by definition, $\tau(a)$ is smaller than or
equal to the average of $\tau$ on the fiber of $a$. So, we only have
to consider the components of $\{\tau\geq \theta_0\}$ which intersect $U$
and there are only finitely many of such components, hence the maximal
value exists.

Let $E$ be the union of irreducible components
of $\{\tau\geq \theta_0\}$ which are not contained in the orbit of
$\Ec_Z$. Since $\theta_0>0$, we know that $E\subset Z$. We want to prove that $E$ is empty.
If $a$ is a generic point in $E$, it does not belong to the
orbit of $\Ec_Z$ and we have $\tau(a)=\theta_0$. If $b$ is a point in
$f^{-1}(a)$, then $b$ is not in the orbit of $\Ec_Z$. Therefore,
$\tau(b)\leq\theta_0$. Since $\tau(a)$ is smaller than or equal to the
average of $\tau$ on $f^{-1}(a)$, we deduce that $\tau(b)=\theta_0$,
and hence $f^{-1}(a)\subset
E$. By induction, we obtain that $f^{-n}(a)\subset E\subset Z$ for every $n\geq
1$. Hence, $a\in \Ec_Z$. This is a contradiction.
\endproof

The following example shows that in general the orbit of $\Ec$ is
not an analytic set. We deduce that in general polynomial-like maps are not
homeomorphically conjugated to restrictions on open sets of
endomorphisms of $\P^k$ (or polynomial maps of $\C^k$ such that the
infinity is attractive) with the same
topological degree.

\begin{example} \rm
Denote by ${\rm D}(a,R)$ the disc of radius $R$ and of center $a$ in $\C$.
Observe that the polynomial $P(z):=6z^2+1$ defines a ramified covering
of degree 2 from $D:=P^{-1}(D(0,4))$ to $D(0,4)$.
The domain $D$ is simply connected and is contained in 
$D(0,1)$. Let
$\psi$ be a bi-holomorphic map between $D(1,2)$ and  $D(0,1)$ such
that $\psi(0)=0$. Define 
$h(z,w):=(P(z),4\psi^m(w))$ with $m$ large enough. This application is
holomorphic and proper from
$W:=D\times D(1,2)$ to 
$V:=D(0,4)\times D(0,4)$. 
Its critical set is given by $zw=0$. 

Define also
$$g(z,w):=10^{-2}\big(\exp(z)\cos(\pi w/2), \exp(z)\sin(\pi w/2)\big).$$ 
One easily check  that $g$ defines a bi-holomorphic map between $W$ and
$U:=g(W)$. Consider now the polynomial-like map
$f:U\rightarrow V$ defined by $f=h\circ g^{-1}$. Its topological
degree is equal to 
$2m$; its critical set $C$ is equal to 
$g\{zw=0\}$. The image of 
$C_1:=g\{z=0\}$ by $f$ is equal to $\{z=1\}$ which is outside $U$. 
The image of 
$g\{w=0\}$ by $f$ is $\{w=0\}\cap V$.
 
The intersection $\{w=0\}\cap U$ contains two components 
$C_2:=g\{w=0\}$ and $C_2':=g\{w=2\}$. They are disjoint because $g$ is
bi-holomorphic. We also have
$f(C_2)= \{w=0\}\cap V$ and $f^{-1}\{w=0\}=C_2$.
Therefore, $\Ec=\{w=0\}\cap V$,
$f^{-1}(\Ec)\subset \Ec$ and $f^{-1}(\Ec)\not =
\Ec\cap U$ since  $f^{-1}(\Ec)$ does not contain $C_2'$. 
The orbit of $\Ec$ is the union of $C_2$ and of the orbit of $C'_2$.
Since $m$ is large, the image of $C_2'$ by 
$f$ is a horizontal curve very close to
$\{w=0\}$. It follows that the orbit of $C_2'$ is a countable union of
horizontal curves close to $\{w=0\}$ and it is not analytic.

It follows that $f$ is not holomorphically conjugate to an
endomorphism (or a polynomial map such that the hyperplane at infinity
is attractive) with the same topological degree. If it were, the
exceptional set would not have infinitely many components in a
neighbourhood of $w=0$.  
\end{example}

\begin{remark}\rm
Assume that $f$ is not with large topological degree but that the series
which defines the ramification current $R$ converges. We can then
construct inverse branches as in Proposition \ref{prop_inverse_ball_polylike}. 
To obtain the same exponential estimates on the diameter of $W_i$, it
is enough to assume that $d_{k-1}<d_t$. In general, we only have that
these diameters tend uniformly to 0
when $n$ goes to infinity. Indeed, we can use the estimate in 
Exercise \ref{exo_d_k-1}. The equidistribution of periodic points and of
negative orbits still holds in this case.
\end{remark}

\bigskip\bigskip

\begin{exercise}
Let $f$ be a polynomial-like map with large topological degree. Show
that there is a small perturbation of $f$, arbitrarily close to $f$, whose
 exceptional set is empty.
\end{exercise}


\section{Properties of the Green measure} \label{section_poly-like-entropy}

Several properties of the equilibrium measure of polynomial-like maps can be proved
using the arguments that we introduced in the case of endomorphisms of
$\P^k$. We have the following result for general polynomial-like maps.

\begin{theorem}
Let $f:U\rightarrow V$ be a polynomial-like map of topological degree
$d_t>1$. Then its equilibrium measure $\mu$ is an invariant measure
of maximal entropy $\log d_t$. Moreover, $\mu$ is K-mixing. 
\end{theorem}
\proof
By Theorem \ref{th_log_jac_polylike}, $\mu$ has no mass on the critical set of
$f$. Therefore, it is an invariant measure of constant Jacobian
$d_t$ in the sense that $\mu(f(A))=d_t\mu(A)$ 
when $f$ is injective on a Borel set $A$. We deduce from Parry's theorem \ref{th_parry} that
$h_\mu(f)\geq \log d_t$. 
The variational principle and Theorem \ref{th_entropy_top_pol_like}
imply that
$h_\mu(f)=\log d_t$. 

We prove the K-mixing property. As in the case of endomorphisms of
$\P^k$, it is enough to show  for $\varphi$ in $L^2(\mu)$ 
that $\Lambda^n(\varphi)\rightarrow \langle\mu,\varphi\rangle$ in $L^2(\mu)$. Since
$\Lambda:L^2(\mu)\rightarrow L^2(\mu)$ is of norm 1, it is enough to
check the convergence for a dense family of $\varphi$. So, we only
have to consider
$\varphi$ smooth. We can also assume that $\varphi$ is
p.s.h. because smooth functions can be written as a difference of p.s.h.
functions. Assume also for simplicity that
$\langle\mu,\varphi\rangle=0$. 

So, the p.s.h. functions $\Lambda^n(\varphi)$ converge to 0
in $L^p_\loc(V)$. By Hartogs' lemma \ref{prop_hartogs_else}, $\sup_U\Lambda^n(\varphi)$ converge to 0.
This and the identity $\langle\mu,\Lambda^n(\varphi)\rangle =\langle \mu,\varphi\rangle=0$ imply that 
$\mu\{\Lambda^n(\varphi)<-\delta\}\rightarrow 0$ for every fixed
$\delta>0$. On the other hand,
by definition of $\Lambda$, $|\Lambda^n(\varphi)|$ is bounded
by $\|\varphi\|_\infty$ which is a constant independent
of $n$. Therefore, $\Lambda^n(\varphi)\rightarrow 0$ in $L^2(\mu)$ and
K-mixing follows.
\endproof

\begin{theorem} 
Let $f$ and $\mu$ be as above. Then the sum of the Lyapounov
  exponents of $\mu$ is at least equal to ${1\over 2}\log d_t$. In particular, $f$
  admits a strictly positive Lyapounov exponent. If $f$ is with large
  topological degree, then $\mu$ is hyperbolic and its Lyapounov
  exponents are at least equal to ${1\over 2}\log(d_t/d_{k-1})$.
\end{theorem}
\proof
By Oseledec's theorem \ref{th_oseledec}, applied in the complex setting, the sum of the Lyapounov
  exponents of $\mu$ (associated to complex linear spaces) is equal to ${1\over 2}\langle\mu,\log
  J\rangle$. Theorem \ref{th_log_jac_polylike} implies that this sum is at least
  equal to ${1\over 2}\log d_t$. The second assertion is proved as in Theorem
  \ref{th_briend_duval} using Proposition \ref{prop_inverse_ball_polylike}.
\endproof

From now on, we only consider maps with large topological degree.
The following result was obtained by Dinh-Dupont in \cite{DinhDupont}.
It generalizes Theorem \ref{th_dinhdupont_pk}.

\begin{theorem} \label{th_dinhdupont_polylike}
Let $f$ be a polynomial-like map with large topological degree as above.
Let $\chi_1,\ldots,\chi_k$ denote
the Lyapounov exponents of the equilibrium measure $\mu$ ordered by
$\chi_1\geq\cdots\geq\chi_k$ and $\Sigma$ their sum. Then the Hausdorff
dimension of $\mu$ satisfies
$${\log d_t\over \chi_1}\leq \dim_H(\mu)\leq 2k-{2\Sigma-\log d_t\over\chi_1}\cdot$$
\end{theorem}

We now prove some stochastic properties of the equilibrium measure.
We first introduce some notions.
Let $V$ be an open subset of $\C^k$ and $\nu$ a probability measure
with compact support in $V$. We consider $\nu$ as a function on  the convex cone
$\PSH(V)$ of p.s.h. functions on $V$, with the $L^1_\loc$-topology.
We say that $\nu$ is {\it PB} if this function is finite,
i.e. p.s.h. functions on $V$ are $\nu$-integrable.
We say that $\nu$ is {\it PC} if it is PB and defines a
continuous functional on $\PSH(V)$. Recall that the weak topology on
$\PSH(V)$ coincides with the $L^p_\loc$ topology for $1\leq p<+\infty$.
In dimension 1, a measure is PB if it has 
locally bounded potentials, a measure is PC if it has locally continuous
potentials. A measure $\nu$ is
 {\it moderate} if for any bounded subset $\Pc$ of
$\PSH(V)$, there are constants  $\alpha>0$ and $A>0$  such that  
$$\langle \nu, e^{\alpha|\varphi|}\rangle\leq A \quad \mbox{for}\quad
\varphi\in \Pc.$$

Let $K$ be a compact subset of $V$.
Define a {\it pseudo-distance} $\dist_{L^1(K)}$ between
$\varphi,\psi$ in $\PSH(V)$
by
$$\dist_{L^1(K)}(\varphi,\psi):=\|\varphi-\psi\|_{L^1(K)}.$$
Observe that if $\nu$ is continuous with respect to $\dist_{L^1(K)}$
then $\nu$ is PC. The following proposition gives a criterion for a
measure to be 
moderate.

\begin{proposition}
If  $\nu$ is H\"older continuous with respect to $\dist_{L^1(K)}$ for
some compact subset $K$ of $V$, then
 $\nu$ is moderate. If $\nu$ is moderate, then 
p.s.h. functions on $V$ are
in $L^p(\nu)$ for every $1\leq p<+\infty$ and $\nu$ has positive
Hausdorff dimension. 
\end{proposition}
\proof
We prove the first assertion. Assume that $\nu$ is H\"older
continuous with respect to $\dist_{L^1(K)}$ for some compact subset $K$ of $V$.
Consider a bounded subset $\Pc$ of $\PSH(V)$. 
The functions in $\Pc$ are uniformly
bounded above on $K$. Therefore, subtracting from these
functions a constant allows to assume that they are negative on
$K$. We want to prove the estimate
$\langle \nu, e^{-\alpha \varphi}\rangle\leq A$ for 
$\varphi\in \Pc$ and for some constants $\alpha,A$.
It is enough to show that $\nu\{\varphi<-M\}\lesssim e^{-\alpha M}$
for some (other) constant $\alpha>0$ and for $M\geq 1$.  
For $M\geq 0$ and $\varphi\in\Pc$, define $\varphi_M:=\max(\varphi,-M)$. We replace 
$\Pc$ by the family of functions $\varphi_M$. This allows to
assume that the family is stable under the operation
$\max(\cdot,-M)$. 
Observe that $\varphi_{M-1}-\varphi_M$ is positive, supported on
$\{\varphi<-M+1\}$, smaller or equal
to 1,  and equal to 1 on
$\{\varphi<-M\}$.  
In order to obtain the above estimate, we only have
to show that 
$\langle\nu,\varphi_{M-1}-\varphi_M\rangle \lesssim e^{-\alpha M}$ for
some $\alpha>0$.

Fix a constant $\lambda>0$ small enough and a constant $A>0$ large enough. 
Since $\nu$ is H\"older continuous and $\varphi_{M-1}-\varphi_M$
vanishes on $\{\varphi>-M+1\}$, we have 
\begin{eqnarray*}
\nu\{\varphi<-M\} & \leq & \langle \nu,\varphi_{M-1}\rangle-\langle\nu,\varphi_M\rangle\leq 
A\ \|\varphi_{M-1}-\varphi_M\|_{L^1(K)}^\lambda \\
& \leq & A\vol\{\varphi\leq -M+1\}^\lambda.
\end{eqnarray*}
On the other hand, since $\Pc$ is a bounded family in $\PSH(V)$, by
Theorem \ref{th_hormander}, we have $\|e^{-\lambda\varphi}\|_{L^1(K)}\leq A$ for
$\varphi\in\Pc$. Hence,
$$\vol\{\varphi\leq -M+1\}\leq A e^{-\lambda (M-1)}.$$
This implies the desired estimate for $\alpha=\lambda^2$ and
completes the proof of the first assertion.

Assume now that $\nu$ is moderate. 
Let $\varphi$ be a p.s.h. function on $V$. Then $e^{\alpha|\varphi|}$
is in $L^1(\nu)$ for some constant $\alpha>0$.
Since $e^{\alpha x}\gtrsim x^p$ for $1\leq
p<+\infty$, we deduce that $\varphi$ is in $L^p(\nu)$. 
For the last assertion in the proposition, it is enough to show that
$\nu(B_r)\leq A r^\alpha$ for any ball $B_r$ of radius $r>0$ where 
$A,\alpha$ are some positive constants. We can assume that $B_r$ is a small
ball centered at a point $a\in K$. Define
$\varphi(z):=\log\|z-a\|$. This function 
belongs to a compact family of p.s.h. functions. Therefore,
$\|e^{-\alpha\varphi}\|_{L^1(\nu)}\leq A$ for some positive constants
$A, \alpha$ independent of $B_r$. Since $e^{-\alpha\varphi}\geq
r^{-\alpha}$ on $B_r$, we deduce that $\nu(B_r)\leq Ar^\alpha$. 
It is well-known that in order to compute the Hausdorff dimension of a set, it
is enough to use only coverings by balls.    
It follows easily that if a Borel set has positive measure, then its
Hausdorff dimension is at least equal to $\alpha$. This completes the
proof. 
\endproof

The following results show that the equilibrium measure of a
polynomial-like map with large topological degree satisfies the
above regularity properties.

\begin{theorem} \label{th_top_PC}
Let $f:U\rightarrow V$ be a polynomial-like map. Then
  the following properties are equivalent:
\begin{enumerate}
\item The map $f$ has large topological degree, i.e. $d_t>d_{k-1}^*$;
\item The measure $\mu$ is PB, i.e. p.s.h. functions on $V$ are integrable
      with respect to $\mu$;
\item The measure $\mu$ is PC, i.e. $\mu$ can be extended to a linear continuous form
      on the cone of p.s.h. functions on $V$;
\end{enumerate}
Moreover, if $f$ is such a map, then there is a constant $0<\lambda<1$
such that
$$\sup_V \Lambda(\varphi)-\langle\mu,\varphi\rangle \leq \lambda
\big[\sup_V \varphi-\langle\mu,\varphi\rangle\big],$$
for $\varphi$  p.s.h. on $V$.
\end{theorem}
\proof
It is clear that $3)\Rightarrow 2)$. 
We show that $1)\Rightarrow 3)$ and $2)\Rightarrow 1)$. 

\medskip
\noindent
$1)\Rightarrow 3)$. Let $\varphi$ be a p.s.h. function on $V$. Let $W$
be a convex open set such that $U\Subset W\Subset V$. For simplicity, assume
that $\|\varphi\|_{L^1(W)}\leq 1$. So, $\varphi$ belongs to a compact
subset of $\PSH(W)$. Therefore, $S:=\ddc \varphi$ has
locally bounded mass in $W$. Define $S_n:=(f^n)_*(S)$. Fix a constant
$\delta>1$ such that $d_{k-1}^*<\delta<d_t$. Condition 1) implies that
$\|(f^n)_*(S)\|_W\lesssim \delta^n$.  
By Proposition
\ref{prop_ddbar_convex_dom}, there are  p.s.h. functions $\varphi_n$ on
$U$ such that $\ddc \varphi_n= S_n$ and $\|\varphi_n\|_U\lesssim
\delta^n$ on $U$. 

Define $\psi_0:=\varphi-\varphi_0$ and
$\psi_n:=f_*(\varphi_{n-1})-\varphi_n$. 
Observe that these functions are pluriharmonic on $U$ and depend
continuously on $\varphi$. Moreover, 
$f_*$ sends continuously p.s.h. functions on $U$ to p.s.h. functions
on $V$. Hence,
$$\|\psi_n\|_{L^1(U)}\leq \|f_*(\varphi_{n-1})\|_{L^1(U)}+
\|\varphi_n\|_{L^1(U)} \lesssim \delta^n.$$
We have
\begin{eqnarray*}
\Lambda^n(\varphi) & = &
\Lambda^n(\psi_0+\varphi_0)=\Lambda^n(\psi_0)+d_t^{-1}\Lambda^{n-1}(f_*(\varphi_0))
\\
& = & \Lambda^n(\psi_0)+ d_t^{-1}\Lambda^{n-1}(\psi_1+\varphi_1)=\cdots\\
& = & \Lambda^n(\psi_0)+d_t^{-1}\Lambda^{n-1}(\psi_1)+\cdots+d_t^{-n+1}\Lambda(\psi_{n-1})+d_t^{-n}
\psi_n+d_t^{-n}\varphi_n.
\end{eqnarray*}
The last term in the above sum converges to 0. 
The above estimate on $\psi_n$ and their pluriharmonicity 
imply, by Proposition \ref{prop_ph_test_polylike}, that the sum 
$$\Lambda^n(\psi_0)+d_t^{-1}\Lambda^{n-1}(\psi_1)+\cdots+d_t^{-n+1}\Lambda(\psi_{n-1})+d_t^{-n}
\psi_n$$
converges uniformly to the finite constant 
$$\langle \mu,\psi_0\rangle + d_t^{-1}\langle \mu,\psi_1\rangle +\cdots+ d_t^{-n}\langle \mu,\psi_n\rangle+\cdots$$
which depends continuously on
$\varphi$. 
We used here the fact that when $\psi$ is pluriharmonic, $\langle
\mu,\psi\rangle$ depends continuously on $\psi$. 
By Theorem \ref{th_polylike_measure}, the above constant is
equal to $\langle\mu,\varphi\rangle$. Consequently, $\mu$ is PC.

\medskip
\noindent
$2)\Rightarrow 1)$. Let $\Fc$ be an $L^1$ bounded family of p.s.h. functions
on a neighbourhood of $\Kc$. We first show that $\langle \mu,\varphi\rangle$ is
uniformly bounded on $\Fc$. Since
$\langle\mu,\Lambda^N(\varphi)\rangle=\langle\mu,\varphi\rangle$, we
can replace $\varphi$ by $\Lambda^N(\varphi)$, with $N$ large enough, in order to assume that
$\Fc$ is a bounded family of p.s.h. functions $\varphi$ on $V$ which are
uniformly bounded above. Subtracting from $\varphi$ a fixed constant
allows to assume that these functions are negative. If  $\langle
\mu,\varphi\rangle$ is not uniformly bounded
on $\Fc$, there are $\varphi_n$ such that $\langle
\mu,\varphi_n\rangle \leq -n^2$. It follows that the series $\sum
n^{-2}\varphi_n$ decreases to a p.s.h. function which is not
integrable with respect to $\mu$. This contradicts that $\mu$ is
PB. We deduce that for any neighbourhood $W$ of $\Kc$, there is a constant $c>0$ such that 
$|\langle \mu,\varphi\rangle|\leq c\|\varphi\|_{L^1(W)}$ for $\varphi$
p.s.h. on $W$.

We now show that there is a constant $0<\lambda<1$ such that
$\sup_V\Lambda(\varphi)\leq \lambda$ if $\varphi$ is a p.s.h. function
on $V$, bounded from above by 1, such that
$\langle\mu,\varphi\rangle=0$. 
This property implies the last assertion in the proposition.
Assume that the property is not satisfied.
Then, there are functions $\varphi_n$ such that
$\sup_V\varphi_n=1$, $\langle\mu,\varphi_n\rangle=0$ and
$\sup_V\Lambda(\varphi_n)\geq 1-1/n^2$.
By definition of $\Lambda$, we have 
$$\sup_U\varphi_n\geq \sup_V\Lambda(\varphi_n)\geq 1-1/n^2.$$ 
The submean value inequality for p.s.h. functions implies that
$\varphi_n$ converge to 1 in $L_\loc^1(V)$. 
On the other hand, we have
$$1=|\langle\mu,\varphi_n-1\rangle|\leq
c\|\varphi_n-1\|_{L^1(W)}.$$
This is a contradiction.

Finally, consider a positive closed $(1,1)$-current $S$ of mass 1 on
$W$. By Proposition \ref{prop_ddbar_convex_dom}, there is a p.s.h. function $\varphi$ on
a neighbourhood of $\overline U$ with bounded $L^1$ norm 
such that $\ddc \varphi=S$. 
The submean inequality for p.s.h functions implies that $\varphi$ is
bounded from above by a constant independent of $S$. 
We can after subtracting from $\varphi$ a constant, assume 
that $\langle \mu,\varphi\rangle=0$. 
The p.s.h. functions
$\lambda^{-n}\Lambda^n(\varphi)$ are bounded above and satisfy $\langle
\mu,\lambda^{-n}\Lambda^n(\varphi)\rangle=\langle\mu,\varphi\rangle=0$.
Hence, they belong to a compact subset of $\PSH(U)$ which is
independent of $S$. If $W'$ is a neighbourhood of $\Kc$ such that
$W'\Subset U$, the mass of 
$\ddc \big[\lambda^{-n}\Lambda^n(\varphi)\big]$ on $W'$ is bounded
uniformly on $n$ and on $S$. Therefore, 
$$\|\ddc (f^n)_*(S)\|_{W'}\leq c\lambda^n d_t^n$$
for some constant $c>0$ independent of $n$ and of $S$. It follows that
$d_{k-1}^*\leq \lambda d_t$. This implies property 1). 
\endproof

\begin{theorem} \label{th_pol_like_holder}
 Let $f:U\rightarrow V$ be a polynomial-like map with large
 topological degree. Let $\Pc$ be a bounded family of p.s.h. functions
 on $V$. Let $K$ be a compact subset of $V$ such that $f^{-1}(K)$ is
 contained in the interior of $K$. Then, the equilibrium measure $\mu$ of $f$ is  H\"older
 continuous on $\Pc$ with respect to $\dist_{L^1(K)}$.
 In particular, this measure is moderate. 
\end{theorem}

Let $\DSH(V)$ denote the space of d.s.h. functions on
$V$, i.e. functions which are differences of p.s.h. functions. They
are in particular in $L^p_\loc(V)$ for every $1\leq p<+\infty$.  
Consider on $\DSH(V)$ the following topology: a sequence 
$(\varphi_n)$ converges to $\varphi$ in $\DSH(V)$ if
$\varphi_n$ converge weakly to $\varphi$ and if we can write
$\varphi_n=\varphi_n^+-\varphi_n^-$ with $\varphi_n^\pm$ in a compact
subset of $\PSH(V)$, independent of $n$.
We deduce from the compactness of bounded sets of p.s.h. functions
that $\varphi_n\rightarrow\varphi$ in all $L^p_\loc(V)$ with $1\leq p<+\infty$. 
Since $\mu$ is PC, it extends
by linearity to a continuous functional on $\DSH(V)$. 

\medskip

\noindent
{\bf Proof of Theorem \ref{th_pol_like_holder}.} 
 Let $\Pc$ be a compact family of
p.s.h. functions on $V$.
We show that $\mu$ is H\"older
continuous on $\Pc$ with respect to $\dist_{L^1(K)}$. 
We claim that
$\Lambda$ is Lipschitz with respect to $\dist_{L^1(K)}$. Indeed, if
$\varphi,\psi$ are in $L^1(K)$, we have for 
the standard volume form $\Omega$  on $\C^k$
$$\|\Lambda(\varphi)-\Lambda(\psi)\|_{L^1(K)}=\int_K
|\Lambda(\varphi-\psi)|\Omega \leq d_t^{-1}\int_{f^{-1}(K)} |\varphi-\psi|
f^*(\Omega).$$
since $f^{-1}(K)\subset K$ and $f^*(\Omega)$ is bounded on
$f^{-1}(K)$, this   implies that
$$\|\Lambda(\varphi)-\Lambda(\psi)\|_{L^1(K)}\leq \const
\|\varphi-\psi\|_{L^1(K)}.$$

Since $\Pc$ is compact, the functions in $\Pc$ are uniformly bounded 
above on $U$. Therefore, replacing
$\Pc$ by the family of $\Lambda(\varphi)$ with
$\varphi\in\Pc$ allows to assume that functions in $\Pc$ are uniformly bounded
above on $V$. On the other hand, since $\mu$ is PC, $\mu$ is bounded
on $\Pc$. Without loss of generality, we can assume that $\Pc$ is the
set of functions $\varphi$ such that $\langle\mu,\varphi\rangle\geq 0$ and
$\varphi\leq 1$. In particular, $\Pc$ is invariant under $\Lambda$.
Let
$\Dc$ be the family of d.s.h. functions 
$\varphi-\Lambda(\varphi)$ with $\varphi\in \Pc$. This is a compact subset of
$\DSH(V)$ which is invariant under $\Lambda$, and we have
$\langle\mu,\varphi'\rangle=0$ for $\varphi'$ in $\Dc$.

Consider a function $\varphi\in\Pc$. Observe that
$\widetilde\varphi:=\varphi-\langle\mu,\varphi\rangle$ is also in
$\Pc$.
Define $\widetilde\Lambda:=\lambda^{-1}\Lambda$ with $\lambda$ the
constant in Theorem \ref{th_top_PC}. 
We deduce from that theorem that
$\widetilde\Lambda(\widetilde\varphi)$ is in $\Pc$. 
Moreover,
$$\widetilde\Lambda\big(\varphi-\Lambda(\varphi)\big)=\widetilde\Lambda\big(
\widetilde\varphi-\Lambda(\widetilde\varphi)\big)=\widetilde\Lambda(
\widetilde\varphi)-\Lambda\big(\widetilde\Lambda(\widetilde\varphi)\big).
$$    
Therefore, $\Dc$ is invariant under
$\widetilde\Lambda$. This is the key point in the proof.
Observe that we can extend $\dist_{L^1(K)}$ to $\DSH(V)$ and 
that $\widetilde\Lambda$ is Lipschitz with respect to this pseudo-distance.

Let $\nu$ be a smooth
probability measure with support in $K$. We have seen that 
$d_t^{-n} (f^n)^*(\nu)$ converge to $\mu$. 
If $\varphi$ is d.s.h. on $V$, then
$$\langle d_t^{-n} (f^n)^*(\nu),\varphi\rangle =
\langle \nu, \Lambda^n(\varphi)\rangle.$$
Define for $\varphi$ in $\Pc$, $\varphi':=\varphi-\Lambda(\varphi)$.
We have
\begin{eqnarray*}
\langle\mu,\varphi\rangle & = & \lim_{n\rightarrow\infty}\langle
\nu,\Lambda^n(\varphi)\rangle \\
& = & 
\langle\nu,\varphi\rangle - \sum_{n\geq 0} \langle
\nu,\Lambda^n(\varphi)\rangle - \langle
\nu,\Lambda^{n+1}(\varphi)\rangle\\
& = & \langle\nu,\varphi\rangle-\sum_{n\geq 0}\lambda^n
\langle \nu, \widetilde\Lambda^n(\varphi')\rangle.
\end{eqnarray*}

Since $\nu$ is smooth with support in
$K$, it is Lipschitz with respect to $\dist_{L^1(K)}$.
We deduce from Lemma \ref{lemma_holder} which is also valid for a pseudo-distance, that
the last series defines a
H\"older continuous function
on $\Dc$. We use here the invariance of $\Dc$
under $\widetilde\Lambda$.
Finally, since the map $\varphi\mapsto\varphi'$ is Lipschitz
on $\Pc$, we conclude that $\mu$ is H\"older continuous on $\Pc$ with
respect to $\dist_{L^1(K)}$.  
\hfill $\square$

\medskip

As in the case of endomorphisms of $\P^k$, we deduce from the above
results the following fundamental estimates on the Perron-Frobenius
operator $\Lambda$.

\begin{corollary} \label{cor_exp_dsh_poly}
Let $f$ be a polynomial-like map with large topological degree as
above. Let $\mu$ be the equilibrium measure and $\Lambda$ the
Perron-Frobenius operator associated to $f$. Let $\Dc$ be a bounded
subset of d.s.h. functions on $V$. There are constants
$c>0$, $\delta>1$ and $\alpha>0$ such that if $\psi$ is in $\Dc$, then
$$\big\langle \mu, e^{\alpha \delta^n |\Lambda^n(\psi)-\langle\mu,\psi\rangle|}\big\rangle \leq c \quad
\mbox{and}\quad \|\Lambda^n(\psi)-\langle\mu,\psi\rangle\|_{L^q(\mu)}\leq cq \delta^{-n}$$
for every $n\geq 0$ and every $1\leq q <+\infty$. 
\end{corollary}

\begin{corollary} \label{cor_exp_holder_poly}
Let $f$, $\mu$, $\Lambda$ be as above. Let $0<\nu\leq 2$ be a
constant. There are constants $c>0$, $\delta>1$ and
$\alpha>0$ such that if $\psi$ is a $\nu$-H{\"o}lder continuous
function on $V$ with $\|\psi\|_{\Cc^\nu}\leq 1$, then 
$$\big\langle \mu, e^{\alpha \delta^{n\nu/2} |\Lambda^n(\psi)-\langle \mu,\psi\rangle|}\big\rangle \leq c
\quad \mbox{and}\quad 
\|\Lambda^n(\psi)-\langle\mu,\psi\rangle\|_{L^q(\mu)}\leq cq^{\nu/2} \delta^{-n\nu/2}$$
for every $n\geq 0$ and every  $1\leq q <+\infty$. Moreover, $\delta$
is independent of $\nu$.
\end{corollary}

The following results are deduced as in the case of endomorphisms of $\P^k$.

\begin{theorem} \label{th_ex_mixing_polylike}
Let $f:U\rightarrow V$ be a polynomial-like map with
  large topological degree and $\mu$ the equilibrium measure of $f$. 
Then $f$ is exponentially mixing. More
  precisely, there is a constant $0<\lambda<1$, such that if $1<p\leq
  +\infty$, we have
$$|\langle \mu,(\varphi\circ
f^n)\psi\rangle -\langle\mu,\varphi\rangle\langle\mu,\psi\rangle| \leq c_p \lambda^n
\|\varphi\|_{L^p(\mu)} \|\psi\|_{L^1(V)}$$
for $\varphi$ in $L^p(\mu)$, $\psi$ p.s.h. on $V$ and $n\geq
0$, where $c_p>0$ is a constant independent of $\varphi,\psi$. If $\nu$ is such that $0\leq \nu\leq 2$, then
there is a constant $c_{p,\nu}>0$ such that 
$$|\langle \mu,(\varphi\circ f^n)\psi\rangle
-\langle\mu,\varphi\rangle\langle\mu,\psi\rangle| 
\leq c_{p,\nu} \lambda^{n\nu/2}
\|\varphi\|_{L^p(\mu)} \|\psi\|_{\Cc^\nu}$$
for $\varphi$ in $L^p(\mu)$, $\psi$ a $\Cc^\nu$ function on
$V$ and $n\geq 0$.
\end{theorem}

The following result  gives the
exponential mixing of any order. It can be extended to 
H{\"o}lder continuous observables using the theory of interpolation 
between Banach spaces.

\begin{theorem} 
Let $f, \mu$ be as in Theorem \ref{th_ex_mixing_polylike} and $r\geq 1$
an integer. Then there are constants $c>0$ and $0<\lambda<1$ such that
$$\Big|\langle \mu, \psi_0(\psi_1\circ f^{n_1})\ldots (\psi_r\circ f^{n_r})\rangle
-\prod_{i=0}^r\langle\mu,\psi_i\rangle \Big|\leq c
\lambda^n\prod_{i=0}^r\|\psi_i\|_{L^1(V)}$$
for $0=n_0\leq n_1\leq \cdots\leq n_r$, $n:=\min_{0\leq i< r} (n_{i+1}-n_i)$ and $\psi_i$ 
p.s.h. on $V$. 
\end{theorem}

As in Chapter \ref{chapter_endo}, we deduce the following result, as a
consequence of Gordin's theorem and of the exponential decay of correlations.

\begin{theorem}
Let $f$ be a polynomial-like map with large topological degree as
above. Let $\varphi$ be a test function which is $\Cc^\nu$ with
$\nu>0$, or is d.s.h. on $V$. Then, either $\varphi$ is a coboundary
or it satisfies the central limit theorem 
with the variance $\sigma>0$ given by 
$$\sigma^2:=\langle\mu,\varphi^2\rangle +2 \sum_{n\geq 1} \langle
\mu,\varphi(\varphi\circ f^n)\rangle.$$ 
\end{theorem}

The following result is obtained as in Theorem \ref{th_deviation}, as
a consequence of the exponential estimates in Corollaries \ref{cor_exp_dsh_poly} and \ref{cor_exp_holder_poly}.

\begin{theorem} \label{th_polylike_deviation}
Let $f$ be  a polynomial-like map with large topological degree as
above.
Then, the equilibrium measure $\mu$ of $f$ satisfies the weak
large deviations theorem for bounded d.s.h. observables and for
H\"older continuous observables. More
precisely, if a function $\psi$ is bounded d.s.h. or H{\"o}lder continuous then for every
$\epsilon>0$ there is a constant $h_\epsilon>0$ such that 
\begin{equation*}
\mu\Big\{ z\in \supp(\mu):\   \big\vert\frac{1}{N}\sum_{n=0}^{N-1}\varphi\circ
f^n(z)- \langle \mu,\varphi\rangle \big\vert >\epsilon \Big\}\leq
e^{-N (\log N)^{-2}h_\epsilon}
\end{equation*}
for all $N$ large enough.
\end{theorem}

\bigskip\bigskip

\begin{exercise}
Assume that for every positive closed $(1,1)$-current $S$ on $V$ we
have $\limsup \|(f^n)_*(S)\|^{1/n}<d_t$. Show that $\mu$ is PB and
deduce that $d_{k-1}^*<d_t$. Hint: write $S=\ddc\varphi$. 
\end{exercise}

\begin{exercise}
Let $\nu$ be a positive measure with compact support in $\C$. Prove
that $\nu$ is moderate if and only if there are positive constants 
$\alpha$ and $c$ such that for every disc $D$ of radius $r$,
$\nu(D)\leq cr^\alpha$. Give an example showing that this condition is
not sufficient in $\C^2$. 
\end{exercise}


\section{Holomorphic families of maps}

In this paragraph, we consider polynomial-like maps
$f_s:U_s\rightarrow V_s$ depending
holomorphically on a parameter $s\in\Sigma$. We will show
that the Green measure $\mu_s$ of $f_s$ depends ``holomorphically'' on
$s$ and then we study the dependence of the Lyapounov exponents on the
parameters. Since the problems are local, we assume for simplicity
that $\Sigma$ is a ball in $\C^l$. Of course, we assume that $U_s$ and
$V_s$ depend continuously on $s$. Observe that if we replace $V_s$ by
a convex open set $V'_s\subset V_s$ and $U_s$ by $f_s^{-1}(V'_s)$ with
$V_s\setminus V_s'$ small enough, the map $f_s$ is still
polynomial-like. So, for simplicity, assume that $V:=V_s$ is
independent of $s$. Let $U_\Sigma:=\cup_s \{s\}\times U_s$. This is an
open set in $V_\Sigma:=\Sigma\times V$. 
Define the holomorphic map
$F:U_\Sigma \rightarrow V_\Sigma$ by $F(s,z):=(s,f_s(z))$. This map is
proper. By continuity, the topological degree $d_t$ of $f_s$ is
independent of $s$. So, the topological degree of $F$ is also
$d_t$. 
Define 
$\Kc_\Sigma:=\cap_{n\geq 0} F^{-n}(V_\Sigma)$. Then $\Kc_\Sigma$ is
closed in $\Uc_\Sigma$. If
 $\pi:\Sigma\times\C^k\rightarrow\Sigma$ is the
canonical projection, then $\pi$ is proper on $\Kc_\Sigma$ and 
$\Kc_s:=\Kc_\Sigma\cap\pi^{-1}(s)$ is the filled Julia set of $f_s$. 

It is not difficult to show that $\Kc_s$ depends upper semi-continuously on
$s$ with respect to the Hausdorff metric on compact sets of $V$. This
means that if $W_{s_0}$ is a neighbourhood of $\Kc_{s_0}$, then
$\Kc_s$ is contained in $W_{s_0}$ for $s$ closed enough to $s_0$.  
We will
see that for maps with large topological degree, $s\mapsto\mu_s$ is
continuous in a strong sense.  However, in general, the Julia set $\Jc_s$, i.e. the support of the
equilibrium measure $\mu_s$, does not depend continuously on $s$.

In our context, the goal is to construct and to study currents which 
measure the bifurcation, i.e. the discontinuity of $s\mapsto
\Jc_s$.
We have the following result due to Pham \cite{Pham}.

\begin{proposition} \label{prop_pham}
Let $(f_s)_{s\in\Sigma}$ be as above. Then, there is a positive closed current $\Rc$
of bidegree $(k,k)$, supported on $\Kc_\Sigma$ such that the slice $\langle
\Rc,\pi,s\rangle$ is equal to the equilibrium measure $\mu_s$
of $f_s$ for $s\in\Sigma$. Moreover, if $\varphi$ is a p.s.h. function
on a neighbourhood of $\Kc_\Sigma$, then the function $s\mapsto
\langle\mu_s,\varphi(s,\cdot)\rangle$ is either equal to $-\infty$ or
is p.s.h. on $\Sigma$. 
\end{proposition}
\proof
Let $\Omega$ be smooth probability measure with compact support in
$V$. Define the positive closed $(k,k)$-current $\Theta$ on
$\Sigma\times V$ by $\Theta:=\tau^*(\Omega)$ where
$\tau:\Sigma\times V\rightarrow V$ is the canonical
projection. Observe that the slice $\langle\Theta,\pi,s\rangle$
coincides with $\Omega$ on $\{s\}\times V$, since $\Omega$ is smooth. Define $\Theta_n:=d_t^{-n}
(F^n)^*(\Theta)$. 
The slice $\langle\Theta_n,\pi,s\rangle$ can be identified with
$d_t^{-n} (f_s^n)^*(\Omega)$ on $\{s\}\times V$. This is a smooth probability
measure which tends to
$\mu_s$ when $n$ goes to infinity. 

Since the problem is local for $s$, we can assume that all the forms
$\Theta_n$ are supported on $\Sigma\times K$ for some
compact subset $K$ of $V$. 
As we mentioned in Appendix \ref{section_intersection}, since these forms have slice mass
1, they belong to a compact family of currents. Therefore, we can
extract a sequence  $\Theta_{n_i}$ which
converges to  some current $\Rc$ with slice mass 1. We want to prove
that $\langle \Rc,\pi,s\rangle =\mu_s$. 

Let $\varphi$ be a smooth p.s.h. function on a neighbourhood of
$\Kc_\Sigma$. So, for $n$ large enough, $\varphi$ is defined on the
support of 
$\Theta_n$ (we reduce $\Sigma$ if necessary).  By slicing theory,
$\pi_*(\Theta_{n_i}\wedge\varphi)$ is equal to the p.s.h. function
$\psi_{n_i}(s):=\langle\Theta_{n_i},\pi,s\rangle(\varphi)$ and
$\pi_*(\Rc\wedge \varphi)$ is equal 
to the p.s.h. function $\psi(s):=\langle\Rc,\pi,s\rangle(\varphi)$ in
the sense of currents. 
By definition of $\Rc$, since $\pi_*$ is continuous on currents
supported on $\Sigma\times K$, $\psi_{n_i}$ converge to $\psi$ in $L^1_\loc(\Sigma)$. On the other hand,
$\langle\Theta_{n_i},\pi,s\rangle$ converge to $\mu_s$. So, the function
$\psi'(s):=\lim \psi_{n_i}(s)=\langle \mu_s,\varphi\rangle$ is equal to  $\psi(s)$
almost
everywhere. Since $\psi_{n_i}$ and $\psi$ are p.s.h., the Hartogs'
lemma implies that $\psi'\leq\psi$. We show the inequality $\psi'(s)\geq
\psi(s)$. 

The function $\psi$ is p.s.h., hence it is strongly upper
semi-continuous. Therefore, there is a sequence
$(s_n)$ converging to $s$ such that $\psi'(s_n)=\psi(s_n)$ and $\psi(s_n)$
converge to $\psi(s)$. Up to extracting a subsequence, we can assume that
$\mu_{s_n}$ converge to some probability measure $\mu_s'$. By
continuity, $\mu_s'$ is totally invariant under $f_s$. We deduce from
Proposition 
\ref{prop_test_ph_polylike} that $\langle \mu_s',\varphi(s,\cdot)\rangle\leq \langle
\mu_s,\varphi(s,\cdot)\rangle$. The first integral is equal to
$\psi(s)$, the second one is equal to $\psi'(s)$. Therefore,
$\psi(s)\leq \psi'(s)$. The identity $\langle \Rc,\pi,s\rangle =\mu_s$
follows.

The second assertion in the proposition is also a consequence of the
above arguments. This is clear when $\varphi$ is smooth. The general
case is deduced using an approximation of $\varphi$ by a decreasing
sequence of smooth p.s.h. functions.
\endproof

Let $\Jac(F)$ denote the Jacobian of $F$ with respect to the standard
volume form on $\Sigma\times \C^k$. Its restriction to
$\pi^{-1}(s)$ is the Jacobian $\Jac(f_s)$ of
$f_s$. Since $\Jac(F)$ is a p.s.h. function, we can apply the previous
proposition and deduce that the function 
$L_k(s):={1\over 2}\langle \mu_s, \log\Jac(f_s)\rangle$ is p.s.h. on
$\Sigma$. Indeed, by Theorem \ref{th_log_jac_polylike}, this function is bounded from below by
${1\over 2}\log d_t$, hence it is not equal to $-\infty$. By Oseledec's theorem
\ref{th_oseledec}, $L_k(s)$ is the sum of the Lyapounov exponents of
$f_s$. We deduce the following result of \cite{DinhSibony1}.

\begin{corollary} \label{cor_sum_lya_poly}
Let $(f_s)_{s\in\Sigma}$ be as above. Then, the sum
  of the Lyapounov exponents associated to the equilibrium measure
  $\mu_s$ of $f_s$ is a p.s.h. function on $s$. In particular, it is
  upper semi-continuous.
\end{corollary}

Pham
defined in \cite{Pham}  {\it the bifurcation $(p,p)$-currents} by $\Bc^p:=(\ddc L_k)^p$ for
$1\leq p \leq \dim\Sigma$. The wedge-product is well-defined since $L_k$ is
locally bounded: it is bounded from below by ${1\over 2}\log d_t$. Very likely, these currents play a crucial role in
the study of bifurcation as we see in the following observation.
Assume that the critical set of $f_{s_0}$ does not
  intersect the filled Julia set $\Kc_{s_0}$ for some $s_0\in\Sigma$.
Since the filled Julia sets $\Kc_s$ vary upper semi-continuously in
  the Hausdorff metric,    $\log\Jac(F)$ is pluriharmonic near $\{s_0\}\times
  \Kc_{s_0}$. It follows that $L_k$
  is pluriharmonic and $\Bc^p=0$ in a neighbourhood of
  $s_0$\footnote{This observation was made by the second author for the family $z^2+c$, with
  $c\in\C$. He showed that the bifurcation measure is the harmonic measure
  associated to the Mandelbrot set \cite{Sibony_O}.}. On the other
hand, using Kobayashi metric, it is easy to show that $f$ is uniformly
hyperbolic on $\Kc_s$ for $s$ close to $s_0$. It follows that 
$\Kc_s=\Jc_s$ and $s\mapsto \Jc_s$ is continuous near $s_0$, see \cite{FornaessSibony6}. 

Note that $L_k$ is equal in the sense of currents to $\pi_*(\log\Jac(F)\wedge \Rc)$,
where $\Rc$ is the current in Proposition \ref{prop_pham}. Therefore,
$\Bc$ can be obtained using the formula
$$\Bc=\ddc\pi_*(\log\Jac(F)\wedge \Rc)= \pi_*([\Cc_F]\wedge \Rc),$$
since $\ddc\log|\Jac(F)|=[\Cc_F]$, the
current of integration on the critical set $\Cc_F$ of $F$. 
We also have the following property of the function $L_k$.

\begin{theorem} \label{th_variation_poly}
Let $(f_s)_{s\in\Sigma}$ be a family of
  polynomial-like maps as above. Assume that $f_{s_0}$ has a large topological
  degree for some $s_0\in\Sigma$. Then $L_k$ is H\"older continuous
  in a neighbourhood of $s_0$. In particular, the bifurcation currents
  $\Bc^p$ are moderate for $1\leq p\leq\dim \Sigma$.
\end{theorem}

Let $\Lambda_s$ denote the Perron-Frobenius operator associated to
$f_s$. For any Borel set $B$, denote by $\Omega_B$
the standard volume form on $\C^k$ restricted to $B$. 
We first prove some preliminary results. 

\begin{lemma} \label{lemma_variation_PB}
Let $W$ be a neighbourhood of the filled Julia set $\Kc_{s_0}$ of
$f_{s_0}$. Then, there is a neighbourhood $\Sigma_0$ of $s_0$ such that
$\langle\mu_s,\varphi\rangle$ depends continuously on $(s,\varphi)$
in $\Sigma_0\times \PSH(W)$. 
\end{lemma}
\proof
We first replace $\Sigma$ by a neighbourhood $\Sigma_0$ of $s_0$ small
enough. So, for every $s\in\Sigma$, the
filled Julia set of $f_s$ is contained in $U:=f_{s_0}^{-1}(V)$ and in
$W$.
We also reduce the size of $V$ in order to assume that $f_s$ is
polynomial-like on a neighbourhood of $\overline U$ with values in a neighbourhood $V'$ of $\overline V$.  
Moreover, since $\langle \mu_s,\varphi\rangle = \langle
\mu_s,\Lambda_s(\varphi)\rangle$ and $\Lambda_s(\varphi)$ depends
continuously on $(s,\varphi)$ in $\Sigma\times \PSH(W)$,
we can replace $\varphi$ by $\Lambda_s^N(\varphi)$
with $N$ large enough and $s\in\Sigma$, in order to assume that
$W=V$. Finally, since $\Lambda_s(\varphi)$ is defined on $V'$, it is
enough to prove the continuity for $\varphi$ p.s.h. on $V$ 
 such that $\varphi\leq 1$ and
$\langle\Omega_U,\varphi\rangle \geq 0$. Denote by $\Pc$ the family of
such functions $\varphi$. Since $\mu_{s_0}$ is PC, we have
$|\langle\mu_{s_0},\varphi\rangle|\leq A$ for some constant
$A\geq 1$ and for $\varphi\in\Pc$. 
Let $\Pc'$ denote the family of p.s.h. functions $\psi$ such that
$\psi\leq 2A$ and $\langle\mu_{s_0},\psi\rangle=0$. 
 The function
$\varphi':=\varphi-\langle\mu_{s_0},\varphi\rangle$ belongs to this
family. Observe that $\Pc'$ is bounded and therefore if $A'\geq 1$ is
a fixed constant large enough, we have 
$|\langle \Omega_U,\psi\rangle|\leq A'$ for $\psi\in\Pc'$.

Fix an integer $N$ large enough. By Theorem \ref{th_top_PC},
$\Lambda_{s_0}^N(\varphi')\leq 1/8$ on $V'$ and 
$|\langle \Omega_U,\Lambda_{s_0}^N(\varphi')\rangle|\leq 1/8$ for
$\varphi'$ as above.
We deduce that 
$2\Lambda_{s_0}^N(\varphi')-\langle \Omega_U,
2\Lambda_{s_0}^N(\varphi')\rangle$ is a function in $\Pc$, smaller
than $1/2$ on $V'$. This
function differs from $2\Lambda_{s_0}^N(\varphi)$ by a constant. So, it is
equal to
$2\Lambda_{s_0}^N(\varphi)-\langle\Omega_U,2\Lambda_{s_0}^N(\varphi)\rangle$.
When $\Sigma_0$ is small enough,  by continuity, 
the operator
$L_s(\varphi):=2\Lambda_s^N(\varphi)-\langle\Omega_U,2\Lambda_s^N(\varphi)\rangle$
preserves $\Pc$ for $s\in\Sigma_0$. Therefore, since $\Lambda_s$ preserves constant
functions, we have
\begin{eqnarray*}
\Lambda_s^{mN}(\varphi) & = & \Lambda_s^{(m-1)N}\big[\langle\Omega_U,\Lambda_s^N(\varphi)\rangle
+ 2^{-1} L_s(\varphi)\big] \\
& = & \langle\Omega_U,\Lambda_s^N(\varphi)\rangle
+2^{-1}\Lambda_s^{(m-1)N}(L_s(\varphi)).
\end{eqnarray*}
By induction, we obtain
\begin{eqnarray*}
\Lambda_s^{mN}(\varphi) & = &
\langle\Omega_U,\Lambda_s^N(\varphi)\rangle +\cdots + 2^{-m+1} \langle
\Omega_U,\Lambda_s^N(L_s^{m-1}(\varphi))\rangle + 2^{-m}
L_s^m(\varphi)\\
& = & \big\langle\Omega_U,\Lambda_s^N\big[\varphi +\cdots + 
2^{-m+1}L_s^{m-1}(\varphi)\big]\big\rangle + 2^{-m}
L_s^m(\varphi)\\
& = & \big\langle d_t^{-N} (f_s^N)^*(\Omega_U),
\varphi +\cdots + 2^{-m+1}L_s^{m-1}(\varphi)\big\rangle + 2^{-m}
L_s^m(\varphi).
\end{eqnarray*}
We deduce from the above property of $L_s$ that the last term converges uniformly to 0 when $m$ goes to
infinity. The sum in the first term converges normally to the
p.s.h. function $\sum_{m\geq 1}
2^{-m+1}L_s^{m-1}(\varphi)$, which depends continuously
on $(s,\varphi)$. Therefore,
$\Lambda_s^{mN}(\varphi)$ converge to a constant which depends
continuously on $(s,\varphi)$.
But we know that the limit is
$\langle \mu_s,\varphi\rangle$.
The lemma follows.
\endproof

Using the same approach as in Theorem \ref{th_pol_like_holder}, we
prove the following result.

\begin{theorem} \label{th_pol_fam_holder}
Let $f_s$, $s_0$ and $W$ be as in
  Theorem \ref{th_variation_poly} and Lemma
  \ref{lemma_variation_PB}. Let $K$ be a compact subset of $W$ such
  that $f_{s_0}^{-1}(K)$ is contained in the interior of $K$. There is a neighbourhood $\Sigma_0$ of $s_0$
  such that if $\Pc$ is a bounded family of p.s.h. functions on $W$,
  then  $(s,\varphi)\mapsto \langle\mu_s,\varphi\rangle$ is
  H\"older continuous on $\Sigma_0\times\Pc$ with respect to
  the pseudo-distance $\dist_{L^1(K)}$ on $\Pc$.
\end{theorem}
\proof
We replace $\Sigma$ by $\Sigma_0$ as in Lemma
\ref{lemma_variation_PB}. 
It is not difficult to check that
$(s,\varphi)\mapsto (s,\Lambda_s(\varphi))$ is locally Lipschitz with respect to
$\dist_{L^1(K)}$. So, replacing $(s,\varphi)$ by $(s,\Lambda_s^N(\varphi))$
with $N$ large enough allows to assume that $W=V$. 
Let $\widehat\Pc$ be the set of $(s,\varphi)$ in $\Sigma\times
\PSH(V)$ such that $\varphi\leq 1$ and $\langle\mu_s,\varphi\rangle\geq 0$. 
By Lemma \ref{lemma_variation_PB}, such functions $\varphi$ belong to
a compact subset of $\PSH(V)$. It is enough to prove that $(s,\varphi)\mapsto \langle\mu_s,\varphi\rangle$ is
  H\"older continuous on $\widehat\Pc$.

Let $\widehat \Dc$ denote the set of $(s,\varphi-\Lambda_s(\varphi))$
with $(s,\varphi)\in\widehat\Pc$. Consider the operator 
$\widehat\Lambda(s,\psi):=(s,\lambda^{-1}\Lambda_s(\psi))$ on $\Dc$ as
in Theorem \ref{th_pol_like_holder} where $\lambda<1$ is a fixed
constant close enough to 1. Theorem \ref{th_top_PC} and the continuity in Lemma
\ref{lemma_variation_PB} imply that $\widehat\Lambda$ preserves
$\widehat\Dc$. Therefore, we only have to follow the arguments in
Theorem \ref{th_pol_like_holder}. 
\endproof

\noindent
{\bf Proof of Theorem \ref{th_variation_poly}.}
We replace $\Sigma$ by a small neighbourhood of $s_0$.
Observe that $\log\Jac(f_s)$, $s\in\Sigma$, is a bounded family of p.s.h. functions
on $U$. By Theorem \ref{th_pol_fam_holder}, it is enough to show that  
$s\mapsto \log\Jac(f_s)$ is H\"older continuous with
respect to $\dist_{L^1(K)}$. 

We also deduce from 
Theorem \ref{th_hormander} that $\langle \Omega_K, e^{\lambda
  |\log\Jac(f_s)|}\rangle \leq A$ for some positive constants $\lambda$ and $A$.
Reducing $V$ and $\Sigma$ allows
to assume that $\Jac(F)$, their derivatives   
and the vanishing order of $\Jac(F)$ are bounded on $\Sigma\times U$ by some
constant $m$.

Fix a constant $\alpha>0$ small enough and a constant $A>0$ large enough.
Define $\psi(s):=\langle\Omega_K,\log\Jac(f_s)\rangle$.  
Consider $s$ and $t$ in $\Sigma$ such
that $r:=\|s-t\|$ is smaller than a fixed small constant. We will compare
$|\psi(s)-\psi(t)|$ with $r^{\lambda\alpha}$ in order to show that
$\psi$ is H\"older continuous with exponent $\lambda\alpha$. Define $S:=\{z\in U,\
\Jac(f_s)<2 r^{2\alpha}\}$. We will bound separately 
$$\langle
\Omega_{K\setminus S},\log\Jac(f_s)-\log\Jac(f_t)\rangle$$ 
and
$$\langle
\Omega_{K\cap S},\log\Jac(f_s)-\log\Jac(f_t)\rangle.$$ 
Note that
$\psi(s)-\psi(t)$ is the sum of the above two integrals.

Consider now the integral on $K\setminus S$. The following estimates
are only valid on $K\setminus S$.  
Since the derivatives of $\Jac(F)$ is
bounded, we have $\Jac(f_t)\geq r^{2\alpha}$.
It follows that the derivatives on $t$ of
$\log\Jac(f_t)$ is bounded by $Ar^{-2\alpha}$. We deduce that
$$|\log\Jac(f_s)-\log\Jac(f_t)|\leq Ar^{1-2\alpha}.$$
Therefore, 
$$|\langle \Omega_{K\setminus
  S},\log\Jac(f_s)-\log\Jac(f_t)\rangle|\leq \|\Omega_{K\setminus S}\|
Ar^{1-2\alpha}\leq r^{\lambda\alpha}.$$

We now estimate the integral on $K\cap S$. Its absolute value is
bounded by 
$$\langle \Omega_{K\cap S},|\log\Jac(f_s)|\rangle
+\langle \Omega_{K\cap S},|\log\Jac(f_t)|\rangle.$$ 
We deduce from the estimate  $\langle \Omega_K, e^{\lambda
  |\log\Jac(f_s)|}\rangle \leq A$ that
$\volume(K\cap S)\leq Ar^{2\lambda\alpha}$. 
Therefore, by Cauchy-Schwarz's inequality, we have
\begin{eqnarray*}
\langle \Omega_{K\cap S},|\log\Jac(f_s)|\rangle & \lesssim & 
\vol(K\cap S)^{1/2} \big\langle \Omega_K,|\log\Jac(f_s)|^2\big\rangle^{1/2} \\
& \lesssim & r^{\lambda\alpha} \langle
\Omega_K,e^{\lambda|\log\Jac(f_s)|}\rangle^{1/2}\lesssim
r^{\lambda\alpha}.
\end{eqnarray*} 
The estimate holds for $f_t$ instead of $f_s$. Hence, $\psi$ is H\"older
continuous. The fact that $\Bc^p$ are moderate follows from Theorem \ref{th_moderate}. 
\hfill $\square$

\medskip

The following result of Pham generalizes Corollary
\ref{cor_sum_lya_poly} 
and allows to define other bifurcation
currents by considering $\ddc L_p$ or their wedge-products \cite{Pham}.

\begin{theorem} Let $(f_s)_{s\in\Sigma}$ be a holomorphic family of polynomial-like
  maps as above. Let $\chi_1(s)\geq \cdots\geq \chi_k(s)$ be the
  Lyapounov exponents of the equilibrium measure $\mu_s$ of
  $f_s$. Then, for $1\leq p\leq k$, the function
$$L_p(s):=\chi_1(s)+\cdots+\chi_p(s)$$
is p.s.h. on $\Sigma$. In particular, $L_p$ is upper semi-continuous.
\end{theorem}
\proof
Observe that $L_p(s)\geq {p\over k} L_k(s)\geq {p\over 2k}\log d_t$. 
We identify the tangent space of $V$ at any point with $\C^k$. So, the
differential $Df_s(z)$ of $f_s$ at a point $z\in U_s$ is a linear
self-map on $\C^k$ which depends holomorphically on $(s,z)$. It induces a
linear self-map on the exterior product $\bigwedge^p\C^k$ that we
denote by $D^pf_s(z)$. This map depends holomorphically on
$(s,z)$. In the standard coordinate system on  $\bigwedge^p\C^k$,
the function $(s,z)\mapsto \log\|D^pf_s(z)\|$ is p.s.h. on
$U_\Sigma$. By Proposition \ref{prop_pham}, the function
$\psi_1(s):=\langle\mu_s,\log\|D^pf_s\|\rangle$ is p.s.h. or equal to
$-\infty$ on $\Sigma$. Define in the same way the functions 
$\psi_n(s):=\langle \mu_s, \log \|D^pf_s^n\|\rangle$ associated to the
iterate $f_s^n$ of $f_s$. We have
$$D^pf_s^{n+m}(z)=D^pf_s^m(f_s^n(z))\circ D^pf_s^n(z).$$ 
Hence,
$$\|D^pf_s^{n+m}(z)\|\leq \|D^pf_s^m(f_s^n(z))\|\ \| D^pf_s^n(z)\|.$$
We deduce using the invariance of $\mu_s$ that
$$\psi_{m+n}(s)\leq \psi_m(s)+\psi_n(s).$$
Therefore, the sequence $n^{-1}\psi_n$ decreases to $\inf_n
n^{-1}\psi_n$. So, the limit is p.s.h. or equal to $-\infty$.
On the other hand, Oseledec's theorem \ref{th_oseledec} implies that
the limit is equal to $L_p(s)$ which is a positive function. It
follows that $L_p(s)$ is p.s.h.
\endproof

Consider now the family $f_s$ of endomorphisms of algebraic degree
$d\geq 2$ of $\P^k$ with $s\in\Hc_d(\P^k)$. We can lift $f_s$ to
polynomial-like maps on $\C^{k+1}$ and apply the above results. The
construction of the bifurcation currents $\Bc^p$ can be obtained
directly using the Green measures of $f_s$. This was done by
Bassanelli-Berteloot in \cite{BassanelliBerteloot}. They also studied
some properties of the bifurcation currents and obtained
nice formulas for that currents in terms of the Green
functions. We also refer to DeMarco, Dujardin-Favre, McMullen and
Sibony \cite{DeMarco,DujardinFavre,McMullen,
  Sibony_O} for results in dimension one.

\bigskip\bigskip

\begin{exercise}
If $f$ is an endomorphism in $\Hc_d(\P^k)$, denote by $L_k(f)$ the sum
of the Lyapounov exponents of the equilibrium measure. Show that
$f\mapsto L_k(f)$ is locally H\"older continuous on
$\Hc_d(\P^k)$. Deduce that the bifurcation currents are
moderate. Hint: use that the lift of $f$ to $\C^{k+1}$ has always a
Lyapounov exponent equal to $\log d$.
\end{exercise}

\begin{exercise} Find a family $(f_s)_{s\in\Sigma}$ such that $\Jc_s$
  does not vary continuously.
\end{exercise}

\begin{exercise}
A family 
$(X_s)_{s\in\Sigma}$ of compact subsets in $V$  is lower semi-continuous at $s_0$ if for
every $\epsilon>0$, $X_{s_0}$ is contained in the $\epsilon$-neighbourhood
of $X_s$ when $s$ is close enough to $s_0$. If
$(\nu_s)_{s\in\Sigma}$ is a continuous family of probability measures
on $V$, show that $s\mapsto\supp(\nu_s)$ is lower semi-continuous. 
If $(f_s)_{s\in\Sigma}$ is a holomorphic family of polynomial-like
maps, deduce that $s\mapsto \Jc_s$ is lower semi-continuous. Show that
if $\Jc_{s_0}=\Kc_{s_0}$, then $s\mapsto \Jc_s$ is continuous at
$s_0$ for the Hausdorff metric.
\end{exercise}

\begin{exercise}
Assume that $f_{s_0}$ is of large topological degree. Let $\delta>0$ be
a constant small enough. Using the
continuity of $s\mapsto\mu_s$, show that if $p_{s_0}$ is a repelling
fixed point in $\Jc_{s_0}$ for $f_{s_0}$, there are repelling fixed
points $p_s$ in $\Jc_s$ for $f_s$, with $|s-s_0|<\delta$, such that
$s\mapsto p_s$ is holomorphic. Suppose $s\mapsto \Jc_s$ is continuous
with respect to the Hausdorff metric. Construct a positive closed
current $\Rc$ supported on $\cup_{|s-s_0|<\delta} \{s\}\times \Jc_s$ with slices
$\mu_s$. Deduce that if $\Jc_{s_0}$ does not contain critical points
of $f_{s_0}$ then $s\mapsto L_k(s)$ is pluriharmonic near $s_0$.
\end{exercise}

\bigskip\bigskip

{\small

\addcontentsline{toc}{section}{Notes}
\noindent
{\bf Notes.} Several results in this chapter
  still hold for larger classes of polynomial-like maps. For example, the
  construction of the equilibrium measure is valid for a manifold $V$
  admitting a smooth strictly p.s.h. function.
The $\ddc$-method was originally introduced for polynomial-like
maps. However, we have seen that it is also effective for
endomorphisms of $\P^k$. In a forthcoming survey, we will show that
the method can be extended to other dynamical systems. 
Several statistical properties obtained in this chapter are new.
}


\begin{appendix}

\chapter{Currents and pluripotential theory}

In this appendix, we recall some basic notions and results
on complex geometry and on currents in the complex setting.
Most of the results are classical and their proofs are not given here.
In constrast, we describe in detail some notions 
in order to help the reader who are not familiar with
complex geometry or currents. The main references for the abstract theory of
currents are \cite{Chemin, deRham, Federer, Schwartz, Warner}. The
reader will find in \cite{Demailly3, Gunning, Hormander1, Lelong,
  Narasimhan} the basics on currents on complex manifolds. We also refer
to \cite{Demailly3, GriffithsHarris, Huybrechts, Voisin} for
the theory of compact
K\"ahler manifolds.


\section{Projective spaces and analytic sets} \label{section_pk}

In this paragraph, we recall the definition of complex projective 
spaces. We then discuss briefly compact K\"ahler manifolds, projective
manifolds and analytic sets.

The complex projective space $\P^k$ is a compact complex manifold of
dimension $k$. It
is obtained as the quotient of $\C^{k+1}\setminus\{0\}$ by the natural
multiplicative action of $\C^*$. In other words, $\P^k$ is the
parameter space of the complex lines passing through $0$ in
$\C^{k+1}$.  The image of a subspace
of dimension $p+1$ of $\C^{k+1}$ is a submanifold of dimension $p$ in
$\P^k$, bi-holomorphic to $\P^p$, and is
called {\it a projective subspace of dimension $p$}. {\it Hyperplanes}
of $\P^k$ are projective
subspaces of dimension $k-1$. The group $\GL(\C,k+1)$ of invertible
linear endomorphisms of $\C^{k+1}$ induces the group $\PGL(\C,k+1)$ of
automorphisms of $\P^k$. It acts transitively on $\P^k$ and sends
projective subspaces to projective subspaces. 

Let $z=(z_0,\ldots,z_k)$ denote the standard coordinates
of $\C^{k+1}$. Consider the equivalence relation: {\it $z\sim z'$ if there
is $\lambda\in\C^*$ such that $z=\lambda z'$}. The projective space $\P^k$ is the
quotient of $\C^{k+1}\setminus\{0\}$ by this relation. We can
cover $\P^k$ by open sets $U_i$ associated to the open sets
$\{z_i\not=0\}$ in $\C^{k+1}\setminus \{0\}$. Each $U_i$ is
bi-holomorphic to $\C^k$ and $(z_0/z_i,\ldots,
z_{i-1}/z_i,z_{i+1}/z_i,\ldots,z_k/z_i)$ is a coordinate system on
this chart. The complement of $U_i$ is the hyperplane defined by $\{z_i=0\}$.
So, $\P^k$ can be considered as 
a natural compactification of $\C^k$. We denote by $[z_0:\cdots:z_k]$
the point of $\P^k$ associated to $(z_0,\ldots, z_k)$. This expression
is {\it the homogeneous coordinates} on $\P^k$. Projective spaces are
compact K\"ahler manifolds. We will describe this notion later.

Let $X$ be a complex manifold of dimension $k$. Let $\varphi$ be a
differential $l$-form on $X$. 
In local holomorphic
coordinates $z=(z_1,\ldots,z_k)$, it can be written as
$$\varphi(z)=\sum_{|I|+|J|=l}\varphi_{IJ} dz_I\wedge d\overline z_J,$$
where $\varphi_{IJ}$ are complex-valued functions,
$dz_I:=dz_{i_1}\wedge\ldots\wedge dz_{i_p}$ if
$I=(i_1,\ldots,i_p)$, and  $d\overline
z_J:=d\overline z_{j_1}\wedge\ldots\wedge d\overline z_{j_q}$ if $J=(j_1,\ldots,j_q)$. 
The {\it conjugate} of
$\varphi$ is 
$$\overline\varphi(z):=\sum_{|I|+|J|=l}\overline \varphi_{IJ} d\overline z_I\wedge d z_J.$$
The form $\varphi$ is real if and only if $\varphi=\overline\varphi$. 

We say that $\varphi$ is a form of {\it of bidegree} $(p,q)$ if
$\varphi_{IJ}=0$ when $(|I|,|J|)\not=(p,q)$. The bidegree does not
depend on the choice of local coordinates.
Let $T_X^\C$ denote the complexification of the tangent bundle of
$X$. The complex structure on $X$ induces a linear endomorphism $\Jc$ on
the fibers of $T_X^\C$ such that $\Jc^2=-\id$. This endomorphism
induces a decomposition of 
$T_X^\C$ into the direct sum of two proper sub-bundles of dimension $k$: {\it the
  holomorphic part} $T_X^{1,0}$ associated to the eigenvalue 1 of
$\Jc$, and the {\it anti-holomorphic part}
$T_X^{0,1}$ associated to the eigenvalue $-1$.
Let $\Omega_X^{1,0}$ and $\Omega_X^{0,1}$ denote the dual bundles
of $T_X^{1,0}$  and  $T_X^{0,1}$. Then,  $(p,q)$-form are sections of
the vector bundle $\bigwedge^p\Omega^{1,0}\otimes \bigwedge^q\Omega^{0,1}$.

If $\varphi$ is a
$(p,q)$-form then the differential $d\varphi$ is the sum of a $(p+1,q)$-form 
and a $(p,q+1)$-form. We then denote by $\partial \varphi$ the part of
bidegree $(p+1,q)$ and   $\dbar \varphi$ the  the part of
bidegree $(p,q+1)$. The operators $\partial$ and $\dbar$ extend
linearly to arbitrary forms $\varphi$. 
The operator $d$ is real, i.e. it
sends real forms to real forms but
$\partial$ and $\dbar$ are not
real. The identity $d\circ d=0$ implies that
$\partial\circ\partial=0$, $\dbar\circ\dbar=0$ and $\partial\dbar+\dbar\partial=0$. 
Define $\dc:={\sqrt{-1}\over2\pi}(\dbar-\partial)$. This operator is real
and satisfies $\ddc = {\sqrt{-1}\over \pi}\ddbar$. 
Note that the above operators 
commute with the pull-back by holomorphic maps. More precisely, if
$\tau:X_1\rightarrow X_2$ is a holomorphic map between complex
manifolds and $\varphi$ is a form on $X_2$ then
$df^*(\varphi)=f^*(d\varphi)$, $\ddc f^*(\varphi)=f^*(\ddc\varphi)$, etc.
Recall that the form $\varphi$ is {\it closed}
(resp. $\partial$-closed, $\dbar$-closed, $\ddc$-closed) if $d\varphi$
(resp. $\partial\varphi$, $\dbar\varphi$, $\ddc\varphi$) vanishes. The
form $\varphi$ is {\it exact} (resp. $\partial$-exact, $\dbar$-exact,
$\ddc$-exact) if it is equal to the differential
$d\psi$ (resp. $\partial \psi$, $\dbar\psi$, $\ddc\psi$) of a form $\psi$. Clearly, exact forms are closed.

A smooth $(1,1)$-form $\omega$ on $X$ is {\it Hermitian} if it
can be written in local coordinates as 
$$\omega(z)=\sqrt{-1}\sum_{1\leq i,j\leq k} \alpha_{ij}(z) dz_i\wedge d\overline z_j,$$ 
where $\alpha_{ij}$ are smooth functions such that the matrix
$(\alpha_{ij})$ is Hermitian. We consider a form $\omega$ such that
the matrix $(\alpha_{ij})$ is
positive definite at every point. It is
strictly positive in the sense that we will introduce later.
If $a$ is a point in
$X$, we can find local coordinates $z$ such that
$z=0$ at $a$ and $\omega$ is equal near 0 to the Euclidean form
$\ddc\|z\|^2$ modulo a term of order $\|z\|$.
The form $\omega$ is always real
and induces a norm on the tangent spaces of $X$. So, it defines a
Riemannian metric on $X$. 
We say that $\omega$ is {\it a K\"ahler form}
if it is a closed positive definite Hermitian form. In this case, one can find  local coordinates $z$ such that
$z=0$ at $a$ and $\omega$ is equal near 0 to 
$\ddc\|z\|^2$ modulo a term of order $\|z\|^2$. So, at the
infinitesimal level, a K{\"a}hler metric is close to the Euclidean
one. This is a crucial property in Hodge
theory in the complex setting.

Consider now a compact complex manifold $X$ of dimension $k$. 
Assume that $X$ is {\it a K\"ahler manifold}, i.e. it admits a
K\"ahler form $\omega$.  
Recall that {\it the de Rham cohomology group} $H^l(X,\C)$ is
the quotient of the space of closed $l$-forms by the
subspace of exact $l$-forms. This complex vector space is of finite dimension.
The real groups  $H^l(X,\R)$ are defined in the
same way using real forms. We  have 
$$H^l(X,\C)= H^l(X,\R)\otimes_\R\C.$$
If $\alpha$ is a closed $l$-form, its class in $H^l(X,\C)$ is denoted
by $[\alpha]$. The group $H^0(X,\C)$ is just the set of constant
functions. So, it is isomorphic to $\C$. The group $H^{2k}(X,\C)$ is also
isomorphic to $\C$. The isomorphism is given by the canonical map
$[\alpha]\mapsto\int_X\alpha$. 
For $l,m$ such that $l+m\leq 2k$, {\it the cup-product}
$$\smile:H^l(X,\C)\times H^m(X,\C)\rightarrow H^{l+m}(X,\C)$$ 
is defined by 
$[\alpha]\smile [\beta]:=[\alpha\wedge\beta]$. The Poincar\'e duality
theorem says that the cup-product is a non-degenerated bilinear form
when $l+m=2k$. So, it defines an isomorphism between $H^l(X,\C)$ and
the dual of $H^{2k-l}(X,\C)$. 

Let $H^{p,q}(X,\C)$, $0\leq p,q\leq k$, denote the subspace of
$H^{p+q}(X,\C)$ generated by the classes
of closed $(p,q)$-forms.
We call $H^{p,q}(X,\C)$ the {\it Hodge cohomology group}. 
Hodge theory shows that
$$H^{l}(X,\C)=\bigoplus_{p+q=l} H^{p,q}(X,\C) \qquad \mbox{and}\qquad
H^{q,p}(X,\C)=\overline{H^{p,q}(X,\C)}.$$ 
This, together with the Poincar{\'e} duality, induces a canonical isomorphism
between $H^{p,q}(X,\C)$ and the dual space of
$H^{k-p,k-q}(X,\C)$. 
Define for $p=q$
$$H^{p,p}(X,\R):=H^{p,p}(X,\C)\cap H^{2p}(X,\R).$$
We have
$$H^{p,p}(X,\C)=H^{p,p}(X,\R)\otimes_\R \C.$$

Recall that {\it the Dolbeault cohomology group} $H^{p,q}_{\dbar}(X)$ is
the quotient of the space of $\dbar$-closed $(p,q)$-forms by the
subspace of $\dbar$-exact $(p,q)$-forms. Observe that a 
$(p,q)$-form is $d$-closed if and only if it is $\partial$-closed and
$\dbar$-closed. Therefore, there is a natural morphism between the 
Hodge and the Dolbeault cohomology groups. Hodge theory asserts that
this is in fact an isomorphism: we have
$$H^{p,q}(X,\C)\simeq H^{p,q}_\dbar(X).$$
The result is a consequence of the following theorem, the so-called
{\it $\ddc$-lemma}, see e.g. \cite{Demailly3, Voisin}.

\begin{theorem} Let $\varphi$ be a smooth $d$-closed $(p,q)$-form on
  $X$. Then $\varphi$ is $\ddc$-exact if and only if it is $d$-exact (or 
  $\partial$-exact or $\dbar$-exact).
\end{theorem}

The projective space $\P^k$ admits a K{\"a}hler form $\omega_\FS$,
called {\it the Fubini-Study form}. It 
is defined on the chart $U_i$ by 
$$\omega_\FS:=\ddc\log \Big(\sum_{j=0}^k \Big|{z_j\over
  z_i}\Big|^2\Big).$$
In other words, if $\pi:\C^{k+1}\setminus\{0\}\rightarrow\P^k$ is the
canonical projection, then $\omega_\FS$ is defined by
$$\pi^*(\omega_\FS):=\ddc\log \Big(\sum_{i=0}^k \big|z_i|^2\Big).$$
One can check that $\omega_\FS^k$ is a probability measure on
$\P^k$. 
The cohomology groups of $\P^k$ are very
simple. We have $H^{p,q}(\P^k,\C)=0$ for $p\not=q$ and
$H^{p,p}(\P^k,\C)\simeq\C$. The groups $H^{p,p}(\P^k,\R)$
and $H^{p,p}(X,\C)$ are generated by the class of $\omega_\FS^p$.
Submanifolds of $\P^k$ are 
K\"ahler, as submanifolds of a K\"ahler manifold.
Chow's theorem says that such a manifold is 
algebraic, i.e. it is the set of common zeros of a finite family of
homogeneous polynomials in $z$. A compact manifold is {\it projective}
if it is bi-holomorphic to a submanifold of a projective
space. Their cohomology groups are in general very rich and difficult
to describe.

A useful result of Blanchard
\cite{Blanchard} says that the blow-up of a compact K{\"a}hler manifold along a
submanifold is always a compact K{\"a}hler manifold. The construction of
the blow-up is as follows. Consider first the case of open sets in  
$\C^k$ with $k\geq 2$. Observe that $\C^k$ is the union of the
complex lines passing through 0 which are parametrized
by the projective space $\P^{k-1}$. {\it The blow-up} $\widehat{\C^k}$ of
$\C^k$ at 0 is obtained by separating these complex lines, that is, we keep
$\C^k\setminus\{0\}$ and replace 0 by a copy of $\P^{k-1}$. More
precisely, if
$z=(z_1,\ldots,z_k)$ denote the coordinates of $\C^k$ and
$[w]=[w_1:\cdots:w_k]$ are homogeneous coordinates of $\P^{k-1}$, then
$\widehat{\C^k}$ is the submanifold of $\C^k\times\P^{k-1}$ defined by
the equations $z_iw_j=z_jw_i$ for $1\leq i,j\leq k$. 
If $U$ is an open set in $\C^k$ containing 0, the blow-up $\widehat U$
of $U$ at 0 is defined by $\pi^{-1}(U)$ where
$\pi:\widehat{\C^k}\rightarrow \C^k$ is the canonical projection.

If $U$ is a neighbourhood of 0 in $\C^{k-p}$, $p\leq k-2$, and $V$ is an open set in
$\C^p$, then the blow-up of $U\times V$ along $\{0\}\times V$ is equal
to $\widehat{U}\times V$. Consider now a submanifold $Y$ of $X$
of dimension $p\leq k-2$. We cover $X$ by charts which either
do not intersect $Y$ or are of the type $U\times V$, where $Y$ is identified with
$\{0\}\times V$. {\it The blow-up} $\widehat{X}$ is obtained by sticking the
charts outside $Y$ with the blow-ups of charts which intersect $Y$. The
natural projection $\pi:\widehat X\rightarrow X$ defines a
bi-holomorphic map between $\widehat X\setminus \pi^{-1}(Y)$ and
$X\setminus Y$. The set $\pi^{-1}(Y)$ is a smooth
hypersurface, i.e. submanifold of codimension 1; it is called {\it the
  exceptional hypersurface}. Blow-up may
be defined using the local ideals of holomorphic functions vanishing on
$Y$. The blow-up of a projective manifold along a submanifold is a
projective manifold.

We now recall some facts on analytic
sets, see \cite{Gunning,Narasimhan}. 
Let $X$ be an arbitrary complex manifold of dimension $k$\footnote{We
  often assume that $X$ is connected for simplicity.}. 
Analytic sets of $X$ can be seen as submanifolds of $X$, possibly with
singularities.  Analytic sets of dimension 0 are locally finite
subsets, those of dimension 1 are (possibly singular) Riemann surfaces. For
example, $\{z_1^2=z_2^3\}$ is an analytic set of $\C^2$ of dimension 1
with a singularity at 0. Chow's theorem holds for analytic sets: any analytic set in $\P^k$ is
the set of common zeros of a finite family of homogeneous
polynomials.

Recall that an {\it analytic
  set} $Y$ of $X$ is locally the set of common zeros of holomorphic
functions: for every point $a\in X$ there is a
neighbourhood $U$ of $a$ and holomorphic functions $f_i$ on $U$ such
that $Y\cap U$ is the intersection of $\{f_i=0\}$. We can choose $U$
so that $Y\cap U$ is defined by a finite family of holomorphic
functions. Analytic sets are closed for the usual topology
on $X$. Local rings of holomorphic functions on $X$ induce local rings
of holomorphic functions on $Y$.
An analytic set $Y$ is {\it irreducible} if it is not a union of
two different non-empty analytic sets of $X$. A general analytic set $Y$
can be decomposed in a unique way into a union of irreducible analytic
subsets $Y=\cup Y_i$, where no component $Y_i$ is contained in another
one. The decomposition is locally finite, that is, given a
compact set $K$ in $X$, only finitely many  $Y_i$ intersect $K$.

Any increasing sequence of
irreducible analytic subsets of $X$ is stationary. A decreasing
sequence $(Y_n)$ of analytic subsets of $X$ is always locally
stationary, that is, for any compact subset $K$ of $X$, the sequence
$(Y_n\cap K)$ is stationary. Here, we do not suppose $Y_n$ irreducible.
The topology on $X$ whose closed sets are exactly the analytic
sets, is called {\it the Zariski topology}. When $X$ is connected, non-empty open Zariski sets
are dense in $X$ for the usual topology. The restriction of the Zariski topology on $X$ to
$Y$ is also called the {\it Zariski topology} of $Y$. When $Y$ is
irreducible, the non-empty Zariski open subsets are also dense in $Y$
but this is not the case for reducible analytic sets.

There is a minimal analytic subset $\sing(Y)$ in $X$ such that $Y\setminus \sing(Y)$
is a (smooth) complex submanifold of $X\setminus \sing(Y)$, i.e. a complex
manifold which is closed and without boundary in $X\setminus
\sing(Y)$. The analytic set $\sing(Y)$ is the {\it singular part} of
$Y$. The {\it regular part} of $Y$ is denoted by $\reg(Y)$; it is
equal to $Y\setminus\sing(Y)$. The manifold $\reg(Y)$ is not
necessarily irreducible; it may have several components.
We call {\it dimension of $Y$}, $\dim(Y)$, the maximum
of the
dimensions of these components; the {\it codimension} $\codim(Y)$ of
$Y$ in $X$ is the integer $k-\dim (Y)$. We say that $Y$ is a proper
analytic set of $X$ if it has positive codimension.
When all the components of $Y$ have the same dimension, we say that
$Y$ is {\it of pure dimension} or {\it of pure codimension}.
When $\sing(Y)$ is non-empty, its
dimension is always strictly smaller than the dimension of $Y$. We can
again decompose $\sing(Y)$ into regular and singular parts. The
procedure can be repeated less than $k$ times and gives a stratification of $Y$
into disjoint complex manifolds.
Note that $Y$ is irreducible if and only if $\reg(Y)$
is a connected manifold. 
The following result is due to Wirtinger.

\begin{theorem}[Wirtinger] \label{th_wirtinger}
Let $Y$ be analytic set of pure dimension
  $p$ of a Hermitian manifold $(X,\omega)$. Then the $2p$-dimensional
  volume of $Y$ on a Borel set $K$ is equal to
$$\vol(Y\cap K)={1\over p!}\int_{\reg(Y)\cap K} \omega^p.$$
Here, the volume is with respect to the Riemannian metric induced by
$\omega$. 
\end{theorem}

Let $D_k$ denote the unit polydisc $\{|z_1|<1,\ldots, |z_k|<1\}$ in $\C^k$.  
The following result describes the local
structure of analytic sets.

\begin{theorem} \label{th_anal_set_local}
Let $Y$ be an analytic set of pure dimension $p$ of
$X$. Let $a$ be a point of $Y$. Then there
is a holomorphic chart $U$ of $X$, bi-holomorphic to $D_k$, with local coordinates
$z=(z_1,\ldots,z_k)$, such that $z=0$ at $a$, $U$ is given by
$\{|z_1|<1,\ldots, |z_k|<1\}$ and the projection $\pi:U\rightarrow
D_p$, defined by
$\pi(z):=(z_1,\ldots ,z_p)$, is
proper on $Y\cap U$. In this case, there is a proper analytic subset $S$
of $D_p$ such that $\pi:Y\cap
U\setminus \pi^{-1}(S)\rightarrow D_p\setminus S$ is a finite covering
and the singularities of $Y$ are contained in $\pi^{-1}(S)$.
\end{theorem}

Recall that a holomorphic map $\tau:X_1\rightarrow X_2$ between complex
manifolds of the same dimension is {\it a covering of degree $d$} if each
point of $X_2$ admits a neighbourhood $V$ such that $\tau^{-1}(V)$ is a
disjoint union of $d$ open sets, each of which is sent bi-holomorphically to
$V$ by $\tau$. Observe the previous theorem also implies that the
fibers of $\pi:Y\cap U\rightarrow D_p$ are finite and contain at most
$d$ points if $d$ is the degree of the covering. We can reduce $U$ in
order to have that $a$ is the unique point in the fiber
$\pi^{-1}(0)\cap Y$. The degree $d$ of the covering depends
on the choice of coordinates and the smallest integer $d$ obtained in this
way is called {\it the multiplicity} of $Y$ at $a$ and is denoted by $\mult(Y,a)$. 
We will see that $\mult(Y,a)$ is the 
Lelong number at $a$ of the positive closed current associated to
$Y$. In other words, if $B_r$ denotes the ball of center $a$ and of radius $r$, then the
ratio between the volume of $Y\cap B_r$ and the volume of a ball of
radius $r$ in $\C^p$ decreases to $\mult(Y,a)$ when $r$ decreases to 0.   

Let $\tau:X_1\rightarrow X_2$ be an open holomorphic map between complex
manifolds of the same dimension. Applying the above result to the
graph of $\tau$, we can show that for any point $a\in X_1$ and for a 
neighbourhood $U$ of $a$ small enough, if $z$ is a generic point in $X_2$
close enough to $\tau(a)$, the number of points in
$\tau^{-1}(z)\cap U$ does not depend on $z$. We call this number {\it
  the multiplicity} or {\it the local topological degree} of
$\tau$ at $a$. We say that $\tau$ 
 is {\it a ramified covering of degree $d$} if $\tau$ is open, proper and
 each fiber of $\tau$ contains exactly $d$ points counted with
 multiplicity. In this case, if $\Sigma_2$ is the set of critical
 values of $\tau$ and  $\Sigma_1:=\tau^{-1}(\Sigma_2)$, then
$\tau:X_1\setminus\Sigma_1\rightarrow X_2\setminus\Sigma_2$ is a
covering of degree $d$.

We recall the notion of analytic space which generalizes complex
manifolds and their analytic subsets.
An {\it analytic space of dimension $\leq p$}
is defined as a complex manifold but a chart is replaced by an
analytic subset of dimension $\leq p$ in an open set of a complex Euclidean
space.
As in the case of analytic subsets, one can decompose analytic
spaces into irreducible components and into regular and 
singular parts. The notions of dimension, of Zariski topology and of
holomorphic maps can be extended to analytic spaces.  The precise definition uses the local ring of holomorphic
functions, see \cite{Gunning, Narasimhan}. An analytic space is {\it normal} if the local ring of
holomorphic functions at every point is integrally closed. This is
equivalent to the fact that for $U$ open in $Z$ holomorphic functions
on $\reg(Z)\cap U$
which are bounded near $\sing(Z)\cap U$, are holomorphic on $U$. In
particular, normal analytic spaces are locally irreducible. A
holomorphic map $f:Z_1\rightarrow Z_2$ between complex spaces is a
continuous map which induces morphisms from local rings of holomorphic
functions on $Z_2$ to the ones on $Z_1$. 
The notions of ramified covering, of multiplicity and of open maps can be
extended to normal analytic spaces.
We have the following useful result where $\widetilde Z$ is called {\it normalization} of $Z$.

\begin{theorem} \label{th_normalization}
Let $Z$ be an analytic space. Then there is a unique, up to a
bi-holomorphic map, normal analytic space $\widetilde Z$ and a finite holomorphic map
$\pi:\widetilde Z\rightarrow Z$ such that
\begin{enumerate}
\item $\pi^{-1}(\reg(Z))$ is a dense Zariski open set of $\widetilde
      Z$ and $\pi$ defines a bi-holomorphic map between $\pi^{-1}(\reg(Z))$ and
      $\reg(Z)$;
\item If $\tau:Z'\rightarrow Z$ is a holomorphic map between 
      analytic spaces, then there is a unique holomorphic map
      $h:Z'\rightarrow\widetilde Z$ satisfying $\pi\circ h=\tau$.
\end{enumerate}
In particular, holomorphic self-maps of $Z$ can be lifted to
holomorphic self-maps of $\widetilde Z$.
\end{theorem}

\begin{example} \rm
Let $\pi:\C\rightarrow\C^2$ be the holomorphic map given by
$\pi(t)=(t^2,t^3)$. This map defines a normalization of the analytic
curve $\{z_1^3=z_2^2\}$ in $\C^2$ which is singular at 0. The
normalization of the analytic set $\{z_1=0\}\cup \{z_1^3=z_2^2\}$ is
the union of two disjoint complex lines. The normalization of a complex curve
(an analytic set of pure dimension 1) is always smooth.
\end{example}

The following desingularization theorem, due to Hironaka, is very
useful.

\begin{theorem} Let $Z$ be an analytic space. Then there is a smooth
  manifold $\widehat Z$, possibly reducible, and a holomorphic map $\pi:\widehat
  Z\rightarrow Z$ such that $\pi^{-1}(\reg (Z))$ is a dense Zariski open
  set of $\widehat Z$ and $\pi$ defines a bi-holomorphic map between
  $\pi^{-1}(\reg (Z))$ and $\reg(Z)$.
\end{theorem}

When $Z$ is an analytic subset of a manifold $X$, then one can obtain
a map $\pi:\widehat X\rightarrow X$ using a sequence of blow-ups
along the singularities of $Z$. The manifold $\widehat Z$ is the
strict transform of $Z$ by $\pi$. The difference with the normalization
of $Z$ is that we do not have the second property in Theorem
\ref{th_normalization} but $\widehat Z$ is smooth.

\bigskip\bigskip

\begin{exercise}
Let $X$ be a compact K{\"a}hler manifold of dimension $k$. 
Show that the Betti number $b_l$, i.e. the dimension of $H^l(X,\R)$, is
even if $l$ is odd and does not vanish if $l$ is even.
\end{exercise}

\begin{exercise}
Let $\Grass(l,k)$ denote the Grassmannian, i.e. the set of linear subspaces of dimension $l$ of
$\C^k$. Show that $\Grass(l,k)$ admits a natural structure of a
projective manifold.
\end{exercise}

\begin{exercise}
Let $X$ be a compact complex manifold of dimension $\geq 2$ and
$\pi:\widehat{X\times X}\rightarrow X\times X$
the blow-up of $X\times X$ along the diagonal $\Delta$. Let $\Pi_1,
\Pi_2$ denote the natural projections from $\widehat{X\times X}$ onto
the two factors $X$ of $X\times X$. Show that $\Pi_1,\Pi_2$ and their
restrictions to $\pi^{-1}(\Delta)$ are submersions. 
\end{exercise}

\begin{exercise} Let $E$ be a finite or countable union of proper analytic
  subsets of a connected manifold $X$. Show that $X\setminus E$ is connected and dense in
  $X$ for the usual topology.
\end{exercise}

\begin{exercise} \label{Exo_push_c0}
Let $\tau:X_1\rightarrow X_2$ be a ramified covering of degree
$n$. Let $\varphi$ be a function on $X_1$. Define
$$\tau_*(\varphi)(z):=\sum_{w\in\tau^{-1}(z)} \varphi(w),$$
where the points in $\tau^{-1}(z)$ are counted with multiplicity. If
$\varphi$ is upper semi-continuous or continuous, show
that $\tau_*(\varphi)$ is upper semi-continuous or continuous
respectively. 
Show that the result still holds
for a general open map $\tau$ between manifolds of the same dimension if
$\varphi$ has compact support in $X_1$. 
\end{exercise}


\section{Positive currents and p.s.h. functions} \label{section_positive}

In this paragraph, we  introduce positive forms, positive
currents and plurisubharmonic functions 
on complex manifolds. The concept of positivity and the notion of 
plurisubharmonic functions are due to Lelong and Oka. The theory has many
applications in complex algebraic geometry and in dynamics.

Let $X$ be a complex manifold of dimension $k$ and $\omega$ a
Hermitian $(1,1)$-form on $X$ which is positive definite at every point. 
Recall that a current $S$ on $X$, of degree $l$ and of dimension $2k-l$, 
is a continuous linear form on the
space $\Dc^{2k-l}(X)$ of smooth $(2k-l)$-forms with compact
support in $X$. Its value on a $(2k-l)$-form $\varphi\in\Dc^{2k-l}(X)$ is
denoted by $S(\varphi)$ or more frequently by $\langle S,\varphi\rangle$.
On a chart, $S$ corresponds to a continuous linear form acting on
the coefficients of $\varphi$. So, it can be represented as an $l$-form with
distribution coefficients. A sequence $(S_n)$ of $l$-currents
converges to an $l$-current $S$ if for every $\varphi\in\Dc^{2k-l}(X)$,
$\langle S_n,\varphi\rangle$ converge to $\langle S,\varphi\rangle$. 
The conjugate of $S$ is the $l$-current $\overline S$ defined by 
$$\langle\overline S,\varphi\rangle:=\overline{\langle
  S,\overline\varphi\rangle},$$
for $\varphi\in\Dc^{2k-l}(X)$. The current $S$ is real if and only if
$\overline S=S$.

The support of $S$ is the smallest closed subset $\supp(S)$ of $X$
such that $\langle S,\varphi\rangle=0$ when $\varphi$ is supported on
$X\setminus \supp(S)$. The current $S$ extends continuously to the
space of smooth forms $\varphi$ such that $\supp(\varphi)\cap\supp(S)$
is compact in $X$. If $X'$ is a complex manifold of dimension $k'$
with $2k'\geq 2k-l$, and if $\tau:X\rightarrow X'$ is a holomorphic map
which is proper on the support of $S$, we can define {\it the
  push-forward} $\tau_*(S)$ of $S$ by $\tau$. This is a current $\tau_*(S)$ of the
same dimension than $S$, i.e. of degree $2k'-2k+l$, which is supported
on $\tau(\supp(S))$, it satisfies
$$\langle \tau_*(S),\varphi\rangle:=\langle S,\tau^*(\varphi)\rangle$$
for $\varphi\in \Dc^{2k-l}(X')$. If $X'$ is a complex manifold of
dimension $k'\geq k$ and if $\tau:X'\rightarrow X$ is a submersion, we
can define {\it the pull-back} $\tau^*(S)$ of $S$ by $\tau$. This is an $l$-current
supported on $\tau^{-1}(\supp(S))$, it satisfies
$$\langle \tau^*(S),\varphi\rangle:=\langle S,\tau_*(\varphi)\rangle$$
for $\varphi\in \Dc^{2k'-l}(X')$. Indeed, since $\tau$ is a
submersion, the current $\tau_*(\varphi)$ is in fact a smooth form
with compact support in $X$; it is given by an integral of $\varphi$
on the fibers of $\tau$.

Any smooth differential $l$-form $\psi$ on $X$ defines a current: it
defines the continuous linear form $\varphi\mapsto
\int_X\psi\wedge\varphi$ on $\varphi\in\Dc^{2k-l}(X)$. So, 
currents extend the notion of differential forms. 
The operators $d,\partial,\dbar$ on differential forms extend to
currents. For example, we have that $dS$ is an $(l+1)$-current defined
by
$$\langle dS,\varphi\rangle:=(-1)^l\langle S, d\varphi\rangle$$
for $\varphi\in\Dc^{2k-l-1}(X)$. One easily check that when $S$ is
a smooth form, the above identity is a consequence of the Stokes' formula. 
We say that $S$ is of 
{\it bidegree} $(p,q)$ and {\it of bidimension} $(k-p,k-q)$ if it
vanishes on forms of bidegree $(r,s)\not=(k-p,k-q)$. The conjugate of
a $(p,q)$-current is of bidegree $(q,p)$. So, if such a current is real,
we have necessarily $p=q$. Note that the push-forward and the
pull-back by holomorphic maps commute with the above operators. They
preserve real currents; the push-forward preserves the bidimension and
the pull-back preserves the bidegree.

There are three notions of positivity which coincide for the bidegrees
$(0,0)$, $(1,1)$, $(k-1,k-1)$ and $(k,k)$. Here, we only use two of
them. They are dual to each other. 
A $(p,p)$-form $\varphi$ is {\it (strongly) positive} if at each
point, it is equal to a combination with positive coefficients of forms
of type
$$(\sqrt{-1}\alpha_1\wedge \overline\alpha_1)\wedge \ldots \wedge
(\sqrt{-1}\alpha_p\wedge \overline\alpha_p),$$
where $\alpha_i$ are $(1,0)$-forms. 
Any $(p,p)$-form can be written as a finite combination of positive $(p,p)$-forms.
For example, in local coordinates $z$, a $(1,1)$-form $\omega$ is written as
$$\omega=\sum_{i,j=1}^k \alpha_{ij} \sqrt{-1}dz_i\wedge d\overline z_j,$$
where $\alpha_{ij}$ are functions. This form is positive if and only
if the matrix $(\alpha_{ij})$ is positive semi-definite at every
point. In local coordinates $z$, the $(1,1)$-form $\ddc\|z\|^2$ is
positive. One can write $dz_1\wedge d\overline z_2$ as a combination
of $dz_1\wedge d\overline z_1$, $dz_2\wedge d\overline z_2$, 
$d(z_1\pm z_2)\wedge d\overline{(z_1\pm z_2)}$ and 
$d(z_1\pm \sqrt{-1}z_2)\wedge d\overline{(z_1\pm \sqrt{-1}z_2)}$.  Hence, positive
forms generate the space of $(p,p)$-forms.

A $(p,p)$-current $S$ is {\it
  weakly positive} if  for
every smooth positive $(k-p,k-p)$-form $\varphi$,
$S\wedge \varphi$ is a positive measure and is {\it positive}
if $S\wedge \varphi$ is a positive measure for
every smooth weakly positive $(k-p,k-p)$-form $\varphi$. 
Positivity implies weak positivity. These properties 
 are preserved
under pull-back by holomorphic submersions and push-forward by 
proper holomorphic maps.
Positive and weakly positive forms and
currents are real. One can consider positive and weakly positive
$(p,p)$-forms as sections of some bundles of salient convex closed
cones which are contained in
the real part of the vector bundle $\bigwedge^p\Omega^{1,0}\otimes
\bigwedge^p\Omega^{0,1}$. 

The wedge-product of a positive current with a positive form is positive.
The wedge-product of a weakly positive current with a positive form is
weakly positive. Wedge-products of weakly positive forms or currents are not
always weakly positive. For real $(p,p)$-currents or
forms $S$, $S'$, we will write $S\geq S'$ and $S'\leq S$ if $S-S'$ is
positive. A current $S$ is {\it negative} if $-S$ is positive. A
$(p,p)$-current or form $S$ is {\it strictly positive} if in local
coordinates $z$, there is a constant $\epsilon>0$ such that $S\geq
\epsilon (\ddc \|z\|^2)^p$. Equivalently, $S$ is strictly positive if
we have locally $S\geq \epsilon\omega^p$ with $\epsilon>0$.

\begin{example} \rm \label{example_current_lelong}
Let $Y$ be an analytic set of pure codimension $p$ of $X$. 
Using the local description of $Y$ near a singularity in Theorem
\ref{th_anal_set_local} and 
Wirtinger's theorem \ref{th_wirtinger}, one can prove 
that the $2(k-p)$-dimensional volume of $Y$
is locally finite in $X$. This allows to define the following
$(p,p)$-current $[Y]$ by
$$\langle [Y],\varphi\rangle:=\int_{\reg(Y)}\varphi$$
for $\varphi$ in $\Dc^{k-p,k-p}(X)$, the space of smooth
$(k-p,k-p)$-forms with compact support in $X$. Lelong proved that this current is
positive and closed \cite{Demailly3, Lelong}. 
\end{example}

If $S$ is a (weakly) positive $(p,p)$-current, it is of order
0, i.e. it extends continuously to the space of continuous forms with
compact support in $X$. In other words, on a chart of $X$, the current
$S$ corresponds to a differential form with measure coefficients. 
We define the mass of $X$ on a Borel set $K$ by 
$$\|S\|_K:=\int_K S\wedge\omega^{k-p}.$$
When $K$ is relatively compact in $X$, we obtain an equivalent
norm if we change the Hermitian metric on $X$. This is a consequence
of the property we mentioned above, which says that $S$ takes values in
salient convex closed cones. Note that the previous mass-norm is
just defined by an integral, which is easier to compute or
to estimate than the usual mass for currents on real
manifolds. 

Positivity implies an important compactness property. As for positive
measures, any family of
positive $(p,p)$-currents with locally uniformly bounded mass, is
relatively compact in the cone of positive $(p,p)$-currents. 
For the current $[Y]$ in Example \ref{example_current_lelong}, by Wirtinger's theorem, the
mass on $K$ is equal to $(k-p)!$ times the volume of $Y\cap K$ with
respect to the considered Hermitian metric. If $S$ is a
negative $(p,p)$-current, its mass is defined by
$$\|S\|_K:=-\int_K S\wedge\omega^{k-p}.$$
The following result is the complex version of the classical support theorem in the
real setting, \cite{Bassanelli, HarveyPolking, Federer}.

\begin{proposition} \label{prop_support_theorem}
Let $S$ be a $(p,p)$-current supported on a smooth complex
submanifold $Y$ of $X$. Let $\tau:Y\rightarrow X$ denote the inclusion
map. Assume that $S$ is $\C$-normal, i.e. $S$ and $\ddc S$ are of
order $0$. Then, $S$ is a current on $Y$. More precisely, there is a
$\C$-normal $(p,p)$-current $S'$ on $Y$ such that $S=\tau_*(S')$. If
$S$ is positive closed and $Y$ is of dimension $k-p$, then $S$ is equal to a
combination with positive coefficients of currents of integration on components of $Y$.
\end{proposition}

The last property holds also when $Y$ is a singular analytic set. 
Proposition \ref{prop_support_theorem} 
applies to positive closed $(p,p)$-currents which play an important
role in complex geometry and dynamics. These currents generalize 
analytic sets of dimension $k-p$, as we have seen in Example
\ref{example_current_lelong}.
They have no mass on Borel sets of $2(k-p)$-dimensional Hausdorff
measure 0. The proposition is used in order to develop a calculus on potentials
of closed currents.

We introduce now the notion of Lelong number for such currents which generalizes the
notion of multiplicity for analytic sets.
Let $S$ be a positive closed $(p,p)$-current on $X$. Consider local
coordinates $z$ on a chart $U$ of $X$ and the local K\"ahler form $\ddc\|z\|^2$. 
Let $B_a(r)$ denote the
ball of center $a$ and of radius $r$ contained in $U$. Then, $S\wedge
(\ddc\|z\|^2)^{k-p}$ is a positive measure on $U$. 
Define for $a\in U$
$$\nu(S,a,r):=\frac{{\|S\wedge (\ddc\|z\|^2)^{k-p}\|_{B_a(r)}}}{{\pi^{k-p}r^{2(k-p)}}}.$$ 
Note that $\pi^{k-p}r^{2(k-p)}$ is $(k-p)!$ times the volume of a ball in $\C^{k-p}$
of radius $r$, i.e. the mass of the current associated to this ball. When $r$ decreases to $0$, 
$\nu(S,a,r)$ is decreasing and the {\it
Lelong number} of $S$ at $a$ is the limit
$$\nu(S,a):=\lim_{r\rightarrow 0} \nu(S,a,r).$$
It does not depend on
the coordinates.
So, we can define the Lelong number for currents on
any manifold. Note that $\nu(S,a)$ is also the mass of the measure $S\wedge
(\ddc\log\|z-a\|)^{k-p}$ at $a$. We will discuss the wedge-product (intersection) of
currents in the next paragraph.

If $S$ is the current of integration on an analytic set
$Y$, by Thie's theorem, $\nu(S,a)$ is equal to
the multiplicity of $Y$ at $a$ which is an integer. 
This implies the following Lelong's inequality:
{\it the Euclidean $2(k-p)$-dimensional volume of $Y$ in a ball
$B_a(r)$ centered at a point $a\in Y$, is at least equal to
${1\over (k-p)!}\pi^{k-p}r^{2(k-p)}$, the volume in $B_a(r)$ of a $(k-p)$-dimensional
linear space passing through $a$}.

Positive closed currents generalize analytic sets but they are much
more flexible. A remarkable fact is that the use of positive closed
currents allows to construct analytic sets. The
following theorem of Siu \cite{Siu} is a beautiful application
of the complex $L^2$ method.

\begin{theorem} \label{th_siu}
Let $S$ be a positive closed $(p,p)$-current on $X$. Then, for $c>0$, 
the level set $\{\nu(S,a)\geq c\}$ of the Lelong number is an analytic
set of $X$, of dimension $\leq k-p$. Moreover, there is a unique
decomposition  $S=S_1+S_2$ where $S_1$ is a locally finite
combination, with positive coefficients, 
of currents of integration on analytic sets of codimension
$p$ and $S_2$ is a positive closed $(p,p)$-current such that
$\{\nu(S_2,z)>0\}$ is a finite or countable union of analytic sets of
dimension $\leq k-p-1$.
\end{theorem}

Calculus on currents is often delicate. However, the theory is
well developped for positive closed $(1,1)$-currents thanks to
the use of plurisubharmonic functions. Note that positive
closed $(1,1)$-currents correspond to hypersurfaces (analytic sets of
pure codimension 1) in complex geometry and working with
$(p,p)$-currents, as with higher codimension analytic sets, is more
difficult.

An upper semi-continuous function $u:X\rightarrow \R\cup\{-\infty\}$,
not identically $-\infty$ on any component of $X$, is {\it
  plurisubharmonic} (p.s.h. for short) if it is subharmonic or
identically $-\infty$ on
any holomorphic  disc in $X$.
Recall that a {\it holomorphic disc} in $X$ is a holomorphic map
$\tau:\Delta\rightarrow X$ where $\Delta$ is the unit disc in
$\C$. One often identifies this holomorphic disc with its image
$\tau(\Delta)$. If $u$ is p.s.h., then $u\circ\tau$ is
subharmonic or identically $-\infty$ on $\Delta$. As for subharmonic
functions, we have the submean inequality: {\it in local
  coordinates, the value at $a$ of a p.s.h. function is smaller or
  equal to the average of the function on a sphere centered at
  $a$}. Indeed, this average increases with the radius of the
sphere. The submean inequality implies the maximum principle: {\it if a
  p.s.h. function on a connected manifold $X$ has a local maximum, it is constant}. 
The semi-continuity implies that p.s.h. functions are locally bounded from above.
A function $v$ is {\it
  pluriharmonic} if $v$ and $-v$ are p.s.h. Pluriharmonic functions
are locally real parts of holomorphic functions, in particular, they
are real analytic. 
The following proposition is of constant use.

\begin{proposition}
A function $u:X\rightarrow \R\cup\{-\infty\}$ is p.s.h. if and only if
the following conditions are satisfied
\begin{enumerate}
\item $u$ is strongly upper semi-continuous, that is, for any subset
  $A$ of full Lebesgue measure in $X$ and for any point $a$ in $X$, we
  have
$u(a)=\limsup u(z)$ when $z\rightarrow a$ and $z\in A$.
\item $u$ is locally integrable with respect to the Lebesgue measure
  on $X$ and $\ddc u$ is a positive closed $(1,1)$-current.
\end{enumerate}
\end{proposition}

Conversely, any positive closed $(1,1)$-current can be locally written
as $\ddc u$ where $u$ is a (local) p.s.h. function. This function is
called {\it a local potential} of the current. Two local potentials
differ by a pluriharmonic function. So, there is
almost a correspondence between positive closed $(1,1)$-currents and
p.s.h. functions. We say that $u$ is {\it strictly p.s.h.} if $\ddc
u$ is strictly positive. The p.s.h. functions are defined at
every point; this is a crucial property in pluripotential
theory. Other important properties of this class of functions are 
some strong compactness properties that we state below.

If $S$ is a positive closed
$(p,p)$-current, one can write locally $S=\ddc U$ with $U$ a
$(p-1,p-1)$-current. We can choose the potential $U$ negative with good
estimates on the mass but the difference of two potentials may be very
singular. The use of potentials $U$ is much more delicate than in the
 bidegree $(1,1)$ case. We state here a useful local estimate,
see e.g. \cite{DinhNguyenSibony2}. 

\begin{proposition} \label{prop_ddbar_convex_dom}
Let $V$ be convex open domain in $\C^k$  and $W$ an open set with
  $W\Subset V$.
Let $S$ be a positive closed $(p,p)$-current on $V$. 
Then there is a  negative $L^1$ form $U$ of bidegree $(p-1,p-1)$ on $W$ such
  that $\ddc U=S$ and $\|U\|_{L^1(W)}\leq
  c\|S\|_V$ where $c>0$ is a constant independent of
  $S$. Moreover, $U$ depends continuously on $S$, where the continuity
  is with respect to the weak topology on $S$ and the $L^1(W)$
  topology on $U$. 
\end{proposition}

Note that when $p=1$, $U$ is equal almost everywhere to a p.s.h.
function $u$ such that $\ddc u=S$. 

\begin{example}\rm
Let $f$ be a holomorphic function on $X$ not identically 0 on any
component of $X$. Then,
$\log|f|$ is a p.s.h. function and we have $\ddc \log|f|=
\sum n_i [Z_i]$ where $Z_i$ are irreducible components of the
hypersurface $\{f=0\}$ and
$n_i$ their multiplicities. The last equation is called {\it
  Poincar{\'e}-Lelong equation}. Locally, the ideal of holomorphic
functions vanishing on $Z_i$ is generated by a holomorphic function
$g_i$ and $f$ is equal to the product of $\prod g_i^{n_i}$ with a
non-vanishing holomorphic function. In some sense, $\log|f|$ is one of
the most singular p.s.h. functions. If $X$ is a ball, 
the convex set generated by such functions is dense in the cone of p.s.h. functions
\cite{Hormander1, Gunning} for the $L^1_\loc$ topology. 
If $f_1,\ldots, f_n$ are
holomorphic on $X$, not identically 0 on a component of $X$, then
$\log(|f_1|^2+\cdots+|f_n|^2)$ is also a p.s.h. function.  
\end{example}

The following proposition is useful in constructing 
p.s.h. functions. 

\begin{proposition} \label{prop_psh_composition}
Let $\chi$ be a function defined on $(\R\cup\{-\infty\})^n$ with
values in $\R\cup\{-\infty\}$, not identically $-\infty$, which is
convex in all variables and increasing in each variable. Let
$u_1,\ldots,u_n$ be p.s.h. functions on $X$. Then
$\chi(u_1,\ldots,u_n)$ is p.s.h. In particular, the function
$\max(u_1,\ldots,u_n)$ is p.s.h.
\end{proposition}

We call {\it complete pluripolar set} the pole set $\{u=-\infty\}$ of
a p.s.h. function and {\it pluripolar set} a subset of a complete
pluripolar one. Pluripolar sets are of Hausdorff dimension $\leq
2k-2$, in particular, they have zero Lebesgue measure. Finite and
countable unions of (locally) pluripolar sets are (locally) pluripolar. In
particular, finite and countable unions of analytic subsets are
locally pluripolar. 

\begin{proposition}
Let $E$ be a closed pluripolar set in $X$ and $u$ a p.s.h. function on
$X\setminus E$, locally bounded above near $E$. Then the extension
of $u$ to $X$ given by
$$u(z):=\limsup_{w\rightarrow z\atop w\in X\setminus E} u(w) \quad
\mbox{for}\quad z\in E,$$ 
is a p.s.h. function.  
\end{proposition}

The following result  describes compactness
properties of p.s.h. functions, see \cite{Hormander1}. 

\begin{proposition} \label{prop_hartogs_else}
Let $(u_n)$ be a sequence of p.s.h. functions on
  $X$, locally bounded from above. Then either it converges locally
  uniformly to $-\infty$ on a component of $X$ or
there is a subsequence
  $(u_{n_i})$  which   converges in $L^p_\loc(X)$ to a p.s.h. function $u$ for every $p$ with
  $1\leq p<\infty$. In the second case, we have $\limsup u_{n_i}\leq
  u$ with equality outside a pluripolar set. Moreover, if $K$ is a compact
  subset of $X$ and if $h$ is a continuous function on $K$ such that
  $u<h$ on $K$, then $u_{n_i}<h$ on $K$ for $i$ large enough.
\end{proposition}

The last assertion is the classical Hartogs' lemma. It suggests the
following notion of convergence introduced in \cite{DinhSibony10}. Let
$(u_n)$ be a sequence of p.s.h. functions 
converging to a p.s.h. function $u$ in $L^1_\loc(X)$. We say that
the sequence $(u_n)$
{\it converges in the
 Hartogs' sense} or {\it is H-convergent} if for any compact subset
$K$ of $X$ there are constants $c_n$ converging to 0 such that
$u_n+c_n\geq u$ on $K$. In this case, Hartogs' lemma implies
that $u_n$ converge pointwise to $u$. If $(u_n)$ decreases to a
function $u$, not identically $-\infty$, then $u$ is p.s.h. and
$(u_n)$ converges in the Hartogs' sense. 
The following result is useful in the calculus with p.s.h. functions.

\begin{proposition}
Let  $u$ be a p.s.h. function on
an open subset $D$ of $\C^k$. Let $D'\Subset D$ be an open set. Then,
there is a sequence of smooth p.s.h. functions $u_n$ on $D'$ which
decreases to $u$.
\end{proposition}

The functions $u_n$ can be obtained 
as the standard convolution of $u$ with some radial function
$\rho_n$ on $\C^k$. The submean inequality for $u$ allows to choose $\rho_n$ so that $u_n$
decrease to $u$.

The following result, see \cite{Hormander2}, may be considered as the
strongest compactness property for p.s.h. functions. The proof can be
reduced to the one dimensional case by slicing.

\begin{theorem} \label{th_hormander}
Let $\Fc$ be a family of p.s.h. functions on $X$ which
  is bounded in $L^1_\loc(X)$. Let $K$ be a compact subset of $X$. 
Then there are constants $\alpha>0$ and $A>0$ such that
$$\|e^{-\alpha u}\|_{L^1(K)}\leq A$$
for every function $u$ in $\Fc$.
\end{theorem}

P.s.h. functions are in general unbounded. However, the last
result shows that such functions are nearly bounded. The above family
$\Fc$ is uniformly  bounded from above on $K$. So, we also have the
estimate 
$$\|e^{\alpha|u|}\|_{L^1(K)}\leq A$$
for $u$ in $\Fc$ and for some (other) constants $\alpha,A$.
More precise estimates can be obtained in terms of the maximal Lelong
number of $\ddc u$ in a neighbourhood of $K$. 

Define the {\it Lelong number} $\nu(u,a)$
of $u$ at $a$ as the Lelong number of $\ddc u$ at $a$. The following
result describes the relation with the singularity of p.s.h. functions
near a pole. We fix here a local coordinate system for $X$.

\begin{proposition}
The Lelong number $\nu(u,a)$ is the supremum of the number $\nu$ such
that the inequality
$u(z)\leq \nu\log \|z-a\|$ holds in a neighbourhood of $a$.
\end{proposition}

If $S$ is a positive closed $(p,p)$-current, the
Lelong number $\nu(S,a)$ can be computed as the mass at $a$ of the
measure $S\wedge (\ddc\log\|z-a\|)^{k-p}$. 
This property allows to prove the following result, due to
Demailly \cite{Demailly3},
which is useful in dynamics.

\begin{proposition}
Let $\tau:(\C^k,0)\rightarrow (\C^k,0)$ be a germ of an open holomorphic
map with $\tau(0)=0$. Let $d$ denote the multiplicity of $\tau$ at
$0$. Let $S$ be a positive closed $(p,p)$-current on a
neighbourhood of $0$. Then, the Lelong number of $\tau_*(S)$ at $0$
satisfies the inequalities
$$\nu(S,0)\leq\nu(\tau_*(S),0)\leq d^{k-p}\nu(S,0).$$
In particular, we have $\nu(\tau_*(S),0)=0$ if and only if $\nu(S,0)=0$. 
\end{proposition}

Assume now that $X$ is a compact K\"ahler manifold and $\omega$ is a
K\"ahler form on $X$. If $S$ is a $\ddc$-closed $(p,p)$-current, we
can, using the $\ddc$-lemma,  define a linear form on $H^{k-p,k-p}(X,\C)$ by
$[\alpha]\mapsto \langle S,\alpha\rangle$. Therefore, the Poincar\'e
duality implies that $S$ is canonically associated to a class $[S]$ in
$H^{p,p}(X,\C)$. If $S$ is real then $[S]$ is in $H^{p,p}(X,\R)$. If
$S$ is positive, its mass $\langle S,\omega^{k-p}\rangle$ depends only
on the class $[S]$. So, the mass of positive $\ddc$-closed currents can be
computed cohomologically. In $\P^k$, the mass of $\omega_\FS^p$ is 1
since $\omega_\FS^k$ is a probability measure. If $H$ is a subspace
of codimension $p$ of $\P^k$, then the current associated to $H$ is of
mass 1 and it belongs to the class $[\omega_\FS^p]$. If $Y$ is an
analytic set of pure codimension $p$ of $\P^k$, {\it the degree}
$\deg(Y)$ of $Y$ is by definition the number of points in its
intersection with a generic projective space of dimension $p$. One can
check that the cohomology class of $Y$ is $\deg(Y)[\omega_\FS^p]$. The
volume of $Y$, obtained using Wirtinger's theorem \ref{th_wirtinger},
is equal to ${1\over p!} \deg(Y)$.

\bigskip\bigskip

\begin{exercise}
With the notation of Exercise \ref{Exo_push_c0}, show that $\tau_*(\varphi)$ is
p.s.h. if $\varphi$ is p.s.h.
\end{exercise}

\begin{exercise}
Using that $\nu(S,a,r)$ is decreasing, show that 
if $(S_n)$ is a sequence of positive closed $(p,p)$-currents on $X$
converging to a current $S$ and
$(a_n)$ is a sequence in $X$ converging to $a$, then $\limsup
\nu(S_n,a_n)\leq \nu(S,a)$.
\end{exercise}

\begin{exercise}
Let $S$ and $S'$ be  positive closed $(1,1)$-currents such that
$S'\leq S$. Assume that the local potentials of $S$ are bounded or
continuous. Show that the local potentials of $S'$ are also bounded or
continuous.
\end{exercise}

\begin{exercise}
Let $\Fc$ be an $L^1_\loc$ bounded family of p.s.h. functions on $X$. Let $K$ be a
compact subset of $X$. Show that $\Fc$ is locally bounded from above and
that there is $c>0$ such that $\|\ddc u\|_K\leq c$ for every $u\in\Fc$. 
Prove that there is a constant $\nu>0$ such that
$\nu(u,a)\leq\nu$ for $u\in\Fc$ and $a\in K$.
\end{exercise} 

\begin{exercise} Let $Y_i$, $1\leq i\leq m$, be analytic sets of pure
  codimension $p_i$ in $\P^k$. Assume $p_1+\cdots+p_m\leq k$. Show
  that the intersection of the $Y_i$'s is a non-empty analytic set of
  dimension $\geq k-p_1-\cdots-p_m$.
\end{exercise}


\section{Intersection, pull-back and slicing}
\label{section_intersection}

We have seen that positive closed currents generalize differential
forms and analytic sets. However, it is not always possible to extend 
the calculus on forms or on analytic sets to currents. We will give here
some results which show how positive closed currents are flexible and
how they are rigid.

The theory of intersection
is much more developed in bidegree $(1,1)$ thanks to the use of their
potentials which are p.s.h. functions. The
case of continuous potentials was considered by
Chern-Levine-Nirenberg \cite{ChernLevineNirenberg}. Bedford-Taylor
\cite{BedfordTaylor} developed a nice theory when the potentials are
locally bounded. The case of unbounded
potentials was considered  by Demailly
\cite{Demailly1} and Forn\ae
ss-Sibony \cite{FornaessSibony2, Sibony2}. We have the following general definition.

Let $S$ be a positive closed $(p,p)$-current on $X$ with $p\leq k-1$. If $\omega$ is
a fixed Hermitian form on $X$ as above, then $S\wedge\omega^{k-p}$
is a positive measure which is called {\it the trace measure} of $S$. In
local coordinates, the coefficients of $S$ are measures, bounded by a
constant times the trace measure. Now, if $u$ is a p.s.h function on
$X$, locally integrable with respect to the trace measure of $S$, then
$uS$ is a current on $X$ and we can define
$$\ddc u\wedge S:=\ddc(uS).$$
Since $u$ can be locally approximated by decreasing sequences of smooth
p.s.h. functions, it is easy to check that the previous intersection
is a positive closed $(p+1,p+1)$-current with support contained in $\supp(S)$. When $u$ is pluriharmonic,
$\ddc u\wedge S$ vanishes identically. So, the intersection depends
only on $\ddc u$ and on $S$. If $R$ is a positive closed
$(1,1)$-current on $X$, one defines $R\wedge S$ as above using local
potentials of $R$. 
In general, $\ddc u\wedge S$ does not depend continuously on $u$ and $S$.
The following proposition is a consequence
of Hartogs' lemma.

\begin{proposition} Let $u^{(n)}$ be p.s.h. functions on $X$ which
  converge in the Hartogs' sense to a p.s.h. function $u$. 
If $u$ is locally integrable with respect to the trace measure of $S$, then
$\ddc u^{(n)}\wedge S$ are well-defined and converge to $\ddc u \wedge S$.
If $u$ is continuous and $S_n$ are positive closed $(1,1)$-currents
converging to $S$, then $\ddc u^{(n)}\wedge S_n$ converge to $\ddc u\wedge
S$. 
\end{proposition}

If $u_1,\ldots, u_q$, with $q\leq k-p$,  are p.s.h. functions, we can define by induction
the wedge-product
$$\ddc u_1\wedge \ldots\wedge \ddc u_q\wedge S$$
when some integrability conditions are satisfied, for example when the
$u_i$ are locally bounded. In particular, if
$u_j^{(n)}$, $1\leq j\leq q$,  are continuous p.s.h. functions
converging locally uniformly to continuous p.s.h. functions $u_j$ and
if $S_n$ are positive closed
converging to $S$, then  
$$\ddc u_1^{(n)}\wedge \ldots\wedge \ddc u_q^{(n)}\wedge S_n
\rightarrow \ddc u_1\wedge \ldots\wedge \ddc u_q\wedge S$$

The following version of the Chern-Levine-Nirenberg inequality is a
very useful result \cite{ChernLevineNirenberg, Demailly3}.

\begin{theorem}
Let $S$ be a positive closed $(p,p)$-current on $X$.
Let $u_1,\ldots, u_q$, $q\leq k-p$, be locally bounded p.s.h. functions on
$X$ and $K$ a compact subset of
$X$. Then there is a constant $c>0$ depending only on $K$ and $X$ such
that if $v$ is p.s.h. on $X$ then
$$\|v\ddc u_1\wedge \ldots \wedge \ddc u_q\wedge S\|_K\leq c
\|v\|_{L^1(\sigma_S)}\|u_1\|_{L^\infty(X)}\ldots \|u_q\|_{L^\infty(X)},$$
where $\sigma_S$ denotes the trace measure of $S$.
\end{theorem}

This inequality implies that p.s.h. functions are integrable
with respect to the current $\ddc u_1\wedge \ldots\wedge \ddc u_q$. We
deduce the following corollary.

\begin{corollary} \label{cor_cln}
Let $u_1,\ldots, u_p$, $p\leq k$, be locally bounded p.s.h. functions on
$X$. Then, the current $\ddc u_1\wedge \ldots\wedge \ddc u_p$
has no mass on locally pluripolar sets, in particular on
proper analytic sets of $X$. 
\end{corollary}

We give now two other regularity
properties of the wedge-product of currents with H{\"o}lder continuous local
potentials. 

\begin{proposition} \label{prop_hol_dim_wedge}
Let $S$ be a positive closed $(p,p)$-current on $X$ and $q$ a positive
integer such that $q\leq k-p$.
Let $u_i$ be H{\"o}lder continuous p.s.h. 
functions of H{\"o}lder exponents
$\alpha_i$ with $0<\alpha_i\leq 1$ and $1\leq i\leq q$. Then, the current
$\ddc u_1\wedge\ldots\wedge\ddc u_q\wedge S$ has no mass on Borel sets
with Hausdorff dimension less than or equal to $2(k-p-q)+\alpha_1+\cdots+\alpha_q$.
\end{proposition}

The proof of this result is given in \cite{Sibony}. It is based on a mass estimate 
on a ball in term of the radius which is a consequence of the
Chern-Levine-Nirenberg inequality.

We say that a positive measure $\nu$ in $X$ is {\it locally moderate} if for
any compact subset $K$ of $X$ and any compact family $\Fc$ of
p.s.h. functions in a neighbourhood of $K$, there are positive
constants $\alpha$ and $c$ such that 
$$\int_K e^{-\alpha u}d\nu \leq c$$
for $u$ in $\Fc$. This notion was introduced in
\cite{DinhSibony1}. 
We say that a positive current is {\it locally moderate} if its trace
measure is locally moderate.
The following result was obtained in \cite{DinhNguyenSibony3}.

\begin{theorem} \label{th_moderate}
Let $S$ be a positive closed $(p,p)$-current on $X$ and $u$ a 
p.s.h. function on $X$. Assume that $S$ is locally moderate
and $u$ is H{\"o}lder continuous. Then the current $\ddc u\wedge S$ is
locally moderate. In particular, wedge-products of positive closed
$(1,1)$-currents with H{\"o}lder continuous local potentials are locally
moderate.
\end{theorem}

Theorem \ref{th_hormander} implies that a measure defined by a smooth
form is locally moderate. Theorem \ref{th_moderate} implies, by induction, that $\ddc u_1\wedge
\ldots\wedge \ddc u_p$ is locally moderate when the p.s.h. functions
$u_j$ are H{\"o}lder continuous. So, using p.s.h. functions as test
functions, the previous currents satisfy similar estimates as smooth
forms do. One may also consider that Theorem \ref{th_moderate}
strengthens \ref{th_hormander} and gives a strong compactness property for
p.s.h. functions. The estimate has many consequences in complex dynamics.

The proof of Theorem \ref{th_moderate} is based on a mass estimate of
$\ddc u\wedge S$ on the sub-level set $\{v<-M\}$ of a p.s.h function
$v$. Some estimates are easily obtained for $u$ continuous using the
Chern-Levine-Nirenberg inequality or for $u$ of class $\Cc^2$. The
case of H{\"o}lder continuous function uses arguments close to the
interpolation between the Banach spaces $\Cc^0$ and $\Cc^2$. However, the non-linearity of the
estimate and the positivity of currents make the problem more subtle.

We discuss now the pull-back of currents by holomorphic maps which are
not submersions. The problem can be considered as a particular case of the
general intersection theory, but we will not discuss this point here. The following
result was obtained in \cite{DinhSibony8}.

\begin{theorem} \label{th_pullback_local}
Let $\tau:X'\rightarrow X$ be an open holomorphic map between complex
manifolds of the same dimension. 
Then the pull-back operator $\tau^*$ on smooth positive closed
$(p,p)$-form can be extended in a canonical way to a continuous
operator on positive closed $(p,p)$-currents $S$ on $X$. 
If $S$ has no mass on a Borel set $K\subset X$, 
then $\tau^*(S)$ has no mass on $\tau^{-1}(K)$. The result also holds
for negative currents $S$ such that $\ddc S$ is positive. 
\end{theorem} 

By canonical way, we mean that the extension is functorial. More
precisely, one can locally approximate $S$ by a sequence of smooth
positive closed forms. The pull-back of these forms converge to some
positive closed $(p,p)$-current which does not depend on the
chosen sequence of forms. This limit defines the pull-back current $\tau^*(S)$. 
The result still holds  when $X'$ is singular. 
In the case of bidegree
$(1,1)$, we have the following result due to  M{\'e}o \cite{Meo}.

\begin{proposition} 
Let $\tau:X'\rightarrow X$ be a holomorphic map between complex
manifolds. Assume that $\tau$ is dominant, that is, the image of $\tau$
contains an open subset of $X$. 
Then the pull-back operator $\tau^*$ on smooth positive closed
$(1,1)$-form can be extended in a canonical way to a continuous
operator on positive closed $(1,1)$-currents $S$ on $X$. 
\end{proposition} 

Indeed, locally we can write $S=\ddc u$ with $u$ p.s.h. The current
$\tau^*(S)$ is then defined by $\tau^*(S):=\ddc (u\circ \tau)$. One can
check that the definition does not depend on the choice of $u$.

The remaining part of this paragraph deals with the slicing of currents.
We only consider a situation used in this course.
Let $\pi: X\rightarrow V$ be a dominant holomorphic map from $X$ to a
manifold $V$ of dimension $l$ and $S$ a current on $X$.
Slicing theory allows to
define the slice $\langle S,\pi,\theta\rangle$ of some currents $S$ on $X$
by the fiber $\pi^{-1}(\theta)$. Slicing theory generalizes the restriction
of forms to fibers. One can also consider it as a generalization of
Sard's and Fubini's theorems for currents or as a special case of
intersection theory: the slice $\langle S,\pi,\theta\rangle$ can
be seen as the wedge-product of $S$ with the current of integration on
$\pi^{-1}(\theta)$. We can consider the slicing of $\C$-flat
currents, in particular, of $(p,p)$-currents such that $S$ and $\ddc S$
are of order 0. The operation preserves  positivity and 
commutes with $\partial$, $\dbar$. 
If $\varphi$ is a smooth form on $X$ then $\langle
S\wedge\varphi,\pi,\theta\rangle = \langle S,\pi,\theta\rangle\wedge
\varphi$. 
Here, we only consider positive closed $(k-l,k-l)$-currents $S$. In this case, 
the slices  $\langle S,\pi,\theta\rangle$ are  positive
measures on $X$ with support in $\pi^{-1}(\theta)$.

Let $y$ denote the coordinates in a chart of $V$ and $\lambda_V:= (\ddc\|y\|^2)^l$ the
Euclidean volume form associated to $y$.
Let $\psi(y)$ be a positive smooth function with compact support
such that $\int\psi\lambda_V=1$. Define
$\psi_\epsilon(y):=\epsilon^{-2l}\psi(\epsilon^{-1} y)$ and 
$\psi_{\theta,\epsilon}(y):=\psi_\epsilon(y-\theta)$. The measures $\psi_{\theta,\epsilon}\lambda_V$
approximate the Dirac mass at $\theta$. For every smooth test function $\Phi$ on $X$,
we have
\begin{equation*} \label{eq_slice}
\langle S,\pi,\theta\rangle (\Phi)=\lim_{\epsilon\rightarrow 0}
\langle S\wedge \pi^*(\psi_{\theta,\epsilon}\lambda_V),\Phi\rangle
\end{equation*}
when $\langle S,\pi,\theta\rangle$ exists.
This property holds for all choice of $\psi$.
Conversely, when the previous limit exists and is independent of $\psi$, 
it defines the measure $\langle S,\pi,\theta\rangle$ and we say 
that $\langle S,\pi,\theta\rangle$ {\it is well-defined}.
The slice $\langle S,\pi,\theta\rangle$ is well-defined for $\theta$
out of a set of Lebesgue measure zero in $V$ and the 
following formula holds for smooth forms $\Omega$ 
of maximal degree with compact support in $V$:
\begin{equation*} \label{eq_slice_bis}
\int_{\theta\in V}\langle S,\pi,\theta\rangle (\Phi)\Omega(\theta) 
 =   \langle S\wedge \pi^*(\Omega),\Phi\rangle.
\end{equation*}
We recall the following result which was obtained in \cite{DinhSibony7}.

\begin{theorem} \label{th_slicing}
Let $V$ be a complex manifold of dimension $l$ and let $\pi$
denote the canonical projection from $\C^k\times V$ onto $V$.
Let $S$ be a positive closed current of bidimension $(l,l)$ on
$\C^k\times V$, supported on $K\times V$ for a given compact subset
$K$ of $\C^k$.
Then the slice $\langle S,\pi,\theta\rangle$ is well-defined for every
$\theta$ in $V$ and is a positive measure whose mass is independent of
$\theta$. 
Moreover, if
$\Phi$ is a p.s.h. function in a neighbourhood of $\supp(S)$,
then the function 
$\theta\mapsto \langle S,\pi,\theta\rangle(\Phi)$ is 
p.s.h. 
\end{theorem}

The mass of $\langle S,\pi,\theta\rangle$ is called the {\it slice
  mass} of $S$. The set of currents $S$ as above with bounded slice
mass is compact for the weak topology on currents. In particular, their
masses are locally uniformly bounded on $\C^k\times V$. 
In general, the slice $\langle
S,\pi,\theta\rangle$ does not depend continuously on $S$ nor on $\theta$. The last
property in Theorem \ref{th_slicing} shows that the dependence on $\theta$
satisfies a semi-continuity property. More generally, we have that
$(\theta,S)\mapsto  \langle S,\pi,\theta\rangle(\Phi)$ is upper
semi-continuous for $\Phi$ p.s.h. We deduce easily from the above definition that 
the slice mass of $S$ depends continuously on $S$.

\bigskip\bigskip

\begin{exercise}
Let $X,X'$ be complex manifolds. Let $\nu$ be a positive measure with
compact support on $X$ such that p.s.h. functions on $X$ are $\nu$-integrable. 
If $u$ is a p.s.h. function on $X\times X'$, show that $x'\mapsto \int
u(x,x')d\nu(x)$ is a p.s.h function on $X'$. Show that if $\nu,\nu'$ are
positive measures on $X, X'$ which are locally moderate, then
$\nu\otimes\nu'$ is a locally moderate measure on $X\times X'$. 
\end{exercise}

\begin{exercise} \label{exercise_lelong_slice}
Let $S$ be a positive closed $(1,1)$-current on the unit ball of
$\C^k$. Let $\pi:\widehat{\C^k}\rightarrow\C^k$ be the blow-up of
$\C^k$ at $0$ and $E$ the exceptional set. 
Show $\pi^*(S)$ is equal to
$\nu[E]+S'$, where $\nu$ is the Lelong number of $S$ at $0$ and $S'$ is
a current without mass on $E$.
\end{exercise}


\section{Currents on projective spaces} \label{section_current_pk}

In this paragraph, we will introduce quasi-potentials of currents, the
spaces of d.s.h. functions, DSH currents and the complex Sobolev space
which are used  as observables in dynamics.
We also introduce  
PB, PC currents and the notion of super-potentials which are crucial
in the calculus with currents in higher bidegree.

Recall that the Fubini-Study form $\omega_\FS$ on $\P^k$ satisfies 
$\int_{\P^k}\omega_\FS^k=1$. 
If $S$ is a positive closed $(p,p)$-current, the mass of $S$ is given by
by $\|S\|:=\langle S,\omega_\FS^{k-p}\rangle$. Since $H^{p,p}(\P^k,\R)$ is generated by
$\omega_\FS^p$, such a current $S$ is cohomologous to
$c\omega_\FS^p$ where $c$ is the mass of $S$. So, 
$S-c\omega_\FS^p$ is exact and the $\ddc$-lemma, which also holds for
currents, implies that there
exists a $(p-1,p-1)$-current $U$, such that $S=c\omega_\FS^p+\ddc U$. 
We call $U$ {\it a quasi-potential} of $S$.
We have in fact the following more precise result \cite{DinhSibony10}.

\begin{theorem}
Let $S$ be a positive closed $(p,p)$-current of mass $1$ in
$\P^k$. Then, there is a negative form $U$ such that $\ddc
U=S-\omega_\FS^p$. For $r,s$ with $1\leq r<k/(k-1)$ and $1\leq
s<2k/(2k-1)$, we have 
$$\|U\|_{L^r}\leq c_r\quad \mbox{and} \quad \|\nabla U\|_{L^s}\leq
c_s,$$
where $c_r,c_s$ are constants independent of $S$. Moreover, $U$
depends linearly and continuously on $S$ with respect to the weak
topology on the currents $S$ and the $L^r$ topology on $U$.  
\end{theorem}

The construction of $U$ uses a kernel constructed in
Bost-Gillet-Soul\'e \cite{BostGilletSoule}. We call $U$ {\it the Green quasi-potential} of
$S$. When $p=1$, two quasi-potentials of $S$ differ by a constant. So, the
solution is unique if we require that $\langle \omega_\FS^k, U\rangle
=0$. In this case, we have a bijective and bi-continuous correspondence
$S\leftrightarrow u$ between positive closed $(1,1)$-currents $S$ and
their normalized quasi-potentials $u$. 

By maximum principle, p.s.h. functions on a compact manifold are
constant. However, the interest of p.s.h. functions is their type of
local singularities.  S.T. Yau introduced in \cite{Yau} the useful
notion of quasi-p.s.h. functions. A {\it quasi-p.s.h. function} is
locally the difference of a p.s.h. function and a smooth one. 
Several properties of quasi-p.s.h. functions can be deduced from
properties of p.s.h. functions.
If $u$ is a quasi-p.s.h. function on $\P^k$ there is a constant
$c>0$ such that $\ddc u\geq -c\omega_\FS$. So, $\ddc u$ is the difference
of a positive closed $(1,1)$-current and a smooth positive closed
$(1,1)$-form: $\ddc u =(\ddc u+c\omega_\FS)-c\omega_\FS$. Conversely, if $S$ is a positive closed $(1,1)$-current
cohomologous to a real $(1,1)$-form $\alpha$, there is a
quasi-p.s.h. function $u$, unique up to a constant, such that $\ddc
u=S-\alpha$. The following proposition is easily obtained using a
convolution on the group of automorphisms of $\P^k$, 
see Demailly \cite{Demailly3} for analogous results on compact
K\"ahler manifolds.

\begin{proposition} Let $u$ be a quasi-p.s.h. function on 
$\P^k$ such that
  $\ddc u\geq -\omega_\FS$. Then, there is a sequence $(u_n)$
  of smooth quasi-p.s.h. functions decreasing to $u$ such that $\ddc u_n\geq
  -\omega_\FS$. In
  particular, if $S$ is a positive closed $(1,1)$-current on $\P^k$, then
  there are smooth positive closed $(1,1)$-forms $S_n$ converging to $S$.
\end{proposition}

A subset $E$ of $\P^k$ is {\it pluripolar} if it is contained in
$\{u=-\infty\}$ where $u$ is a quasi-p.s.h. function. It is {\it
  complete pluripolar} if there is a quasi-p.s.h. function $u$ such
that $E=\{u=-\infty\}$. It is easy to check that analytic sets are
complete pluripolar and that a countable union of pluripolar sets is
pluripolar. The following capacity is close to a notion of capacity 
introduced by H. Alexander in \cite{Alexander}. The interesting point
here is that our definition extends to general compact
K\"ahler manifold \cite{DinhSibony6}. We will see that the same idea
allows to define the capacity of a current.
Let $\Pc_1$ denote the set of quasi-p.s.h. functions $u$ on $\P^k$
such that $\max_{\P^k} u=0$. 
 {\it The capacity} of a Borel set $E$ in $\P^k$  is
$$\capacity(E):=\inf_{\varphi\in\Pc_1}\exp \Big(\sup_E u\Big).$$
The Borel set $E$ is pluripolar if and only if $\capacity(E)=0$.
It is not difficult to show that when the volume of $E$ tends to the
volume of $\P^k$ then $\capacity(E)$ tends to 1.

The space of d.s.h. functions (differences of quasi-p.s.h. functions) and the complex Sobolev space of
functions on compact K{\"a}hler manifolds were introduced by
the authors in \cite{DinhSibony6, DinhSibony11}. They satisfy strong compactness properties and  are
invariant under the action of holomorphic maps. 
Using them as test functions, permits to obtain several results in
complex dynamics.

A function on $\P^k$ is called {\it d.s.h.} if it is equal outside a
pluripolar set to the difference of two quasi-p.s.h. functions. We
identify two d.s.h. functions if they are equal outside a pluripolar
set. Let $\DSH(\P^k)$ denote the space of d.s.h. functions on $\P^k$. We
deduce easily from properties of p.s.h. functions that $\DSH(\P^k)$ is
contained in $L^p(\P^k)$ for $1\leq p<\infty$. If $u$ is d.s.h. then $\ddc
u$ can be written as the difference of two positive closed
$(1,1)$-currents which are cohomologous. 
Conversely, if $S^\pm$ are positive
closed $(1,1)$-currents of the same mass, then there is a
d.s.h. function $u$, unique up to a constant, such that $\ddc
u=S^+-S^-$. 

We introduce several equivalent norms on $\DSH(\P^k)$. 
Define
$$\|u\|_\DSH:=|\langle \omega_\FS^k,u\rangle|+\min\|S^\pm\|,$$
where the minimum is taken over positive closed $(1,1)$-currents
$S^\pm$ such that $\ddc
u=S^+-S^-$. The term $|\langle \omega_\FS^k,u\rangle|$
may be replaced by $\|u\|_{L^p}$
  with $1\leq p<\infty$; we then obtain equivalent norms. The space of d.s.h. functions endowed with
  the above norm is a Banach space. However, we will use on this space
  {\it a weaker topology}: we say that a sequence $(u_n)$ converges to $u$
  in $\DSH(\P^k)$ if $u_n$ converge to $u$ in the sense of currents and
  if $(u_n)$ is bounded with respect to $\|\cdot\|_\DSH$. Under the
  last condition on the DSH-norm,
the convergence in the sense of currents of $u_n$ is equivalent to the
convergence in $L^p$ for $1\leq p<\infty$. We have the following
proposition \cite{DinhSibony6}.

\begin{proposition}
Let $u$ be a d.s.h. function on $\P^k$ such that $\|u\|_\DSH\leq 1$. Then there
are negative quasi-p.s.h. function $u^\pm$ such that $u=u^+-u^-$,
$\|u^\pm\|_\DSH\leq c$ and $\ddc u^\pm\geq -c\omega_\FS$, where $c>0$ is a
constant independent of $u$. 
\end{proposition}

A positive measure on $\P^k$ is said to be {\it PC} \footnote{In dimension
  1, the measure is PC if and only if its local Potentials are
  Continuous.} if it can be extended to a continuous linear form on
$\DSH(\P^k)$. Here, the continuity is with respect to the weak
topology on d.s.h. functions.
A positive measure is {\it PB} \footnote{In dimension
  1, the measure is PB if and only if its local Potentials are
  Bounded.} if quasi-p.s.h. functions are integrable with respect to
this measure. PB measures have no mass on pluripolar sets and
 d.s.h. functions are integrable with respect to such measures. PC measures are always PB. Let $\mu$ be a
non-zero PB positive measure on $X$. Define
$$\|u\|_\mu:=|\langle \mu,u\rangle|+\min\|S^\pm\|,$$
with $S^\pm$ as above. We have the following useful property \cite{DinhSibony6}.

\begin{proposition} \label{prop_dsh_norm_eq}
The semi-norm $\|\cdot\|_\mu$ is in fact a norm on $\DSH(\P^k)$ which is 
equivalent to $\|\cdot\|_\DSH$.
\end{proposition}

One can extend the above notions to currents but the definitions are
slightly different. Let $\DSH^p(\P^k)$
denote the space generated by negative 
$(p,p)$-currents $\Phi$ such that $\ddc\Phi$ is the difference of two
positive closed $(p+1,p+1)$-currents.  A DSH $(p,p)$-current, i.e. a
current in $\DSH^p(\P^k)$, is not an $L^1$ form in general. Define the
$\|\Phi\|_\DSH$-norm of a negative current $\Phi$ in $\DSH^p(\P^k)$ by 
$$\|\Phi\|_\DSH:=\|\Phi\|+\min \|\Omega^\pm\|,$$ 
where $\Omega^\pm$ are positive closed such that
$\ddc\Phi=\Omega^+-\Omega^-$. For a general $\Phi$ in $\DSH^p(\P^k)$
define
$$\|\Phi\|_\DSH:=\min(\|\Phi^+\|_\DSH+\|\Phi^-\|_\DSH),$$
where $\Phi^\pm$ are negative currents in $\DSH^p(X)$ such that
$\Phi=\Phi^+-\Phi^-$. We also consider on this space {\it the weak
topology}: a sequence $(\Phi_n)$ converges to $\Phi$ in $\DSH^p(\P^k)$ if
it converges to $\Phi$ in the sense of currents and if
$(\|\Phi_n\|_\DSH)$ is bounded. Using a convolution on the group of
automorphisms of $\P^k$, we can show 
that smooth forms are dense in $\DSH^p(\P^k)$. 

A positive closed $(p,p)$-current $S$ is called {\it PB} if there is a
constant $c>0$ such that $|\langle S,\Phi\rangle|\leq c\|\Phi\|_\DSH$
for any real smooth $(k-p,k-p)$-form $\Phi$. The current $S$ is {\it
  PC} if it can be extended to a linear continuous form on
$\DSH^{k-p}(\P^k)$. The continuity is with respect to the weak
topology we consider on $\DSH^{k-p}(\P^k)$. PC currents
are PB. We will see that these notions correspond to currents with
bounded or continuous super-potentials. As a consequence of Theorem
\ref{th_pullback_local}, we have the following useful result.

\begin{proposition} Let $f:\P^k\rightarrow\P^k$ be a holomorphic
  surjective map. Then, the operator $f^*$ on smooth forms has a
  continuous extension $f^*:\DSH^p(\P^k)\rightarrow\DSH^p(\P^k)$. If $S$
  is a current on $\DSH^p(\P^k)$ with no mass on a Borel set $A$, then
  $f^*(S)$ has no mass on $f^{-1}(A)$. 
\end{proposition}

Another useful functional space is {\it the complex Sobolev space} $W^*(\P^k)$. Its
definition uses the complex structure of $\P^k$.
In dimension one, $W^*(\P^1)$ coincides
with the Sobolev space $W^{1,2}(\P^1)$ of real-valued functions in $L^2$ with gradient in
$L^2$. In higher dimension, $W^*(\P^k)$ is the space of
functions $u$ in $W^{1,2}(\P^k)$ such that $i\partial u\wedge \dbar u$
is bounded by a positive closed $(1,1)$-current $\Theta$. We define
$$\|u\|_{W^*}:=|\langle \omega_\FS^k,u\rangle|+\min\|\Theta\|^{1/2}$$
with $\Theta$ as above, see \cite{DinhSibony11, Vigny3}.
By Sobolev-Poincar{\'e} inequality, the term $|\langle \omega_\FS^k,u\rangle|$ may be replaced by
$\|u\|_{L^1}$ or $\|u\|_{L^2}$; we then obtain equivalent norms.
{\it The weak topology} on $W^*(\P^k)$ is defined as in the case of
d.s.h. functions: a sequence $(u_n)$ converges in $W^*(\P^k)$ to a function
$u$ if it converges to $u$ in the sense of currents and if
$(\|u_n\|_{W^*})$ is bounded. 
A positive measure $\mu$ is {\it WPC} if it can be extended to a
linear continuous form on $W^*(\P^k)$. If $u$ is a strictly negative
quasi-p.s.h. function, one can prove that $\log(-u)$ is in
$W^*(\P^k)$. This allows to show that WPC measures have no mass on
pluripolar sets.

In the rest of the paragraph, we will introduce the notion of
super-potentials associated to positive closed
$(p,p)$-currents.
They are  canonical functions defined on an infinite dimensional spaces and
are, in some sense, quasi-p.s.h. functions there. Super-potentials were
introduced by the authors in order to replace ordinary quasi-p.s.h. functions
which are used as quasi-potentials for currents of bidegree $(1,1)$. The
theory is satisfactory in the case of projective spaces
\cite{DinhSibony10} and can be easily extended to homogeneous
manifolds.

Let $\Cc_{k-p+1}(\P^k)$ denote the convex set of positive closed currents of
bidegree $(k-p+1,k-p+1)$ and of
mass 1, i.e. currents cohomologous to $\omega_\FS^{k-p+1}$. 
Let $S$ be a positive closed $(p,p)$-current on $\P^k$. We
assume for simplicity that $S$ is of mass 1; the general case can be
deduced by linearity. The super-potential\footnote{The super-potential
  we consider here corresponds to the super-potential of mean 0 in
  \cite{DinhSibony10}. The other super-potentials differ from $\Uc_S$ by 
  constants.}
$\Uc_S$ of $S$ is a function
on $\Cc_{k-p+1}(\P^k)$ with values in $\R\cup\{-\infty\}$. Let $R$ be a
current in $\Cc_{k-p+1}(\P^k)$ and $U_R$ a potential of
$R-\omega_\FS^{k-p+1}$. Subtracting from $U_R$ a constant times
$\omega_\FS^{k-p}$ allows to have $\langle U_R,\omega_\FS^p\rangle=0$.
We say that $U_R$ is a quasi-potential of mean 0 of $R$. 
Formally, i.e. in the case where $R$ and $U_R$ are smooth, the value of $\Uc_S$ at $R$ is defined by
$$\Uc_S(R):=\langle S,U_R\rangle.$$
One easily check using Stokes' formula that formally if $U_S$ is a
quasi-potential of mean 0 of $S$, then $\Uc_S(R)=\langle
U_S,R\rangle$. Therefore, the previous definition does not depend on
the choice of $U_R$ or $U_S$. By definition, we have
$\Uc_S(\omega_\FS^{k-p+1})=0$. 
Observe also that when $S$ is smooth, the above definition makes sense
for every $R$ and $\Uc_S$ is a continuous affine function on $\Cc_{k-p+1}(\P^k)$.
It is also clear that if $\Uc_S=\Uc_{S'}$, then $S=S'$. 
The following theorem  allows to define
$\Uc_S$ in the general case. 

\begin{theorem} The above function $\Uc_S$, which is defined on smooth
  forms $R$ in $\Cc_{k-p+1}(\P^k)$, can be extended to an affine function on
  $\Cc_{k-p+1}(\P^k)$ with values in $\R\cup\{-\infty\}$ by
$$\Uc_S(R):=\limsup \Uc_S(R'),$$
where $R'$ is smooth in $\Cc_{k-p+1}(\P^k)$ and converges to $R$. 
We have $\Uc_S(R)=\Uc_R(S)$. Moreover,
there are smooth positive closed $(p,p)$-forms $S_n$ of mass $1$ and
constants $c_n$ converging to $0$ such that $\Uc_{S_n}+c_n$ decrease
to $\Uc_S$. In particular, $\Uc_{S_n}$ converge pointwise to $\Uc_S$.  
\end{theorem}

For bidegree $(1,1)$, there is a unique quasi-p.s.h. function $u_S$ such that
$\ddc u_S=S-\omega_\FS$ and $\langle \omega_\FS^k, u_S\rangle=0$. If $\delta_a$
denotes the Dirac mass at $a$, we have
$\Uc_S(\delta_a)=u_S(a)$. Dirac masses are extremal elements in
$\Cc_k(\P^k)$. The super-potential $\Uc_S$ in this case is just the affine
extension of $u_S$, that is, we have for any probability measure
$\nu$:
$$\Uc_S(\nu)=\int \Uc_S(\delta_a)d\nu(a)=\int u_S(a) d\nu(a).$$

The function $\Uc_S$ extends the action $\langle S,\Phi\rangle$ on
smooth forms $\Phi$ to $\langle S,U\rangle$ where $U$ is a quasi-potential
of a positive closed current.
Super-potentials satisfy analogous properties as
quasi-p.s.h. functions do. They are upper semi-continuous and bounded
from above by a universal constant. Note that we consider here the
weak topology on $\Cc_{k-p+1}(\P^k)$. We  have 
the following version of the Hartogs' lemma.

\begin{proposition} Let $S_n$ be positive closed $(p,p)$-currents of
  mass $1$ on $\P^k$ converging to $S$. Then for every continuous
  function $\Uc$ on $\Cc_{k-p+1}$ with $\Uc_S<\Uc$, we have
  $\Uc_{S_n}<\Uc$ for $n$ large enough. In particular, $\limsup
  \Uc_{S_n}\leq \Uc_S$. 
\end{proposition}

We say that $S_n$ converge to $S$ {\it in the Hartogs' sense} if
$S_n$ converge to $S$ and if there are constants $c_n$ converging to 0
such that $\Uc_{S_n}+c_n\geq \Uc_S$. If $\Uc_{S_n}$ converge uniformly
to $\Uc_S$, we say that $S_n$ {\it converge SP-uniformly} to $S$.

One can check that PB and PC
currents correspond to currents of bounded or continuous
super-potential. In the case of bidegree $(1,1)$, they correspond to currents with
bounded or continuous quasi-potential. We say that $S'$ is {\it more
diffuse than} $S$ if $\Uc_{S'}-\Uc_S$ is bounded from below. So, PB
currents are more diffuse than any other currents.

In order to prove the above results and to work with super-potentials,
we have to consider a geometric structure on $\Cc_{k-p+1}(\P^k)$. In a weak sense,
$\Cc_{k-p+1}(\P^k)$ can be seen as a space of infinite dimension which
contains many ``analytic'' sets of finite dimension that we call {\it
  structural varieties}. Let $V$ be a complex manifold and $\Rc$ a
positive closed current of bidegree $(k-p+1,k-p+1)$ on $V\times \P^k$. Let $\pi_V$ denote
the canonical projection map from $V\times\P^k$ onto $V$. One can prove
that the slice $\langle
\Rc,\pi_V,\theta\rangle$ is defined for $\theta$ outside a locally
pluripolar set of $V$. Each slice can be identified with a positive
closed $(p,p)$-current $R_\theta$ on $\P^k$. Its mass does not depend
on $\theta$. So, multiplying $\Rc$ with a constant, we can assume
that all the $R_\theta$ are in $\Cc_{k-p+1}(\P^k)$. The map
$\tau(\theta):=R_\theta$ or the family $(R_\theta)$ is called {\it
 a structural variety} of $\Cc_{k-p+1}(\P^k)$. The restriction of $\Uc_S$ to
this structural variety, i.e. $\Uc_S\circ\tau$, is locally a
d.s.h. function or identically $-\infty$. 
When the structural variety is nice enough, this restriction is
quasi-p.s.h. or identically $-\infty$. In practice, we often use
some special structural discs parametrized by $\theta$ in the unit
disc of $\C$. They are obtained by convolution of a given current $R$ with
a smooth probability measure on the group $\PGL(\C,k+1)$ of
automorphisms of $\P^k$. 

Observe that since the correspondence $S\leftrightarrow \Uc_S$ is
$1:1$, the compactness on positive closed currents should
induce some compactness on super-potentials. We have the following result.

\begin{theorem} \label{th_hormander_sp}
Let $W\subset \P^k$ be an open set and $K\subset W$ a compact
set. Let $S$ be a current in $\Cc_p(\P^k)$ 
with support in $K$ and $R$  a current in $\Cc_{k-p+1}(\P^k)$. Assume that
the restriction of $R$ to $W$ is a bounded form.
Then, the super-potential $\Uc_S$ of $S$ satisfies 
$$|\Uc_S(R)|\leq
A\big(1+\log^+\|R\|_{\infty,W}\big)$$
where $A>0$ is a constant independent of $S$, $R$ and $\log^+:=\max(0,\log)$.
\end{theorem}

This result can
be applied to $K=W=\P^k$ and can be
considered as a version of the exponential estimate in Theorem \ref{th_hormander}.
Indeed, the weaker estimate $|\Uc_S(R)|\lesssim
1+\|R\|_\infty$ is easy to obtain. It corresponds to the $L^1$
estimate on the quasi-p.s.h. function $u_S$ in the case of bidegree
$(1,1)$. 

Using the analogy with the bidegree $(1,1)$ case, we define {\it the
  capacity} of a current $R$ as
$$\capacity(R):=\inf_S\exp\big(\Uc_S(R)-\max \Uc_S\big).$$
This capacity describes rather the regularity of $R$: an $R$ with
big capacity is somehow more regular.
Theorem \ref{th_hormander_sp} implies that 
$\capacity(R)\gtrsim \|R\|_\infty^{-\lambda}$
for some universal constant $\lambda>0$. 
This property is close to the capacity estimate for Borel sets in term
of volume.

Super-potentials allow to develop a theory of intersection of currents
in higher bidegree. Here, the fact that $\Uc_S$ has a value at every
point (i.e. at every current $R\in\Cc_{k-p+1}(\P^k)$) is crucial. Let $S$, $S'$ be positive closed currents of
bidegree $(p,p)$ and $(p',p')$ with $p+p'\leq k$. We assume for simplicity that their
masses are equal to 1. We say that $S$ and $S'$ are {\it wedgeable} if
$\Uc_S$ is finite at $S'\wedge \omega_\FS^{k-p-p'+1}$. This property
is symmetric on $S$ and $S'$. 
If $\widetilde S, \widetilde S'$ are more diffuse than
$S,S'$ and if $S,S'$ are wedgeable, then  $\widetilde S,
\widetilde S'$ are wedgeable.

Let $\Phi$ be a real smooth
form of bidegree $(k-p-p',k-p-p')$. Write $\ddc
\Phi=c(\Omega^+-\Omega^-)$ with $c\geq 0$ and $\Omega^\pm$ positive
closed of mass 1.  If $S$ and $S'$ are wedgeable, define
the current $S\wedge S'$ by
$$\langle S\wedge S',\Phi\rangle := \langle S',\omega_\FS^p\wedge\Phi\rangle +
c\Uc_S(S'\wedge \Omega^+)-c\Uc_S(S'\wedge \Omega^-).$$
A simple computation shows that the definition coincides with the usual
wedge-product when $S$ or $S'$ is smooth. One can also prove that the
previous definition does not depend on the choice of $c$, $\Omega^\pm$
and is symmetric with respect to $S,S'$. If $S$ is of bidegree
$(1,1)$, then $S,S'$ are wedgeable if and only if the quasi-potentials
of $S$ are integrable with respect to the trace measure of $S'$. In
this case, the above definition coincides with the definition in
Section \ref{section_intersection}. We have the following general result.

\begin{theorem} \label{th_wedge_product_pp}
Let $S_i$ be positive closed currents of bidegree
  $(p_i,p_i)$ on $\P^k$ with $1\leq i\leq m$ and $p_1+\cdots+p_m\leq
  k$. Assume that for $1\leq i\leq m-1$, $S_i$ and
  $S_{i+1}\wedge\ldots\wedge S_m$ are wedgeable. Then, this condition
  is symmetric on $S_1,\ldots, S_m$. The wedge-product $S_1\wedge
  \ldots\wedge S_m$ is a positive closed current of mass
  $\|S_1\|\ldots\|S_m\|$ supported on $\supp(S_1)\cap\ldots\cap\supp(S_m)$.
It depends linearly on each variable and is
  symmetric on the variables. If $S_i^{(n)}$ converge to $S_i$ in the
  Hartogs' sense, then the $S_i^{(n)}$ are wedgeable and 
$S_1^{(n)}\wedge \ldots \wedge S_m^{(n)}$ converge in the Hartogs'
sense to $S_1\wedge \ldots\wedge S_m$.
\end{theorem}

We deduce from this result that wedge-products of PB currents are PB. One
can also prove that wedge-products of PC currents are PC.
If $S_n$ is defined by analytic sets, they are wedgeable if the
intersection of these analytic sets is of codimension
$p_1+\cdots+p_m$. In this case, the intersection in the sense of
currents coincides with the intersection of cycles, i.e. is equal to the
current of integration on the intersection of cycles where we count
the multiplicities. We have the following criterion of 
wedgeability which contains the case of cycles.

\begin{proposition}
Let $S,S'$ be positive closed currents on $\P^k$ of bidegrees $(p,p)$
and $(p',p')$. Let $W$, $W'$ be open sets such that $S$ restricted to
$W$ and $S'$ restricted to $W'$ are  bounded forms. Assume that $W\cup
W'$ is $(p+p')$-concave in the sense that there is a positive closed
smooth form of bidegree $(k-p-p'+1,k-p-p'+1)$ with compact support in $W\cup
W'$. Then $S$ and $S'$ are wedgeable.
\end{proposition}

The following result can be deduced from Theorem
\ref{th_wedge_product_pp}.

\begin{corollary}
Let $S_i$ be positive closed $(1,1)$-currents 
on $\P^k$ with $1\leq i\leq p$.
Assume that for $1\leq i\leq p-1$, $S_i$ admits a quasi-potential
which is integrable with respect to the trace measure of 
  $S_{i+1}\wedge\ldots\wedge S_p$. Then, this condition
  is symmetric on $S_1,\ldots, S_p$. The wedge-product $S_1\wedge
  \ldots\wedge S_p$ is a positive closed $(p,p)$-current of mass
  $\|S_1\|\ldots\|S_p\|$ supported on $\supp(S_1)\cap\ldots\cap\supp(S_p)$.
It depends linearly on each variable and is
  symmetric on the variables. If $S_i^{(n)}$ converge to $S_i$ in the
  Hartogs' sense, then the $S_i^{(n)}$ are wedgeable and 
$S_1^{(n)}\wedge \ldots \wedge S_p^{(n)}$ converge to $S_1\wedge \ldots\wedge S_p$.
\end{corollary}

We discuss now currents with H\"older continuous super-potential and
moderate currents.
The space $\Cc_{k-p+1}(\P^k)$ admits a natural distances $\dist_\alpha$,
with $\alpha>0$, defined by 
$$\dist_\alpha(R,R'):=\sup_{\|\Phi\|_{\Cc^\alpha}\leq 1} |\langle
R-R',\Phi\rangle|,$$
where $\Phi$ is a smooth $(p-1,p-1)$-form on $\P^k$. 
The norm $\Cc^\alpha$ on $\Phi$ is the sum of the $\Cc^\alpha$-norms
of its coefficients for a fixed atlas of $\P^k$. 
The topology
associated to $\dist_\alpha$ coincides with the weak topology. Using
the theory of interpolation between Banach spaces \cite{Triebel}, we obtain for
$\beta>\alpha>0$ that
$$\dist_\beta\leq\dist_\alpha\leq
c_{\alpha,\beta}[\dist_\beta]^{\alpha/\beta}$$
where $c_{\alpha,\beta}>0$ is a constant. So, a function on
$\Cc_{k-p+1}(\P^k)$ is H{\"o}lder continuous with respect to $\dist_\alpha$ if
and only if it is H{\"o}lder continuous with respect to
$\dist_\beta$. The following proposition is useful in dynamics.

\begin{proposition} The wedge-product of positive closed currents
on $\P^k$ with H{\"o}lder continuous super-potentials has a H{\"o}lder continuous
  super-potential.  
Let $S$ be a positive closed $(p,p)$-current with a H\"older
continuous super-potential. Then, the Hausdorff
dimension of $S$ is strictly larger than $2(k-p)$. Moreover, $S$ is moderate,
i.e. for any bounded family $\Fc$ of d.s.h. functions on $\P^k$, there
are constants $c>0$ and $\alpha>0$ such that 
$$\int e^{\alpha|u|}d\sigma_S\leq c$$
for every $u$ in $\Fc$, where $\sigma_S$ is the trace measure of $S$. 
\end{proposition}

\bigskip\bigskip

\begin{exercise} \label{exo_cap_vol}
Show that there is a constant $c>0$ such that 
$$\capacity(E)\geq
\exp(-c/\volume(E)).$$ 
Hint: use the compactness of $\Pc_1$ in $L^1$. 
\end{exercise}

\begin{exercise} \label{exercise_sum_dsh}
Let $(u_n)$ be a sequence of d.s.h. functions such that $\sum
\|u_n\|_\DSH$ is finite. Show that $\sum u_n$ converge pointwise out of a
pluripolar set to a d.s.h. function. Hint: write $u_n=u_n^+-u_n^-$
with $u_n^\pm\leq 0$, $\|u_n^\pm\|_\DSH\lesssim \|u_n\|_\DSH$ and
$\ddc u_n^\pm\geq -\|u_n\|_\DSH\omega_\FS$. 
\end{exercise}

\begin{exercise}
If $\chi$ is a convex increasing function on $\R$ with bounded
derivative and $u$ is a d.s.h. function, show that $\chi\circ u$ is
d.s.h. If $\chi$ is Lipschitz and $u$ is in $W^*(\P^k)$, show that $\chi\circ
u$ is in $W^*(\P^k)$. Prove that bounded d.s.h. functions are in
$W^*(\P^k)$. 
Show that $\DSH(\P^k)$ and  $W^*(\P^k)$ are stable under the max and min operations.
\end{exercise}

\begin{exercise}
Let $\mu$ be a non-zero positive measure which is WPC. Define 
$$\|u\|_\mu^*:=|\langle \mu, u\rangle|+\min\|\Theta\|^{1/2}$$
with $\Theta$ as above. Show that $\|\cdot\|_\mu^*$ defines a norm
which is equivalent
to $\|\cdot\|_{W^*}$.
\end{exercise}

\begin{exercise}
Show that the capacity of $R$ is positive  if and only if $R$ is PB.
\end{exercise}

\begin{exercise}
Let $S$ be a positive closed $(p,p)$-current of mass $1$ with positive
Lelong number at a point $a$. Let
$H$ be a hyperplane containing $a$ such that $S$ and $[H]$ are
wedgeable. Show that the Lelong number of $S\wedge [H]$ at $a$
is the same if we consider it as a current on $\P^k$
or on $H$. If $R$ is a positive closed current of bidimension
$(p-1,p-1)$ on $H$, show that $\Uc_S(R)\leq \Uc_{S\wedge [H]}(R)+c$
where $c>0$ is a constant independent of $S,R$ and $H$. Deduce that PB
currents have no positive Lelong numbers.
\end{exercise}

\begin{exercise} Let $K$ be a compact subset in $\C^k\subset \P^k$. Let
  $S_1,\ldots,S_p$ be positive closed $(1,1)$-currents on
  $\P^k$. Assume that their quasi-potentials are bounded on
  $\P^k\setminus K$. Show that $S_1,\ldots,S_p$ are wedgeable. Show
  that the wedge-product $S_1\wedge\ldots\wedge S_p$ is continuous for
  Hartogs' convergence.
\end{exercise}

\end{appendix}


\

\noindent
T.-C. Dinh, UPMC Univ Paris 06, UMR 7586, Institut de
Math{\'e}matiques de Jussieu, F-75005 Paris, France. {\tt
  dinh@math.jussieu.fr}, {\tt http://www.math.jussieu.fr/$\sim$dinh} 

\

\noindent
N. Sibony,
Universit{\'e} Paris-Sud, Math{\'e}matique - B{\^a}timent 425, 91405
Orsay, France. {\tt nessim.sibony@math.u-psud.fr} 
\end{document}